\documentclass{article}
\usepackage{amsfonts, color}
\usepackage{textgreek}
\usepackage[dvipsnames]{xcolor}
\usepackage{geometry}
\usepackage{enumitem}
\usepackage[normalem]{ulem}
\usepackage{relsize}
\usepackage{comment}
\usepackage{lineno}
\usepackage{xstring}
\usepackage{braket}
\usepackage{amstext, amsthm, amssymb, amsmath, mathtools}
\numberwithin{equation}{section}
\usepackage{mathrsfs}
\usepackage{bbm}
\usepackage{graphicx}
\usepackage{float}
\usepackage[hypertexnames=false]{hyperref}
\usepackage[capitalize,nameinlink,noabbrev]{cleveref}
\usepackage{xr}
\makeatletter
\newcommand*{\addFileDependency}[1]{%
  \typeout{(#1)}%
  \@addtofilelist{#1}%
  \IfFileExists{#1}{}{\typeout{No file #1.}}%
}
\makeatother
\newcommand*{\myexternaldocument}[2][]{%
  \externaldocument[#1]{#2}%
  \addFileDependency{#2.tex}%
  \addFileDependency{#2.aux}%
}
\myexternaldocument[M-]{main_metric}
\usepackage[style=alphabetic, maxbibnames=99, isbn=false]{biblatex}
\usepackage[colorinlistoftodos,textsize=small]{todonotes}
\usepackage{bookmark}
\usepackage[utf8]{inputenc}
\makeatletter\def\@todonotes@font{\rmfamily}\makeatother
\presetkeys{todonotes}{color=yellow!50}{}

\usepackage{subcaption}
\usepackage{adjustbox}
\usepackage{extpfeil}
\usepackage{caption}
\usepackage{longtable}
\usepackage{booktabs}
\usepackage{array}
\usepackage{tikz}
\usepackage{tikz-cd}
\usepackage{bm}
\LoadFontDefinitionFile{U}{stmry}
\DeclareFontShape{U}{stmry}{b}{n}{<->ssub * stmry/m/n}{}
\usepackage{stackengine}
\usepackage{thmtools}
\usepackage{soul}

\setstackEOL{\\}
\usetikzlibrary{calc,decorations.pathreplacing}

\usepackage{crossreftools}

\newcommand{\Rext}{\mathbb{R} \cup \{+\infty\}}

\newcommand{\diff}{\mathrm{d}}

\newcommand{\law}{\mathcal{L}}

\newcommand{\Pp}{\mathcal{P}}
\newcommand{\R}{\mathbb{R}}

\newcommand{\EE}{{\mathbb{E}}}

\newcommand{\NN}{{\mathbb N}} 
\newcommand{\Q}{\mathbb{Q}}

\newcommand{\ACt}{\ensuremath{AC_{\mathrm{II}}}}

\NewDocumentCommand{\Tang}{ O{\varepsilon} }{
  \mathcal{T}_{#1} \mathcal{S}
}

\DeclareMathOperator*{\argmin}{arg\,min}
\newcommand{\limsupTauTS}{\limsup_{\substack{T \searrow s,\, t \searrow s, \\ T \geq t}}}

\usepackage{xcolor}

\newif\ifshowcomments
\showcommentstrue

\newcommand{\FundingLogos}{%
  \raisebox{0pt}{\includegraphics[height=1.5cm]{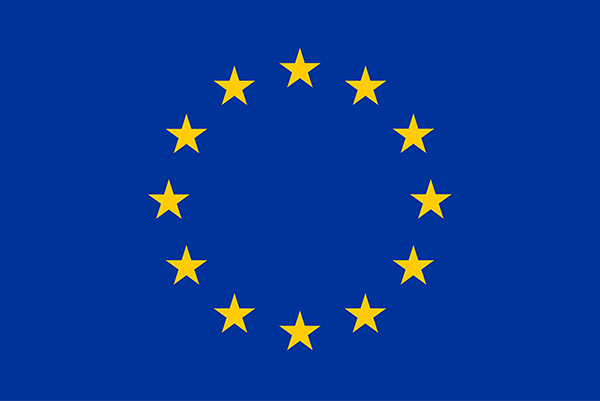}}%
  \hspace{1em}%
  \raisebox{0pt}{\includegraphics[height=1.5cm]{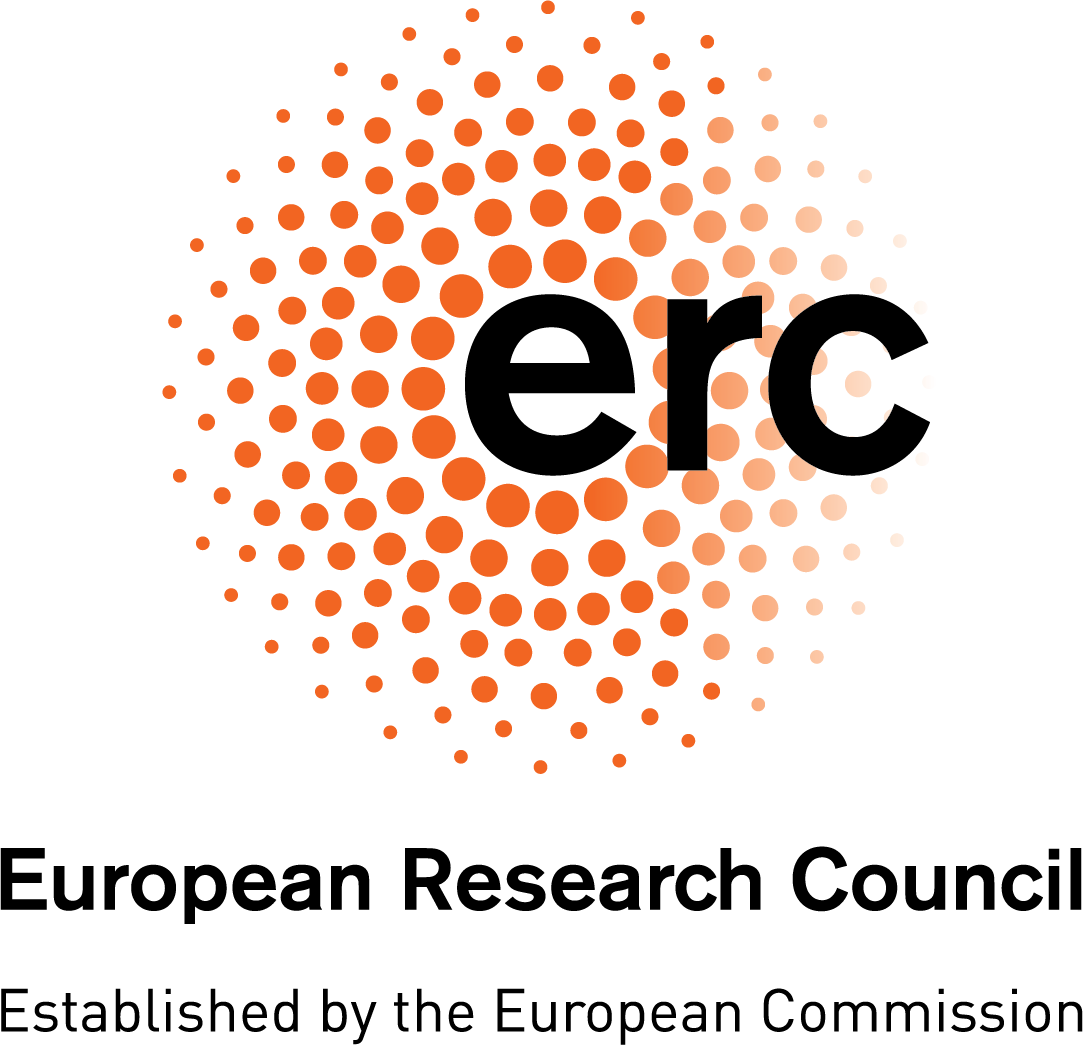}}%
}

\usepackage{tcolorbox}
\tcbset{colback=white,boxrule=0.5pt,arc=2pt,left=4pt,right=4pt,top=4pt,bottom=4pt}

\theoremstyle{plain}
\newtheorem{theorem}{Theorem}[section]
\newtheorem{lemma}[theorem]{Lemma}
\newtheorem{proposition}[theorem]{Proposition}
\newtheorem{corollary}[theorem]{Corollary}
\theoremstyle{definition} 
\newtheorem{claim}{Claim}[section] 

\theoremstyle{definition}
\newtheorem{definition}{Definition}
\newtheorem{assumption}{Assumption}

\newtheoremstyle{boldremark}%
  {0.5\topsep}{0.5\topsep}
  {\normalfont}
  {}
  {\bfseries}
  {.}
  {.5em}
  {}
\theoremstyle{boldremark}
\newtheorem{remark}{Remark}
\theoremstyle{remark}
\newtheorem{example}{Example}

\usepackage{authblk}

\let\oldproof\proof
\let\oldendproof\endproof

  \newlist{claimlist}{enumerate}{1}
  \setlist[claimlist]{label= (\arabic*)}
  \crefname{claimlisti}{Claim}{Claims}
  \Crefname{claimlisti}{Claim}{Claims}
  \newlist{conditionlist}{enumerate}{1}
  \setlist[conditionlist]{label= (\arabic*),nosep}
  \crefname{conditionlisti}{Condition}{Conditions}
  \Crefname{conditionlisti}{Condition}{Conditions}
\crefname{theorem}{Theorem}{Theorems}
\crefname{lemma}{Lemma}{Lemmas}
\crefname{corollary}{Corollary}{Corollaries}
\crefname{section}{Section}{Sections}
\crefname{proposition}{Proposition}{Propositions}
\crefalias{lemma}{lemma}
\crefalias{corollary}{corollary}
\crefalias{proposition}{proposition}
\crefname{definition}{Definition}{Definitions}
\crefname{remark}{Remark}{Remarks}
\crefname{conjecture}{Conjecture}{Conjectures}
\crefname{example}{Example}{Examples}
\crefname{annoyingproof}{Proof}{Proofs}
\crefname{assumption}{Assumption}{Assumptions}
\crefname{algorithm}{Algorithm}{Algorithms}

\let\ORGhypersetup\hypersetup
\protected\def\hypersetup{\ORGhypersetup}
\title{A Lions Lift Framework for Gradient and Newtonian Evolutions in Wasserstein Spaces}

\author[1,2]{Massimo Fornasier\thanks{\texttt{massimo.fornasier@cit.tum.de}}}

\author[1,2]{Jona Klemenc\thanks{\texttt{jona.klemenc@tum.de}}}
\author[1,2]{Alessandro Scagliotti\thanks{\texttt{scag@ma.tum.de}}}

\affil[1]{Department of Mathematics, CIT School, Technical University of Munich, Munich, Germany}
\affil[2]{Munich Center for Machine Learning (MCML), Munich, Germany}

\usepackage[algoruled, linesnumbered]{algorithm2e}
\crefname{AlgoLine}{line}{lines}
\crefname{algocf}{algorithm}{algorithms}
\makeatletter
\let\oldnl\nl
\newcommand{\nonl}{\renewcommand{\nl}{\let\nl\oldnl}}
\makeatother
\tikzset{degil/.style={
                decoration={markings,
                mark= at position 0.5 with {
                \node[transform shape] (tempnode) {$\backslash$};
                }
                },
                postaction={decorate}
}
}
\crefname{enumi}{point}{points}
\Crefname{enumi}{Point}{Points}
\newcounter{property}

\newlist{properties}{enumerate}{1}
\setlist[properties]{
  label=\textbf{Property \arabic*:},
  ref=\arabic*,
  leftmargin=*,
}

\crefname{propertiesi}{property}{properties}
\Crefname{propertiesi}{Property}{Properties}

\newcommand{\propertyitem}{%
  \refstepcounter{property}%
  \item
}
\makeatletter
\newcommand{\extrabold}[1]{\mathpalette\extrabold@{#1}}
\newcommand{\extrabold@}[2]{%
  \ooalign{%
    $\m@th#1\bm{#2}$\cr
    \hidewidth$\m@th#1\kern0.02em\bm{#2}$\hidewidth\cr
    \hidewidth$\m@th#1\kern-0.02em\bm{#2}$\hidewidth\cr
    \hidewidth$\m@th#1\raise0.06ex\hbox{$\m@th#1\bm{#2}$}$\hidewidth\cr
    \hidewidth$\m@th#1\lower0.02ex\hbox{$\m@th#1\bm{#2}$}$\hidewidth\cr
  }%
}
\makeatother
\newcommand{\bgamma}{\extrabold{\gamma}}
\newcommand{\norm}[1]{\left\lVert#1\right\rVert}
\newcommand{\Leb}{\mathrm{Leb}}
\makeatletter

\newcommand\ful[2][\dimexpr2.45\fontdimen5\textfont2\relax]{%
  \mathpalette{\ful@{#1}}{#2}%
}

\newcommand\ful@[3]{%
  \vbox to0pt{%
    \vss
    \hrule height .7pt
    \vskip-\normalbaselineskip
    \kern#1
    \hbox{%
      \rlap{%
        \color{white}%
        \smash[t]{$\m@th#2#3$}%
      }%
      \smash[t]{$\m@th#2#3$}%
    }%
  }%
}

\makeatother
\newcommand{\LtwoR}{L^2_{\R}}
\newcommand{\LtwoRful}{\ful{L}^2_{\R}}

\begin{document}
\maketitle
\begin{abstract}
  To analyze evolutions driven by functionals over the space of probability measures, one can opt for the extrinsic Lions lift approach or the intrinsic Wasserstein viewpoint with its abstract metric structure. While the latter has been a successful playground for the study of first-order gradient flows, it classically relies on $\lambda$-convexity, and has eluded attempts of extensions to second-order Newtonian-like dynamics. Motivated by these obstructions, we show that the extrinsic Lions lift has the full capabilities to carry out the first-order gradient-flow program, which is classically done intrinsically in the Wasserstein space. As illustrated in an example, this paves the way for a framework of gradient evolutions beyond the $\lambda$-convex case, which is classically demanded for energies in the Wasserstein space. Furthermore, our interpretation of the Lions lift enables the formulation of a notion of second-order Newtonian curves in the Wasserstein phase space. Such curves can be constructed through a second-order JKO-like scheme and, when applied to linear and interacting energies, naturally lead to Vlasov kinetic equations and existence of solutions thereof.
\end{abstract}
\noindent {\it Keywords:} Lions lift, Wasserstein space, gradient flows, 
JKO scheme, Newtonian dynamics, Wasserstein phase space, Vlasov equation, interaction energies
\\

\noindent {\it MSC:} Primary 49Q22, 35Q83; Secondary 49J27, 35R06, 47J35, 49M25, 82C40\\

\noindent
\textit{AI use disclosure:} 
All concepts, research questions and proofs in the main text were developed by the authors. 
AI tools were employed for bug checking, copy-editing tasks and for the coding of the figures.
For more details, see \cref{tab:ai-disclosure}.
\tableofcontents

\section{Introduction}
To analyze gradient flows on the space of probability measures on $\R^d$ with respect to its Euclidean structure on $\R^d$, two distinct perspectives have evolved over the past decades.

On the one hand, we find the analytic perspective of metric gradient flows, which extends the notion of gradient flows to general metric spaces as \emph{curves of maximal slope} (see \cite{ambrosioGradientFlows2008}). Curves of maximal slope are characterized solely by identities between scalar values, as opposed to the ODE
\begin{equation} \label{eq:intro_ode}
  \dot u = - \nabla E,
\end{equation}
which describes ordinary gradient flows and is an identity settled in the tangent space, which is generally unavailable in metric spaces. More precisely, one first observes that for smooth functionals $E$ over some Hilbert space $\mathcal{H}$ and absolutely continuous curves $u\colon [0, T] \to \mathcal{H}$, we always have, for $s, t \in [0, T]$ with $s<t$,
\begin{equation} \label{eq:intro_coms_inequlality_hilbert}
  E\big(u(s)\big) - E\big(u(t)\big) \leq \frac{1}{2} \int_s^t |\dot{u}(\sigma)|^2 + |\nabla E\big(u(\sigma)\big)|^2 \, \diff \sigma,
\end{equation}
where equality holds if and only if $u$ fulfills \eqref{eq:intro_ode} for almost every $\sigma \in [s, t]$. For a general metric space $\mathcal{S}$, $|\dot{u}(\sigma)|$ may be replaced by
\begin{equation*}
  |u'|(\sigma) \coloneqq \limsup_{\tilde{\sigma} \to \sigma} \frac{d_{\mathcal{S}}(u(\sigma), u(\tilde{\sigma}))}{|\sigma - \tilde{\sigma}|}
\end{equation*}
and a suitable replacement $g$ for $|\nabla E(\cdot)|$ is formulated such that \cref{eq:intro_coms_inequlality_hilbert} is fulfilled for any absolutely continuous curve $u$; such a replacement $g$ is called an upper gradient. 
We then say that an absolutely continuous trajectory is a curve of maximal slope if equality holds\footnote{For brevity, we gloss over the technicalities involving weak and strong upper gradients in this introduction, see \cite[Def. 1.2.1, Def. 1.2.2]{ambrosioGradientFlows2008}; the exposition given here is slightly simplified.}. 
A prime candidate for such an upper gradient is the metric slope
\begin{equation*}
  |\partial E|(x) = \limsup_{y \to x} \frac{(E(x) - E(y))^+}{d(x, y)}
\end{equation*}
and the relaxed metric slope
\begin{equation*}
  |\partial^- E|(x) = \inf_{y_i \to x}\liminf_{i \to \infty} |\partial E|(y_i).
\end{equation*}
Furthermore, this analytic framework comes with an algorithmic machinery to \emph{construct} such curves of maximal slope using the minimizing movement scheme, introduced in~\cite{degiorgi1993new}: Fixing $U_0$ and for $\tau > 0$, we recursively choose $U_\tau^{n+1} \in J_{\tau, \mathcal{S}}[U_\tau^{n}] = \argmin_{y \in \mathcal{S}}[E(y) + \frac{d(U_\tau^n, y)^2}{2\tau}]$. Then, setting $\bar{U}_\tau(t) \coloneqq U_\tau^{\lceil t/\tau \rceil}$, we have that under certain topology and regularity assumptions\footnote{Proving these regularity assumptions is an application-specific challenge and is often done using ad-hoc tools. Under the assumptions of $\lambda$-convexity, much of this ad-hoc work simplifies.}, there exists a limit curve $\bar{u}\colon [0, +\infty) \to \mathcal{S}$ such that $\bar{U}_\tau(t) \to \bar{u}(t)$ along a subsequence, for all $t$, and this limiting curve $\bar{u}$ is a curve of maximal slope. 
Even more, we can take as a metric space $\mathcal{S}$ the Wasserstein space $(\Pp_2, W_2)$ and combine this framework with the analysis of the tangent space on $(\Pp_2, W_2)$ spearheaded by Otto~\cite{OttoGeometry2001}.
This standpoint has unveiled a dissipative structure underlying---and allowing the analysis of---a whole class of relevant PDEs.

On the other hand, the probabilistic viewpoint developed by Lions around 2007 (see \cite{cardaliaguetNotesMeanField,carmonaProbabilisticTheoryMean2018}) in the context of mean field games approaches energies $E\colon \Pp_2 \to \R$ by composing them with the function $\law\colon L^2(\Omega; \R^d) \to \Pp_2$, where $\Omega$ is some fixed and ``large enough'' probability space. 
The resulting function $\tilde{E}\colon L^2(\Omega; \R^d) \to \R$ is called the \emph{Lions lift} (or \emph{Lagrangian lift}) of $E$.
Under the assumption that $\tilde{E}$ is Fréchet differentiable at some $X \in L^2 \coloneqq L^2(\Omega; \R^d)$, $\nabla \tilde{E}(X)$ can be written as $\nabla \tilde{E}(X) = \eta \circ X$ for some $\eta \in L^2_{\law(X)}(\R^d; \R^d)$; moreover, for $X' \in L^2(\Omega; \R^d)$ such that $\law(X) = \law(X') \eqqcolon \nu$, also $\nabla \tilde{E}(X') = \eta \circ X'$. This $\eta$ is called the L-derivative of $E$ at $\nu$, and coincides with the Wasserstein gradient of $E$ at $\nu$ in the sense of Otto (for more details, see \cite{gangbo2019differentiability,jimenezDynamicalSystemsHamilton2023}). 
Thus, under this assumption of $\tilde{E}$ being Fréchet differentiable---called L-differentiability of $E$---, the probabilistic viewpoint and the analytic viewpoint are closely intertwined. 
In practice, however, L-differentiability is a strong property to demand, and the analytic framework of curves of maximal slope has been so successful precisely because it works in the nonsmooth regime.
We also recall that the problem of lifting the Wasserstein Hessian has been addressed in \cite{chassagneux2022probabilistic}.
Further bridges between the Lions lift and the nonsmooth analysis in $\Pp_2$ already exist: For example, the Bellman equation of McKean-Vlasov optimal control is treated in the viscosity sense by lifting it from $\Pp_2$ to the space of random variables \cite{PhamDynProgramming,MasterBellmanEquation2022}. 
Furthermore, in \cite{cavagnariLagrangianApproachTotally2026}, the authors have, for a generalization of gradient flows along $\lambda$-convex functionals, shown an equivalence between those flows in $\Pp_2$ and the mirrored ones taking place in the space of random variables.
We also mention the contribution \cite{cavagnari2022lagrangian} on the analysis (and the equivalence) of these different standpoints for multi-agent optimal control problems.
While recent years have seen a growing number of such bridges, the two viewpoints have largely been pursued in parallel: the maximal-slope theory is worked out and applied intrinsically on the Wasserstein space, while the extrinsic, Hilbertian picture on the space of random variables is seldom exploited for more specific applications.

This dichotomy is most unfortunate for two reasons. First, the notions of metric slope and derivative, upper gradients and curves of maximal slope---all from the analytic framework of metric gradient flows---mirror beautifully via the law map $\law$ between the $L^2$-space of random variables and the Wasserstein space ($\Pp_2, W_2$). Second, the roughness of Wasserstein spaces prevents us from applying off-shelf the theory of Newtonian motion in metric spaces developed in \cite{fornasierNewtonianMotionCompanion}; the Lions lift machinery holds promise for adressing this issue.

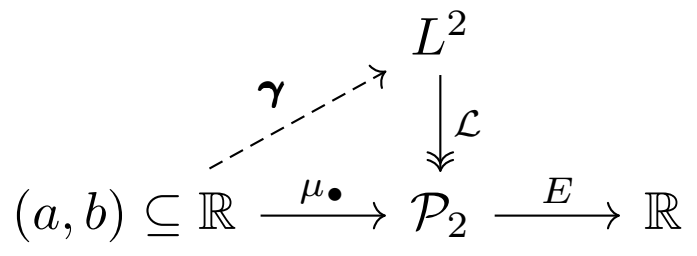
\begin{figure}[ht]
  \centering
\begin{tikzcd}
                                                                              & L^2 \arrow[d, "\law", two heads] &            \\
{(a, b) \subseteq \mathbb{R}} \arrow[ru, "\bgamma", dashed] \arrow[r, "\mu_{\bullet}"] & \Pp_2 \arrow[r, "E"]         & \mathbb{R}
\end{tikzcd}
\caption{Lifting curves of maximal slope from $\Pp_2$ to $L^2$.}\label{fig:lifting}
\end{figure}

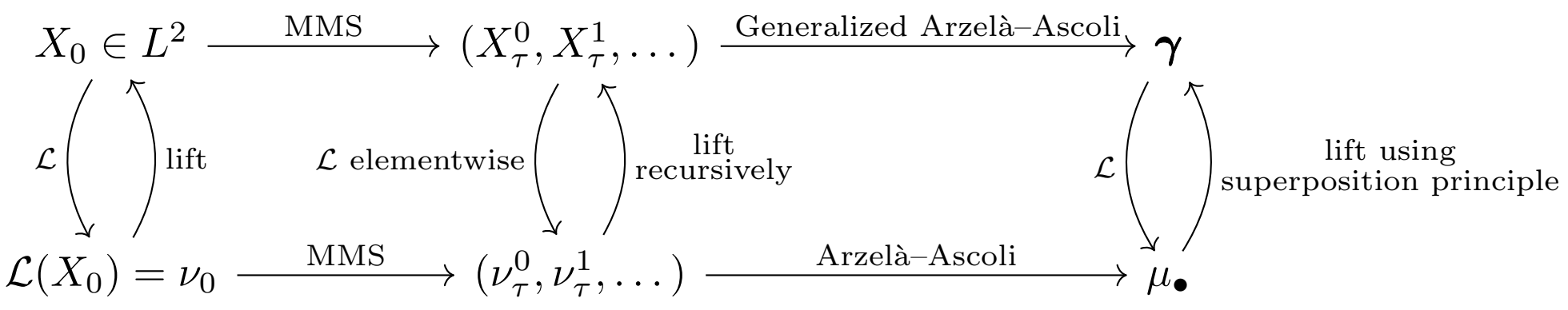
\begin{figure}[H]
  \centering
\begin{tikzcd}
X_0 \in L^2 \arrow[dd, "\law"', bend right] \arrow[rr, "\text{MMS}"]      &  & {(X_\tau^0, X_\tau^1, \dots)} \arrow[dd, "\law \text{ elementwise}"', bend right] \arrow[rrrr, "\text{Generalized Arzelà--Ascoli}"] &  && &\bgamma \arrow[dd, "\law"', bend right]                                         \\
                                                                                             &  &                                                                                                                                  &  &&&                                                                                   \\
\law(X_0) = \nu_0 \arrow[uu, "\text{lift}"', bend right] \arrow[rr, "\text{MMS}"] &  & {(\nu_\tau^0, \nu_\tau^1, \dots)} \arrow[uu, "\substack{\text{lift}\\\text{recursively}}"', bend right] \arrow[rrrr, "\text{Arzelà--Ascoli}"]          &&&  & \mu_\bullet \arrow[uu, "\substack{\text{lift using}\\\text{superposition principle}}"', bend right]
\end{tikzcd}
\caption{Completing the correspondence: Every part of the construction of curves of maximal slope as in \cite{ambrosioGradientFlows2008} can be carried out equivalently in $\Pp_2$ or $L^2$. For the recursive lift from $(\nu_\tau^0, \nu_\tau^1, \dots)$ to $(X_\tau^0, X_\tau^1, \dots)$, one needs $\Omega$ ``big enough'' to argue as in \cref{rmk:infimum_is_achieved}. Alternatively, one can use a countable joining argument (see \cite[Lemma 5.4.3]{ambrosioGradientFlows2008}) and lift jointly. In this case, one gets only the equality $\law(X_\tau^0) = \law(X_0)$, instead of the stronger equality $X_\tau^0 = X_0$.}\label{fig:completing_mirror}
\end{figure}

To make the first point more precise---we will expand on this in \cref{sec:lifting_gradient_problems}---, for any functional $E\colon \Pp_2 \to \R$, we have that $|\partial E|(\law(X)) = |\partial (E \circ \law)|(X)$ for any $X \in L^2$; furthermore, by the superposition principle (see \cite{lisiniCharacterizationAbsolutelyContinuous2007,ambrosioGradientFlows2008}), every Wasserstein-$\mathrm{AC}^2$ curve $\mu_\bullet\colon [0, T] \to \Pp_2$ can be lifted to its \emph{horizontal lift}, i.e., an absolutely continuous curve $\bgamma\colon [0, T] \to L^2$ so that $\law \circ \bgamma = \mu_\bullet$ and such that $|\mu'|(t) = |\bgamma'|(t)$ for almost all $t \in [0, T]$. Even further, it is easy to show that any functional $g\colon \Pp_2 \to [0, +\infty]$ is an upper gradient for the functional $E$ if and only if $g \circ \law$ is an upper gradient for the Lions lift $\tilde{E}$. Lastly, we obtain the correspondence, or mirroring
\begin{equation}\label{eq:correspondence_coms}
\left\{\begin{array}{c}\text{Horizontal lifts of curves of }\\
\text{ maximal slope of } E \text{ in } \Pp_2\end{array}\right\}
\leftrightarrow
\left\{\begin{array}{c}\text{Curves of maximal slope}\\
\text{of } E \circ \law \text{ in } L^2\end{array}\right\}.
\end{equation}
In \cref{sec:topology-problem}, we make rigorous and explore an unsettling crack in this mirror. Namely, while we can \emph{lift} any curve of maximal slope in $\Pp_2$ to a curve of maximal slope in $L^2$, it seems difficult to use the algorithmic machinery to \emph{construct} those curves in $L^2$ directly: The Arzelà--Ascoli theorem---a key part of the aforementioned machinery---is not available in $L^2$, due to topological obstructions, namely the lack of a topology on $L^2$ in which bounded sets are precompact, the distance is lower semicontinuous, and $\law$ is narrowly continuous. We address this issue by uncovering a weaker property that we call \emph{A-A-readiness}, which suffices for an Arzelà--Ascoli-like argument; finally, it turns out that $L^2(\Omega; \R^d)$ is A-A-ready for $\Omega$ ``large enough''. For such $\Omega$, we can furthermore employ the minimizing movement scheme in $L^2$ instead of in $\Pp_2$, completing the correspondence and making every step in the analytic framework of curves of maximal slope available equivalently in $\Pp_2$ and $L^2$ (see \cref{fig:completing_mirror}). \\
As the analytic framework of curves of maximal slope in $\Pp_2$ is employed mostly with PDEs in mind, we spend \cref{sec:PDE_gradient_flows} chasing a description of PDEs solved by curves of maximal slope in $L^2$ driven by the Lions lift of a functional. We emphasize that we carry out this analysis without explicitly using the tangent structure of $\Pp_2$; instead, leveraging the Hilbert space structure of $L^2$, we provide a very weak description of those PDEs. Owing to the correspondence \eqref{eq:correspondence_coms}, this result yields a weaker characterization of solutions of PDEs fulfilled by curves of maximal slope in $\Pp_2$, which is, to the best of our knowledge, novel. \\
In \cref{sec:examples}, we present two examples showing the flexibility gained by being able to run all arguments in either $\Pp_2$ or in $L^2$, seamlessly switching back and forth. The first example to challenge the Lions correspondence is rather pedagogical, showing that we can use the weak characterization of PDEs even for energies defined through densities, for instance an entropy functional, which do not admit a direct expression in terms of the Lions lift. In the second example, we leverage the flexibility to derive, in a heretofore untackled regime, the existence of curves of maximal slope and their corresponding PDE for a nonlocal attraction-repulsion potential with singular interactions. This functional is known not to be $\lambda$-convex (see \cite{HertrichGrafBeinertSteidl2024}), so that the results of~\cite{ambrosioGradientFlows2008} do not apply off-the-shelf. \\
Finally, \cref{sec:newtonian_curves} concerns the quest to go beyond first-order evolutions driven by functionals on Wasserstein space. There, we build on the Newtonian framework introduced by us in~\cite{fornasierNewtonianMotionCompanion} whose pillars are analogous to the analytic program of curves of maximal slope---$\mathrm{AC}$-curves, the minimizing movement scheme and compactness thereof---, but whose direct application is obstructed by the ``excessive roughness'' of $\Pp_2$. Here, we leverage the Lions lift. Combining this Newtonian framework with the novel compactness arguments in $L^2$ we developed in \cref{sec:a-a-w-o-lsc,sec:a-a-readiness-of-l2}, we manage to generate curves of Newtonian motion in $L^2$, driven by the Lions lift of functionals on $\Pp_2$. Furthermore, the Newtonian framework comes with a variational characterization; we show that this variational characterization can be intrinsically checked on the Wasserstein phase space $\ful{\Pp}_2 \coloneqq \Pp_2(\R^d \times \R^d)$, thus proposing an intrinsic notion of \emph{Wasserstein-Newtonian curves}, which in turn yields another correspondence, or mirroring
  \begin{equation*}
\left\{\begin{array}{c}\text{Horizontal lifts of Wasserstein-}\\
\text{ Newtonian curves in } \ful{\Pp}_2\end{array}\right\} 
  \leftrightarrow
  \left\{\begin{array}{c}\text{Curves of Newtonian motion}\\
  \text{of } E \circ \law \text{ in } L^2\end{array}\right\}.
  \end{equation*}
Again, at this point, we face an unfortunate crack in the mirror: While we are able to \emph{characterize} Wasserstein-Newtonian curves in $\ful{\Pp}_2$ intrinsically, to construct the MMS-like approximating sequence, we are forced to lift our data and functional to $L^2$. To sidestep this detour, we propose an algorithm intrinsic to $\ful{\Pp}_2$, providing the mirroring of the MMS-like approximating sequence. 
\begin{figure}[H]
  \centering
\begin{tikzcd}
(X_0, V_0) \in L^2 \arrow[dd, "\ful{\law}"', bend right] \arrow[rr, "\text{MMS-like}"]      &  & {((X_{\tau, 0}, V_{\tau, 0}), (X_{\tau, 1}, V_{\tau, 1}), \dots)} \arrow[dd, "\ful{\law} \text{ elementwise}"', bend right] \arrow[rr, "\substack{\text{Compactness}\\\text{argument}}"] &  &(\bgamma, \dot{\bgamma}) \arrow[dd, "\ful{\law}"', bend right]                                         \\
                                                                                             &  &                                                                                                                                  &  &                                                                                   \\
\law((X_0, V_0)) = \ful{\nu}_0 \arrow[uu, "\text{lift}"', bend right] \arrow[rr, "\text{\cref{alg:wasserstein-newtonian-seq}}"] &  & {(\ful{\nu}_{\tau, 0}, \ful{\nu}_{\tau, 1}, \dots)} \arrow[uu, "\substack{\text{lift}\\\text{using}\\\text{Thm. \labelcref{thm:algorithm_in_p_2}}}"', bend right]    \arrow[rr, dashed]     &  & \ful{\mu}_\bullet \arrow[uu, "\substack{\text{lift using}\\\text{ a second-order}\\\text{superposition}\\\text{principle}}"{swap, yshift=6mm}, bend right]
\end{tikzcd}
\caption{The correspondence between $\Pp_2$ and $L^2$ for the Newtonian case. While we do not explicitly construct the dashed arrow, the convergences in the compactness argument in $L^2$ pass down to the laws, yielding the convergence below.}\label{fig:completing_mirror_newtonian}
\end{figure}
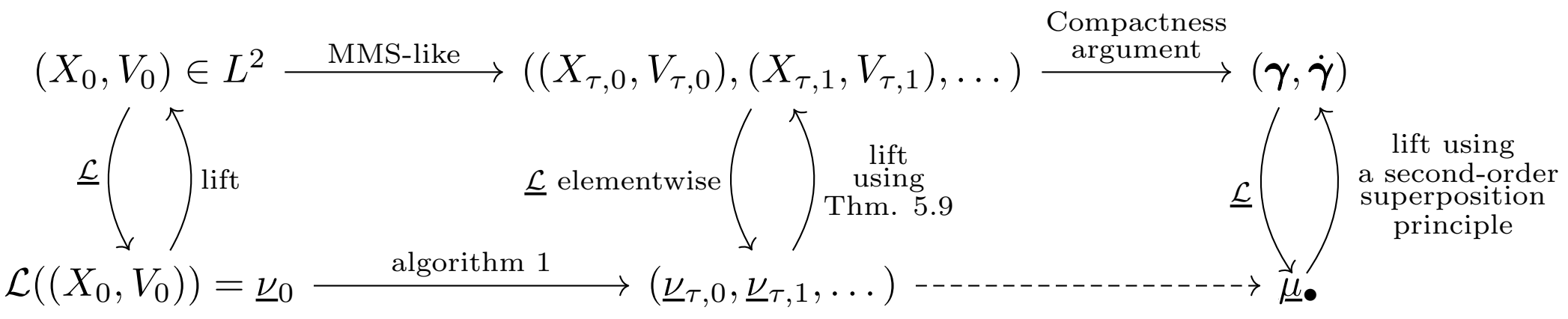

The work done in \cref{sec:newtonian_curves} relates to two prior strains of research: On one hand, in~\cite{ambrosio2008hamiltonian}, the authors considered an explicit discrete-time scheme for Hamiltonian systems---of which Newtonian dynamics are a special case---on the Wasserstein phase space; they showed that under mild conditions on the Hamiltonian, the scheme converges to a continuous curve in $\ful{\Pp}_2$ fulfilling a Vlasov equation. 
On the other hand, Brigati et al. analyzed in \cite{brigatiKineticOptimalTransport2025} those curves in $\ful{\Pp}_2$ (called \emph{physical curves} therein) that form the candidate class for Wasserstein-Newtonian curves, analogously to the class $W^{2,2}$ in a Hilbert space, and its metric generalization \ACt (considered in \cite{fornasierNewtonianMotionCompanion}). 
In contrast to \cite{ambrosio2008hamiltonian}, our scheme is implicit, and we show as our major contribution the convergence to a curve fulfilling the intergal second-order variational characterization (see \cref{rmk:recall_curve _Newtonian} and \cref{eq:def_variational_Newtonian}) and fulfilling the Valsov equation a posteriori (as in \cref{prop:abstract_pde_newtonian}, consistently with \cite{ambrosio2008hamiltonian}). 
One of the tools employed in this variational characterization is precisely the framework of physical curves introduced in \cite{brigatiKineticOptimalTransport2025}.

\paragraph{Methodological Note.} Other more established literature, such as \cite{ambrosioGradientFlows2008} and  \cite{carmonaProbabilisticTheoryMean2018}, addresses evolutions by adopting strictly separate standpoints, which is the purely measure theoretical one in \cite{ambrosioGradientFlows2008} and the probabilistic point of view in \cite{carmonaProbabilisticTheoryMean2018}.
In our case, in order to connect the perspective of metric gradient flows in $\Pp_2$ and that of the Lions lift on $L^2_{\mathbb R}$, we must leverage both fronts at the same time. While this juggling requires careful handling, it is precisely what is needed to overcome the dichotomy mentioned earlier.

\subsection{Contributions}
In \cref{sec:lifting_gradient_problems}, we settle some results concerning the Lions lift, part of which have been fruitfully employed in mean field games (see, e.g., \cite{chassagneux2022probabilistic}) and optimal control (see, e.g., \cite{averboukh2025pontryagin,jimenezDynamicalSystemsHamilton2023} and the recent works \cite{capuani2027generalized,bonnet2026exponential}). 
The core contributions of the present paper are the following:
\begin{enumerate}
  \item In \cref{sec:topology-problem}, we detail the topological obstruction that prevents the direct use of classical Arzelà--Ascoli arguments for constructing curves of maximal slope in $L^2$ driven by the Lions lift of narrowly lower semicontinuous functionals on $\Pp_2$.
  \item In \cref{sec:a-a-w-o-lsc}, we amend this issue by distilling an essential property that allows the usage of a relaxed Arzelà--Ascoli argument. Furthermore, this generalized version can be used as a drop-in replacement for constructing curves of maximal slope. We call the property \emph{A-A-readiness}.
  \item In \cref{sec:a-a-readiness-of-l2}, we prove that for the probability space
  \begin{equation*}
    (\Omega_\R, \mathcal{F}_\R, \mathbb{P}_\R) \coloneqq \Big([0,1]^\R, \bigotimes_\R \mathcal{B}\big([0, 1]\big), \bigotimes_\R \Leb|_{[0, 1]} \Big),
  \end{equation*}
  the space $\LtwoR \coloneqq L^2(\Omega_\R; \R^d)$ is A-A-ready.
  \item In \cref{sec:PDE_gradient_flows}, we apply a geometric argument in the proof of \cref{lem:technical-coms}, where we show that, given a functional $F\colon \mathcal{H} \to \R$ defined over a Hilbert space $\mathcal{H}$ and a curve of maximal slope $u\colon [0, T] \to \mathcal{H}$, we have
  \begin{equation*}
      \langle \dot{u}(t), v \rangle_{\mathcal{H}} = - \frac{\diff}{\diff \tau} F(u(t) + \tau v) \Big|_{\tau = 0}
    \end{equation*}
    for almost every $t \in [0, T]$ and for any direction $v \in \mathcal{H}$ such that the right hand side exists.
    Applying this proposition to the Hilbert space $\mathcal{H} = L^2$, we obtain that a curve of maximal slope $\bgamma$ in $L^2$ for the Lions lift of a functional $E$ fulfills the weak PDE
    \begin{equation*}
      \frac{\diff}{\diff t} \int \varphi_t(x) \,\diff \mu_t(x) = \int \partial_t \varphi_t(x) \,\diff \mu_t(x) - \frac{\diff}{\diff s} E\big((T_s^{\varphi_t})_\# \mu_t\big) \Big|_{s = 0},
    \end{equation*}
    where $\mu_t = \law(\bgamma_t)$ and $T_s^{\varphi_t}(x) = x + s \nabla_x \varphi_t(x)$. By virtue of the correspondence \cref{eq:correspondence_coms}, we obtain that the same equation is fulfilled for any curve of maximal slope $\mu_\bullet$ for the functional $E$. This unlocks a straightforward but rigorous characterization of curves of maximal slope in $\Pp_2$, moving beyond the usual $\lambda$-convexity assumption for the driving functional.
  \item In \cref{sec:examples}, we apply the theory from \cref{sec:wasserstein_as_quotient,sec:PDE_gradient_flows} to two exemplary energies on $\Pp_2$. The first example (\cref{sec:ex_energy_through_density}) is not aimed at providing novel existence proofs, but rather to show that the framework of \cref{sec:wasserstein_as_quotient,sec:PDE_gradient_flows} works even in the setting of energies defined through densities, for instance an entropy functional, for which the Lions lift does not admit a direct expression. Conversely, for the second example (\cref{sec:ex_nonlocal}), we exploit the correspondence with the Lions lift to show the existence of curves of maximal slope for the nonlocal attraction-repulsion transport potential
  \begin{equation*}
    E\colon \Pp_2 \to \R, \qquad  \mu \mapsto -\frac{1}{2} \iint K_s(x-y) \,\diff (\mu-\nu)(x) \, \diff (\mu-\nu)(y)
  \end{equation*}
  with power kernels $K_s(z) = |z|^{s+1}$ and $s \in [0, 1)$, for a sufficiently regular reference signed measure $\nu$; gradient flows of such a potential have previously been the object of study in \cite{RosenzweigSlepcevWang2026} and \cite{fornasierNonlocalAttractionrepulsionTransport2026}. We furthermore use the tools from \cref{sec:PDE_gradient_flows} to prove that such curves of maximal slope actually solve in the weak sense the heuristically expected PDE.
  \item In \cref{sec:characterizing_wasserstein_newtonian_curves}, we descend the variational characterization of curves of Newtonian motion in $L^2$ to the phase space $\ful{\Pp}_2$, describing \emph{Wasserstein-Newtonian curves}, i.e., those trajectories in $\ful{\Pp}_2$ which arise as the image of curves of Newtonian motion in $L^2$. 
  One of the ingredients of this variational characterization is the notion of second metric derivative and a metric generalization of $W^{2,2}$-curves. 
  In \cite{fornasierNewtonianMotionCompanion}, we construct such a generalization for metric spaces with extendable geodesics, called \ACt-curves. 
  In $\ful{\Pp}_2$, it turns out that the avatar of $L^2$-valued $W^{2,2}$-curves is given by the concept of ``physical curves'' from \cite{brigatiKineticOptimalTransport2025}; we recall this notion in \cref{sec:physical_connection}. 
   By introducing as well avatars of the other ingredients of the variational characterization of curves of Newtonian motion, namely tangent curves and upper gradient-Hessian pairs, we obtain a characterization of Wasserstein-Newtonian curves which is fully intrinsic to $\ful{\Pp}_2$.
  \item In \cref{sec:constructing_newtonian_curves_in_l2}, we connect the Newtonian framework from~\cite{fornasierNewtonianMotionCompanion} with the compactness arguments developed in \cref{sec:a-a-readiness-of-l2} to construct curves of Newtonian motion in $\LtwoR$ for the Lions lift of functionals on $\Pp_2$. We furthermore provide a construction of the laws of such curves completely within $\ful{\Pp}_2$.
  \item In \cref{sec:identifying_newtonian}, we show how under certain regularity assumptions, the law of curves of Newtonian motion in $\LtwoR$ fulfill a Vlasov kinetic equation.
  \item In \cref{sec:example_newtonian}, we apply our theory to construct curves of Newtonian motion along regular interaction kernels.
  \item In Appendix~\ref{sec:generalized_newtonian}, we use ad-hoc arguments to relax some of the regularity assumptions employed in the example of \cref{sec:example_newtonian}.
  \item In Appendix~\ref{sec:proof-of-frechet-derivative-l2-w2} and~\ref{sec:proof-of-lipschitz-cont}, we collect proofs for some technical lemmas.
\end{enumerate}

\newpage
\subsection*{General Notations}
In the table below, we collect a few notations that we will use throughout the manuscript. 

\begingroup
\setlength{\LTpre}{4pt}
\setlength{\LTpost}{4pt}
\renewcommand{\arraystretch}{1.1}

\begin{longtable}{@{} p{0.30\textwidth} p{0.65\textwidth} @{}}
\toprule
\textbf{Symbol} & \textbf{Description} \\
\midrule
\endfirsthead
\multicolumn{2}{l}{\emph{(continued)}} \\
\toprule
\textbf{Symbol} & \textbf{Description} \\
\midrule
\endhead
\midrule
\multicolumn{2}{r}{\emph{(continued on next page)}} \\
\endfoot
\bottomrule
\endlastfoot

\multicolumn{2}{@{}l}{\textit{Spaces and sets}} \\[1pt]
$\mathbb{R}^d$; $\mathbb{R}^{2d}$
  & State space; phase space $\mathbb{R}^d \times \mathbb{R}^d$ with variables $(x, v)$. \\
$\Leb^d$;  $\Leb|_{[0,1]}$
  & Lebesgue measure on $\R^d$; Lebesgue measure on $[0,1]$. \\
$(\Omega, \mathcal{F}, \mathbb{P})$
  & Generic base probability space; unless stated otherwise, we assume it to admit a random variable $Z$ such that $\law(Z) = \Leb|_{[0, 1]}$. \\
$L^2$
  & $L^2(\Omega; \R^d)$, random variables with finite second moment. \\
$(\Omega_\R, \mathcal{F}_\R, \mathbb{P}_\R)$
  & Specific choice of a probability space where $\Omega_\R = [0, 1]^{\R}$, $\mathcal{F}_\R = \bigotimes_\R \mathcal{B}([0, 1])$ and $\mathbb{P}_\R = \bigotimes_\R \Leb|_{[0,1]}$. \\
$\LtwoR$
  & $L^2(\Omega_\R; \R^d)$; A-A-ready by \cref{cor:l2_is_aa}. \\
$\ful{L}^2$, $\LtwoRful$
  & Phase-space of random variables $L^2(\Omega; \R^{2d}) \cong L^2 \times L^2$ and $L^2(\Omega_\R; \R^{2d}) \cong \LtwoR \times \LtwoR$, respectively. \\
$\Pp_2$
  & $\Pp_2(\R^d)$, probability measures with finite second moment,
    with the Wasserstein distance $W_2$. \\
$\ful{\Pp}_2$
  & Wasserstein phase space $\Pp_2(\R^d \times \R^d)$. \\
$\Pi(\mu, \nu)$, $\{\pi_x\}_{x}$
  & Couplings of $\mu$ and $\nu$; disintegration of $\pi$ with respect to the first marginal. \\
$\mathcal{M}$, $\ful{\mathcal{M}}$
  & Finite signed Radon measures on $\R^d$, resp.\ $\R^{2d}$. \\
$C_0^2$, $\ful{C}_0^2$
  & $C^2$-functions on $\R^d$, resp.\ $\R^{2d}$, vanishing at infinity together with their first and
    second derivatives; $(C_0^2)'$, $(\ful{C}_0^2)'$ are the duals;
    see \cref{fig:all_spaces,fig:all_spaces_phase_space}. \\
$\mathrm{AC}^2\big((a,b); \mathcal{S}\big)$
  & Absolutely continuous curves in a metric space $\mathcal{S}$ with square-integrable metric
    derivative; see \cite[Def.~1.1.1]{ambrosioGradientFlows2008}. \\
\ACt{}
  & Second-order absolutely continuous curves; see \cite[\cref*{M-def:AC2_curve}]{fornasierNewtonianMotionCompanion}. \\[4pt]

\multicolumn{2}{@{}l}{\textit{Lifts, projections, and related maps}} \\[1pt]
$\law\colon L^2 \to \Pp_2$
  & Law map, $\law\colon X \mapsto X_\# \mathbb{P}$; surjective and $1$-Lipschitz (\cref{sec:wasserstein_as_quotient}). \\
$\ful{\law}$
  & Phase-space law map $\ful{L}^2 \to \ful{\Pp}_2$. \\
$\tilde{E} \coloneqq E \circ \law$
  & Lions lift of a functional $E\colon \Pp_2 \to \Rext$. \\
$\bgamma$
  & Horizontal lift of a curve $\mu_\bullet$ in $\Pp_2$: $\law \circ \bgamma = \mu_\bullet$ and
    $|\bgamma'| = |\mu'|$ a.e.; see \cref{sec:lifting_gradient_problems}.\\
$\eta_\nu$
  & L-derivative of $E$ at $\nu \in \Pp_2$, so that $\nabla \tilde{E}(X) = \eta_{\law(X)} \circ X$. \\
$G$, $\ful{G}$
  & $G \coloneqq r \circ b \circ \law\colon L^2 \to (C_0^2)'$; $\ful{G}$ is its phase-space analogue. \\
$p_x$, $\ful{\mu}^x$
  & Projection $(x, v) \mapsto x$; state marginal $\ful{\mu}^x \coloneqq (p_x)_\# \ful{\mu}$. \\

\multicolumn{2}{@{}l}{\textit{Topologies and convergence}} \\[1pt]
$\sigma_{\mathrm{narrow}}$, $\sigma_w$
  & Narrow topology on $\Pp_2$; weak topology on $L^2$. \\
$\law^{-1}(\sigma_{\mathrm{narrow}})$, $\law^{-1}(\sigma_{W_2})$
  & Pullbacks along $\law$ of the narrow, resp.\ $W_2$-topology;
    $\ful{\law}^{-1}(\sigma_{\mathrm{narrow}})$ is defined analogously. \\
A-A-ready
  & Relaxed Arzelà--Ascoli property; see \cref{def:a-a-readiness,prop:aa-ready-allows-aa}. \\
$\xrightarrow{Y}$
  & Young convergence; see \cref{thm:young_compactness_from_rindler,prop:countable_young}. \\
$\mathrm{lsc}(F)$, $\mathrm{lsc}_\sigma(F)$
  & Lower semicontinuous envelope of $F$ (with respect to the topology $\sigma$). \\[4pt]

\multicolumn{2}{@{}l}{\textit{Curves and their derivatives}} \\[1pt]
$\mu_\bullet$, $X_\bullet$, $\ful{\mu}_\bullet$
  & Curves in $\Pp_2$, $L^2$, resp.\ $\ful{\Pp}_2$, with values $\mu_t$, $X_t$, $\ful{\mu}_t$. \\
$|\mu'|(t)$
  & Metric derivative:
    $\limsup_{s \to t} d(\mu_s, \mu_t)/|s - t|$. \\
$\dot{u}$, $\ddot{u}$
  & First and second derivatives of curves in a Hilbert space. \\
$\ful{d}_T$
  & Kinetic transport cost of \cite{brigatiKineticOptimalTransport2025};
    see \cref{sec:characterizing_wasserstein_newtonian_curves}. \\
Physical curve
  & Curve $\ful{\mu}_\bullet$ in $\ful{\Pp}_2$ with
    $\ful{d}_{t-s}(\ful{\mu}_s, \ful{\mu}_t) \leq \int_s^t m(r)\,\diff r$, $m \in L^2$;
    see \cref{def:horizontal_lift_physical} (from \cite[Def.~1.3]{brigatiKineticOptimalTransport2025}) and \cref{cor:lifting_physical_curves}. \\
$|\ful{\mu}''|(t)$
  & $\limsup_{h \downarrow 0} \ful{d}_h(\ful{\mu}_t, \ful{\mu}_{t+h})/h$;
    see \cref{def:horizontal_lift_physical}. \\
$V_t^{\bgamma}$, $V_t^{\ful{\mu}}$
  & Tangent geodesic to $\bgamma$ in $L^2$; its image
    $V_t^{\ful{\mu}}(\sigma) = [(x, v) \mapsto x + (\sigma - t) v]_\# \ful{\mu}_t$ in $\Pp_2$;
    see \cref{sec:wasserstein-newtonian-curves}. \\[4pt]

\multicolumn{2}{@{}l}{\textit{Energies, slopes, and upper gradient-Hessian pairs}} \\[1pt]
$E$
  & Energy functional $E\colon \Pp_2 \to \Rext$. \\
$|\partial E|(x)$
  & Metric slope: $\limsup_{y \to x} (E(x) - E(y))^+/d(x, y)$. \\
$|\partial^- E|(x)$
  & Relaxed slope: $\inf_{y_i \to x} \liminf_{i} |\partial E|(y_i)$. \\
$J_{\tau, \mathcal{S}}[x]$
  & Proximal map in a generic metric space $(\mathcal{S}, d_{\mathcal{S}})$: $\argmin_{y} \big\{E(y) + \tfrac{1}{2\tau} d_{\mathcal{S}}^2(x, y)\big\}$;
    $J_{\tau, \LtwoR}$ and $J_{\tau, \Pp_2}$ as in \cref{rmk:mms_downstairs_upstairs}. \\
$\partial^+ F$, $D^+F$, $\partial^0 F$
  & Superdifferential; Fréchet superdifferential;
    Fréchet subdifferential of minimal norm. \\
$H_0^{\tilde{E}}(\xi)$
  & Limit Hessian of $\tilde{E}$ along a geodesic $\xi$ in $L^2$; see \cref{eq:hessian_0}. \\
$H_0^E(\ful{\nu})$
  & $\frac{\diff^2}{\diff \sigma^2} E\big([(x, v) \mapsto x + \sigma v]_\# \ful{\nu}\big)\big|_{\sigma = 0}$
    for $\ful{\nu} \in \ful{\Pp}_2$; see \cref{sec:wasserstein-newtonian-curves}. \\
$(g, H)$
  & Upper gradient-Hessian pair; see \cref{def:upg_wasserstein}, and
    \cite[\cref*{M-def:upper_gradient}]{fornasierNewtonianMotionCompanion} for the metric setting. \\
Wasserstein-Newtonian curve
  & See \cref{def:wasserstein_newtonian_curve}. \\[4pt]

\multicolumn{2}{@{}l}{\textit{Example energies}} \\[1pt]
$F(x, \rho)$
  & Integrand of an energy defined on densities; see \cref{eq:formula_internal_energy}. \\
$K_s$
  & Power kernel $K_s(z) = |z|^{s+1}$, $s \in [0, 1)$; see \cref{eq:definition_of_potential}. \\
$\nu = f \Leb^d$
  & Reference signed measure; see \cref{eq:nu}. \\
$E_1$, $E_2$
  & Repulsive and attractive parts of the energy \eqref{eq:definition_of_potential}. \\
$V$
  & Background potential $K_s * \nu$ in \cref{sec:ex_nonlocal};
    external potential in \cref{sec:example_newtonian}. \\
$K$
  & Symmetric interaction kernel; see \cref{eq:newtonian_potential_definition,ass:kernel_and_potential}. \\[4pt]

\multicolumn{2}{@{}l}{\textit{Discrete schemes}} \\[1pt]
$\tau$
  & Time step size. \\
$U_\tau^n$
  & $n$-th iterate of the minimizing movement scheme; see \cref{eq:mms}. \\
$\bar{U}_\tau$, $\tilde{U}_\tau$
  & Piecewise constant and De Giorgi interpolations; see \cref{eq:definition_de_giorgi}. \\
$|U'_\tau|$
  & Discrete metric derivative $d(U_\tau^i, U_\tau^{i-1})/\tau$ on $(\tau(i-1), \tau i]$. \\
$U_{\tau, k}$
  & $k$-th iterate of a Newtonian sequence; see \cref{eq:newtonian_sequence}. \\
$U'_{\tau, k}$
  & Discrete velocity $(U_{\tau, k} - U_{\tau, k-1})/\tau$;
    see \cref{eq:L2_newtonian_seq_velocities}. \\
$Z_{\tau, k}$
  & Predictor point $U_{\tau, k} + \tau U'_{\tau, k}$; see \cref{eq:wasserstein_uprime}. \\
$u_\tau$, $\tilde{u}_\tau$, $\dot{u}_\tau$
  & Piecewise constant, De Giorgi, and velocity interpolations;
    see \cref{eq:de_giorgi_interpolant}. \\
$\mathcal{V}^{\tau, \varepsilon}$, $\tilde{\mathcal{V}}^{\tau, \varepsilon}$, $\bar{\mathcal{V}}^{\tau, \varepsilon}$
  & Approximate tangents; see \cref{eq:tilde_v_wasserstein,eq:bar_v_wasserstein}. \\
$|u_\tau^{(1)}|$, $|u_\tau^{(2)}|$
  & Discrete first- and second-order metric derivatives;
    see \cref{eq:second_derivative_metric}. \\
$|G_\tau|$
  & Discrete slope of the De Giorgi interpolation; see \cref{eq:G_tau}. \\
$\ful{M}_{\tau, i}$, $\nu_i$
  & Iterates and intermediate measures of \cref{alg:wasserstein-newtonian-seq}. \\[4pt]

\multicolumn{2}{@{}l}{\textit{Miscellaneous}} \\[1pt]
$T_\# \mu$
  & Pushforward of $\mu$ under a measurable map $T$. \\
$\delta_x$, $\chi_A$
  & Dirac measure at $x$; indicator function of $A$. \\
$\diff$
  & Differential operator ($\diff t$, $\diff x$, etc.). \\
$\mathbb{N},\, \mathbb{Q},\, \R,\, \R_+,\, \Rext$
  & Usual sets of numbers,
    $\R_+ \coloneqq [0, +\infty)$. \\
\end{longtable}
\endgroup

\section{The Wasserstein Space as a Quotient of \texorpdfstring{$L^2$}{L²}} \label{sec:wasserstein_as_quotient}

Let us start by showing how the Wasserstein space $\big(\Pp_2(\R^d), W_2 \big)$ can be constructed as a quotient space of $L^2(\Omega; \R^d)$ for a suitable $\Omega$. 
Consider a probability space $(\Omega, \mathcal{F}, \mathbb{P})$ which admits a random variable $Z$ such that $\law(Z) = \Leb|_{[0, 1]}$, e.g.,\ $(\Omega, \mathcal{F}, \mathbb{P}) = \big([0, 1], \mathcal{B}([0, 1]), \Leb|_{[0, 1]}\big)$ (the unit interval with the Borel $\sigma$-algebra and the Lebesgue measure), and set $L^2 \coloneqq L^2(\Omega; \R^d)$. 
The map $\law\colon L^2 \to \Pp_2(\R^d)$ defined as $ X \mapsto \law(X) \coloneqq  X_\# \mathbb{P}$ mapping each random variable to its law is surjective\footnote{The fact that we can realize any probability measure as the distribution of a random variable is, e.g., given in \cite[Proposition 10.7.6]{bogachevMeasureTheory2007}, for the case $\Omega = [0, 1]$ equipped with the Lebesgue measure (in that Proposition, take $X = \{\star\}$). Precomposing with $Z$ gives the result. The result in \cite[Proposition 10.7.6]{bogachevMeasureTheory2007} works for any codomain which is a standard Borel space.\label{footnote:realize_prop_measure}}, as we may realize each probability measure with finite second moment as the law of a suitable random variable on $(\Omega, \mathcal{F}, \mathbb{P})$. 
Thus, as sets, we have the equality $\Pp_2(\R^d) = L^2/\sim$, where the equivalence relation $\sim$ is defined by $X \sim Y$ if and only if $\law(X) = \law(Y)$. What is left to show is that the metric structure of $\big(\Pp_2(\R^d), W_2 \big)$ coincides with the quotient metric induced by $L^2$ via the map $\law$, i.e.,\ that we have
\begin{equation} \label{eq:basic_inf_result}
  W_2(\mu, \nu) = \inf\left\{\|X - Y\|_{L^2(\Omega; \R^d)} : \law(X) = \mu, \law(Y) = \nu\right\}.
\end{equation}
We recall that the Wasserstein distance $W_2$ is defined as
\begin{equation*}
    W_2(\mu, \nu) \coloneqq 
    \inf \left\{ 
    \sqrt{\int_{\R^d\times \R^d} |x-y|^2 \, \diff \pi(x,y)} : \pi( \cdot \times \R^d) = \mu, \ \pi(\R^d\times \cdot)=\nu
    \right\}
\end{equation*}
for every $\mu, \nu \in \Pp_2(\R^d)$.
To show \cref{eq:basic_inf_result}, we notice that the inequality ``$\leq$'' is immediate, as for any $X, Y$ such that $\law(X) = \mu$ and $\law(Y) = \nu$, the joint law $\pi = \law(X, Y)$ is a coupling of $\mu$ and $\nu$. For the other direction, consider an optimal coupling $\pi$ of $\mu$ and $\nu$ and let $(X, Y)$ be a random variable on $(\Omega, \mathcal{F}, \mathbb{P})$ with law $\pi$. Then, by construction, we have $\law(X) = \mu$ and $\law(Y) = \nu$, and hence the optimality of $\pi$ yields the desired result.
When there is no possibility of ambiguity, we write $\Pp_2$ in place of $\Pp_2(\R^d)$.

\subsection{Lifting Gradient Problems from \texorpdfstring{$\Pp_2$}{P₂} to \texorpdfstring{$L^2$}{L²}} \label{sec:lifting_gradient_problems}
From the superposition principle (see, e.g., \cite{lisiniCharacterizationAbsolutelyContinuous2007} or~\cite[Theorem~8.2.1 and Theorem~8.3.1]{ambrosioGradientFlows2008}), we observe immediately\footnote{This simple consequence of the superposition principle has been observed before, e.g. in \cite[Theorem 7.1]{bertucciTangentSpaceWasserstein2025} or \cite[Theorem 2.1 (i)]{jimenezDynamicalSystemsHamilton2023}.} that we can lift absolutely continuous curves from $\Pp_2$ to $L^2$. 
Indeed, let $\mu_{\bullet}\colon  [0, T] \to \Pp_2$ be an absolutely continuous curve. Then, by \cite[Lemma 1.1.4]{ambrosioGradientFlows2008}, we may choose $r \colon [0, L] \to [0, T]$ Lipschitz and strictly monotone increasing with absolutely continuous inverse $q \coloneqq r^{-1}$ such that $\bar{\mu}_{\bullet} \coloneqq \mu_{\bullet} \circ r$ is Lipschitz and in particular $\mathrm{AC}^2$. Then, by the superposition principle, there exists a probability measure $\boldsymbol \eta$ on $AC^2([0, L]; \R^d)$ such that $(e_t)_\# \boldsymbol \eta = \bar{\mu}_t$ for all $t \in [0, L]$ and such that
\begin{equation*}
  \int_0^L | \bar{\mu}' |^2(t) \, \diff t = \iint_0^L |\gamma'_t|^2 \, \diff t \, \diff \boldsymbol{\eta}(\gamma),
\end{equation*}
where we emphasize that, on the left hand side, we integrate the metric derivative in the sense of~\cite{ambrosioGradientFlows2008}, while on the right hand side, we integrate the norm of the classical derivative---which, in the case of absolutely continuous curves in a Hilbert space, coincides with the metric derivative. \\
Setting $X$ to be a measurable function from $\Omega$ to $AC^2([0, L]; \R^d)$ such that $\law(X) = \boldsymbol{\eta}$---which we can do, as $AC^2([0, L]; \R^d)$ is a standard Borel space\footref{footnote:realize_prop_measure}---, we can define the curve $\bar{\bgamma}\colon [0, L] \to L^2$ as $t\mapsto \bar{\bgamma}_t = X(\cdot)(t)$. 
By construction, we have that $\law(\bar{\bgamma}_t) = \bar{\mu}_t$ for all $t \in [0, L]$. Next, for every $s, t \in [0, L]$ such that $s \leq t$, we calculate:
\begin{align*}
  \|\bar{\bgamma}_t - \bar{\bgamma}_s\|_{L^2} &= \left(\int |X(\omega)(t) - X(\omega)(s)|^2\, \diff \mathbb{P}(\omega)\right)^{\frac{1}{2}} \leq \left(\int \left|\int_s^t |X(\omega)'(r)|\,\diff r \right|^2\, \diff \mathbb{P}(\omega)\right)^{\frac{1}{2}} \\
  &\overset{(*)}{\leq} \int_s^t \left( \int |X(\omega)'(r)|^2\, \diff \mathbb{P}(\omega) \right)^{\frac{1}{2}} \, \diff r,
\end{align*}
where the application of Minkowski's integral inequality in $(*)$ requires $(\omega, t) \mapsto |X(\omega)'(t)|$ to be jointly measurable in $(\omega, t)$; we get this by the fact that $X$ is jointly measurable---it is measurable in $\omega$ and continuous in $t$---and $(\omega, t) \mapsto |X(\omega)'(t)|$ is the $\limsup$ of difference quotients. 
Since
\begin{equation} \label{eq:estimate_horizontal_lift_velocity}
  \int_0^L \|X(\cdot)'(r)\|^2_{L^2} \, \diff r = \int \int_0^L |\gamma'_t|^2 \,\diff t \, \diff \boldsymbol{\eta}(\gamma) = \int_0^L | \bar{\mu}' |^2(t) \, \diff t,
\end{equation}
we have that $\bar{\bgamma}$ is in $\mathrm{AC}^2$, and $\| \dot{\bar{\bgamma}}_t \|_{L^2} \leq \|X(\cdot)'(t)\|_{L^2}$ for almost every $t \in [0, L]$. 
Thus, we further see from \eqref{eq:estimate_horizontal_lift_velocity} that $\int_0^L \| \dot{\bar{\bgamma}}_t \|_{L^2}^2 \, \diff t \leq \int_0^L | \bar{\mu}' |^2(t) \, \diff t$. 
As we always have that $|\bar{\mu}'|(t) \leq \|\dot{\bar{\bgamma}}_t\|_{L^2}$ by $\law\colon L^2(\Omega; \R^d) \to (\Pp_2, W_2)$ being 1-Lipschitz, we conclude from the above that $|\bar{\mu}'|(t) = \|\dot{\bar{\bgamma}}_t\|_{L^2}$ for a.e.~$t \in [0, L]$, i.e., that $\bar{\bgamma}$ is a horizontal lift of $\bar{\mu}_{\bullet}$. 
Then, since $\bar{\bgamma}$ is Lipschitz and $q$ is absolutely continuous, $\bgamma \coloneqq \bar{\bgamma} \circ q$ is absolutely continuous as well, and a direct calculation shows that $|\bgamma'|(t) = |\bar{\bgamma}'|(q(t)) q'(t)$. 
Likewise, $|\mu'|(t) = |\bar{\mu}'|(q(t)) q'(t)$, and so we get $|\bgamma'|(t) = |\mu'|(t)$. 
Since $\law(\bgamma_t) = \law(\bar{\bgamma}_{q(t)}) = \bar{\mu}_{q(t)} = \mu_t$, we have that $\bgamma$ is a horizontal lift of $\mu_{\bullet}$.\\
To lift gradient problems, what is left to prove is that for a functional $E\colon\Pp_2 \to \R \cup \{+\infty\}$, we have that $|\partial E| (\law(X)) = |\partial (E \circ \law)| (X)$ for all $X \in L^2(\Omega; \R^d)$. 
Indeed, this equality would imply that $\bgamma$ is a curve of maximal slope of $E \circ \law$ if and only if $\mu_{\bullet}$ is a curve of maximal slope of $E$. One direction follows directly from $\law$ being 1-Lipschitz. 
More specifically, let $X_n \xrightarrow{n \to \infty} X$ in $L^2$ be a sequence such that
\begin{equation*}
    |\partial (E \circ \law)| (X) = \limsup_{Y\to_{L^2} X} \frac{\big( E \circ \law(X) - E \circ \law(Y) \big)^+}{\| Y -X \|_{L^2}} = \limsup_{n\to \infty} \frac{\big( E \circ \law(X) - E \circ \law(X_n) \big)^+}{\| X_n -X \|_{L^2}}.
\end{equation*}
Then (discarding those $n$ for which $\law(X_n) = \law(X)$, since those do not contribute to the $\limsup$ above), we get
\begin{equation*}
\begin{split}
    |\partial (E \circ \law)| (X) &= \limsup_{n\to \infty} \frac{\big( E\big( \law(X)\big) - E\big( \law(X_n)\big) \big)^+}{\| X_n -X \|_{L^2}} \leq \limsup_{n\to \infty} \frac{\big( E\big( \law(X)\big) - E\big( \law(X_n)\big) \big)^+}{W_2\big( \law(X_n),\law(X) \big)} \\
    &\leq |\partial E|(\law(X)).
\end{split}
\end{equation*}
For the other direction, consider a sequence $\mu_n \xrightarrow{n \to \infty} \law(X)$ in $\Pp_2$ along which the $\limsup$ in the definition of $|\partial E|(\law(X))$ is achieved. 
If we could lift this sequence to a sequence $X_n \xrightarrow{n \to \infty} X$ in $L^2$ such that $\law(X_n) = \mu_n$ and such that we can control $\|X_n - X\|_{L^2}$ by $W_2(\mu_n, \law(X))$, we would be done. More precisely, assume the following claim:
\begin{claim} \label{claim:convergence_lifts}
Let $X \in L^2$, let $\mu_n \xrightarrow{n \to \infty} \law(X)$ in $\Pp_2$, and let $(\varepsilon_n)$ be any sequence of positive real numbers. Then, there exists a sequence $X_n \to X$ in $L^2$ such that $\law(X_n) = \mu_n$ and $\|X_n - X\|_{L^2} \leq W_2(\mu_n, \law(X)) + \varepsilon_n$.
\end{claim}
Given the above claim, we would have
\begin{equation*}
  |\partial E| (\law(X)) = \limsup_{n \to \infty} \frac{\big(E\big(\law(X)\big) - E(\mu_n) \big)^+}{W_2(\mu_n, \law(X))} \leq \limsup_{n \to \infty} \frac{\big( E \circ \law(X) - E \circ \law(X_n) \big)^+}{\|X_n - X\|_{L^2}} \leq |\partial (E \circ \law)| (X),
\end{equation*}
thus concluding the proof of the equality $|\partial E| (\law(X)) = |\partial (E \circ \law)| (X)$.
To prove the claim, we will instead prove the following stronger result:
\begin{proposition} \label{prop:wasserstein_distance_as_inf_of_l2_distance}
  Let $(\Omega, \mathcal{F}, \mathbb{P})$ be a probability space which admits a random variable $Z$ such that $\law(Z) = \Leb|_{[0, 1]}$, and set $L^2 \coloneqq L^2(\Omega; \R^d)$. Furthermore, let $X \in L^2$, and let $\mu \in \Pp_2$. Then, we have
  \begin{equation} \label{eq:wasserstein_distance_as_inf_of_l2_distance}
    W_2(\law(X), \mu) = \inf\left\{\|X - Y\|_{L^2(\Omega; \R^d)} \mid \law(Y) = \mu\right\}.
  \end{equation}
\end{proposition}
We stress that \cref{eq:wasserstein_distance_as_inf_of_l2_distance} is significantly stronger than the more standard \cref{eq:basic_inf_result}, since one of the two random variables is fixed. 
\begin{proof}
  The inequality ``$\leq$'' is immediate, as for any $Y$ such that $\law(Y) = \mu$, the joint law $\pi = \law(X, Y)$ is a coupling of $\law(X)$ and $\mu$. 
For the other direction, we rely on \cite[Lemma 5.23]{carmonaProbabilisticTheoryMean2018}, which states that for any atomless standard probability space $(\Omega', \mathcal{F}', \mathbb{P}')$ and $Z, Z' \in L^2(\Omega'; \R^d)$ such that $\law(Z) = \law(Z')$ and for any $\varepsilon > 0$, we can find two measurable measure preserving mappings $\iota, \iota^{-1}$ (i.e., $\iota_\# \mathbb{P}' = \mathbb{P}'$, and likewise for $\iota^{-1}$) from $\Omega'$ onto itself such that 
\begin{equation*}
\mathbb{P}'\Big(\big\{\omega \in \Omega' \mid (\iota \circ \iota^{-1})(\omega) = (\iota^{-1} \circ \iota)(\omega) = \omega\big\}\Big) = 1, \quad \text{and }\mathbb{P}'\Big(|Z - Z' \circ \iota^{-1}| \leq \varepsilon\Big) = 1.
\end{equation*}
We will apply this lemma with $\Omega' = \R^d \times [0, 1]$, $\mathcal{F}'$ being the Borel $\sigma$-algebra and $\mathbb{P}' \coloneqq (X, Z)_\# \mathbb{P}$; note that this probability space is atomless, since the pushforward along the projection onto the last coordinate is the Lebesgue measure. 
We fix $\varepsilon > 0$ and first use \cref{eq:basic_inf_result} to pick $X', Y' \in L^2(\Omega'; \R^d)$ such that $\|X' - Y'\|_{L^2(\Omega'; \R^d)} \leq  W_2(\law(X), \mu) + \varepsilon$ and $\law_{\mathbb{P}'}(X') = \law_{\mathbb{P}}(X)$ and $\law_{\mathbb{P}'}(Y') = \mu$.
Setting $p_1\colon \R^d \times [0, 1] \to \R^d$ to be the projection onto the first component---note that $\law_{\mathbb{P}'}(p_1) = \law_{\mathbb{P}}(X)$---, we use \cite[Lemma 5.23]{carmonaProbabilisticTheoryMean2018} to choose $\iota$ such that $\mathbb{P}'(|p_1 - X' \circ \iota^{-1}| \leq \varepsilon) = 1$. 
We set $Y \coloneqq Y' \circ \iota^{-1} \circ (X, Z)$---note that $\law_{\mathbb{P}}(Y) = \law_{\mathbb{P}'}(Y') = \mu$---and calculate
\begin{align*}
  \|X - Y\|_{L^2(\Omega; \R^d)} &= \|p_1 \circ (X, Z)  - Y' \circ \iota^{-1} \circ (X, Z)\|_{L^2(\Omega; \R^d)} = \|p_1  - Y' \circ \iota^{-1}\|_{L^2(\Omega'; \R^d)} \\
  &= \|p_1\circ \iota - Y'\|_{L^2(\Omega'; \R^d)} \leq \|p_1 \circ \iota - X'\|_{L^2(\Omega'; \R^d)} + \|X' - Y'\|_{L^2(\Omega'; \R^d)} \\
  &\leq W_2(\law(X), \mu) + 2 \varepsilon.
\end{align*}
Since $\varepsilon$ was arbitrary, this concludes the proof.
\end{proof}
We summarize the previous discussion in the following lemma.
\begin{lemma}\label{lem:summarize_sec_1}
  Let $(\Omega, \mathcal{F}, \mathbb{P})$ be a probability space which admits a random variable $Z$ such that $\law(Z) = \Leb|_{[0, 1]}$, set $L^2 \coloneqq L^2(\Omega; \R^d)$, and let $\law\colon L^2 \to \Pp_2$ be the function mapping a random variable to its law. \\
  Then, given a functional $E\colon \Pp_2 \to \Rext$, for every $g\colon \Pp_2 \to [0, +\infty]$ we have that $g$ is a weak/strong upper gradient of $E$ if and only if $g \circ \law$ is a weak/strong upper gradient of $E \circ \law$. 
  Furthermore, for any $Y \in L^2$, we have that $|\partial E|(\law(Y)) = |\partial (E \circ \law)|(Y)$. Lastly, whenever $g \circ \law$ is a weak upper gradient for $E \circ \law$ (or equivalently, $g$ is a weak upper gradient for $E$), then a curve $\bgamma\colon (a, b) \to L^2$ is a curve of maximal slope for $E \circ \law$ with respect to $g \circ \law$ if and only if (1.) $\law \circ \bgamma$ is a curve of maximal slope for $E$ with respect to $g$ and (2.) $\bgamma$ is a horizontal lift of $\law \circ \bgamma$.
\end{lemma}
\begin{proof}
  From the discussion above, we have that $|\partial E|(\law(Y)) = |\partial (E \circ \law)|(Y)$. Next, we fix $E$ and $g$ as in the lemma. Then, we have
  \begin{align}
    \forall \mu_\bullet \in AC((a, b); \Pp_2)&\colon \begin{cases}
    g \circ \mu_\bullet \text{ is Borel } \\ |E(\mu_t) - E(\mu_s)| \leq \int_s^t g(\mu_r) |\mu'|(r) \, \diff r  \quad \forall a < s \leq t < b
    \end{cases} \label{eq:sup_equivalence_1}\\
    \Leftrightarrow \forall \bgamma \in AC((a, b); L^2)&\colon \begin{cases}
    g \circ \law \circ \bgamma \text{ is Borel } \\ |E(\law(\bgamma_t)) - E(\law(\bgamma_s))| \leq \int_s^t g(\law(\bgamma_r)) |\bgamma'|(r) \, \diff r \quad \forall a < s \leq t < b,
    \end{cases} \label{eq:sup_equivalence_2}
  \end{align}
  where the implication $\Rightarrow$ follows from the fact that for $\bgamma \in AC((a, b); L^2)$, we have that $\law \circ \bgamma \in AC((a, b); \Pp_2)$ and $|(\law \circ \bgamma)'| \leq |\bgamma'|$ by $\law$ being 1-Lipschitz; the implication $\Leftarrow$ follows likewise by picking a horizontal lift. 
  Since \cref{eq:sup_equivalence_1} and \cref{eq:sup_equivalence_2} are just the definitions of $g$ and $g \circ \law$ being strong upper gradients (see \cite[Def. 1.2.1]{ambrosioGradientFlows2008}), we have shown that $g$ is a strong upper gradient if and only if $g \circ \law$ is a strong upper gradient; the argument for weak upper gradients (see \cite[Def. 1.2.2]{ambrosioGradientFlows2008}) is analogous. 
  Finally, let $\bgamma\colon (a, b) \to L^2$ such that $\bgamma \in \mathrm{AC}_{\mathrm{loc}}$ and $g$ be a weak upper gradient for $E$. Then, by \cite[Def. 1.3.2]{ambrosioGradientFlows2008}, $\bgamma$ is a curve of maximal slope for $E \circ \law$ with respect to $g \circ \law$ if and only if $E \circ \law \circ \bgamma$ is $L^1$-almost everywhere equal to a nonincreasing map $\varphi$ and
  \begin{equation} \label{eq:defining_inequ_coms}
    \varphi'(t) \leq -\frac{1}{2} |\bgamma'|^2(t) - \frac{1}{2} g^2(\law(\bgamma_t)), \qquad \text{for almost every } t \in (a, b).
  \end{equation}
  However, we also have that 
  \begin{equation*}
    -\frac{1}{2} |\bgamma'|^2(t) - \frac{1}{2} g^2(\law(\bgamma_t)) \leq -\frac{1}{2} |(\law \circ \bgamma)'|^2(t) - \frac{1}{2} g^2(\law(\bgamma_t)) \leq - g(\law(\bgamma_t)) |(\law \circ \bgamma)'|(t) \leq \varphi'(t),
  \end{equation*}
  where the first inequality comes from $\law$ being 1-Lipschitz and is an equality almost everywhere if and only if $\bgamma$ is a horizontal lift of $\law \circ \bgamma$, while the last inequality comes from $g$ being a weak upper gradient and uses the fact that $(g \circ \law \circ \bgamma) |(\law \circ \bgamma)'| \in L^1_{\mathrm{loc}}$. 
  This latter claim follows from either $\bgamma$ or $\law \circ \bgamma$ being a curve of maximal slope.
  Thus, we see that $\bgamma$ is a curve of maximal slope of $E\circ\law$ with respect to $g\circ\law$ if and only if $\bgamma$ is a horizontal lift of $\law \circ \bgamma$ and $\law\circ \bgamma$ is a curve of maximal slope of $E$ with respect to $g$.
\end{proof}

To close this section, we remark on continuity properties of functionals on $\Pp_2$.

\begin{remark} \label{rmk:continuity_of_functionals}
  From \cref{claim:convergence_lifts}, we can immediately deduce that for any functional $G\colon \Pp_2 \to \R$, $W_2$-(lower/upper semi-)continuity is equivalent to the $L^2$-(lower/upper semi-)continuity of $G \circ \law$. 
  Likewise, for the lower semicontinuous envelope, we have that $\mathrm{lsc}(G) \circ \law = \mathrm{lsc}(G \circ \law)$; this follows from the same argument and using the definition $\mathrm{lsc}(G)(\mu) = \inf_{\nu_i \to \mu} \liminf_{i} G(\nu_i)$.
\end{remark}

\subsection{Constructing Curves of Maximal Slope in Metric Spaces as per \texorpdfstring{\cite{ambrosioGradientFlows2008}}{[AGS08]}}\label{sec:constructing_curves_of_maximum_slope}
Now that we have established the correspondence \eqref{eq:correspondence_coms} between horizontal lifts of curves of maximal slope of a functional $E\colon\Pp_2 \to \R \cup \{+\infty\}$ in $\Pp_2$ and curves of maximal slope of $E \circ \law$ in $L^2$, a natural question to ask is whether we can construct curves of maximal slope in $L^2$ directly, as the latter is (geometrically speaking) a nice space. 
However, constructing curves of maximal slope in $L^2$ directly comes with topological obstructions.
In this section, we sketch the procedure of how curves of maximal slope in some metric space $(\mathcal{S}, d_{\mathcal{S}})$ are constructed (for the full construction, we refer to \cite{ambrosioGradientFlows2008}) to lay bare the difficulties that are posed by the choice $\mathcal{S}=L^2$ when following through this program. 
Throughout, we will assume an energy $E \colon \mathcal{S} \to \R \cup \{+\infty\}$ that fulfills the coercivity assumption as in \cite[2.1b]{ambrosioGradientFlows2008}, i.e., such that there exist $\tau^\star > 0$ and $x^\star \in \mathcal{S}$ such that $\inf_{y \in \mathcal{S}} E(y) + \frac{1}{2\tau^\star} d_{\mathcal{S}}^2(x^\star, y) > -\infty$. 
Furthermore, for $\tau > 0$ and $x \in \mathcal{S}$, we define the proximal point operator
\begin{equation} \label{eq:def_proximal_metric}
  J_{\tau, \mathcal{S}}[x] \coloneqq \argmin_{y \in \mathcal{S}}\left\{E(y) + \frac{1}{2\tau} d_{\mathcal{S}}^2(x, y)\right\}.
\end{equation}
\subsubsection{Step 1: The Minimizing Movement Scheme}\label{sec:step_1}
Given any $U_0 \in \mathcal{S}$ and $\tau \in (0, \tau^\star)$, we iteratively construct
\begin{equation}\label{eq:mms}\tag{MMS}
  \begin{split}
  U_\tau^{0} &\coloneqq U_0, \\
  U_\tau^{i+1} &\in J_{\tau, \mathcal{S}}[U_\tau^{i}].
  \end{split}
\end{equation}
The coercivity condition on $E$ together with the fact that $\tau < \tau^\star$ ensures at each step that $\inf_{y \in \mathcal{S}}\left\{E(y) + \frac{1}{2\tau} d_{\mathcal{S}}^2(U_\tau^{i}, y)\right\} > -\infty$; to ensure existence of the minimizers (i.e., that $J_{\tau, \mathcal{S}}[U_\tau^{i}] \neq \emptyset$), we use the direct method of the calculus of variations.
To this end, it suffices to have a topology $\sigma$ on $\mathcal{S}$ such that (1.) $E$ and $d_{\mathcal{S}}$ are $\sigma$-sequentially lower semicontinuous, and that (2.) the sublevel sets of $d_{\mathcal{S}}(\cdot, U_\tau^{i})$ are sequentially precompact in $\sigma$.\footnote{The fact that $\tau < \tau^\star$ together with the coercivity condition ensures that each sublevel set of the function $y \mapsto E(y) + \frac{1}{2\tau}d_\mathcal{S}^2(U_\tau^{i}, y)$ is contained in a sublevel set of $d_{\mathcal{S}}(\cdot, U_\tau^{i})$.}

\subsubsection{Step 2: Interpolations and A Priori Bounds}

For a minimizing movement sequence $(U_\tau^{i})_{i \in \NN}$, i.e., one generated as in \cref{eq:mms}, one can define the two interpolations (see \cite[Def. 2.0.2, Def. 3.2.1]{ambrosioGradientFlows2008})
\begin{alignat}{2}
  \bar{U}_\tau(t) &= U_\tau^{i+1}, &\quad &\text{for } t \in (\tau i, \tau (i+1)] \label{eq:constant_interpolation}\\
  \tilde{U}_\tau(t) &\in J_{\delta, \mathcal{S}}[U_\tau^{i}], &\quad &\text{for } t \in (\tau i, \tau (i+1)], \text{ where } \delta \coloneqq t - \tau i,\label{eq:definition_de_giorgi}
\end{alignat}
where we emphasise that $\tilde{U}$ (called the De Giorgi interpolation) is not uniquely defined by \eqref{eq:definition_de_giorgi}, and that the existence of a De Giorgi interpolation is guaranteed by the existence of a topology $\sigma$ which is compatible with $E$ and $d_{\mathcal{S}}$ as in \cref{sec:step_1}; we furthermore define the function (see \cite[(2.3.4)]{ambrosioGradientFlows2008})
\begin{equation*}
  |U'_\tau|(t) = \frac{d_{\mathcal{S}}(U_\tau^{i}, U_\tau^{i-1})}{\tau}, \qquad \text{if } t \in (\tau (i-1), \tau i].
\end{equation*}
By mostly algebraic calculations, one can, for each $S, T > 0$,  and $\tau \in (0, \tau^\star)$ small enough, and each minimizing movement sequence $(U_\tau^{i})_{i \in \NN}$ such that $E(U_\tau^{0}) \leq S$ and $d_{\mathcal{S}}^2(U_\tau^{0}, x^\star) \leq S$, derive the a priori bounds (see \cite[Lemma 3.2.2, (3.2.3)]{ambrosioGradientFlows2008})
\begin{gather}
  E(U_\tau^{i}) - E(U_\tau^{j}) \geq \int_{\tau i}^{\tau j} \frac{1}{2} |U'_{\tau}|^2(t) + \frac{1}{2} |\partial E|^2\big(\tilde{U}_{\tau}(t)\big)\, \diff t\, \qquad \text{ for all } i, j \in \NN \text{ such that } i < j \label{eq:a_priori_ags_1}\\
  d_{\mathcal{S}}^2(U_\tau^{n}, x^\star) \leq C, \quad \sum_{j=1}^n\frac{d_{\mathcal{S}}^2(U_\tau^{j}, U_\tau^{j-1})}{2 \tau} \leq C \qquad \text{ for all } 1 \leq n \leq \left\lfloor \frac{T}{\tau} \right\rfloor \label{eq:a_priori_ags_2}, \\
  d_{\mathcal{S}}^2(\tilde{U}_\tau(t), \bar{U}_\tau(t)) \leq C \tau \qquad \text{ for all } t \in [0, T], \label{eq:a_priori_ags_3}
\end{gather}
where $C$ depends only on $T$ and $S$.

\subsubsection{Step 3: Compactness and Arzelà--Ascoli}
Let us assume we have fixed a point $U_0$, are given a positive sequence $(\tau_j)_{j \in \NN}$ converging to $0$ and that, for each $\tau_j$, we have constructed a minimizing movement sequence $(U_{\tau_j}^{i})_{i \in \NN}$.
If there exists a topology $\sigma$ such that (1.) $d_{\mathcal{S}}$ is $\sigma$-sequentially lower semicontinuous, (2.) $\sigma$ is weaker than the topology induced by $d_{\mathcal{S}}$, and (3.) the ball $B_{\sqrt{C}}(x^\star) \coloneqq \{ x\in \mathcal{S} : d_{\mathcal{S}}(x,x^\star)\leq \sqrt{C} \}$ is $\sigma$-sequentially precompact, we can use an Arzelà--Ascoli argument (see \cite[Prop. 3.3.1]{ambrosioGradientFlows2008}) to obtain a subsequence (without relabeling) $(\tau_j)_{j \in \NN}$, an absolutely continuous limit curve $u$ and a function $A \in L^2(0, T)$ such that $|U'_{\tau_j}| \rightharpoonup A$ weakly in $L^2(0, T)$, $|u'| \leq A$ almost everywhere, and $\bar{U}_{\tau_j}(t) \xrightharpoonup{\sigma} u(t)$, for all $t \in [0, T]$. 
If $E$ is sequentially $\sigma$-lower semicontinuous, we can use \cref{eq:a_priori_ags_3} to take the limit in \cref{eq:a_priori_ags_1} and obtain
\begin{equation*}
  E(u(0)) - E(u(t)) \geq \int_0^t \frac{1}{2} |u'|^2(r) + \frac{1}{2} (\mathrm{lsc}_{\sigma}(|\partial E|))^2(u(r)) \, \diff r.
\end{equation*}
Whenever $\mathrm{lsc}_{\sigma}(|\partial E|)$ is a strong upper gradient of $E$, we can conclude that $u$ is a curve of maximal slope.

\subsubsection{Requirements on the Topology and \texorpdfstring{$\Pp_2$}{P₂}}

In the three steps above, we saw that the topology $\sigma$ entered into the construction of curves of maximal slope in two places: When ensuring that $J_{\tau, \mathcal{S}}[x]$ is nonempty for $\tau$ small enough and any $x \in \mathcal{S}$, and when using the Arzelà--Ascoli argument to construct the limit curve. 
For both passages, a topology $\sigma$ fulfilling the following four properties suffices:
\begin{properties}[label=\textbf{(S\arabic*)}, ref=S\arabic*]
  \propertyitem\label{prop:general_weaker_than_distance} The topology $\sigma$ is weaker than the one induced by $d_{\mathcal{S}}$.
  \propertyitem\label{prop:general_lsc} The distance $d_{\mathcal{S}}$ is sequentially lower semicontinuous with respect to $\sigma$.
  \propertyitem The balls of $d_{\mathcal{S}}$ are sequentially precompact with respect to $\sigma$.\label{prop:general_seq_precompact}
  \propertyitem The energy functional $E$ is sequentially $\sigma$-lower semicontinuous.\label{prop:E_lsc}
\end{properties}

For the choice $(\mathcal{S}, d_{\mathcal{S}}) = (\Pp_2, W_2)$, the usual candidate for $\sigma$ is $\sigma_{\mathrm{narrow}}$, the narrow topology: $\sigma_{\mathrm{narrow}}$ fulfills \cref{prop:general_weaker_than_distance,prop:general_lsc,prop:general_seq_precompact}, and many functionals $E\colon \Pp_2 \to \R \cup \{+\infty\}$ of interest are narrowly lower semicontinuous, fulfilling \cref{prop:E_lsc}.

\subsection{The Topology Problem of \texorpdfstring{$L^2$}{L²} for Arzelà--Ascoli}\label{sec:topology-problem}
If we want to construct curves of maximal slope as above for the choice $\mathcal{S} = L^2$, and for energies of the form $\tilde{E} = E \circ \law$, it is natural to replace \cref{prop:E_lsc} by asking for $\law\colon (L^2, \sigma) \to (\Pp_2, \sigma_{\mathrm{narrow}})$ to be continuous: Then, \cref{prop:E_lsc} is fulfilled for any functional of the form $\tilde{E} = E \circ \law$ with $E$ narrowly lower semicontinuous. In place of \cref{prop:general_lsc,prop:general_seq_precompact,prop:general_weaker_than_distance,prop:E_lsc}, we thus end up with the following desirable properties of $\sigma$:
\begin{properties}[label=\textbf{(L\arabic*)}, ref=L\arabic*]
  \propertyitem\label{prop:weaker_than_distance} The topology $\sigma$ is weaker than the one induced by $d_{L^2}$.
  \propertyitem\label{prop:lsc} The distance $d_{L^2}$ is sequentially lower semicontinuous with respect to $\sigma$.
  \propertyitem The balls of $d_{L^2}$ are sequentially precompact with respect to $\sigma$.\label{prop:seq_precompact}
  \propertyitem The projection $\law: (L^2, \sigma) \to (\Pp_2, \sigma_\text{narrow})$ is continuous.\label{prop:cont_projections}
\end{properties} 
Note that the weak topology $\sigma_w$ fulfills \cref{prop:weaker_than_distance,prop:lsc,prop:seq_precompact} but not \cref{prop:cont_projections}, while the initial topology $\sigma_{\law} \coloneqq \law^{-1}(\sigma_\text{narrow})$ induced by $\law$ fulfills \cref{prop:weaker_than_distance,prop:seq_precompact,prop:cont_projections} but not \cref{prop:lsc}. 
Thus, a first candidate for $\sigma$ could be the topology $\sigma(\sigma_w, \sigma_{\law})$ generated by the union of the two topologies. However, this construction fails, because there exist sequences $X_i$ in $L^2$ for which there does not exist a \emph{joint} limit point in the two topologies, so that $\sigma(\sigma_w, \sigma_{\law})$ does not fulfill \cref{prop:seq_precompact}. 
To see this unfortunate fact, let us choose a random variable $Z\colon \Omega \to \R$ such that $\law(Z) = \Leb|_{[0,1]}$ and let us consider the sequence $X_i$ defined as
\begin{equation*}
  X_i(\omega) = \left\{\begin{array}{ll}
    1 & \text{if } Z(\omega) \in [2 j / 2^{i+1}, (2 j + 1) / 2^{i+1}) \text{ for some } j = 0, \ldots, 2^i - 1, \\
    -1 & \text{otherwise.}
  \end{array}\right.
\end{equation*}
Then, it is immediate to see that $X_i \rightharpoonup 0$ in $L^2$, while $\law(X_i) = \frac{1}{2}\delta_{-1} + \frac{1}{2}\delta_1$ for all $i$, and hence any limit point $X$ of $X_i$ in the topology induced by $\sigma_{\law}$ must also fulfill $\law(X) = \frac{1}{2}\delta_{-1} + \frac{1}{2}\delta_1$; both these observations also hold for every subsequence. 
Thus, no subsequence can have a joint limit, and so no subsequence converges in $\sigma(\sigma_w, \sigma_{\law})$. 
Even worse, this elementary example proves that a topology $\sigma$ fulfilling \cref{prop:lsc,prop:seq_precompact,prop:cont_projections} cannot exist: Otherwise, by \cref{prop:seq_precompact}, there would exist a subsequence $X_{i_k}$ converging in $\sigma$ to some $Y$, and by \cref{prop:cont_projections}, $\mathbb{P}(Y = -1) = \mathbb{P}(Y = 1) = \frac{1}{2}$. Taking $\bar{Y} = -Y$, we have by \cref{prop:lsc} that $\liminf_{k \to \infty} d_{L^2}(X_{i_k}, \bar{Y}) \geq d_{L^2}(Y, \bar{Y}) = 2$, which a simple calculation\footnote{Taking $j, k \in \NN$ such that $j \neq k$, we calculate $d_{L^2}(X_{i_j}, \bar{Y})^2 + d_{L^2}(X_{i_k}, \bar{Y})^2 = 4 \big(\mathbb{P}(X_{i_j} \neq \bar{Y}) + \mathbb{P}(X_{i_k} \neq \bar{Y})\big) = 4\big(\mathbb{P}(X_{i_j} \neq \bar{Y} \lor X_{i_k} \neq \bar{Y}) + \mathbb{P}(X_{i_j} \neq \bar{Y} \land X_{i_k} \neq \bar{Y})\big) \leq 4\big(1 + \mathbb{P}(X_{i_k} = X_{i_j})\big) = 6$, and thus we have either $d_{L^2}(X_{i_j}, \bar{Y}) \leq \sqrt{3} < 2$ or $d_{L^2}(X_{i_k}, \bar{Y}) \leq \sqrt{3} < 2$.} shows is impossible.

\subsection{Arzelà--Ascoli without Lower Semicontinuity}\label{sec:a-a-w-o-lsc}
At this point, we find ourselves in an unsatisfying position: While we know that first-order gradient problems in $\Pp_2$ may be lifted to first-order gradient problems in $L^2$, we cannot use the classical Arzelà--Ascoli framework to construct solutions in $L^2$, as we know that a suitable topology does not exist. \\
A solution to this problem lies in relaxing the requirements in the statement of the Arzelà--Ascoli theorem; a possible relaxation is the following, which we introduce in a general metric space:

\begin{definition} \label{def:a-a-readiness}
  Let $(\mathcal{S}, d_{\mathcal{S}})$ be a complete metric space, and let $\sigma$ be a pre-Hausdorff\footnote{Meaning that two points being $T_0$-separable implies that they are also $T_2$-separable, see, e.g.,~\cite{stinePreHausdorffSpaces2008a}.} topology on $\mathcal{S}$ which is weaker than the topology induced by $d_{\mathcal{S}}$. 
  The metric space $(\mathcal{S}, d_{\mathcal{S}})$ is said to be \emph{A-A-ready} w.r.t.~$\sigma$ if, for any countable index set $I$ and a collection of sequences $\{(x_n^m)_{n \in \mathbb{N}}\}_{m \in I}$ in $\mathcal{S}$ such that for all fixed $m$, the sequence $(x_n^m)_{n \in \mathbb{N}}$ is bounded in $n$, there exists a subsequence $n_k$ (independent of $m$) and limit points $\bar{x}^m \in \mathcal{S}$ such that:
  \begin{enumerate}
    \item For all $m_1, m_2 \in I$, we have
    \begin{equation} \label{eq:aa_property_liminf}
      \liminf_{k \to \infty} d_{\mathcal{S}}(x_{n_k}^{m_1}, x_{n_k}^{m_2}) \geq d_{\mathcal{S}}(\bar{x}^{m_1}, \bar{x}^{m_2}).
    \end{equation}
    \item For all $m \in I$, we have
    \begin{equation} \label{eq:aa_property_limit}
      x_{n_k}^m \xrightarrow{\sigma} \bar{x}^m, \qquad \text{as } k \to \infty.
    \end{equation}
    \item \label{item:respects_strong_limits} For all $m \in I$ such that $d_{\mathcal{S}}(x_n^m, \tilde{x}^m) \to 0$ as $n \to \infty$ for some $\tilde{x}^m \in \mathcal{S}$, we have that $\bar{x}^m = \tilde{x}^m$. 
  \end{enumerate}
\end{definition}

\begin{remark}
    We observe that \cref{item:respects_strong_limits} in \cref{def:a-a-readiness} does not follow automatically from the fact that the topology $\sigma$ is weaker than the one induced by the distance $d_{\mathcal{S}}$, as a pre-Hausdorff topology---as $\sigma$ actually is---need not have unique limits.
    For this reason, \cref{item:respects_strong_limits} enforces that the $\sigma$-limit selected by the subsequence agrees with the strong limit.
\end{remark}

With the above definition, we are in a position to prove the following version of Arzelà--Ascoli:
\begin{proposition} \label{prop:aa-ready-allows-aa}
  Let $\sigma$ be a pre-Hausdorff topology on a complete metric space $(\mathcal{S}, d_{\mathcal{S}})$ such that $\mathcal{S}$ is A-A-ready w.r.t.~$\sigma$ (see \cref{def:a-a-readiness}).\\
  Let $I$ be a (potentially unbounded) interval in $\R$, $x \in \mathcal{S}$ and let $u_n\colon I \to \mathcal{S}$ be curves such that $\sup_{n \in \NN} d_{\mathcal{S}}(u_n(t), x) < \infty$ for all $t \in I$ and
  \begin{equation} \label{eq:aa-ready-allows-aa-equicontinuity-assumption}
    \limsup_{n \to \infty} d_{\mathcal{S}}(u_n(s), u_n(t)) \leq \delta(s, t) \quad \text{for all } s, t \in I,
  \end{equation}
  for a (symmetric) function $\delta\colon I \times I \to [0, \infty)$ such that
  \begin{equation} \label{eq:AA_modus_continuity_is_modus_continuity}
    \lim_{(s, t) \to (r, r)} \delta(s, t) = 0 \quad \text{for all } r \in I \setminus \mathcal{C},
  \end{equation}
  where $\mathcal{C}$ is an (at most) countable subset of $I$. \\
  Then, there exists a subsequence $(n_k)_{k\in \NN}$ and a limit curve $\bar{u}\colon I \to \mathcal{S}$ which is $d_{\mathcal{S}}$-continuous in $I \setminus \mathcal{C}$, such that:
  \begin{itemize}
    \item For all $s, t \in I \cap (\Q \cup \mathcal{C})$, we have that $d_{\mathcal{S}}(\bar{u}(s), \bar{u}(t)) \leq \delta(s, t)$. In particular, if $\delta$ is upper semicontinuous on $I \times I$, then this estimate holds for all $s, t \in I$.
    \item For all $t \in I$, we have that $u_{n_k}(t) \xrightarrow{\sigma} \bar{u}(t)$ as $k\to \infty$.
    \item For all $t \in I$ such that $d_{\mathcal{S}}\big( u_n(t), \tilde{x}_t \big)\to 0$  as $n\to\infty$ for some $\tilde{x}_t \in \mathcal{S}$, we have that $\bar{u}(t) = \tilde{x}_t$.
  \end{itemize}
\end{proposition}
\begin{proof}
  The proof is a simple adaptation of the classical proof of Arzelà--Ascoli. We first introduce $\mathcal{D} \coloneqq \{t \in I \mid u_n(t) \xrightarrow{d_{\mathcal{S}}} \tilde{x}_t \text{ for some } \tilde{x}_t \in \mathcal{S}\}$ and choose a countable set $\mathfrak{C} \subseteq I$ such that
  \begin{itemize}
      \item[--] $I \cap \Q \subseteq \mathfrak{C}$;
      \item[--] $\mathcal{C} \subseteq \mathfrak{C}$ (recall that $\mathcal{C}$ is introduced in \cref{eq:AA_modus_continuity_is_modus_continuity});
      \item[--] $\mathfrak{C} \cap \mathcal{D}$ is dense in $\mathcal{D}$.
  \end{itemize}
  We apply the A-A-readiness property to the countable collection of sequences $\{(u_n(s))_{n \in \mathbb{N}}\}_{s \in \mathfrak{C}}$, obtaining a subsequence $n_k$ and limit points $\hat{u}(s)$ for all $s \in \mathfrak{C}$. By \cref{eq:aa_property_liminf,eq:aa-ready-allows-aa-equicontinuity-assumption}, we thus have that $d_{\mathcal{S}}(\hat{u}(s_1), \hat{u}(s_2)) \leq \delta(s_1, s_2)$ for all $s_1, s_2 \in \mathfrak{C}$. We may thus extend $\hat{u}$ to a function $\bar{u}\colon I \to \mathcal{S}$ by $d_{\mathcal{S}}$-continuity so that $\bar{u}$ is $d_{\mathcal{S}}$-continuous on $I \setminus \mathcal{C}$---and thus, in particular $\sigma$-continuous, since $\sigma$ is weaker than the topology induced by $d_{\mathcal{S}}$.
  What is left to prove is the $\sigma$-convergence of $u_{n_k}(t)$ to $\bar{u}(t)$ for all $t \in I \setminus \mathfrak{C}$. To this end, we fix such a $t$, and note that $\bar{u}$ is continuous in $t$. We will show that any subsequence $n_k'$ of $n_k$ contains a further subsequence $n_k''$ such that $u_{n_k''}(t)$ $\sigma$-converges to $\bar{u}(t)$. To proceed, we fix such a subsequence $n_k'$. Then, we can choose a sequence $(s_i)_{i \in \NN}$ in $\mathfrak{C}$ with $s_i \to t$ and apply the A-A-readiness property to the countable collection of sequences $\{(u_{n_k'}(t))_{k \in \mathbb{N}}\} \cup \{(u_{n_k'}(s_i))_{k \in \mathbb{N}}\}_{i \in \NN}$, so that we obtain a subsequence $n_k''$ and limit points $\tilde{u}(t)$, $\tilde{u}(s_i)$ such that
  \begin{align}
    u_{n_k''}(t) \xrightarrow[k \to \infty]{\sigma }\tilde{u}(t), \qquad u_{n_k''}(s_i) \xrightarrow[k \to \infty]{\sigma} \tilde{u}(s_i), \quad &\text{for all } i \in \NN; \label{eq:AA_subsequence_convergences}\\
    \limsup_{k \to \infty} d_{\mathcal{S}}(u_{n_k''}(s_i), u_{n_k''}(t)) \geq d_{\mathcal{S}}(\tilde{u}(s_i), \tilde{u}(t)), \quad &\text{for all } i \in \NN. \label{eq:AA_subsequence_distances}
  \end{align}
  \cref{eq:AA_subsequence_distances} together with \cref{eq:aa-ready-allows-aa-equicontinuity-assumption,eq:AA_modus_continuity_is_modus_continuity} in particular implies that $\tilde{u}(s_i) \xrightarrow{d_{\mathcal S}} \tilde{u}(t)$, and because $\sigma$ is weaker than the topology induced by $d_{\mathcal S}$, we further have $\tilde{u}(s_i) \xrightarrow{\sigma} \tilde{u}(t)$.
  Since $\sigma$ is pre-Hausdorff, we can deduce from \cref{eq:AA_subsequence_convergences} via
  \begin{equation*}
    \bar{u}(s_i) \xleftarrow[k \to \infty]{\sigma} u_{n_k''}(s_i) \xrightarrow[k \to \infty]{\sigma} \tilde{u}(s_i) \qquad \forall i \in \NN
  \end{equation*}
  that $\bar{u}(s_i)$ and $\tilde{u}(s_i)$ are topologically indistinguishable in $\sigma$, for all $i \in \NN$. Thus, $\tilde{u}(s_i) \xrightarrow{\sigma} \tilde{u}(t)$ implies $\bar{u}(s_i) \xrightarrow{\sigma} \tilde{u}(t)$ and from
  \begin{equation*}
    \bar{u}(t) \xleftarrow[i \to \infty]{\sigma} \bar{u}(s_i) \xrightarrow[i \to \infty]{\sigma} \tilde{u}(t),
  \end{equation*}
  we can in turn deduce that $\bar{u}(t)$ is topologically indistinguishable from $\tilde{u}(t)$ in $\sigma$. 
  Finally, from the first limit in \cref{eq:AA_subsequence_convergences}, we deduce that also $u_{n_k''}(t) \xrightarrow{\sigma }\bar{u}(t)$. \\
  Lastly, we have to show that $\bar{u}$ respects strong limits. For $s \in \mathfrak{C} \cap \mathcal{D}$, this comes directly from \cref{item:respects_strong_limits} of \cref{def:a-a-readiness}; for general points $t \in \mathcal{D}$, we use the fact that $\mathfrak{C} \cap \mathcal{D}$ is dense in $\mathcal{D}$ and again apply a double-limit argument as above.
\end{proof}

\begin{remark}\label{rmk:proximal_point_operator_aa}
  As explained in \cref{sec:constructing_curves_of_maximum_slope}, when constructing curves of maximal slope, the topology $\sigma$ enters in two places: (1.) when ensuring that $J_{\tau, \mathcal{S}}[y]$ is nonempty, for $\tau$ small enough and for all $y \in \mathcal{S}$, and (2.) when using an Arzelà--Ascoli argument. While \cref{prop:aa-ready-allows-aa} is aimed at solving (2.), let us remark that whenever $E$ is sequentially $\sigma$-lower semicontinuous and coercive in the sense of \cite[2.1b]{ambrosioGradientFlows2008}, then (1.) is also resolved by $(\mathcal{S}, d_{\mathcal{S}})$ being A-A-ready w.r.t.~$\sigma$, by arguing similarly to the direct method of the calculus of variations: Fixing $y$ and taking a minimizing sequence $(z_n)_{n \in \NN}$ of the functional $F\colon z \mapsto E(z) + \frac{1}{2\tau} d_{\mathcal{S}}^2(y, z)$, we can invoke the A-A-readiness\footnote{As explained in \cref{sec:constructing_curves_of_maximum_slope}, because $E$ is coercive, every sublevel of $F$ is contained in a sublevel of $d(\cdot, y)$, for $\tau$ small enough.} with the choice $\{(x_n^m)_{n \in \NN}\}_{m \in \{1, 2\}} = \{(y, y, y, \dots), (z_1, z_2, z_3, \dots)\}$ to obtain a subsequence $(z_{n_k})_{k \in \NN}$ and limit points $\bar{y}, \bar{z}$. \Cref{item:respects_strong_limits} of \cref{def:a-a-readiness} guarantees that $\bar{y} = y$, and \cref{eq:aa_property_liminf} ensures that $\liminf_{k \to \infty} d_{\mathcal{S}}(y, z_{n_k}) \geq d_{\mathcal{S}}(y, \bar{z})$. Together with the sequential $\sigma$-lower semicontinuity of $E$, we deduce that $\liminf_{k \to \infty}F(z_{n_k}) \geq F(\bar{z})$, i.e., that $\bar{z} \in J_{\tau, \mathcal{S}}[y]$.
\end{remark}

\begin{remark}\label{rmk:generalized_aa_yields_coms}
    \cref{prop:aa-ready-allows-aa} allows us to construct curves of maximal slope on a complete metric space $(\mathcal{S}, d_{\mathcal{S}})$, given that $\mathcal{S}$ is A-A-ready with respect to a suitable topology $\sigma$. More precisely, let us assume that we can find a topology $\sigma$ such that $\mathcal{S}$ is A-A-ready with respect to $\sigma$, and such that the functional $E$ as well as $|\partial E|$ are lower semicontinuous with respect to $\sigma$, and let us further assume that $|\partial E|$ is a strong upper gradient. Then we may use \cref{prop:aa-ready-allows-aa} as a drop-in replacement for the Arzelà--Ascoli argument \cite[Proposition 3.3.1]{ambrosioGradientFlows2008}, and argue as in \cref{rmk:proximal_point_operator_aa} for $J_{\tau, \mathcal{S}}[x]$ being nonempty, such that the machinery of \cite{ambrosioGradientFlows2008}, as sketched in \cref{sec:constructing_curves_of_maximum_slope}---and in particular the conclusions of \cite[Theorem 2.3.1, Theorem 2.3.3]{ambrosioGradientFlows2008}---goes through essentially unchanged. In the next section, we will see that $L^2$---for a large enough base space $\Omega$---is A-A-ready with respect to $\law^{-1}(\sigma_{\text{narrow}})$, which allows us to construct curves of maximal slope for the composition $E \circ \law$, where $E$ is lower semicontinuous with respect to the narrow topology.
  
\end{remark}
\subsection{A-A-Readiness of \texorpdfstring{$\LtwoR$}{L²\_ℝ}} \label{sec:a-a-readiness-of-l2}

Now, we aim to apply \cref{prop:aa-ready-allows-aa} and to proceed as in \cref{rmk:generalized_aa_yields_coms} for an energy $\tilde{E}$ defined on $L^2\coloneqq L^2(\Omega;\R^d)$ and having the form $\tilde{E} \coloneqq E \circ \law$. 
More precisely, we are interested in the case where the energy $E\colon \Pp_2 \to \Rext$ is narrowly lower semicontinuous (so that $\tilde{E}$ is $\law^{-1}(\sigma_{\mathrm{narrow}})$-lower semicontinuous, where $\law^{-1}(\sigma_{\mathrm{narrow}})$ denotes the pullback of the narrow topology along the map $\law\colon L^2 \to \Pp_2$). 
Hence, we have to prove A-A-readiness of $L^2$ with respect to $\law^{-1}(\sigma_{\text{narrow}})$\footnote{We note that the pullback of a pre-Hausdorff topology is again pre-Hausdorff (\cite[Theorem 1.6]{stinePreHausdorffSpaces2008a}).}.
To establish this fact for $L^2(\Omega; \R^d)$, we need $\Omega$ to be `suitably large', and our arguments will rely on Young measures.

Before we show a countable extension to the classical main theorem about Young measures, we state the standard version for reference.
\begin{theorem}[Young 1937-1942, cited from \cite{rindlerCalculusVariations2018}, Theorem 4.1]\label{thm:young_compactness_from_rindler}
  Let $(\Omega, \mathcal{F}, \mathbb{P})$ be a  complete non-atomic standard probability space. Let $(V_j) \subseteq L^p(\Omega; \R^N)$ be a norm-bounded sequence, where $p \in [1,\infty]$. Then, there exists a subsequence (not explicitly labeled) and a family of probability measures
  \begin{equation*}
    (\nu_\omega)_{\omega \in \Omega} \subseteq \Pp(\R^N),
  \end{equation*}
  called the ($L^p$-)Young measure generated by the (sub)sequence $(V_j)$, such that the following assertions are true:
  \begin{enumerate}
    \item The family $(\nu_\omega)_\omega$ is weakly$^*$-measurable, i.e., for all Carathéodory integrands $f\colon \Omega \times \R^N \to \R$, the compound function
    \begin{equation*}
      \omega \mapsto \langle f(\omega, \bullet), \nu_\omega \rangle \coloneqq \int_{\R^N} f(\omega, x) \, \diff \nu_\omega(x), \qquad \omega \in \Omega
    \end{equation*}
    is $\mathcal{F}$-measurable. \label{item:young_limit_1}
    \item If $p \in [1, \infty)$, it holds that 
    \begin{equation*}
      \int_\Omega \int_{\R^N} |x|^p \, \diff \nu_\omega(x) \, \diff \mathbb{P}(\omega) < \infty.
    \end{equation*}\label{item:young_limit_2}
    \item For all Carathéodory integrands $f \colon \Omega \times \R^N \to \R$ with the property that the family $(f(\cdot, V_j(\cdot)))_j$ is uniformly $L^1$-bounded and equiintegrable, it holds that
    \begin{equation} \label{eq:defining_property_young_measure}
      f(\cdot, V_j(\cdot)) \rightharpoonup \Big( \omega \mapsto \int_{\R^N} f(\omega, x) \, \diff \nu_\omega(x) \Big) \text{ in } L^1(\Omega; \R),
    \end{equation}
    where $\rightharpoonup$ denotes weak convergence, i.e., for all $\varphi \in L^\infty(\Omega; \R)$, we have the convergence
    \begin{equation*}
      \int_\Omega f(\omega, V_j(\omega)) \cdot \varphi(\omega) \, \diff \mathbb{P}(\omega) \to \int_\Omega \int_{\R^N} f(\omega, x) \cdot \varphi(\omega) \, \diff \nu_\omega(x) \, \diff \mathbb{P}(\omega). 
    \end{equation*}\label{item:young_limit_3}
  \end{enumerate}
\end{theorem}

For any norm-bounded sequence $(V_j)_{j \in \NN} \subseteq L^p(\Omega; \R^N)$, where $p \in [1,\infty]$, we call a \emph{Young limit of $(V_j)$} any family $(\nu_\omega)_{\omega \in \Omega}$ such that \cref{item:young_limit_1,item:young_limit_2,item:young_limit_3} in \cref{thm:young_compactness_from_rindler} hold, and we denote this convergence by $V_j \xrightarrow{Y} \nu$.

\begin{remark}
  In the context of \cite{rindlerCalculusVariations2018}, $\Omega$ is assumed to be a bounded Lipschitz domain of $\R^d$; however, by choosing a mod-$0$ isomorphism\footnote{I.e., a measure-preserving function that is, when restricted to some full-measure set, a bijection with measurable inverse onto the full-measure image.} between a nonatomic standard probability space and $([0,1], \Leb)$ (existing by \cite[Theorem 9.2.2]{bogachevMeasureTheory2007}), we may instead assume $(\Omega, \mathcal{F}, \mathbb{P})$  to be any complete atomless standard probability space.
  We also recall that in the context of \cite{rindlerCalculusVariations2018}, a Carathéodory integrand $f$ is a function that is Lebesgue measurable (or $\mathcal{F}$-measurable, after generalizing to nonatomic standard probability spaces) in the first argument, fixing the second argument, and continuous in the second argument, for almost every fixed first argument (see \cite[Lemma 2.4]{rindlerCalculusVariations2018}).
\end{remark}

In the following, we describe a countable extension of Young measures, i.e., Young measures taking values in $\R^{\NN}$. To achieve this, we will build a tower of compatible Young measures taking values in $\R^N$, for each $N$ finite, and then use the Kolmogorov extension theorem to extend to $\R^{\NN}$.
In the next statement, $\R^{\NN}$ is endowed with the product topology and the corresponding Borel $\sigma$-algebra, which in turn coincides with the product $\sigma$-algebra $\mathcal B(\R)^{\otimes\NN}$. Finally, we observe that $L^0(\Omega; \R^\NN)$ denotes the space of measurable functions valued in $\R^\NN$.

\begin{proposition} \label{prop:countable_young}
  Let $(\Omega, \mathcal{F}, \mathbb{P})$ be a complete non-atomic standard probability space, and let $(Z_n)_{n \in \mathbb{N}}$ be a sequence in $L^0(\Omega; \R^\NN)$ such that for all $N \in \NN$, $\|p_{[N]} \circ Z_n\|_{L^2(\Omega; \R^N)} \leq C_N$ for all $n$ and some constant $C_N > 0$, where $p_{[N]}\colon \R^\NN \to \R^N$ is the projection onto the first $N$ coordinates. Then, there exists a subsequence $n_k$ and a function $\nu: \Omega \to \mathcal{P}(\R^\mathbb{N})$ such that:
  \begin{enumerate}
    \item For all $A \subseteq \R^\mathbb{N}$ measurable, the function $\omega \mapsto \nu_\omega(A)$ is measurable.
    \item For all $N \in \NN$, we have the Young convergence $p_{[N]} \circ Z_{n_k} \xrightarrow{Y} (p_{[N]})_\# \nu$, where the pushforward $(p_{[N]})_\# \nu$ is defined as the pointwise pushforward $\omega \mapsto (p_{[N]})_\# \nu_\omega$.
  \end{enumerate}
\end{proposition}

\begin{proof}
  By a diagonal argument, we can find a subsequence $n_k$ such that $p_{[N]} \circ Z_{n_k} \xrightarrow{Y} \nu_N$ for some Young measure $\nu_N \colon \Omega \to \mathcal{P}(\R^N)$ for all $N \in \NN$. 
  Since Young measures are essentially unique\footnote{If both $\nu^{(1)}$ and $\nu^{(2)}$ are a Young limit of some sequence $(V_j)_{j \in \NN}$ in $L^p(\Omega; \R^N)$ for some $p \in [1, \infty]$, then $\nu^{(1)}_\omega=\nu^{(2)}_\omega$ for almost every $\omega \in \Omega$.}, we have that there exists a full-measure set $\Omega' \subseteq \Omega$ such that for all $\omega \in \Omega'$ and all $N_1 \leq N_2$, we have that $(p_{[N_1]})_\# \nu_{N_2, \omega} = \nu_{N_1, \omega}$. 
  By applying the Kolmogorov extension theorem \cite[Theorem 7.7.1]{bogachevMeasureTheory2007} to each $\omega \in \Omega'$, we obtain a limit function $\nu: \Omega' \to \mathcal{P}(\R^\NN)$ such that $(p_{[N]})_\# \nu_\omega = \nu_{N, \omega}$ for all $N$ and all $\omega \in \Omega'$. Extending $\nu$ to $\Omega \setminus \Omega'$ in any arbitrary measurable way, we obtain $\nu$ such that property 2 above is fulfilled. For property 1, we use the monotone class theorem twice: We first note that for all $N \in \NN$ and $A \subseteq \R^N$ Borel, the function $\omega \mapsto \nu_{N, \omega}(A)$ is measurable: for open $A$, this holds by monotonically approximating $\chi_A$ by continuous functions, and then we generalize with the functional monotone class theorem. Thus, (1) holds for finite unions of finite cylinders. Moreover, for $(A_n)_{n \in \NN}$ an increasing sequence of measurable sets in $\R^{\NN}$ such that $\omega \mapsto \nu_\omega(A_n)$ is measurable for all $n$, we have that also $\omega \mapsto \nu_\omega(\bigcup_{n \in \NN} A_n) = \lim_{n \to \infty} \nu_\omega(A_n)$ is measurable; the analogue holds for decreasing sequences of measurable sets. Thus, by the monotone class theorem, $\omega \mapsto \nu_\omega(A)$ is measurable for all measurable $A$.
\end{proof}
As a next step, we note that we can realize the limit measure by a random variable. We use a standard result, slightly rephrased from \cite{bogachevMeasureTheory2007}:
\begin{proposition}[{\cite[Prop. 10.7.6]{bogachevMeasureTheory2007}}]\label{prop:young_measurable_realization}
  Let $(\Omega, \mathcal{F}, \mathbb{P})$ be a complete non-atomic standard probability space. Let $\nu\colon \Omega \to \mathcal{P}(\R^\NN)$ be a probability-measure-valued function such that for all $A \subseteq \R^\NN$ measurable, the function $\omega \mapsto \nu_\omega(A)$ is measurable. Then, there exists a $\mathcal{F} \otimes \mathcal{B}([0,1])$-measurable map $X\colon \Omega \times [0, 1] \to \R^\NN$ such that for $\mathbb{P}$-a.e.~$\omega \in \Omega$, we have $\law(X(\omega, \cdot)) = \nu_\omega$, where the law is computed with respect to the Lebesgue measure $\Leb$.
\end{proposition}
The above two propositions suggest a way forward to the A-A-readiness of $L^2$ with respect to $\law^{-1}(\sigma_{\text{narrow}})$ (we make this more precise in \cref{prop:l2_is_aa_prep,cor:l2_is_aa}): Given a countable collection of sequences $\{(X_n^m)_{n \in \NN}\}_{m \in \NN}$ (we may assume $I = \NN$ in \cref{def:a-a-readiness}), we define the sequence $(Z_n)_{n \in \NN}$ with $Z_n = (X_n^1, X_n^2, \ldots)$ in $L^0(\Omega; (\R^d)^\NN)$. Identifying $(\R^d)^\NN$ with $\R^\NN$, we can use \cref{prop:countable_young} to obtain a Young measure $\nu\colon \Omega \to \Pp((\R^d)^{\NN})$ such that, along a joint subsequence (without relabeling), $p_{[N]} \circ Z_n \xrightarrow{Y} (p_{[N]})_{\#} \nu$. Next, we could use \cref{prop:young_measurable_realization} to obtain a random variable $\bar{X}\colon \Omega \times [0, 1] \to (\R^d)^{\NN}$ realizing $\nu$, i.e., such that $\law(\bar{X}(\omega, \cdot)) = \nu_\omega$ for almost every $\omega \in \Omega$, and set $\bar{X}^m\colon \Omega \times [0, 1] \to \R^d$ to be the $m$-th component of this random variable, for each $m \in \NN$. Then, \cref{eq:defining_property_young_measure} expressed in terms of these realizing variables becomes
\begin{equation*}
  f(\cdot, X_n^{m_1}(\cdot), \dots, X_n^{m_N}(\cdot)) \rightharpoonup \left(\omega \mapsto \int_{[0, 1]} f(\omega, \bar{X}^{m_1}(\omega, \omega'), \dots, \bar{X}^{m_N}(\omega, \omega')) \, \diff \Leb_{[0, 1]}(\omega')\right),
\end{equation*}
where the convergence is weak in $L^1(\Omega; \R)$.
Testing with $f(\omega, x_1, x_2) = \min(|x_1 - x_2|, M)^2$ and using a truncation argument, we obtain
\begin{equation*}
  \liminf_{n \to \infty} d_{L^2(\Omega; \R^d)}(X_{n}^{m_1}, X_{n}^{m_2}) \geq d_{L^2(\Omega \times [0, 1]; \R^d)}(\bar{X}^{m_1}, \bar{X}^{m_2});
\end{equation*}
testing instead with $f(\omega, x) = g(x)$ for $g\colon \R^d \to \R$ a bounded continuous function, we obtain the narrow convergence of the laws, i.e., that $\law(X_n^{m}) \xrightarrow{\sigma_{\mathrm{narrow}}} \law(\bar{X}^m)$ as $n \to \infty$, for all $m \in \NN$. These two properties (nonexpansion and convergence of the laws) are very close to what is required for A-A-readiness of $L^2$ w.r.t. $\law^{-1}(\sigma_{\mathrm{narrow}})$. 

However, the ``limiting random variables'' $(\bar{X}^{m})_{m \in \NN}$ require the larger base space $\Omega \times [0, 1]$. Here, the choice of $\Omega$ comes in: By making $\Omega$ large enough, we ensure that we can always absorb the extra component $[0, 1]$ into $\Omega$, without modifying any of the random variables in $\{(X_n^m)_{n \in \NN}\}_{m \in \NN}$. 
We make this idea rigorous in the following:

\begin{proposition} \label{prop:l2_is_aa_prep}
  Let $\Omega_\R = [0, 1]^\R$ endowed with the product Borel $\sigma$-algebra $\mathcal{F}_\R \coloneqq \bigotimes_\R \mathcal{B}\big([0, 1]\big)$ and the product measure $\mathbb{P}_\R = \bigotimes_\R \Leb|_{[0, 1]}$, where $\Leb|_{[0, 1]}$ is the Lebesgue measure on $[0, 1]$. Then, for any countable collection of sequences $\{(X_n^m)_{n \in \mathbb{N}}\}_{m \in \NN}$ in $\LtwoR = L^2(\Omega_\R; \R^d)$ such that for all $m \in \NN$, the sequence $(X_n^m)_{n \in \NN}$ is bounded in $n$, there exists a subsequence $n_k$, independent of $m$, and limit points $\bar{X}^m \in \LtwoR$ such that, for any $N \in \NN$,
  \begin{equation*}
    \law\Big((X_{n_k}^1, \dots, X_{n_k}^N)\Big) \xrightarrow{\text{narrow}} \law\Big((\bar{X}^1, \dots, \bar{X}^N)\Big).
  \end{equation*}
  In other words, for any continuous bounded function $f\colon (\R^d)^N \to \R$, we have that
  \begin{equation} \label{eq:aa_prep_limit_equality}
    \lim_{k \to \infty} \EE\Big[f(X_{n_k}^1, \dots, X_{n_k}^N)\Big] = \EE\Big[f(\bar{X}^1, \dots, \bar{X}^N)\Big].
  \end{equation}
  Furthermore, we can choose $(\bar{X}^m)_{m \in \NN}$ such that, whenever $X_n^m \to \tilde{X}^m$ as $n \to \infty$ \emph{strongly} in $\LtwoR$ for some $m \in \NN$, then $\bar{X}^m = \tilde{X}^m$ in $\LtwoR$.
\end{proposition}

Before we prove \cref{prop:l2_is_aa_prep}, let us show how it implies the desired A-A-readiness.
\begin{corollary} \label{cor:l2_is_aa}
  Let $\Omega_\R = [0, 1]^\R$ endowed with the product Borel $\sigma$-algebra $\bigotimes_\R \mathcal{B}\big([0, 1]\big)$ and the product measure $\mathbb{P}_\R = \bigotimes_\R \Leb|_{[0, 1]}$, where $\Leb|_{[0, 1]}$ is the Lebesgue measure on $[0, 1]$. Then, $\LtwoR = L^2(\Omega_\R; \R^d)$ is A-A-ready with respect to $\law^{-1}(\sigma_\text{narrow})$.
\end{corollary}

\begin{proof}
  We have to show that the limit sequence whose existence is guaranteed by \cref{prop:l2_is_aa_prep} fulfills \cref{eq:aa_property_liminf,eq:aa_property_limit}; we start with \cref{eq:aa_property_liminf}. For $x, y \in \R^d$ and $M > 0$, we define $d_M(x, y) \coloneqq \min\{|x - y|, M\}$, so that $d_M(x, y) \nearrow |x-y|$ as $M \to \infty$. 
  Let us consider a countable index set $I$ and a collection of sequences $\{(X_n^m)_{n \in \mathbb{N}}\}_{m \in I}$ in $\LtwoR$ such that for all fixed $m$, the sequence $(X_n^m)_{n \in \NN}$ is bounded in $n$. By adding constant sequences if $I$ is finite, we may assume $I = \NN$. We use \cref{prop:l2_is_aa_prep} to obtain a subsequence $n_k$ fulfilling \cref{eq:aa_prep_limit_equality}.
  By applying \cref{eq:aa_prep_limit_equality}, with $f\colon (\R^d)^N \to \R$ being $(x^1, \dots, x^N) \mapsto d_M(x^{m_1}, x^{m_2})^2$, we see that
  \begin{equation} \label{eq:10}
    \lim_{k \to \infty} \EE\Big[d_M(X_{n_k}^{m_1}, X_{n_k}^{m_2})^2\Big] = \EE\Big[d_M(\bar{X}^{m_1}, \bar{X}^{m_2})^2\Big]
  \end{equation}
  for all $m_1, m_2 \in \NN$, where we take the expected value over $\Omega_\R$. \Cref{eq:aa_property_liminf} now follows:
  \begin{align}
    d_{\LtwoR}^2(\bar{X}^{m_1}, \bar{X}^{m_2}) &= \EE\Big[|\bar{X}^{m_1} - \bar{X}^{m_2}|^2\Big] = \EE\Big[\lim_{M \to +\infty} d_M(\bar{X}^{m_1}, \bar{X}^{m_2})^2\Big] \notag \\
    &= \lim_{M \to +\infty} \EE\Big[d_M(\bar{X}^{m_1}, \bar{X}^{m_2})^2\Big] = \sup_{M > 0} \lim_{k \to \infty }\EE\Big[d_M(X_{n_k}^{m_1}, X_{n_k}^{m_2})^2\Big] \label{eq:12}\\
    &\leq \liminf_{k \to \infty} \sup_{M > 0} \EE\Big[d_M(X_{n_k}^{m_1}, X_{n_k}^{m_2})^2\Big] = \liminf_{k \to \infty} \lim_{M \to \infty} \EE\Big[d_M(X_{n_k}^{m_1}, X_{n_k}^{m_2})^2\Big] \\
    &= \liminf_{k \to \infty} \EE\Big[\lim_{M \to \infty} d_M(X_{n_k}^{m_1}, X_{n_k}^{m_2})^2\Big] = \liminf_{k \to \infty} \EE\Big[|X_{n_k}^{m_1} - X_{n_k}^{m_2}|^2\Big], \label{eq:14}
  \end{align}
  where the first equalities in \cref{eq:12,eq:14} follow from monotone convergence.  \\
  To show \cref{eq:aa_property_limit}, we test narrow convergence by picking an arbitrary continuous bounded function $g\colon \R^d \to \R$. Then, again setting $f\colon (\R^d)^N \to \R$ to be $(x^1, \dots, x^N) \mapsto g(x^m)$, we get from \cref{eq:aa_prep_limit_equality} that
  \begin{equation*}
    \lim_{k \to \infty}\int g(x)\,\diff \big(\law(X_{n_k}^m)\big)(x) = \lim_{k \to \infty}\EE\Big[g(X^m_{n_k})\Big] = \EE\Big[g(\bar{X}^m)\Big] = \int g(x)\,\diff \big(\law(\bar{X}^m)\big)(x);
  \end{equation*}
  lastly, \cref{item:respects_strong_limits} of \cref{def:a-a-readiness} follows from the last sentence of \cref{prop:l2_is_aa_prep}.
\end{proof}
We now present the proof of \cref{prop:l2_is_aa_prep}. Since it involves a number of objects, we list them in the table below.
\begin{table}[H]
  \centering
  \renewcommand{\arraystretch}{1.25}
  \begin{tabular}{@{}lp{0.28\textwidth}p{0.50\textwidth}@{}}
    \toprule
    Object & Lives in & Meaning \\
    \midrule
    $X_n^m$ & $\LtwoR$ & $n$-th element of the $m$-th given sequence \\
    $Z_n$ & $L^0(\Omega_\R; (\R^d)^\NN) \cong L^0(\Omega_\R; \R^\NN)$ & $(X_n^1, X_n^2, \ldots)$ \\
    $\mathfrak{p}_\alpha$ & $\mathfrak{p}_\alpha\colon \Omega_\R \to [0,1]$ & projection onto the coordinate $\alpha$ \\
    $\mathfrak{p}_{\mathcal{B}}$ & $\mathfrak{p}_{\mathcal{B}} \colon \Omega_\R \to [0,1]^{\mathcal{B}}$ & projection onto the coordinates in $\mathcal{B} \subseteq \R$ \\
    $\mathcal{A}_n$, $\mathcal{A}$ & countable subsets of $\R$ & coordinates on which $Z_n$, resp.\ all $Z_n$, depend \\
    $\mathcal{A}'$ & countable subset of $\R$ &  $\mathcal{A} \cup \{\beta\}$, where $\beta \in \R \setminus \mathcal{A}$\\
    $\bar{Z}_n$ & $[0,1]^{\mathcal{A}} \to \R^\NN$ & $Z_n = \bar{Z}_n \circ \mathfrak{p}_{\mathcal{A}}$ \\
    $\mathbb{P}_{\mathcal{A}}$ & $\Pp([0,1]^{\mathcal{A}})$ & $\mathfrak{p}_{\mathcal{A}\,\#}\mathbb{P}_\R = \bigotimes_{\mathcal{A}} \Leb|_{[0, 1]}$ \\
    $\nu$ & $[0,1]^{\mathcal{A}} \to \Pp(\R^\NN)$ & Young limit of $\bar{Z}_{n_k}$ from \cref{prop:countable_young} \\
    $Y$ & $[0,1]^{\mathcal{A}'} \to \R^\NN$ & realization of $\nu$ from \cref{prop:young_measurable_realization} \\
    $\underline{p}_m$ & $(\R^d)^\NN \to \R^d$ & projection onto the $m$-th $\R^d$-factor \\
    $p_{[N]}$ & $\R^\NN \to \R^N$ & projection onto the first $N$ coordinates \\
    $\bar{X}^m$ & $\LtwoR$ & $\underline{p}_m \circ Y \circ \mathfrak{p}_{\mathcal{A}'}$, the limit point of $(X_{n}^m)_{n \in \NN}$ \\
    \bottomrule
  \end{tabular}\label{tab:objects_l2_is_aa_prep}
\end{table}
\begin{proof}[Proof of \cref{prop:l2_is_aa_prep}]
  Let a countable collection of sequences $\{(X_n^m)_{n \in \NN}\}_{m \in \NN}$ in $\LtwoR$ be given, such that for all $m$, $\|X_n^m\|_{\LtwoR} \leq C_m$ for all $n$ and some constant $C_m > 0$. 
  We define the sequence $Z_n = (X_n^1, X_n^2, \ldots)$ in $L^0(\Omega_\R; (\R^d)^\NN)$, and identify $(\R^d)^\NN$ with $\R^\NN$. Because each $Z_n$ is measurable with respect to the product $\sigma$-algebra $\bigotimes_\R \mathcal{B}\big([0, 1]\big)$, we have that there exists\footnote{Since the $\sigma$-algebra $\mathcal{B}([0,1]^{\NN}) = \bigotimes_{\NN}\mathcal{B}([0, 1])$ is countably generated, so is $Z_n^{-1}(\mathcal{B}([0,1]^{\NN}))$. For each of the countable generators $\{B_i\}_{i \in \NN}$, we have that $B_i$ lies in the $\sigma$-algebra generated by contably many projections $\{\mathfrak{p}_\alpha\}_{\alpha \in \tilde{\mathcal{A}}_{n, i}}$ (since $B_i \in \bigotimes_\R \mathcal{B}\big([0, 1]\big)$). Taking $\mathcal{A}_n = \bigcup_{i \in \NN} \tilde{\mathcal{A}}_{n, i}$, we obtain the result.} a countable subset $\mathcal{A}_n \subseteq \R$ such that $Z_n$ is measurable with respect to the $\sigma$-algebra generated by the projections $\mathfrak{p}_\alpha\colon \Omega_\R \to [0, 1]$ for $\alpha \in \mathcal{A}_n$.
  From this, it follows from the Doob-Dynkin lemma\footnote{See, e.g., \cite[Theorem 2.12.3]{bogachevMeasureTheory2007}, where it is stated for random variables with the codomain $\R$. By choosing a Borel isomorphism, it readily generalizes for the codomain being any nonatomic standard Borel space.} that the function $Z_n$ factors through the projection $\mathfrak{p}_{\mathcal{A}_n}\colon \Omega_\R \to [0, 1]^{\mathcal{A}_n}$, i.e., has the form $Z_n = \tilde{Z}_n \circ \mathfrak{p}_{\mathcal{A}_n}$ for some $\tilde{Z}_n\colon [0,1]^{\mathcal{A}_n} \to \R^\NN$ measurable. 
  By taking $\mathcal{A} = \bigcup_{n \in \NN} \mathcal{A}_n$, we deduce that all $Z_n$ factor through the projection $\mathfrak{p}_{\mathcal{A}}\colon \Omega_\R \to [0, 1]^{\mathcal{A}}$ as some $\bar{Z}_n\colon [0,1]^{\mathcal{A}} \to \R^{\NN}$, where $\mathcal{A}$ is countable. We can also assume that $\mathcal{A}$ is nonempty. 
  Choosing some index $\beta \in \R \setminus \mathcal{A}$ and setting $\mathcal{A}' \coloneqq \mathcal{A} \cup \{\beta\}$, we thus find ourselves in the following situation:\\
\begin{center}
\begin{tikzcd}
\Omega_\R \arrow[rrrd, "Z_n"] \arrow[d, "\mathfrak{p}_{\mathcal{A} \cup \{\beta\}}"]                               &  &  &                       \\
{[0,1]^\mathcal{A} \times [0,1] \cong [0,1]^{\mathcal{A'}}} \arrow[d, "\mathfrak{p}_A"] &  &  & \mathbb{R}^\mathbb{N} \\
{[0,1]^\mathcal{A}} \arrow[rrru, "\bar{Z}_n", dashed]                                                 &  &  &                      
\end{tikzcd}    
\end{center}
Next, setting $\mathbb P_{\mathcal{A}}\coloneqq \mathfrak{p}_{\mathcal{A}\, \#}\mathbb{P}_\R$, we apply \cref{prop:countable_young} to the countable collection of sequences $\bar{Z}_n$ in $L^0([0,1]^{\mathcal{A}}, \R^\NN)$, obtaining a subsequence $n_k$ and a function $\nu\colon [0,1]^{\mathcal{A}} \to \mathcal{P}(\R^\NN)$ such that for all $M \in \NN$, we have the Young convergence $p_{[M]} \circ \bar{Z}_{n_k} \xrightarrow{Y} (p_{[M]})_\# \nu$. 
By \cref{prop:young_measurable_realization}, we can furthermore find a measurable\footnote{Since \cref{prop:young_measurable_realization} assumes completeness of the $\sigma$-algebra of the domain, we may first have to replace $Y$ (resp.\ $\tilde{Y}$ in the proof of \cref{lem:lift_coupling}) by a $\Leb^{\mathcal{A} \cup \{\beta\}}$-a.e.\ equal function to get measurability with respect to the Borel $\sigma$-algebra in the domain (as opposed to the Lebesgue completion thereof).\label{footnote:replace_random_variable_to_get_measurable}} function $Y\colon [0,1]^{\mathcal{A}} \times [0, 1] \cong [0,1]^{\mathcal{A}'} \to \R^\NN$ such that for a.e.~$\omega_{\mathcal{A}} \in [0,1]^{\mathcal{A}}$, we have $\law(Y(\omega_{\mathcal{A}}, \cdot)) = \nu_{\omega_{\mathcal{A}}}$. Reidentifying $\R^\NN$ with $(\R^d)^\NN$ and letting $\underline{p}_n\colon (\R^d)^\NN \to \R^d$ be the projection onto the $n$-th $\R^d$-factor, we set $\bar{X}^m \coloneqq \underline{p}_m \circ Y \circ \mathfrak{p}_{\mathcal{A}'}$. 
To prove \cref{eq:aa_prep_limit_equality}, we let $f\colon (\R^d)^N \to \R$ be a bounded continuous function and use the Young convergence of $p_{[d \cdot N]} \circ \bar{Z}_{n_k}$ to $(p_{[d \cdot N]})_\# \nu$:
  \begin{align*}
    \lim_{k \to \infty} \EE\Big[f(X^1_{n_k}, \dots, X^N_{n_k})\Big] &= \lim_{k \to \infty} \int_{[0,1]^{\mathcal{A}}} f\Big(p_{[d \cdot N]}\big( \bar{Z}_{n_k}(\omega_{\mathcal{A}})\big)\Big)\, \diff \mathbb{P}_\mathcal{A}(\omega_{\mathcal{A}}) \\
    &= \int_{[0,1]^{\mathcal{A}}} \int_{\R^{\mathbb N}} f \Big(p_{[d \cdot N]}\big( x\big)\Big) \,\diff  \nu_{\omega_\mathcal{A}}(x)\, \diff \mathbb{P}_\mathcal{A}(\omega_{\mathcal{A}}) \\
    &= \int_{[0,1]^{\mathcal{A}}} \int_{[0,1]} f \Big(p_{[d \cdot N]}\big( Y(\omega_{\mathcal{A}}, \omega_{\beta})\big)\Big) \,\diff  \Leb(\omega_\beta)\, \diff \mathbb{P}_\mathcal{A}(\omega_{\mathcal{A}}) \\
    &= \int_{\Omega_\R } f \Big(\bar{X}^1(\omega), \dots, \bar{X}^N(\omega)\Big) \, \diff \mathbb{P}_\R(\omega) = \EE\Big[f(\bar{X}^1, \dots, \bar{X}^N)\Big].
  \end{align*}
  Thus, we have in particular that $\law(X_{n_k}^i) \xrightarrow{\text{narrow}} \law(\bar{X}^i)$ for all $i \in \NN$, and hence $\bar{X}^i \in \LtwoR$ follows from lower semicontinuity of the second moment under narrow convergence. 
  Finally, the fact that the limits coincide whenever some sequence $(X_n^m)_{n \in \NN}$ converges strongly follows from the essential uniqueness of the Young limit.
\end{proof}

\begin{remark}\label{rmk:infimum_is_achieved}
  By applying \cref{prop:l2_is_aa_prep} with $(X_n^1)_{n \in \NN} = (X, X, \dots)$ and $(X_n^2)_{n \in \NN} = (Y_1, Y_2, \dots)$, with $(Y_i)_{i \in \NN}$ a minimizing sequence of the infimum in \cref{eq:wasserstein_distance_as_inf_of_l2_distance}, we can see that the infimum in \cref{eq:wasserstein_distance_as_inf_of_l2_distance} is achieved. In other words, given any $X \in \LtwoR$ and $\mu \in \Pp_2$, there exists $Y \in \LtwoR$ with $\law(Y) = \mu$ such that the joint law $\law((X, Y)) \in \Pi(\law(X), \mu)$ is an optimal coupling.
  In the following, we prove a refinement of that statement: Given any $X \in \LtwoR$, $\mu \in \Pp_2$, and \emph{any} coupling $\pi \in \Pi(\law(X), \mu)$, we can find $Y \in \LtwoR$ such that $\law((X, Y)) = \pi$. Furthermore, we can choose $Y$ to be conditionally independent of any countable number of additional random variables.
\end{remark}

\begin{lemma} \label{lem:lift_coupling}
  Let $\Omega_\R = [0, 1]^\R$ endowed with the product Borel $\sigma$-algebra $\bigotimes_\R \mathcal{B}\big([0, 1]\big)$ and the product measure $\mathbb{P}_\R = \bigotimes_\R \Leb|_{[0, 1]}$, where $\Leb|_{[0, 1]}$ is the Lebesgue measure on $[0, 1]$. Moreover, let $X \in \LtwoR$, and let $\mu \in \Pp_2$; finally, let $\pi \in \Pi(\law(X), \mu)$. Then, there exists $Y \in \LtwoR$ such that $\law\big((X, Y)\big) = \pi$. Furthermore, given a countable collection $\{Z_i\}_{i \in I}$ of random variables $Z_i\colon \Omega_\R \to \R^d$, we can choose $Y$ such that, for all $i \in I$, $\mathbb{P}_\R(Y \in \cdot \mid Z_i=z, X=x) = \mathbb{P}_\R(Y \in \cdot \mid X=x)$\footnote{Here, $\mathbb{P}_\R(Y \in \cdot \mid Z_i=z, X=x)$ denotes a regular conditional distribution, i.e., a disintegration of $\law((Y, Z_i, X))$ along $(z, x)$; see \cite[Chapter 8]{kallenbergFoundationsModernProbability2021}. Hence, from $\law((X, Y)) = \pi$, we automatically have that $\mathbb{P}_\R(Y \in \cdot \mid X=x) = \pi_x$, a disintegration of $\pi$ along $x$.} for $\law\big((Z_i, X)\big)$-almost every $(z, x) \in \R^d \times \R^d$.
\end{lemma}
\begin{proof}
  As in the proof of \cref{prop:l2_is_aa_prep}, we pick a countable index set $\mathcal{A}$ and $\beta \in \R \setminus \mathcal{A}$ such that $X$ and all $Z_i$ factor through $\mathcal{A}$ as $\tilde{X}$ and $\tilde{Z}_i$, respectively, and set $\mathcal{A}' \coloneqq \mathcal{A} \cup \{\beta\}$. Using \cref{prop:young_measurable_realization}, we can pick $\tilde{Y}\colon [0, 1]^{\mathcal{A}'} \to \R^d$ measurable\footref{footnote:replace_random_variable_to_get_measurable} such that for $\mathbb{P}_{\mathcal{A}}$-a.e.~$\omega_{\mathcal{A}} \in [0, 1]^\mathcal{A}$ (with $\mathbb{P}_{\mathcal{A}} = \mathfrak{p}_{\mathcal{A}\, \#} \mathbb{P}_\R$ as in the proof of \cref{prop:l2_is_aa_prep}), we have $\law(\tilde{Y}(\omega_{\mathcal{A}}, \cdot)) = \pi_{\tilde{X}(\omega_{\mathcal{A}})}$; we recall that by $(\pi_x)_{x \in \R^d}$, we denote the disintegration of $\pi$ with respect to the first component. Setting $Y \colon \Omega_\R \to \R^d$ to be the composition of $\tilde{Y}$ with the projection $\mathfrak{p}_{\mathcal{A}'}\colon \Omega_\R \to [0, 1]^{\mathcal{A}'}$, we first show that $\law\big((X, Y)\big) = \pi$. To this end, note that we have
  \begin{align*}
    \law\big((X, Y)\big) &= \int_{[0, 1]^{\mathcal{A}}} \int_{[0,1]} \delta_{\tilde{X}(\omega_{\mathcal{A}})} \otimes \delta_{\tilde{Y}(\omega_{\mathcal{A}}, \omega_\beta)} \, \diff \Leb(\omega_\beta) \, \diff \mathbb{P}_{\mathcal{A}}(\omega_{\mathcal{A}}) \\
    &= \int_{[0, 1]^{\mathcal{A}}} \delta_{\tilde{X}(\omega_{\mathcal{A}})} \otimes \pi_{\tilde{X}(\omega_{\mathcal{A}})} \, \diff \mathbb{P}_{\mathcal{A}}(\omega_{\mathcal{A}}) = \int_{\R^d} \delta_x \otimes \pi_x \, \diff \law(X)(x) = \pi.
  \end{align*}
  Next, we consider the function $\nu\colon (\R^d \times \R^d \times \mathcal{B}(\R^d)) \to [0, 1]$ which maps $(z, x, A)$ to $\pi_{x}(A)$ (disregarding the argument $z$). Then, for each $z, x$, the function $\nu(z, x, \bullet)$ is a probability measure, while for fixed $A \in \mathcal{B}(\R^d)$, the function $\nu(\bullet, \bullet, A)$ is $\mathcal{B}(\R^d) \otimes \mathcal{B}(\R^d)$-measurable. Furthermore, we have, for $i \in I$ and $A \subseteq \R^d$, $B \subseteq \R^{2d}$ both Borel measurable,
  \begin{align*}
    \mathbb{P}_\R(Y \in A, (Z_i, X) \in B) &= \int_{[0, 1]^{\mathcal{A}}} \int_{[0,1]} \chi_A(\tilde{Y}(\omega_{\mathcal{A}}, \omega_\beta)) \chi_B\big((\tilde{Z}_i(\omega_{\mathcal{A}}), \tilde{X}(\omega_{\mathcal{A}}))\big) \, \diff \Leb(\omega_\beta) \, \diff \mathbb{P}_{\mathcal{A}}(\omega_{\mathcal{A}}) \\
    &=\int_{[0, 1]^{\mathcal{A}}} \pi_{\tilde{X}(\omega_{\mathcal{A}})}(A) \chi_B\big((\tilde{Z}_i(\omega_{\mathcal{A}}), \tilde{X}(\omega_{\mathcal{A}}))\big) \, \diff \mathbb{P}_{\mathcal{A}}(\omega_{\mathcal{A}}) \\
    &=\int_B \pi_x(A) \, \diff \law\big((Z_i, X)\big)(z, x) = \int_B \nu(z, x, A) \, \diff \law\big((Z_i, X)\big)(z, x).
  \end{align*}
  Thus, $\nu$ fulfills the defining property\footnote{The three properties shown above define $\mathbb{P}_\R(Y \in \cdot \mid (Z_i, X))$ essentially uniquely; see either \cite[Theorem 3.4/Theorem 8.5]{kallenbergFoundationsModernProbability2021} or \cite[Theorem 5.3.1]{ambrosioGradientFlows2008}.} of $\mathbb{P}_\R(Y \in \cdot \mid (Z_i, X))$, which shows that $\mathbb{P}_\R(Y \in \cdot \mid Z_i=z, X=x) = \mathbb{P}_\R(Y \in \cdot \mid X=x)$ for $\law\big((Z_i, X)\big)$-almost every $(z, x) \in \R^d \times \R^d$.
\end{proof}
\begin{remark}\label{rmk:mms_downstairs_upstairs}
  For $X \in \LtwoR$ and $\tau > 0$, and some functional $E\colon \Pp_2(\R^d) \to \Rext$, we have that $\law(J_{\tau, \LtwoR}[X]) = J_{\tau, \Pp_2}[\law(X)]$, where $J_{\tau, \LtwoR}[X] \coloneqq \argmin_{Y \in \LtwoR}(E(\law(Y)) + \frac{1}{2\tau}\|Y - X\|^2_{\LtwoR})$ and $J_{\tau, \Pp_2}[\law(X)] \coloneqq \argmin_{\nu \in \Pp_2}(E(\nu) + \frac{1}{2\tau} W_2^2(\law(X), \nu))$. This follows directly from the fact that as in \cref{rmk:infimum_is_achieved}, the infimum in \cref{eq:wasserstein_distance_as_inf_of_l2_distance} is achieved for $\Omega_\R=[0,1]^\R$. Thus, we can see that for $X \in \LtwoR$, running the MMS scheme~\eqref{eq:mms} in $\Pp_2$ starting at $\law(X)$ is equivalent to running the MMS scheme in $\LtwoR$ starting at $X$ (and considering the law $\law$ of the resulting sequence).
\end{remark}
\begin{remark}
  The space $\LtwoR$ is not separable. In fact, there is no hope in finding non-atomic $(\Omega, \mathcal{F}, \mathbb{P})$ such that $L^2(\Omega; \R^d)$ is separable and $\law^{-1}(\sigma_{\text{narrow}})$ is A-A-ready as per \cref{def:a-a-readiness}. To see this, we claim that for any topology $\sigma$ on some separable Hilbert space $\mathcal{H}$, $\mathcal{H}$ being A-A-ready with respect to $\sigma$ implies that for any norm bounded sequence $(x_n)_{n \in \NN}$, we can find a $\sigma$-limit $\bar{x}$ along some subsequence $n_k$ such that $\liminf_{k \to \infty}\|y - x_{n_k}\| \geq \|y - \bar{x}\|$, for all $y \in \mathcal{H}$---something we proved to be impossible in \cref{sec:topology-problem} for $\sigma=\law^{-1}(\sigma_{\text{narrow}})$. To see the claim, we pick a countable dense subset $\{y^m\}_{m \in \NN}$. Then, we can combine any norm-bounded sequence $(x_n)_{n \in \NN}$ with the constant sequences $y^m_n = y^m$  for all $n \in \NN$ and for all $m \in \NN$, and apply the A-A-readiness property to the resulting countable collection $\{(y^m_n)_{n \in \NN}\}_{m \in \NN} \cup \{(x_n)_{n \in \NN}\}$ of sequences to yield a subsequence and limit points $\bar{x}$ and $\bar{y}^m$; since the constant sequences converge strongly, we have that $\bar{y}^m = y^m$. Then, we get from \cref{eq:aa_property_liminf} that $\liminf_{k \to \infty} \|y^m - x_{n_k}\| \geq \|y^m - \bar{x}\|$; since $y^m$ is dense, this extends to all $y \in \mathcal{H}$. \\
  On the other hand, it is \cref{item:respects_strong_limits} of \cref{def:a-a-readiness}---the $\sigma$-limits respecting strong limits---which forces non-separability of the $L^2$ space. Without requiring \cref{item:respects_strong_limits}, we could have chosen $\Omega = [0,1]$: For any countable collection $\{(X_n^m)_{n \in \NN}\}_{m \in \NN}$ of sequences of $\R^d$-random variables on $[0, 1]$ with bounded second moments, it suffices to apply the A-A-readiness property to $\LtwoR$ for the collection of sequences $\{(X_n^m \circ \mathfrak{p}_0)_{n \in \NN}\}_{m \in \NN}$, for $\mathfrak{p}_0\colon [0, 1]^\R \to [0, 1]$ the projection to the $0$-th index. Then, the resulting limit points $(\bar{X}^m)_{m \in \NN}$ jointly factor through $[0, 1]^\mathcal{A}$ for some countable $\mathcal A$ as $\tilde{X}_m \colon [0, 1]^\mathcal{A} \to \R^d$, and choosing any mod-$0$-isomorphism $F\colon [0,1] \to [0,1]^\mathcal{A}$, we have that $\law_{\Leb|_{[0, 1]}}(\tilde{X}_m \circ F) = \law_{\otimes_{\R}\Leb|_{[0, 1]}}(\bar{X}^m)$. Since also $\law_{\Leb|_{[0, 1]}}(X_n^m) = \law_{\otimes_{\R}\Leb|_{[0, 1]}}(X_n^m \circ \mathfrak{p}_0)$, we get that $X_{n_k}^m \to \tilde{X}_m \circ F$ in $\law^{-1}(\sigma_{\text{narrow}})$. Likewise, both $F$ and $\mathfrak{p}_0$ preserve the $L^2$-distance, so also \cref{eq:aa_property_liminf} holds. \\
  Here, we decided to demand \cref{item:respects_strong_limits}, since it allows the limit curve in the Arzelà--Ascoli \cref{prop:aa-ready-allows-aa} to, for example, preserve joint starting points. Furthermore, we use that the limit respects strong limits in \cref{sec:newtonian_curves}.
\end{remark}

\section{Identifying PDEs of Gradient Flows in \texorpdfstring{$L^2$}{L²}}\label{sec:PDE_gradient_flows}

One of the properties of Wasserstein-$\mathrm{AC}^2$ curves $\mu_\bullet\colon (a, b) \to \Pp_2$ is the fact that they fulfill weakly a continuity equation
\begin{equation}\label{eq:CE}\tag{CE}
  \partial_t \mu_t + \nabla \cdot (v_t \mu_t) = 0, \qquad \text{in } \R^d \times (a, b),
\end{equation}
where $v\colon (x, t) \mapsto v_t(x) \in \R^d$ is a Borel velocity field satisfying
\begin{equation*}
  \int_a^b\int_{\R^d} |v_t(x)|^2 \, \diff \mu_t(x) \, \diff t < +\infty,
\end{equation*}
see \cite[Theorem 8.3.1]{ambrosioGradientFlows2008}.
In addition, whenever $\mu_\bullet$ is a curve of maximal slope with respect to the metric slope $|\partial E|$ for sufficiently regular $E$ (see \cite[Definition 10.3.9]{ambrosioGradientFlows2008}), the vector field $v$ observed in \cref{eq:CE} is given as $v_t = - \partial^0 E(\mu_t)$, where $\partial^0$ denotes the minimal-norm (in the sense of $L^2(\mu_t)$) element of the reduced Wasserstein subdifferential (see \cite[Theorem 11.1.3]{ambrosioGradientFlows2008}). In that setting, the continuity equation~\eqref{eq:CE} becomes
\begin{equation} \label{eq:CE_COMS}
  \partial_t \mu_t - \nabla \cdot (\partial^0 E(\mu_t) \, \mu_t) = 0, \qquad \text{in } \R^d \times (a, b).
\end{equation}
The connection of Wasserstein-curves of maximal slope with PDEs of the form~\eqref{eq:CE_COMS} is one of the major drivers of the success of the related theory.

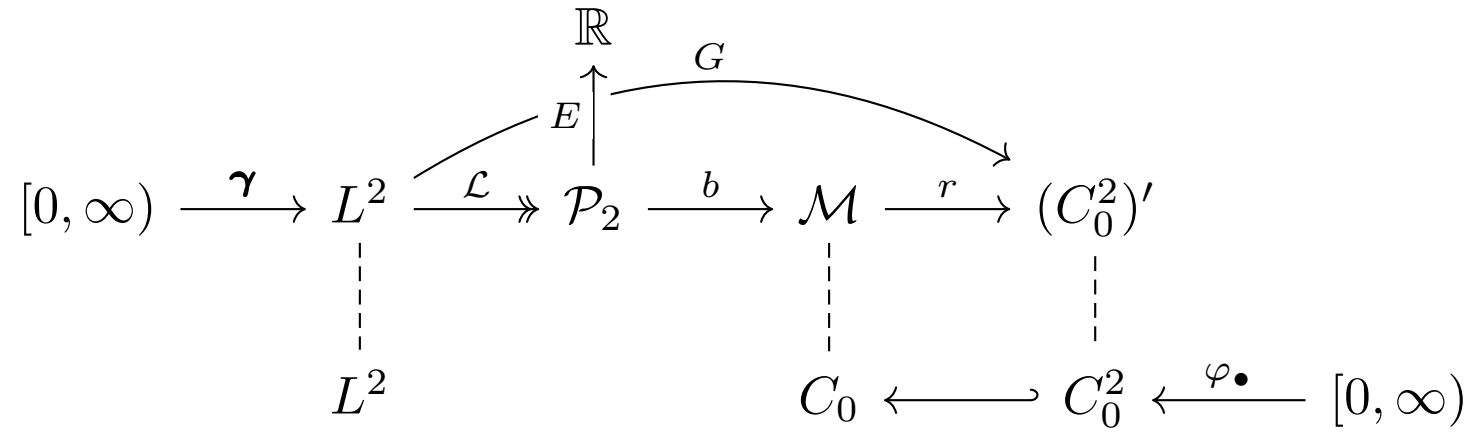
\begin{figure}[H]
\begin{center}
\begin{tikzcd}
                             &                                                                                      & \mathbb{R}                           &                                                         &                                     &                                     \\
{[0, \infty)} \arrow[r, "\bgamma"] & L^2 \arrow[r, "\law", two heads] \arrow[d, no head, dashed] \arrow[rrr, "G", bend left] & \Pp_2 \arrow[r, "b"] \arrow[u, "E"{fill=white, inner sep=2pt}, crossing over] & \mathcal{M} \arrow[r, "r"] \arrow[d, no head, dashed] & (C_0^2)' \arrow[d, no head, dashed] &                                     \\
                             & L^2                                                                                  &                                      & C_0                                                   & C_0^2 \arrow[l, hook']              & {[0, \infty)} \arrow[l, "\varphi_\bullet"']
\end{tikzcd}
\end{center}
\caption{The relations of the relevant spaces. Here, the vertical dashed lines denote duality pairings. Moreover, $\mathcal{M}$ denotes the space of finite signed Radon measures on $\R^d$, $b$ is the canonical inclusion and $r$ the restriction of test functions. While we will not make use of the intermediate spaces $\mathcal{M}$ and $C_0$, we still show them for the complete conceptual picture.}\label{fig:all_spaces}
\end{figure}

If we are given a curve of maximal slope $\bgamma$ in $L^2$ for an energy $\tilde{E}\colon L^2 \to \R \cup \{+\infty\}$ of the form $\tilde{E} = E \circ \law$, we can deduce through the correspondence~\eqref{eq:correspondence_coms} that $\mu_\bullet = \law \circ \bgamma$ fulfills weakly the PDE~\eqref{eq:CE_COMS}, whenever $E$ is regular in the sense of~\cite[Definition 10.3.9]{ambrosioGradientFlows2008}. 
This derivation of the PDE fulfilled by $\mu_\bullet$ is, from the Lions lift perspective, not completely transparent: First, such a derivation does not connect the vector field $v_\bullet$ in \cref{eq:CE} to the absolutely continuous curve $\bgamma$, as its connection to the ``velocity'' of $\bgamma$ remains elusive. 
Second, we may wonder if it is strictly necessary to go through the Wasserstein theory, including the theory of Wasserstein subgradients and regularity therein, to obtain \cref{eq:CE_COMS}. Third, by deriving a PDE fulfilled by $\mu_\bullet$ through the Lions lift, we might hope to lower the regularity assumptions imposed on $E$. 
In this section, we derive a PDE fulfilled by the laws $\mu_\bullet$ of curves of maximal slope in $L^2$ through the Lions lift, deriving a weaker form of \cref{eq:CE_COMS} (in \cref{prop:weak-formulation-abstract}), and, on our way, revealing the connection between $v_\bullet$ in~\eqref{eq:CE} and $\bgamma$ (in \cref{rmk:CE_X}).
To this end, we do not rely on the Fr\'echet derivative of $\tilde{E} = E \circ \law$, so as to avoid imposing too strong regularity assumptions on $E$. Instead, we will identify the PDE fulfilled by $\mu_\bullet$ by testing it against functions $\varphi \in C_c^\infty(\R^{d+1})$. We interpret such a test function as an absolutely continuous curve $\varphi_\bullet\colon [0, \infty) \to C_0^2(\R^d)$ defined as $t \mapsto \varphi_t(\cdot)$, where $C_0^2(\R^d)$ denotes the space of twice continuously differentiable functions which, together with their first and second derivatives, vanish at infinity, equipped with the norm $\|f\|_{C^2} = \|f\|_\infty + \|\nabla f\|_\infty + \|\nabla^2 f\|_\infty$. Note that $C_0^2(\R^d)$ is a Banach space, and that the inclusion $C_0^2(\R^d) \to C_0(\R^d)$ induces a map on the dual spaces (see \cref{fig:all_spaces}).
Our main idea is to use the product rule in the duality pairing $((C_0^2)', C_0^2)$ of the two curves $G \circ \bgamma$ and $\varphi_\bullet$, where $G \coloneqq r \circ b \circ \law$,  i.e., we use the fact that 
\begin{equation*}
  \frac{\diff}{\diff t} \langle G(\bgamma_t), \varphi_t \rangle_{(C_0^2)' \times C_0^2} = \langle \frac{\diff}{\diff t} ( G(\bgamma_t)), \varphi_t \rangle_{(C_0^2)' \times C_0^2} +  \langle G(\bgamma_t), \partial_t \varphi_t \rangle_{(C_0^2)' \times C_0^2},
\end{equation*}
where $\langle \cdot, \cdot \rangle_{(C_0^2)' \times C_0^2}$ denotes the duality pairing in $(C_0^2)' \times C_0^2$. 
To this end, we first have to identify the Fréchet derivative of the map $G \colon L^2 \to (C_0^2)'$. We defer the proof of the following lemma to \cref{sec:proof-of-frechet-derivative-l2-w2}, where we make precise the following formal computation:
\begin{equation*}
  \langle DG(X)[H], \psi \rangle_{(C_0^2)' \times C_0^2} 
  =
  DG(X)[H](\psi) = \lim_{s \to 0} \frac{1}{s} \int \big(\psi(X + s H) - \psi(X)\big) \, \diff \mathbb{P} = \int \langle \nabla \psi(X), H \rangle_{\R^d} \, \diff \mathbb{P},
\end{equation*}
for all $X, H \in L^2$ and all $\psi \in C^2_0$.

\begin{lemma} \label{prop:frechet-derivative-L2-W2}
  The map $G\coloneqq r \circ b \circ \law\colon L^2 \to (C_0^2)'$ is $1$-Lipschitz and Fréchet differentiable, and the Fréchet derivative is given by the function $DG(X)\colon L^2 \to (C_0^2)'$, defined as
  \begin{equation} \label{eq:frechet-derivative-L2-W2}
  \begin{split}
    DG(X)[H] \colon  C_0^2 & \to \R \\
    \psi &\mapsto \langle DG(X)[H], \psi \rangle_{(C_0^2)' \times C_0^2} =  \int \big\langle \nabla \psi(x), Z(x) \big\rangle_{\R^d} \, \diff(\law(X))(x) \\
    &  \qquad \qquad \qquad \qquad \qquad \qquad= \int \big\langle \nabla \psi\big(X(\omega)\big), H(\omega) \big\rangle_{\R^d} \, \diff \mathbb{P}(\omega),
  \end{split}
  \end{equation}
  where $Z \colon \R^d \to \R^d$ is any Borel map such that $Z(X(\omega)) =  \mathbb{E}[H | X](\omega)$ for $\mathbb{P}$-a.e.~$\omega$. In particular, $DG(X)$ has the preadjoint
  \begin{equation} \label{eq:pre-adjoint-L2-W2}
    DG(X)_*\colon C_0^2 \to L^2, \qquad DG(X)_*[\psi] = \Big[ \omega \mapsto \nabla \psi\big(X(\omega)\big) \Big].
  \end{equation}
\end{lemma}

  To allow us to describe the PDE fulfilled by curves of maximal slope in $L^2$, we need the following abstract proposition about curves of maximal slope in Hilbert spaces.
  
  \begin{proposition} \label{lem:technical-coms}
    Let $\mathcal{H}$ be a Hilbert space, let $F\colon \mathcal{H} \to \R \cup \{+\infty\}$ be a functional and let $u\colon [0, \infty) \to \mathcal{H}$ be a curve of maximal slope for $F$ with respect to the metric slope $|\partial F|$, defined using the metric induced by Hilbert space norm. Then, for almost every $t \in [0, \infty)$, the following holds: For any $v \in \mathcal{H}$ such that $\frac{\diff}{\diff \tau} F(u(t) + \tau v) \big|_{\tau = 0}$ exists in $\R$ (and $v = v_t$ may depend on $t$), we have
    \begin{equation} \label{eq:hilbert_space_coms}
      \langle \dot{u}(t), v \rangle_{\mathcal{H}} = - \frac{\diff}{\diff \tau} F(u(t) + \tau v) \Big|_{\tau = 0}.
    \end{equation}
  \end{proposition}
  
  \begin{remark}
    The requirement of the limit $\frac{\diff}{\diff \tau} F(u(t) + \tau v) \big|_{\tau = 0}$ existing \emph{in $\R$}---as opposed to existing in the extended reals---is not restrictive, as it can only happen for a negligible set of times $t$ that $\frac{\diff}{\diff \tau} F(u(t) + \tau v) \big|_{\tau = 0} \in \{-\infty, +\infty\}$. Indeed, the limit having infinite absolute value forces $|\partial F|(u(t)) = +\infty$, which is false for almost every $t$, since $u$ is a curve of maximal slope with respect to $|\partial F|$.
  \end{remark}

  \begin{remark}
    \Cref{lem:technical-coms} does not provide a characterization of curves of maximal slope in general. Indeed, let us consider the example of $F\colon \R^2 \to \R$ defined by 
    \begin{equation*}
      F(x) \coloneqq \begin{cases}
        0 \qquad &\text{if } |x|=1\\
        1 \qquad &\text{otherwise}.
      \end{cases}
    \end{equation*} 
    Then, taking any curve $u\colon [0, +\infty) \to \R^2$ such that $|u(t)|=1$ for all $t \in [0, +\infty)$, the limit on the right hand side of \cref{eq:hilbert_space_coms} does not exist for any $t \in [0, +\infty)$ and any $v \in \R^2$, except if $v=0$, in which case \cref{eq:hilbert_space_coms} is trivially fulfilled. Thus, any such curve fulfills the conclusion of \cref{lem:technical-coms}, but only the constant curves are curves of maximal slope with respect to the metric slope $|\partial F|$.
  \end{remark}

  Before we start with the rigorous proof of \cref{lem:technical-coms}, we explain the broad strategy. We provide a proof by contradiction, so we assume that there exists a set $I$ of times $t$ of nonzero (outer) measure for which there exists some direction $v_t$ with $\|v_t\| = 1$ contradicting \cref{eq:hilbert_space_coms}, and the limit on the right hand side of \cref{eq:hilbert_space_coms} existing in $\R$; by possibly replacing $v_t$ by $-v_t$, we may assume that the equality in \cref{eq:hilbert_space_coms} is replaced by a strict inequality ``$<$''. The proof then proceeds in two stages: First, we use a countable saturation argument: Writing 
  \begin{equation*}
    I = \bigcup_{T'_1 < T'_2 \in \Q} \bigg\{t \in I\colon \langle \dot{u}(t), v_t \rangle_{\mathcal{H}} < T'_1 < T'_2 < - \frac{\diff}{\diff \tau} F(u(t) + \tau v_t) \bigg|_{\tau = 0} \bigg\},
  \end{equation*}
  we deduce from the countable subadditivity of the Lebesgue measure the existence of some set $I' \subseteq I$ with nonzero (outer) measure and of $T_1, T_2 \in \R$ such that for all $t \in I'$, we even have
  \begin{equation*}
    \langle \dot{u}(t), v_t \rangle_{\mathcal{H}} < T_1 < T_2 < - \frac{\diff}{\diff \tau} F(u(t) + \tau v_t) \bigg|_{\tau = 0}.
  \end{equation*}
  The second part is more geometric and heavily leans on the Hilbert space structure (see \cref{fig:sketch-geometric-arg}): 
  Here, we first use the uncountability of $I'$ to select an accumulation point $t^* \in I'$ and a sequence $t_\varepsilon$ in $I'$ such that $t_\varepsilon \to t^*$ as $\varepsilon \to 0$. 
  Then we set, for each $\varepsilon$,  $w_{t_\varepsilon}^s = u(t_\varepsilon) + s v_{t_\varepsilon} - u(t^*)$. For $t_\varepsilon$ close to $t^*$, we then have that
  \begin{equation*}
    w_{t_\varepsilon}^0 = u(t_\varepsilon) - u(t^*) \approx (t_\varepsilon - t^*) \dot{u}(t_\varepsilon) \approx (t_\varepsilon - t^*) \dot{u}(t^*).
  \end{equation*}
  Using that $\frac{\diff}{\diff s} \|w_{t_\varepsilon}^s\|\, \Big|_{s=0} = \frac{1}{\| w_{t_\varepsilon}^0\|} \langle w_{t_\varepsilon}^0, v_{t_\varepsilon} \rangle$, and linearizing we thus have, for small $s > 0$,
  \begin{equation}  \label{eq:proof_sketch_technical_coms}
  \begin{split}
    \|\dot{u}(t^*)\| \, \big(\|w_{t_\varepsilon}^s\| - & \|w_{t_\varepsilon}^0\|\big) \approx s \langle \dot{u}(t_\varepsilon), v_{t_\varepsilon} \rangle_{\mathcal{H}} < s \, T_1 < s \, T_2 \\
    &< s \,\left( -\frac{\diff}{\diff \tau} F(u(t_\varepsilon) + \tau v_{t_\varepsilon}) \bigg|_{\tau = 0}\right) \approx F(u(t^*) +  w_{t_\varepsilon}^0) - F(u(t^*) + w_{t_\varepsilon}^s).
\end{split}
  \end{equation}
  Since $u$ is a curve of maximal slope with respect to $|\partial F|$, we lastly have that $\|\dot{u}(t^*)\| = |\partial F|(u(t^*))$ and $F(u(t^*)) - F(u(t^*) + w_{t_\varepsilon}^0) \approx |\partial F|(u(t^*)) \|w_{t_\varepsilon}^0\|$. Combining with \cref{eq:proof_sketch_technical_coms}, we get
  \begin{equation*}
    F(u(t^*)) - F(u(t^*) + w_{t_\varepsilon}^s) \approx |\partial F|(u(t^*)) \|w_{t_\varepsilon}^0\| + s \, T_2 > |\partial F|(u(t^*)) \|w_{t_\varepsilon}^0\| + s\, T_1 \approx |\partial F|(u(t^*)) \|w_{t_\varepsilon}^s\|.
  \end{equation*}
  By controlling the errors in the approximations above and sending $\varepsilon$ and $s$ to $0$ (with compatible convergence speeds), we construct a contradiction to the definition of $|\partial F|$ in the point $u(t^*)$.
  
\definecolor{curvegreen}{RGB}{20,72,66}
\definecolor{inkred}{RGB}{166,18,45}
\definecolor{inkblue}{RGB}{30,27,130}
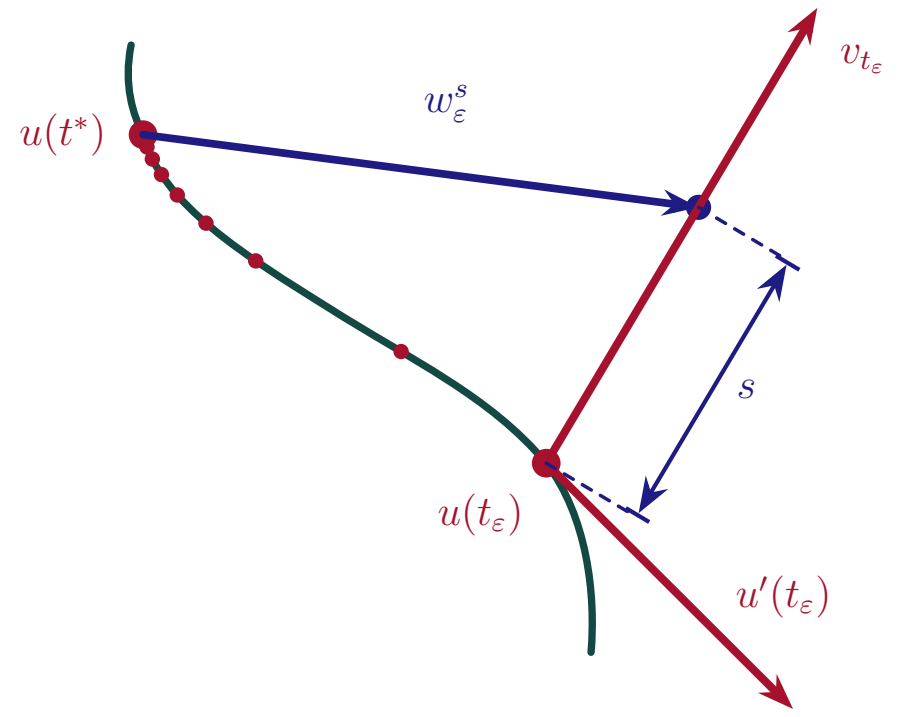
\begin{figure}[H]
    \centering
  \resizebox{0.4\textwidth}{!}{
  \begin{tikzpicture}[
    line cap=round,
    line join=round,
    >={Stealth[length=3.6mm,width=2.6mm]},
    every node/.style={font=\large}
]

    \coordinate (P0) at (-0.12,5.75);
    \coordinate (A) at (0,4.85);
    \coordinate (P1) at (1.85,3.12);
    \coordinate (B) at (4.05,1.55);
    \coordinate (Pend) at (4.50,-0.35);

    \coordinate (C) at (5.58,4.12);
    \coordinate (Vtip) at ($(B)!1.78!(C)$);
    \coordinate (Utip) at (6.53,-0.92);

    \path (A)
        .. controls (0.18,4.15) and (1.05,3.62)
        .. coordinate[pos=0.06] (d1)
           coordinate[pos=0.12] (d2)
           coordinate[pos=0.20] (d3)
           coordinate[pos=0.31] (d4)
           coordinate[pos=0.47] (d5)
           coordinate[pos=0.70] (d6)
        (P1);

    \path (P1)
        .. controls (2.70,2.58) and (3.48,2.23)
        .. coordinate[pos=0.30] (d7)
           coordinate[pos=0.63] (d8)
           coordinate[pos=0.90] (d9)
        (B);

    \path (B)
        .. controls (4.46,1.02) and (4.55,0.23)
        .. coordinate[pos=0.56] (d10) (Pend);

    \draw[curvegreen,line width=2.1pt]
        (P0)
        .. controls (-0.20,5.35) and (-0.10,5.05) .. (A)
        .. controls (0.18,4.15) and (1.05,3.62) .. (P1)
        .. controls (2.70,2.58) and (3.48,2.23) .. (B)
        .. controls (4.46,1.02) and (4.55,0.23) .. (Pend);

    \foreach \p in {d1,d2,d3,d4,d5,d6,d7}
        \fill[inkred] (\p) circle[radius=2.2pt];

    \fill[inkred] (A) circle[radius=4.1pt];
    \fill[inkred] (B) circle[radius=4.1pt];

    \draw[inkblue,line width=2.2pt,->]
        (A) -- (C);
    \fill[inkblue] (C) circle[radius=3.6pt];
    \node[inkblue,above=6pt] at (3.05,4.64) {$w^s_\varepsilon$};

    \draw[inkred,line width=2.3pt,->] (B) -- (Vtip);
    \draw[inkred,line width=2.3pt,->] (B) -- (Utip);

    \coordinate (Bo) at ($ (B)!30pt!-90:(C) $);
    \coordinate (Co) at ($ (C)!30pt! 90:(B) $);

    \draw[inkblue, dashed, line width=1pt] (B) -- (Bo);
    \draw[inkblue, dashed, line width=1pt] (C) -- (Co);

    \draw[inkblue, line width=1.2pt, |<->|]
        (Bo) -- (Co)
        node[midway, right=3pt] {$s$};
    \node[inkred,left=7pt] at (A) {$u(t^{*})$};
    \node[inkred,below left=5pt and 3pt] at (B) {$u(t_{\varepsilon})$};
    \node[inkred,below right=7pt and 3pt] at (Vtip) {$v_{t_{\varepsilon}}$};
    \node[inkred,above right=-3pt and 8pt] at (5.55,-0.02) {$u'(t_{\varepsilon})$};
\end{tikzpicture}
  }
  \caption{A sketch of the geometric argument in the proof of \cref{lem:technical-coms}.}
  \label{fig:sketch-geometric-arg}
\end{figure}
\begin{proof}[{Proof of \cref{lem:technical-coms}}]
  We will prove the claim by contradiction, so we start by assuming that there exists a non-negligible set $I \subseteq [0, \infty)$ such that for all $t \in I$, there exists $v_t$ such that $\|v_t\| = 1$,
  \begin{equation*}
    \langle \dot{u}(t), v_t \rangle_{\mathcal{H}} < - \frac{\diff}{\diff \tau} F(u(t) + \tau v_t) \Big|_{\tau = 0}.
  \end{equation*}
  where both sides exist; by changing the sign of $v_t$ if necessary, we can ensure the direction of the inequality.
  Because $u$ is a 2-curve of maximal slope (see~\cite[Definition 1.3.2]{ambrosioGradientFlows2008}), and because the metric slope is a weak upper gradient (see~\cite[Theorem 1.2.5]{ambrosioGradientFlows2008}), as in~\cite[Remark 1.3.3]{ambrosioGradientFlows2008}, we may assume the existence of a non-increasing map $\varphi$ such that $F \circ u = \varphi$ almost everywhere, and
  \begin{equation*}
    |u'|^2(t) = |\partial F|^2(u(t)) = - \varphi'(t) \qquad \text{for a.e. } t \in [0, \infty).
  \end{equation*}
  First, note that by intersecting $I$ with full-measure sets, we may assume that for all $t \in I$, the following hold:
  \begin{enumerate}
    \item $u$ is differentiable at $t$;
    \item $F(u(t)) = \varphi(t)$;
    \item $\varphi$ is differentiable at $t$;
    \item\label{item:technical-proof-4} $|u'|^2(t) = |\partial F|^2(u(t)) = - \varphi'(t)$.
  \end{enumerate}
  Next, we will use a countable saturation argument to find a subset $J \subseteq I$, two numbers $T_1 < T_2$ and $\delta > 0$ such that
  \begin{equation} \label{eq:coms-technical-after-saturation}
    \frac{\langle u(t + \tau) - u(t), v_t \rangle}{\tau} < T_1 < T_2 < \frac{F(u(t)) - F(u(t) + \tau \cdot v_t)}{\tau} \qquad \forall t \in J, \tau \in (-\delta, \delta).
  \end{equation}
  To this end, we first choose $T_1, T_2$ by noting that
  \begin{equation*}
    I = \bigcup_{T'_1 < T'_2 \in \Q} \underbrace{\bigg\{t \in I\colon \langle \dot{u}(t), v_t \rangle_{\mathcal{H}} < T'_1 < T'_2 < - \frac{\diff}{\diff \tau} F(u(t) + \tau v_t) \bigg|_{\tau = 0} \bigg\}}_{\eqqcolon I_{T'_1, T'_2}}.
  \end{equation*}
  Since the union over the countable $I_{T'_1, T'_2}$ has non-zero outer measure, there must exist $T_1 < T_2 \in \Q$ such that $I_{T_1, T_2} \eqqcolon J'$ has non-zero outer measure as well\footnote{Here, we resort to the outer measure because the sets $I_{T'_1, T'_2}$ are not guaranteed to be measurable. The saturation argument works because the outer measure is countably subadditive.}. Next, note that for each $t \in J'$, we have
  \begin{equation*}
    \lim_{\tau \to 0} \frac{\langle u(t + \tau) - u(t), v_t \rangle}{\tau} < T_1 < T_2 < \lim_{\tau \to 0}\frac{F(u(t)) - F(u(t) + \tau \cdot v_t)}{\tau}.
  \end{equation*}
  Thus, for each $t \in  J'$, there exists $\delta_t \in (0, t)$ such that
  \begin{equation*}
    \frac{\langle u(t + \tau) - u(t), v_t \rangle}{\tau} < T_1 < T_2 < \frac{F(u(t)) - F(u(t) + \tau \cdot v_t)}{\tau} \qquad \forall \tau \in (-\delta_t, \delta_t).
  \end{equation*}
  Finally, we repeat the countable union argument by noting that
  \begin{equation*}
    J' = \bigcup_{0 < \delta' \in \Q} \underbrace{\Big\{t \in J'\colon \delta_t > \delta'\Big\}}_{\eqqcolon J'_{\delta'}}.
  \end{equation*}
  As before, there must exist some $\delta$ such that $J'_{\delta} \eqqcolon J$ has positive outer measure. \\
  In particular, $J$ is uncountable, hence there exists a sequence $(t_i)_{i \in \NN}$ in $J$ and $t^*$ in $J$ such that $t_i \searrow t^*$ as $i \to \infty$. By \cref{item:technical-proof-4} above, we have that $|u'|(t^*) \cdot |\partial F|(u(t^*)) = -\varphi'(t^*)$ and $|u'|(t^*) = |\partial F|(u(t^*)) \eqqcolon T > 0$\footnote{The case where $T = 0$ immediately yields a contradiction, as in this case the directional derivative of $F$ at $t^*$---if it exists---cannot be anything other than $0$, and likewise $\langle \dot{u}(t^*), v\rangle = \langle 0, v \rangle = 0$ for all $v \in \mathcal{H}$.} and thus, for each $\varepsilon > 0$, there exists $\alpha_\varepsilon \in (0, \varepsilon]$ such that, for every $\tau \in (0, \alpha_\varepsilon)$, we have
  \begin{align}
    \frac{\varphi(t^*) - \varphi(t^* + \tau)}{\tau} \cdot (1 + \varepsilon) \geq \frac{\|u(t^*) - u(t^*+\tau)\|}{\tau} \cdot \sup_{\substack{x \in B_\tau(u(t^*))\\x \neq u(t^*)}}\frac{F(u(t^*)) - F(x)}{\|u(t^*) - x\|}, \label{eq:technical-product} \\
    \frac{\|u(t^* + \tau) - u(t^*)\|}{\tau} \leq (1 + \varepsilon) T, \qquad\frac{\|u(t^* + \tau) - u(t^*)\|}{\tau} \geq \frac{1}{1 + \varepsilon} T.  \label{eq:technical-limit-stuff}
  \end{align}
  Likewise, for each $\epsilon$, there exists $\beta_\epsilon$ such that
  \begin{equation*}
    \frac{F(u(t^*)) - F(x)}{\|u(t^*) - x\|} \leq (1 + \epsilon) T, \qquad \forall x \in B_{\beta_\epsilon}(u(t^*)) \setminus \{u(t^*)\}.
\end{equation*}
  Rewriting \cref{eq:technical-product}, we get
  \begin{equation*}
    \frac{\varphi(t^*) - \varphi(t^* + \tau)}{\|u(t^*) - u(t^*+\tau)\|} \cdot (1 + \varepsilon) \geq \sup_{x \in B_\tau(u(t^*))}\frac{F(u(t^*)) - F(x)}{\|u(t^*) - x\|} \geq T, \qquad \forall \tau \in (0, \alpha_\varepsilon).
  \end{equation*}
  Next, for each $\varepsilon > 0$, we may pick some $t^\varepsilon \in J$ from the sequence $(t_i)_{i \in \NN}$ such that $t^\varepsilon \in (t^*, t^* + \alpha_\varepsilon) \cap (t^*, t^* + \delta)$, and we set $\Delta_\varepsilon \coloneqq t^\varepsilon - t^* > 0$. By using \cref{eq:coms-technical-after-saturation} with $\tau=-\Delta_\varepsilon$, we have that
  \begin{equation*}
    \frac{\langle u(t^*) - u(t^\varepsilon), v_{t^\varepsilon} \rangle}{-\Delta_\varepsilon} = \frac{\langle   u(t^\varepsilon) - u(t^*), v_{t^\varepsilon} \rangle}{\Delta_\varepsilon} < T_1,
  \end{equation*}
  and thus, setting $w_\varepsilon^s \coloneqq u(t^\varepsilon)+ s v_{t^\varepsilon} - u(t^*) $, we deduce
  \begin{equation*}
    \|w_\varepsilon^s\|^2 = \|w_\varepsilon^0\|^2 + 2s \langle u(t^\varepsilon) - u(t^*), v_{t^\varepsilon} \rangle + \|s v_{t^\varepsilon}\|^2 \leq \|w_\varepsilon^0\|^2 + 2s T_1 \Delta_\varepsilon + s^2.
  \end{equation*}
  Moreover, we calculate, for $s < \delta$,
  \begin{equation*}
    F\big( u(t^\varepsilon) \big) - F\big( u(t^*) + w_\varepsilon^s \big) = F\big( u(t^\varepsilon) \big) - F\big( u(t^\varepsilon) + s v_{t^\varepsilon} \big) \geq T_2 \cdot s.
  \end{equation*}
  Thus, we have, for any $\epsilon > 0$, $\varepsilon > 0$ and $s \in (0, \delta)$ such that $\|w_\varepsilon^0\| + s < \beta_\epsilon$, $s |T_2| \leq T^2 \Delta_\varepsilon \frac{1}{(1 + \epsilon)^2 }$ and $\varepsilon \leq \epsilon$,
  \begin{equation} \label{eq:some_inequality_techical_lem}
    \begin{split}
    (1+\epsilon) T &\geq \frac{F\big( u(t^*) \big) - F\big( u(t^*) + w_\varepsilon^s\big )}{\|w_\varepsilon^s\|} = \frac{\Big(F(u(t^*)) - F\big( u(t^\varepsilon) \big) \Big) + \Big( F\big( u(t^\varepsilon) \big) - F\big( u(t^*) + w_\varepsilon^s \big)\Big)}{\|w_\varepsilon^s\|} \\
    &\geq \frac{T\frac{\|w_\varepsilon^0\|}{1+\epsilon} + s T_2}{\sqrt{\|w_\varepsilon^0\|^2 + 2sT_1 \Delta_\varepsilon + s^2}} \geq \frac{T^2 \Delta_\varepsilon \frac{1}{(1+\epsilon)^2} + s T_2}{\sqrt{T^2 \Delta_\varepsilon^2(1 + \epsilon)^2 + 2 s T_1 \Delta_\varepsilon + s^2}}
    \end{split}
  \end{equation}
  Squaring and rearranging the above, we obtain
  \begin{align*}
    (1 + \epsilon)^4T^4 \Delta_\varepsilon^2 + (1 + \epsilon)^2 T^2 2 s T_1 \Delta_\varepsilon + (1+\epsilon)^2T^2s^2 \geq T^4 \Delta_\varepsilon^2 \frac{1}{(1+\epsilon)^4} + 2 T^2 \Delta_\varepsilon \frac{1}{(1+\epsilon)^2} s T_2 + s^2 T_2^2
  \end{align*}
  Setting $a_\epsilon \coloneqq (1 + \epsilon)^4 - \frac{1}{(1+\epsilon)^4}$, $b_\epsilon \coloneqq \frac{T_2}{(1 + \epsilon)^2} - (1 + \epsilon)^2 T_1$, and $c_\epsilon \coloneqq (1+\epsilon)^2 T^2 - T_2^2$, we may rearrange the above to the following inequality
  \begin{equation*}
    T^4 \Delta_\varepsilon^2 a_\epsilon - 2 s T^2 \Delta_\varepsilon b_\epsilon + s^2 c_\epsilon \geq 0.
  \end{equation*}
  We aim to show that, in our situation, the above inequality cannot hold for all $\epsilon > 0$, $\varepsilon > 0$ and $s \in (0, \delta)$ such that $\|w_\varepsilon^0\| + s < \beta_\epsilon$, $s |T_2| \leq T^2 \Delta_\varepsilon \frac{1}{(1 + \epsilon)^2}$ and $\varepsilon \leq \epsilon$. To this end, we first rewrite the inequality (note that $\Delta_\varepsilon > 0$ by construction) to obtain
  \begin{equation*}
    T^4 a_\epsilon - \frac{s}{\Delta_\varepsilon} 2 T^2 b_\epsilon + \bigg(\frac{s}{\Delta_\varepsilon}\bigg)^2 c_\epsilon \geq 0.
  \end{equation*}
  Since $a_\epsilon \to 0$, $b_\epsilon \to T_2 - T_1 > 0$ and $c_\epsilon \to T^2 - T_2^2$, we may first choose $\epsilon$ small enough such that the quadratic equation $x \mapsto T^4 a_\epsilon - x 2T^2 b_\epsilon + x^2  c_\epsilon$ is negative in some interval $(x_{-}, x_{+})$ with $x_{-} |T_2| < \frac{T^2}{(1 + \epsilon)^2}$; we then choose $x^\star \in (x_{-}, x_{+})$ such that $x^\star |T_2| \leq \frac{T^2}{(1 + \epsilon)^2}$. As furthermore $\Delta_\varepsilon \to 0$ as $\varepsilon \to 0$, we may next choose $\varepsilon$ small enough such that $\|w_\varepsilon^0\| + \Delta_\varepsilon x^\star < \beta_\epsilon$, $\Delta_\varepsilon x^\star < \delta$ and $\varepsilon \leq \epsilon$. Then, setting $s \coloneqq \Delta_\varepsilon x^\star$, we have that $s |T_2| \leq T^2 \Delta_\varepsilon\frac{1}{(1 + \epsilon)^2}$ (hence we are in the regime where \cref{eq:some_inequality_techical_lem} holds and our calculations are valid), and we have found a contradiction.
\end{proof}
    We are now ready to prove the main result of this section:
  \begin{proposition}\label{prop:weak-formulation-abstract}
    Let $E\colon\Pp_2 \to \R \cup \{+\infty\}$, let $\bgamma\colon [0, \infty) \to L^2$ be a curve of maximal slope for the composition $\tilde{E} \coloneqq E \circ \law$ with respect to the metric slope $|\partial (E \circ \law)|$. For $\varphi \in C_c^\infty(\R^{d+1})$, $t \in [0, \infty)$, and $s \in \R$, we define the map $T_s^{\varphi_t} \colon \R^d \to \R^d$ as
    \begin{equation*}
      T_s^{\varphi_t}(x) = x + s \nabla_x \varphi_t(x),
    \end{equation*}
    and the measure $\mu_t \coloneqq \law\big(\bgamma_t\big)$. Then, for a.e.~$t \in [0, \infty)$ and every $\varphi \in C_c^\infty(\R^{d+1})$, we have
    \begin{equation} \label{eq:weak_formulation_pde}
      \frac{\diff}{\diff t} \int \varphi_t(x) \,\diff \mu_t(x) = \int \partial_t \varphi_t(x) \,\diff \mu_t(x) - \frac{\diff}{\diff s} E\big((T_s^{\varphi_t})_\# \mu_t\big) \Big|_{s = 0},
    \end{equation}
    whenever the right hand side is well defined. Equivalently,
    \begin{equation*}
      \frac{\diff}{\diff t}\EE_{\mathbb{P}}[\varphi_t(\bgamma_t)] = \EE_{\mathbb{P}}[\partial_t \varphi_t(\bgamma_t)] - \frac{\diff}{\diff s} \tilde{E}\big(T_s^{\varphi_t} \circ \bgamma_t\big) \Big|_{s = 0},
    \end{equation*}
    whenever the right hand side is well defined.
  \end{proposition}
  \begin{remark}
    In the above, we can already see that whenever we can express $\frac{\diff}{\diff s} E\big((T_s^{\varphi_t})_\# \mu_t\big) \Big|_{s = 0}$ in terms of derivatives of $\varphi$ integrated along $\mu_t$, we identify $\law \circ \bgamma$ as the weak solution of a certain PDE\@. In particular, the weak formulation of the term~\eqref{eq:weak_formulation_pde} becomes evident for energies $E$ acting on $\mu$ in an integral manner, as in \cref{sec:ex_energy_through_density}.
  \end{remark}
  
    \begin{proof}[Proof of \cref{prop:weak-formulation-abstract}]
      As we announced earlier, we will use the product rule in the duality pairing $((C_0^2)', C_0^2)$ of the two curves $G \circ \bgamma$ and $\varphi_\bullet$. By \cref{prop:frechet-derivative-L2-W2}, $G \circ \bgamma$ is absolutely continuous (as the composition of a Lipschitz continuous map with an absolutely continuous curve), and we have for a.e.~$t$:
      \begin{align}
        \frac{\diff}{\diff t} \int \varphi_t(x) \,\diff \mu_t(x) &=  \frac{\diff}{\diff t} \langle G(\bgamma_t), \varphi_t \rangle_{(C_0^2)' \times C_0^2} \label{eq:continuity equation_-1}\\
        &= \langle DG(\bgamma_t)[\dot{\bgamma}_t], \varphi_t \rangle_{(C_0^2)' \times C_0^2} +  \langle G(\bgamma_t), \partial_t \varphi_t \rangle_{(C_0^2)' \times C_0^2} \label{eq:continuity equation_0}\\
        &= \big\langle \dot{\bgamma}_t, DG(\bgamma_t)_*[\varphi_t] \big\rangle_{L^2} + \int \partial_t \varphi_t(x) \,\diff \mu_t(x) \label{eq:continuity equation}\\
        &\overset{(*)}{=} - \frac{\diff}{\diff s} E\Big(\law \big( \bgamma_t + s \, DG(\bgamma_t)_*[\varphi_t] \big) \Big) \Big|_{s = 0} + \int \partial_t \varphi_t(x) \,\diff \mu_t(x), \notag
      \end{align}
      where we used \cref{lem:technical-coms} in $(*)$.
      Next, we have to identify $\law(\bgamma_t + s \, DG(\bgamma_t)_*[\varphi_t])$. To this end, we choose any $\psi \in C_b^0(\R^d)$ and calculate:
      \begin{align*}
        \int \psi(x) \,\diff ((\bgamma_t + s \, DG(\bgamma_t)_*[\varphi_t] )_\# \mathbb{P})(x) &= \int \psi(\bgamma_t(\omega) + s \nabla_x \varphi_t(\bgamma_t(\omega))) \,\diff \mathbb{P}(\omega) \\
        &= \int \psi(T_s^{\varphi_t}(x)) \,\diff ((\bgamma_t)_\# \mathbb{P})(x) = \int \psi(x) \,\diff \big((T_s^{\varphi_t})_\# \mu_t\big)(x),
      \end{align*}
      which concludes the proof.
    \end{proof}
    
    The following remark aims at closing this part by reconnecting to the discussion about the continuity equation reported at the beginning of the present section.
    
    \begin{remark}\label{rmk:CE_X}
      Using \cref{eq:continuity equation_-1,eq:continuity equation_0,eq:continuity equation,eq:pre-adjoint-L2-W2}---which do not depend on $\bgamma$ being a curve of maximal slope---, we get
      \begin{equation*}
        \frac{\diff}{\diff t}\EE_{\mathbb{P}}[\varphi_t(\bgamma_t)] = \EE_{\mathbb{P}}\left[\nabla_x \varphi_t(\bgamma_t) \cdot \dot{\bgamma}_t\right] + \EE_{\mathbb{P}}[\partial_t \varphi_t(\bgamma_t)],
      \end{equation*}
      for arbitrary $\varphi \in C_c^\infty(\R^{d+1})$.
      Thus, choosing any velocity field $v_\bullet\colon \R^d \to \R^d$ such that $\EE\left[\dot{\bgamma}_t \mid \bgamma_t\right] = v_t \circ \bgamma_t$, we equivalently have
      \begin{equation*}
      \begin{split}
        \frac{\diff}{\diff t} \int \varphi_t(x) \,\diff \mu_t(x) &= \big\langle \dot{\bgamma}_t, \nabla \varphi_t \circ \bgamma_t\rangle_{L^2} + \int \partial_t \varphi_t(x) \,\diff \mu_t(x) \\
        &= \int \nabla \varphi_t(x) \cdot v_t(x) \,\diff \mu_t(x) + \int \partial_t \varphi_t(x) \,\diff \mu_t(x), 
        \end{split}
      \end{equation*}
      which is precisely the weak formulation of the continuity equation
      \begin{equation*}
        \partial_t \mu_t + \nabla \cdot (\mu_t v_t) = 0.
      \end{equation*}
    \end{remark}

\section{The Mirroring in Action: Examples} \label{sec:examples}

In this section, we provide two examples showing the flexibility gained by switching between $\Pp_2$ and $L^2$: one pedagogical sanity check, for energies defined through densities, and one showing novel results for the nonlocal attraction-repulsion transport potential.

\subsection{Example I: Energies Defined on Densities} \label{sec:ex_energy_through_density}
In \cref{prop:weak-formulation-abstract}, we have shown an abstract formula that can be used to identify weak formulations of PDEs associated to curves of maximal slope in $L^2$. 
To make this more concrete and as a sanity check, we explore how this theory applies to the case where the functional $E\colon\Pp_2 \to \R \cup \{+\infty\}$ is defined on densities, i.e., there exists a function $F\colon \R^d \times [0, \infty) \to \R \cup \{+\infty\}$ such that for any $\mu \in \Pp_2$, $E$ is defined as
\begin{equation} \label{eq:formula_internal_energy}
  E(\mu) = \begin{cases}
    \int F\big(x, \rho(x)\big) \, \diff x, & \text{if } \mu \ll \Leb^d, \ \frac{\diff \mu}{\diff \Leb^d} = \rho , \\
    +\infty, & \text{otherwise}.
  \end{cases}
\end{equation}
 Our goal in this example is to show that the Lions lift yields precisely the same results as the purely metric construction of gradient flows in $\Pp_2$, even when $E$ cannot be directly reformulated as an integral functional of random variables in $L^2_{\mathbb R}$. 
 To this end, following the path paved by \cref{lem:technical-coms} and \cref{prop:weak-formulation-abstract}, we have to compute the derivative $\frac{\diff}{\diff s} E\big((T_s^{\varphi_t})_\# \mu_t\big) \big|_{s = 0}$.
 \\
 
 We start by computing the density of $(T_s^{\varphi_t})_\# \mu_t$, where we use the following standard result (see, e.g., \cite[Theorem 3.7.1]{bogachevMeasureTheory2007} or \cite[Lemma 5.5.3]{ambrosioGradientFlows2008}) without providing a proof:
 
\begin{proposition}
  Let $\mu = \rho \, \Leb^d$ be a measure absolutely continuous w.r.t.~the Lebesgue measure, and let $T\colon \R^d \to \R^d$ be a diffeomorphism. Then, the pushforward measure $T_\# \mu$ is absolutely continuous w.r.t.\ the Lebesgue measure, with density
  \begin{equation*}
    \rho_T(x) = \frac{\rho(T^{-1}(x))}{|\det(DT(T^{-1}(x)))|}.
  \end{equation*}
\end{proposition}

Note that, importantly in our case, $T_s^{\varphi_t}$ is indeed a diffeomorphism for $s$ small enough, since $\varphi$ is smooth and compactly supported. 
Furthermore, for $s\to 0$, we have 
$$|\det DT_s^{\varphi_t}(x)| = |\det(I + s \nabla_x^2 \varphi_t(x))| = 1 + s \Delta \varphi_t(x) + O(s^2) \eqqcolon \kappa_s(x),$$
where $I$ is the identity matrix and the coefficients of the $O(s^2)$ error have a global upper bound independent of $x$. 
Thus, we have, for $s$ small enough
\begin{align*}
  E\big((T_s^{\varphi_t})_\# \mu_t\big) &= \int F\Bigg(x, \frac{\rho_t\big((T_s^{\varphi_t})^{-1}(x)\big)}{1 + s \Delta \varphi_t((T_s^{\varphi_t})^{-1}(x)) + O(s^2)}\Bigg) \diff x \\
  &= \int \underbrace{F\Bigg(T_s^{\varphi_t}(x), \frac{\rho_t\big(x \big)}{1 + s \Delta \varphi_t(x) + O(s^2)}\Bigg) \underbrace{(1 + s \Delta \varphi_t(x) + O(s^2))}_{= \kappa_s(x)}}_{\eqqcolon\mathcal{F}_s(x)} \diff x,
\end{align*}
where $F\colon \R^d \times [0, \infty) \to \R \cup \{+\infty\}$ is the function involved in the definition of the energy $E$ (see \cref{eq:formula_internal_energy}).
Finally, if we can differentiate under the integral sign, we have
\begin{equation*}
  \frac{\diff}{\diff s} E\big((T_s^{\varphi_t})_\# \mu_t\big) \Big|_{s = 0} = \int \frac{\diff}{\diff s} \mathcal{F}_s(x) \Big|_{s = 0} \diff x.
\end{equation*}
Hence, we calculate the derivative inside the integral:
\begin{equation} \label{eq:derivative_of_internal_energy}
  \begin{split}
  \frac{\diff}{\diff s} \mathcal{F}_s(x) = &\nabla_x F\left(T_s^{\varphi_t}(x), \frac{\rho_t(x)}{\kappa_s(x)}\right) \cdot (\partial_s T_s^{\varphi_t}) \kappa_s(x) + F\left(T_s^{\varphi_t}(x), \frac{\rho_t(x)}{\kappa_s(x)}\right) \partial_s \kappa_s(x) \\
  &- \partial_\rho F\left(T_s^{\varphi_t}(x), \frac{\rho_t(x)}{\kappa_s(x)}\right) \frac{\rho_t(x)}{\kappa_s(x)} \partial_s \kappa_s(x),
  \end{split}
\end{equation}
and thus we obtain:
\begin{equation*}
  \frac{\diff}{\diff s} \mathcal{F}_s(x) \Big|_{s = 0} = \nabla_x F(x, \rho_t(x)) \cdot \nabla_x \varphi_t(x) + \Big(F(x, \rho_t(x)) - \rho_t(x) \partial_\rho F(x, \rho_t(x))\Big) \Delta \varphi_t(x).
\end{equation*}
To justify the differentiation under the integral, we demand the following assumption:

\begin{assumption} \label{ass:differentiation_below_integral}
  The function $F\colon \R^d \times [0, \infty) \to \R \cup \{+\infty\}$ (see \cref{eq:formula_internal_energy}) is finite, continuously differentiable on $\R^d \times (0, \infty)$, and $F(x, 0) = 0$ for all $x \in \R^d$. 
  For every $\rho \in L^1(\R^d)$ such that $\int_{\R^d} |F(x, \rho(x))| \diff x < +\infty$, we have that
  \begin{equation*}
    \int_{\R^d} |\rho(x) \partial_\rho F(x, \rho(x))| \diff x < +\infty.
  \end{equation*}
  Furthermore, there exists $C > 0$ such that, for all $x, y \in \R^d$ with $|x - y| \leq 1$, all $r > 0$ and all $\lambda \in [\frac{1}{2}, 2]$, we have that 
  \begin{equation*}
    |\nabla_x F(x, r)| \leq C(1 + r), \qquad |F(y, \lambda r)| + |\lambda r \, \partial_\rho F(y, \lambda r)| \leq C\big(|F(x, r)| + r + 1\big).
  \end{equation*}
\end{assumption}

Whenever \cref{ass:differentiation_below_integral} is met, taking the derivative inside the integral is justified by dominated convergence whenever $\int |F(x, \rho_t(x))| \diff x < +\infty$. To show this, we want to apply dominated convergence to the finite differences; by the mean value theorem, it suffices to dominate the right hand side in \cref{eq:derivative_of_internal_energy} by an $L^1$ function on $\R^d \times [-\delta, \delta]$, for some $\delta > 0$. Outside the compact support $K=\mathrm{supp}(\varphi_t)$ of $\varphi_t$, $\mathcal{F}_s(x) = F(x, \rho_t(x))$ does not depend on $s$, so the right hand side in \cref{eq:derivative_of_internal_energy} is $0$ there. Furthermore, for $\delta$ small enough, we have that $\kappa_s(x) \in [\frac{1}{2}, 2]$, that $|T_{s}^{\varphi_t}(x) - x| \leq 1$, and that $\kappa^{-1}_s(x)$, $\partial_s \kappa_s(x)$ and $\delta_s T_s^{\varphi_t}(x)$ are bounded by some $D$, uniformly over $(x, s) \in \R^d \times [-\delta, \delta]$. Thus, using \cref{ass:differentiation_below_integral}, we can domiate the right hand side in \cref{eq:derivative_of_internal_energy} by
\begin{align}
  \left|\frac{\diff}{\diff \sigma} \mathcal{F}_\sigma(x) \Big|_{\sigma=s}\right| \leq \chi_K \, C(1 + \rho_t(x) D) D^2 + C(|F(x, \rho_t(x))| + D + 1),
\end{align}
which is in $L^1$ because $\rho_t$ and $|F(\cdot, \rho_t(\cdot)|$ are, and $K$ is compact.
\\
Putting everything together, we have shown the following:
\begin{proposition} \label{prop:weak-formulation-densities}
  Let $E\colon\Pp_2 \to \R \cup \{+\infty\}$ be a functional defined on densities as in \cref{eq:formula_internal_energy} and fulfilling \cref{ass:differentiation_below_integral}, and let $\bgamma\colon [0, \infty) \to L^2$ be a curve of maximal slope for $E \circ \law$ with respect to the metric slope $|\partial (E \circ \law)|$. 
  For any $t \in [0, \infty)$, we define the measure $\mu_t \coloneqq \law\big(\bgamma_t\big)$, and we assume that $\mu_t = \rho_t \, \Leb^d$ for some density $\rho_t$. 
  Then, for a.e.~$t \in [0, \infty)$, we have
  \begin{equation*}
  \begin{split}
    \frac{\diff}{\diff t} \int \varphi_t(x) \rho_t(x) \diff x &= \int \partial_t \varphi_t(x) \rho_t(x) \diff x - \int \nabla_x F(x, \rho_t(x)) \cdot \nabla_x \varphi_t(x) \diff x \\
    &\qquad - \int \Big(F(x, \rho_t(x)) - \rho_t(x) \partial_\rho F(x, \rho_t(x))\Big) \Delta_x \varphi_t(x) \diff x,
    \end{split}
  \end{equation*}
  for all $\varphi \in C^\infty_c(\R^{d+1})$.
\end{proposition}
\begin{example} \label{ex:bolzmann-coms}
  Let $F(x, \rho) = \rho \log \rho$, so that $E$ is the negative Boltzmann entropy. Then, we have $\nabla_x F(x, \rho) = 0$ and $\rho \partial_\rho F(x, \rho) = \rho \log \rho + \rho$. Thus, $F$ fulfills \cref{ass:differentiation_below_integral} and by \cref{prop:weak-formulation-densities}, any curve of maximal slope $\mu_\bullet$ for $E$ satisfies for a.e.~$t$:
  \begin{equation*}
    \frac{\diff}{\diff t} \int \varphi_t(x) \rho_t(x) \diff x = \int \partial_t \varphi_t(x) \rho_t(x) \diff x + \int \rho_t(x) \Delta \varphi_t(x) \diff x, \qquad \forall \varphi \in C_c^\infty(\R^{d+1}),
  \end{equation*}
  which is the weak formulation of the heat equation $\partial_t \rho_t = \Delta \rho_t$.
\end{example}

These results, as a sanity check, do confirm that the PDE description of gradient flow in $L^2_{\mathbb R}$ according to \cref{lem:technical-coms} and \cref{prop:weak-formulation-abstract} reproduces and is fully compatible with  the well-known metric approach to gradient flows in $\Pp_2$ as in \cite{ambrosioGradientFlows2008}.

\subsection{Example II: The Nonlocal Attraction-Repulsion Transport Potential} \label{sec:ex_nonlocal}
In this section, we will use the framework developed in earlier sections to investigate the $W_2$-gradient flow along the nonlocal attraction-repulsion transport potential $E\colon \Pp_2 \to \R$ defined as 
\begin{equation} \label{eq:definition_of_potential}
  E(\mu) \coloneqq -\frac{1}{2} \iint K_s(x-y) \, \diff (\mu-\nu)(x) \, \diff (\mu-\nu)(y).
\end{equation}
with power kernels
\begin{equation*}
  K_s(z) = |z|^{s+1}, \qquad s \in [0, 1),
\end{equation*}
and a reference signed measure
\begin{equation} \label{eq:nu}
  \nu = f \, \Leb^d, \qquad f \in L^\infty(\R^d) \cap L^1(\R^d), \qquad \int_{\R^d} |y|^{s+1} |f(y)| \, \diff y < +\infty.
\end{equation}
In the special case of $\nu = 0$, \cref{eq:definition_of_potential} reduces to the repulsion potential
\begin{equation*}
  \mu \mapsto -\frac{1}{2}  \iint K_s(x-y) \,\diff \mu(x) \, \diff \mu(y)
\end{equation*}

The energy \eqref{eq:definition_of_potential} has its roots in Riesz potential theory \cite{Riesz1938} and rests on
Schoenberg's result that $|x-y|^q$ is conditionally negative definite for $0<q<2$
\cite{Schoenberg1938a,Schoenberg1938b}, which makes it nonnegative and embeds it in a Hilbert
space as long as $\nu$ is also a probability measure; Deny and Fuglede later gave it a Fourier-analytic, negative-Sobolev representation
\cite{Deny1950,Fuglede1960}. It was independently rediscovered by Sz\'ekely as the basis of
energy statistics \cite{SzekelyRizzoBakirov2007,SzekelyRizzo2013,SzekelyRizzo2023}, and in
machine learning as the {\it Maximum Mean Discrepancy} (MMD) with a Riesz-type kernel \cite{Gretton2006,Gretton2012}, the two
viewpoints being shown equivalent in \cite{Sejdinovic2013}. Around 2010 it reappeared in image
processing for electrostatic halftoning \cite{Schmaltz2010,Teuber2011}, whose empirical
minimizers were proved consistent in \cite{FornasierHutter2016}; we refer to
\cite{ColasantoFocardiFornasierMattesini2026} for a detailed historical account.

The formal Wasserstein gradient flow of $E$ is the continuity equation
\begin{equation} \label{eq:WGF}
  \partial_t \mu_t + \nabla \cdot (\mu_t v_t) = 0,
  \qquad
  v_t = -\nabla \frac{\delta E}{\delta \mu}(\mu_t) = \nabla K_s * (\mu_t - \nu),
\end{equation}
that is,
\begin{equation*}
  \partial_t \mu_t = -\nabla \cdot \Big( \mu_t \int_{\R^d} \nabla K_s(\cdot - y) \,\diff (\mu_t - \nu)(y) \Big),
\end{equation*}
where
\begin{equation*}
  \nabla K_s(z) = (s+1)\,|z|^{s-1} z \quad (z \neq 0), \qquad \nabla K_s(0) := 0.
\end{equation*}
The energy $E$ is in general \emph{not} $\lambda$-convex on $\Pp_2$ as was noted in \cite{HertrichGrafBeinertSteidl2024,HertrichWaldAltekrugerHagemann2024}. The lack of $\lambda$-convexity prevents a straightforward application of the construction of curves of maximal slope, i.e., the gradient flow, in $\Pp_2$ via the minimizing movement scheme in purely metric form as in \cite{ambrosioGradientFlows2008}. 
Yet, recently the papers \cite{RosenzweigSlepcevWang2026,fornasierNonlocalAttractionrepulsionTransport2026} obtain well-posedness results for \eqref{eq:WGF} in the class of Lagrangian solutions with bounded  moments  by means of suitable (square radius) regularizations of the potentials $K_s$,  uniform a priori moment estimates, and compactness arguments. 

In what follows, we show that existence of distributional solutions can be  recovered in terms of the usual construction by minimizing movement scheme, by leveraging the Lions lift. As already mentioned above, $\lambda$-convexity on $\Pp_2$ is not available, while, through the Lions lift, we gain the Hilbert space structure of $L^2_{\mathbb{R}}$, and properties and notions there, such as  more classical convexity and subdifferentials, which do allow us to establish more easily existence of gradient flow lines. 
Indeed, it turns out that $- \tilde E \coloneqq -  E \circ \law\colon L^2 \to \R$ is a $C^1$ perturbation of a convex function, as we will clarify in this section. 
Eventually, we will prove the following theorem:

\begin{theorem} \label{thm:cusp}
  For every $X_0 \in \LtwoR$, there exists a curve $\bgamma\colon [0, +\infty) \to \LtwoR$, which is the $\law^{-1}(\sigma_{W_2})$-limit of piecewise constant interpolations $(\bar{U}_{\tau_i})_{i \in \NN}$ of MMS sequences\footnote{\label{footnote:recall_mms}We recalled the piecewise constant interpolation in \cref{eq:constant_interpolation} and the MMS scheme in~\eqref{eq:mms}; here, we use the cases $(\mathcal{S}, d) = (\LtwoR, d_{\|\cdot\|})$ (in \cref{thm:cusp}) and $(\mathcal{S}, d) = (\Pp_2, W_2)$ (in \cref{cor:cusp}).} and has the following properties:
  \begin{enumerate}
    \item $\bgamma$ is a curve of maximal slope for the functional $\tilde{E}$ with respect to the strong upper gradient $\mathrm{lsc}(|\partial \tilde{E}|)$, which factors through $\law$.
    \item $\mu_\bullet \coloneqq \law \circ \bgamma$ fulfills the following PDE:
  \begin{equation} \label{eq:PDE_0}
        \partial_t \mu_t + \nabla \cdot(\mu_t v_t) = 0, \text{ where } v_t(x) =  \int_{y \neq x} (s+1) (x-y)|x-y|^{s-1} \, \diff \mu_t(y) - \nabla V(x),
  \end{equation}
  in distributional sense, where $V = K_s \ast \nu$ is the background potential.
  \end{enumerate}
\end{theorem}

\begin{remark} \label{rmk:notonlyprob}
    It is perhaps interesting to highlight that the distributional existence result of \cref{thm:cusp} holds for every given signed measure $\nu$ complying with \eqref{eq:nu}, while the well-posedness results of Lagrangian solutions in \cite{RosenzweigSlepcevWang2026,fornasierNonlocalAttractionrepulsionTransport2026} require that $\nu$ is a probability measure.
\end{remark}

By the equivalence in \cref{lem:summarize_sec_1} we can directly deduce from \cref{thm:cusp} the following:
\begin{corollary}\label{cor:cusp}
  For every $\mu_0 \in \Pp_2(\R^d)$, there exists a curve $\mu_\bullet\colon [0, +\infty) \to \Pp_2(\R^d)$, which is the $W_2$-limit of piecewise constant interpolations $(\mu^{\tau_i}_\bullet)_{i \in \NN}$ of MMS sequences\footref{footnote:recall_mms} and has the following properties:
  \begin{enumerate}
    \item $\mu_\bullet$ is a curve of maximal slope for the functional $E$ with respect to the strong upper gradient $\mathrm{lsc}(|\partial E|)$.
    \item $\mu_\bullet$ fulfills the following PDE:
  \begin{equation} \label{eq:PDE_0_meas}
        \partial_t \mu_t + \nabla \cdot(\mu_t v_t) = 0, \text{ where } v_t(x) =  \int_{y \neq x} (s+1) (x-y)|x-y|^{s-1} \, \diff \mu_t(y) - \nabla V(x),
  \end{equation}
  in the distributional sense, where $V = K_s \ast \nu$ is the background potential.
  \end{enumerate}
\end{corollary}
We will split our analysis of \cref{eq:definition_of_potential} into two parts: Analyzing the repulsion part $E_1$ and analyzing the attraction part $E_2$, i.e., 
\begin{align*}
  E(\mu) &= -\frac{1}{2} \iint K_s(x-y) \,\diff (\mu-\nu)(x) \, \diff (\mu-\nu)(y) \\
  &= \underbrace{-\frac{1}{2} \iint K_s(x-y) \,\diff \mu(x) \, \diff \mu(y)}_{E_1(\mu)} + \underbrace{\iint K_s(x-y) \, \diff \nu(y) \diff \mu(x)}_{E_2(\mu)} - \underbrace{\frac{1}{2} \iint K_s(x-y) \,\diff \nu(x) \, \diff \nu(y)}_{\text{constant in }\mu}. 
\end{align*}
We will ignore the part that is constant in $\mu$, since it does not matter for the curve of maximal slope in $\mu$.

\subsubsection{Analysis of the Repulsive Potential}

As usual, we analyze $\tilde{E}_1 \coloneqq E_1 \circ \law$. In particular, we do that by factorizing it into the composition $\tilde{E}_1 = F \circ \circleddash$, with
\begin{equation} \label{eq:composition_E1}
\begin{tikzcd}
L^2(\Omega; \R^d) \arrow[rr, "\circleddash"] &  & L^2(\Omega \times \Omega; \R^d) \arrow[rr, "F"] &  & \R,
\end{tikzcd}
\end{equation}
where the linear map $\circleddash$ is given by $X \mapsto \big[(\omega, \omega') \mapsto X(\omega) - X(\omega')\big]$, and the function $F$ is given by $Y \mapsto - \frac{1}{2}\EE[K_s(Y)]$, where the expectation is taken with respect to the product measure $\mathbb{P} \times \mathbb{P}$. Since $F$ is concave on $L^2(\Omega \times \Omega; \R^d)$ and $\circleddash$ is linear, we immediately deduce that $\tilde{E}_1$ is concave as well; as such, we turn to analyzing its supergradient. 
For this analysis, the composition \eqref{eq:composition_E1} is useful because we may apply the chain rule $\partial^+(F \circ \circleddash)(X) = \circleddash^* \, \partial^+ F(\circleddash(X))$ for every $X \in L^2$, where $\circleddash^*$ denotes the adjoint of $\circleddash$ (see, e.g., \cite[Corollary 16.52]{bauschkeConvexAnalysisMonotone2017}).

\begin{lemma} \label{lemma:supergradient}
  For all $X \in L^2$ and all $\xi \in \partial^+ \tilde{E}_1(X)$, we have that\footnote{For brevity, we leave the expression ``There exists a representative $g\colon \Omega \to \R$ of $\xi$ such that [...]'' implicit.} for all $\omega \in \Omega$,
  \begin{equation*}
    |\xi(\omega)| \leq 2 (1 + |X(\omega)| + \|X\|_{L^2}).
  \end{equation*}
\end{lemma}
\begin{proof}
  For every fixed $Y \in L^2(\Omega \times \Omega; \R^d)$, we have that every $\zeta \in \partial^+ F(Y)$ lies almost surely in the graph of $\big(\partial^+ (-\frac{1}{2} K_s)\big)(Y)$. 
  Moreover, since
  \begin{equation*}
    \partial^- K_s(x) = \begin{cases}
      \{(s+1) x |x|^{s-1}\}, \qquad &x \neq 0, \\
      \{0\}, \qquad &s\neq 0, x=0, \\
      \overline{B_1(0)}, \qquad &s=0, x=0,
    \end{cases}
  \end{equation*}
  we have for any $(\omega, \omega') \in \Omega \times \Omega$ the estimate $|\zeta(\omega, \omega')| \leq 1 + |Y(\omega, \omega')|$. 
  To obtain an estimate on the superdifferentials of $F \circ \circleddash$, we need to  compute the adjoint $\circleddash^*$. To this end, we pick arbitrary $X \in L^2$ and $Y \in L^2(\Omega \times \Omega; \R^d)$ and calculate
  \begin{align*}
    \langle \circleddash(X), Y \rangle_{L^2(\Omega \times \Omega; \R^d)} &= \iint \langle X(\omega) - X(\omega'), Y(\omega, \omega') \rangle \, \diff \mathbb{P}^2(\omega, \omega') \\
    &= \int \langle X(\omega), \int Y(\omega, \omega') \, \diff \mathbb{P}(\omega') \rangle \, \diff \mathbb{P}(\omega) - \int \langle X(\omega'), \int Y(\omega, \omega') \, \diff \mathbb{P}(\omega) \rangle \, \diff \mathbb{P}(\omega') \\
    &= \int \langle X(\omega), \int Y(\omega, \omega') - Y(\omega', \omega) \, \diff \mathbb{P}(\omega') \rangle \, \diff \mathbb{P}(\omega) \\
    &= \langle X, \left[\omega \mapsto \int Y(\omega, \omega') - Y(\omega', \omega) \, \diff \mathbb{P}(\omega')\right] \rangle_{L^2(\Omega; \R^d)},
  \end{align*}
  and see that $\circleddash^*(Y) = \left[\omega \mapsto \int Y(\omega, \omega') - Y(\omega', \omega) \, \diff \mathbb{P}(\omega')\right]$. Finally, for $X \in L^2$ and $\xi \in \partial^+ F(\circleddash(X))$, we estimate
  \begin{align*}
    |(\circleddash^*(\xi))(\omega)| &\leq \int |\xi(\omega, \omega')| + |\xi(\omega', \omega)| \, \diff \mathbb{P} \leq \int 2 (1 + |X(\omega) - X(\omega')|) \, \diff \mathbb{P}(\omega') \\
    &\leq 2 (1 + |X(\omega)| + \|X\|_{L^2})
  \end{align*}
  and we conclude the proof.
\end{proof}
Having analyzed the superdifferential of $\tilde{E_1}$, we now turn to our analysis of the attractive potential.

\subsubsection{Analysis of the Attractive Potential}

In this part, we focus on the attractive potential $E_2$ by analyzing the composition $\tilde{E}_2 = E_2 \circ \law$ defined by
\begin{equation*}
  \tilde{E}_2\colon X \mapsto \int_{\Omega} (K_s \ast \nu)(X(\omega)) \, \diff \mathbb{P}(\omega).
\end{equation*}
To eventually show that $\tilde{E}_2$ is $C^1$, we will examine the smoothness of the convolution $K_s \ast \nu$.

\begin{lemma}\label{lem:linear_potential_globally_lipschitz}
  Let $\nu$ have the form as in \cref{eq:nu}. Then, $K_s \ast \nu$ is in $W^{2,\infty}_{\mathrm{loc}}(\R^d)$ with $\nabla^2 (K_s \ast \nu) \in L^\infty(\R^d; \R^{d \times d})$ and thus, $\nabla (K_s \ast \nu)$ is globally Lipschitz continuous.
\end{lemma}
\begin{proof}
We will first do the calculation for $s > 1 - d$. We have that $|\nabla^2 K_s| \leq (s+1) (2-s) |z|^{s-1}$. 
We can thus estimate the distributional Hessian, recalling that $f = \frac{\diff \nu}{\diff \Leb^d}$,
  \begin{align*}
    |\nabla^2 (K_s \ast \nu)(x)| &= |((\nabla^2 K_s) \ast \nu)(x)| \\
    &\leq \int_{B_1(x)} |\nabla^2 K_s(x-y)| |f(y)| \, \diff y + \int_{B_1^c(x)} |\nabla^2 K_s(x-y)| |f(y)| \, \diff y \\
    &\leq (s+1)(2-s) \Big(\|f\|_{L^\infty} \int_{B_1(0)} |z|^{s-1} \, \diff z + \|f\|_{L^1}\Big).
  \end{align*}
  The only remaining case is $d=1$, $s=0$, in which case we have that $\nabla^2 K_s = 2 \delta_0$ and $|\nabla^2(K_s * \nu)(x)| \leq 2 \|f\|_{L^\infty}$.
\end{proof}

We are now in a position to prove that $\tilde{E}_2$ is $C^1$.

\begin{corollary}\label{cor:e2_c1}
  The functional $\tilde{E}_2$ is of class $C^1$, with\footnote{We furthermore have $\nabla (K_s \ast \nu) = (\nabla K_s) \ast \nu$.} $\nabla \tilde{E}_2(X) = \big[\omega \mapsto \nabla (K_s \ast \nu)(X(\omega))\big]$.
\end{corollary}
\begin{proof}
  Since $K_s \ast \nu$ is in $W^{2, \infty}_{\mathrm{loc}}$ with $\nabla^2 (K_s \ast \nu) \in L^\infty$, we have that the Gâteaux derivative is given by $\nabla \tilde{E}_2(X) = \big[\omega \mapsto \nabla (K_s \ast \nu)(X(\omega))\big]$. In particular, we have, for $X_1, X_2 \in L^2$, the estimate
  \begin{equation*}
    \|\nabla \tilde{E}_2(X_1) - \nabla \tilde{E}_2(X_2)\|_{L^2} \leq \|\nabla^2 (K_s \ast \nu)\|_{L^\infty} \|X_1 - X_2\|_{L^2},
  \end{equation*}
  upgrading the Gâteaux derivative to a Fréchet derivative and showing that $\tilde{E}_2 \in C^1$.
\end{proof}

\subsubsection{A \texorpdfstring{$C^1$}{C¹}-Perturbation of a Concave Functional}
In this section, we will analyze the curves of maximal slope along $\tilde{E}$ by using the fact that $\tilde{E}_1$ is concave, while $\tilde{E}_2$ is $C^1$.

\begin{lemma}\label{lem:main_lemma}
  The lower semicontinuous envelope $\mathrm{lsc}(|\partial \tilde{E}|)$ of the metric slope is a strong upper gradient for $\tilde{E}$. Furthermore, for any curve $\bgamma$ of maximal slope with respect to $\mathrm{lsc}(|\partial \tilde{E}|)$, the law $\mu_\bullet = \law \circ \bgamma$ fulfills the PDE~\eqref{eq:PDE_0}.
\end{lemma}

\newcommand{\mm}{\mathord{-}\mkern-1.5mu}
\begin{proof}
  Since $-\tilde{E}$ is a $C^1$-perturbation of a convex lower semicontinuous functional, we have that $|\partial \mm\tilde{E}| = \|\partial^0\mm\tilde{E}(\bullet)\|$ is itself lower semicontinuous and a strong upper gradient for $\mm\tilde{E}$, where $\partial^0$ denotes the Fréchet subdifferential of minimal norm (see \cite[Corollary 1.4.5]{ambrosioGradientFlows2008}; lower semicontinuity follows from \cite[(1.4.17), (1.4.18)]{ambrosioGradientFlows2008} and reflexivity of $L^2$). 
  Furthermore, $\tilde{E}$ being a $C^1$-perturbation of a concave functional, we have that 
  $$|\partial \tilde{E}|(Y) \geq \sup_{\zeta \in \partial^+ \tilde{E}(Y)} \|\zeta\|_{L^2} \geq |\partial \mm \tilde{E}|(Y),$$
  and thus, also $\mathrm{lsc}(|\partial \tilde{E}|) \geq |\partial \mm\tilde{E}|$. 
  As $|\partial \mm \tilde{E}|$ is a strong upper gradient for $\mm\tilde{E}$, we deduce that also $\mathrm{lsc}(|\partial \tilde{E}|)$ is a strong upper gradient for $\mm\tilde{E}$, and hence also for $\tilde{E}$. 
  To show the second part of the Lemma (i.e., the fulfillment of the PDE~\eqref{eq:PDE_0}), we would like to invoke \cref{prop:weak-formulation-abstract}; however, \cref{prop:weak-formulation-abstract} requires a curve of maximal slope with respect to the metric slope $|\partial \tilde{E}|$, while we are discussing curves of maximal slope with respect to $\mathrm{lsc}(|\partial \tilde{E}|)$.
  For this reason, we need to recover a curve of maximal slope with respect to the metric slope. To this end, let us take any curve\footnote{The existence of such curves is not known at this point, but will be shown in the proof of \cref{thm:cusp} below.} of maximal slope $\bgamma\colon [0, T] \to L^2$ for $\tilde{E}$ with respect to the strong upper gradient $\mathrm{lsc}(|\partial \tilde{E}|)$. 
  Then, the time-reversed $\bar{\bgamma}\colon t \mapsto \bgamma_{T-t}$ is a curve of maximal slope for $\mm\tilde{E}$ with respect to $\mathrm{lsc}(|\partial \tilde{E}|)$, and hence also with respect to $|\partial \mm \tilde{E}|$. 
  At this point, we can use \cref{prop:weak-formulation-abstract} to determine the weak PDE that is fulfilled by $\bar{\bgamma}$. 
  Therefore, taking $\varphi \in C^\infty_c(\R^d)$ and $\mu \in \Pp_2(\R^d)$, we compute
  \begin{align*}
    \frac{\diff}{\diff h} E_1((T_h^\varphi)_\# \mu)\Big|_{h = 0} &= - \lim_{h \to 0}\frac{1}{2}  \iint \frac{1}{h} \left[ K_s\Big(x - y + h \big(\nabla \varphi(x) - \nabla \varphi(y)\big)\Big) - K_s(x - y) \right] \, \diff \mu(y) \, \diff \mu(x) \\
    &= -\frac{1}{2} \iint_{y \neq x} \langle \nabla K_s(x-y), \nabla \varphi(x) - \nabla \varphi(y) \rangle \, \diff \mu(y) \, \diff \mu(x) \\
    &= -\iint_{y \neq x} \langle \nabla K_s(x-y), \nabla \varphi(x) \rangle \, \diff \mu(y) \diff \mu(x) \\
    &= -\int \langle \nabla \varphi(x), \int_{y \neq x} \nabla K_s(x - y) \, \diff \mu(y) \rangle \, \diff \mu(x),
  \end{align*}
  where we used the antisymmetry of $\nabla K_s$ in the third equality, and used dominated convergence to move the limit into the integral, as explained below.
  Indeed, since $\varphi$ is supported on some $B_R(0) \subseteq \R^d$ and $|\nabla \varphi|\leq C$ for some $C>0$, and since $|\nabla K_s(z)|=(s+1)|z|^s$, by the mean value theorem we have
  \begin{equation*}
  \begin{split}
      \left|\frac{K_s\Big(x - y + h \big(\nabla \varphi(x) - \nabla \varphi(y)\big)\Big) - K_s(x - y)}{h}\right|
      &\leq 2C (s+1) \big( |x-y| + 2C\big)^s \\
      &\leq 2C(s+1)\big(|x| + |y| + 2C \big)^s
  \end{split}
  \end{equation*}
  for every $h \in [-1, 1]$, and the right-hand side is $(\mu \otimes \mu)$-integrable because $\mu \in \Pp_2$ and $s < 1$. 
  For $E_2$, we may use \cref{cor:e2_c1} to immediately deduce that 
  \begin{equation*}
  \frac{\diff}{\diff h} E_2((T_h^\varphi)_\#\mu)\big|_{h = 0} = \int \langle \nabla \varphi(x), \nabla (K_s \ast \nu)(x) \rangle \, \diff \mu(x).      
  \end{equation*}
  We now apply \cref{prop:weak-formulation-abstract} to $\bar{\bgamma}$ and $-\tilde{E}$, where we just computed that, for $\varphi \in C^\infty_c(\R^d)$ and $\mu \in \Pp_2(\R^d)$, we have 
  \begin{equation*}
    \frac{\diff}{\diff h} (-E) \big((T_h^\varphi)_\# \mu\big) \Big|_{h = 0} =\int \Big\langle \nabla \varphi(x), \left[ \int_{y \neq x} \nabla K_s(x - y) \, \diff \mu(y) \right] - \nabla(K_s * \nu) (x) \Big\rangle  \, \diff \mu(x),
  \end{equation*}
   and we deduce that $\bar{\bgamma}$ fulfills the PDE
  \begin{equation*}
        \partial_t \mu_t + \nabla \cdot(\mu_t v_t) = 0, \text{ where } v_t(x) =  - \int_{y \neq x} (s+1) (x-y)|x-y|^{s-1} \, \diff \mu_t(y) + \nabla V(x),
  \end{equation*}
  with $V = K_s \ast \nu$, which is the time-reverse of the PDE~\eqref{eq:PDE_0}; thus, $\bgamma$ fulfills~\eqref{eq:PDE_0}.
\end{proof}

At this point, we are almost ready to prove the main theorem. 
Indeed, after the construction of solutions $(U_\tau^n)_{n \in \NN}$ in $\LtwoR$ for the minimizing movement scheme, we can use compactness arguments as in \cref{rmk:generalized_aa_yields_coms}.
In this way, we obtain that there exists a sequence $(\tau_i)_{i\in \NN}$ with $\tau_i\searrow0$ such that the curves $t \mapsto \bar{U}_{\tau_i}(t)$ (i.e., the piecewise constant interpolations of $(U_{\tau_i}^n)_{n \in \NN}$) fulfill $\bar{U}_{\tau_i}(t) \to \bgamma_t$ for all $t \in [0, +\infty)$ in the sense of the  $\law^{-1}(\sigma_\text{narrow})$ topology, for some limit curve $\bgamma \in \mathrm{AC}^2_{\mathrm{loc}}([0, +\infty); \LtwoR)$. 
Then, we would like to conclude that $\bgamma$ is a curve of maximal slope with respect to $\mathrm{lsc}(|\partial \tilde{E}|)$. 
To say this, however, we would need to upgrade the pointwise convergence $\bar{U}_{\tau_i}(t) \to \bgamma_t$ from $\law^{-1}(\sigma_{\text{narrow}})$ to be strong in $\LtwoR$---or, since $|\partial \tilde{E}|$ factors through $\law$, in \footnote{I.e., we need to upgrade the statement ``\emph{For all $t \in [0, +\infty)$ and all $\psi \in C_b(\R^d)$, we have $\EE_{\mathbb{P}_\R}[\psi(\bar{U}_\tau(t))] \to \EE_{\mathbb{P}_\R}[\psi(\bgamma_t)]$''} to the statement ``\emph{For all $t \in [0, +\infty)$, we have that $W_2(\law(\bar{U}_\tau(t)), \law(\bgamma_t)) \to 0$}''.} $\law^{-1}(\sigma_{W_2})$, as per \cref{rmk:continuity_of_functionals}. 
To this end, we show that the $\bar{U}_\tau(t)$ have uniformly integrable second moments.

\begin{lemma} \label{lemma:uniformly_integrable_sec_moments}
  Let $X_0 \in L^2$, let $(\tau_i)_{i \in \NN}$ be a sequence of step sizes in $\R_+$ such that $\tau_i \to 0$, let $(U_{\tau_i}^n)_{n \in \NN}$ be solutions for~\eqref{eq:mms} such that $U_{\tau_i}^0 = X_0$ for all $i \in \NN$, and let $\bar{U}_{\tau_i}$ denote the corresponding piecewise constant interpolations as in \cref{eq:constant_interpolation}. 
  Then, for each $t \in [0, T]$, the family  $(\bar{U}_{\tau_i}(t))_{i \in \NN}$ has uniformly integrable second moments.
\end{lemma}

Before we prove the above lemma, we quickly show an auxiliary result, which holds for general energies:
\begin{lemma} \label{lemma:frechet_superdifferential}
  Let $\mathcal{H}$ be a Hilbert space, $F\colon \mathcal{H} \to \R$ be a functional. Furthermore, let $X, Y \in \mathcal{H}$ such that
  \begin{equation*}
    Y \in \argmin_{Z \in \mathcal{H}} F(Z) + \frac{|Z - X|_\mathcal{H}^2}{2 \tau}.
  \end{equation*}
  Then, the Fr\'echet superdifferential $D^+F$ has at most one element at $Y$. Furthermore, if $D^+F(Y)$ is nonempty, say $D^+F(Y) = \{\xi\}$, then $Y - X = -\tau \xi$, under the usual Riesz identification.
\end{lemma}
\begin{proof}
  On the one hand, by optimality, for every $V \in \mathcal{H}$ we have that 
  $$\langle Y-X, V \rangle \geq \tau \limsup_{h \downarrow 0} \frac{F(Y) - F(Y + h V)}{h}.$$ 
  On the other hand, for every element $\xi \in D^+F(Y)$, we also have that 
  $$\limsup_{h \downarrow 0} \frac{F(Y) - F(Y + h V)}{h} \geq - \langle \xi, V \rangle.$$ 
  Combining those and replacing $V$ by $-V$ (notice that both inequalities hold for every $V \in \mathcal{H}$) we get that $\langle Y - X, V \rangle = -\tau \langle \xi, V \rangle$ for every $V \in \mathcal{H}$, which yields $Y - X = -\tau \xi$, for every $\xi \in D^+F(Y)$, which finishes the proof.
\end{proof}

\begin{proof}[Proof of \cref{lemma:uniformly_integrable_sec_moments}]
  The goal of this proof is to find some dominating $Z \in L^2(\Omega; \R)$ such that $Z(\omega) \geq |\bar{U}_{\tau_i}(t)(\omega)|$ for all $\omega \in \Omega$, uniformly over $i \in \NN$ and $t \in [0, T]$. To this end, we first note that by combining \cref{lemma:supergradient,lem:linear_potential_globally_lipschitz}, we see that there exists some $C > 0$ such that for every $X \in L^2$ and every Fréchet superdifferential $\xi \in \partial^+\tilde{E}$, we have that $|\xi(\omega)| \leq C (1 + |X(\omega)| + \|X\|_{L^2})$. 
  Furthermore, we can use the a priori bounds from \cite[Lemma 3.2.2]{ambrosioGradientFlows2008} to find some bound $D \geq \|U_{\tau_i}^{n}\|_{L^2}$, uniformly over $n \in \NN$ such that $\tau n \leq T$. 
  Next, we see that\footnote{Again, for brevity, we decide to drop the phrase ``There exists a representative of $U_{\tau_i}^{n+1}$ such that for all $\omega \in \Omega$ [...]''.}, for $i$ fixed and $n \in \NN$ such that $1 \leq n \leq \lfloor T / {\tau_i} \rfloor$, we have that
  \begin{align*}
    |U_{\tau_i}^{n}(\omega) - U_{\tau_i}^{n-1}(\omega)| &\overset{(1)}{=} \tau_i |\xi(\omega)|
    \overset{(2)}{\leq} {\tau_i} C(1 + |U_{\tau_i}^{n}(\omega)| + \|U_{\tau_i}^{n}\|_{L^2})
    \\&\overset{(3)}{\leq} {\tau_i} C' (1 + |U_{\tau_i}^{n}(\omega)|),
  \end{align*}
  where we denote by $U_{\tau_i}^{n}$ the $n$-th step of the MMS with stepsize ${\tau_i}$, $\xi$ is the unique---by \cref{lemma:frechet_superdifferential}---element of the superdifferential of $\tilde{E}$ at $U_{\tau_i}^{n}(\omega)$. Moreover, we also used \cref{lemma:frechet_superdifferential} in $(1)$, the estimate on $|\xi(\omega)|$ from above in $(2)$ and the uniform bound on $\|U_{\tau_i}^{n}\|_{L^2}$ in $(3)$, for some $C'$ independent of $i$ and $n \in \NN$ such that $1 \leq n \leq \lfloor T / {\tau_i} \rfloor$. If we combine the triangle inequality and the last estimate, we see that
  \begin{equation*}
    |U_{\tau_i}^{n}(\omega)| \leq \frac{1}{1-\tau_i C'} (|U_{\tau_i}^{n-1}(\omega)| + \tau_i C'),
  \end{equation*}
  from which we iteratively get that 
  \begin{align*}
    |U_{\tau_i}^{n}(\omega)| &
    \leq \frac{1}{(1-\tau_i C')^{n}} |X_0(\omega)| +  \sum_{k=1}^{n} \frac{\tau_i C'}{(1-\tau_i C')^k} 
    \leq \frac{1}{(1-\tau_i C')^{n}} (|X_0(\omega)| + C' T)  \\
    &\leq (|X_0(\omega)| + C' T) \left(\frac{1}{1- C'\tau_i}\right)^{T/\tau_i} 
    \leq (|X_0(\omega)| + C' T) (1 + e^{C'T}),
  \end{align*}
  where the last estimate holds for $\tau_i$ small enough, uniformly over $n$ and $\omega \in \Omega$. Thus, for $\tau_i$ small enough, we have that for all $n \in \NN$ such that $0 \leq n \leq \lfloor T / \tau_i \rfloor$, $U_{\tau_i}^{n}$ is dominated by $(|X_0| + C'T) (1 + e^{C'T}) \in L^2$, from which we deduce uniformly integrable second moments.
\end{proof}

We are now ready to prove \cref{thm:cusp}.
\begin{proof}[Proof of \cref{thm:cusp}]
  Let us first show that for any $\tau > 0$ small enough, we can construct an MMS sequence $(U_\tau^n)_{n \in \NN}$ starting in $X_0$ (i.e., fulfilling~\eqref{eq:mms}; see \cref{sec:constructing_curves_of_maximum_slope}). 
  To this end, we first argue that for each $U_\tau^n$, we can find $\mu_{\tau, n+1} \in \Pp_2(\R^d)$ such that
  \begin{equation*}
      \mu_{\tau, n+1} \in J_{\tau, \Pp_2}[\law(U_\tau^n)] = \argmin_{\rho \in \Pp_2}\big(E(\rho) + \frac{1}{2\tau} W_2^2(\law(U_\tau^n), \rho)\big).
  \end{equation*} 
  To this end, recalling the decomposition $E=E_1 +E_2$, let us first note that $E_1$ is narrowly continuous on $W_2$-bounded subsets of $\Pp_2$; this follows from considering a narrowly converging sequence $\rho_n$ bounded in $W_2$ and applying \cite[Lemma 5.1.7]{ambrosioGradientFlows2008} with $\mu_n = \rho_n \otimes \rho_n$---which is itself narrowly converging and has bounded second moments, and thus uniformly integrable $(s+1)$-moments---and $f(x, y) = K_s(x-y)$---which has $(s+1)$-growth. 
  Likewise, $E_2$ is narrowly continuous on $W_2$-bounded subsets of $\Pp_2$ because the $(s+1)$-moment bound of $\nu$ and $(s+1)$-growth of $K_s$ give $(s+1)$-growth of $K_s * \nu$; furthermore, $K_s * \nu$ is continuous by \cref{lem:linear_potential_globally_lipschitz}. Thus, also $E= E_1 + E_2$ is narrowly continuous on $W_2$-bounded subsets. 
  Moreover, $E= E_1 + E_2$ is also clearly coercive since $s+1 < 2$ (in the sense of \cite[2.1b]{ambrosioGradientFlows2008}\footnote{I.e., such that there exist $\tau^\star > 0$ and $\mu^\star \in \Pp_2$ such that $\inf_{\mu \in \Pp_2} E_1(\mu) + E_2(\mu) + \frac{1}{2\tau^\star} W_2^2(\mu^\star, \mu) > -\infty$}) and by the direct method (using lower semicontinuity of $W_2^2(\law (U_\tau^n), \bullet)$), we see that such $\mu_{\tau, n+1}$ exists. Therefore, we conclude that $J_{\tau, \Pp_2}[\law(U_\tau^n)]$ is nonempty. 
  By \cref{rmk:mms_downstairs_upstairs}, also $J_{\tau, \LtwoR}[U_\tau^n]$ is nonempty, and we can pick $U_\tau^{n+1} \in J_{\tau, \LtwoR}[U_\tau^n]$. Proceeding recursively, we can construct the required MMS sequence $(U_\tau^n)_{n \in \NN}$. Next, we choose a sequence $(\tau_i)_{i \in \NN}$ converging to $0$ and, for each $i \in \NN$, an MMS sequence $(U_{\tau_i}^n)_{n \in \NN}$ with stepsize $\tau_i$. 
  By the standard a priori bounds (see \cite[Lemma 3.2.2]{ambrosioGradientFlows2008}), for each $T > 0$, there exists a uniform bound $C$ such that $\|U_{\tau_i}^n\|_{\LtwoR} \leq C$, for all $i \in \NN$ and all $n \leq \lfloor T / \tau_i \rfloor$; furthermore, we have the standard $AC$-like bounds for piecewise constant interpolations (again, see \cite[Lemma 3.2.2]{ambrosioGradientFlows2008}). 
  Thus, we may proceed as in \cref{rmk:generalized_aa_yields_coms}, using \cref{cor:l2_is_aa}, to get that there exists a subsequence $\tau_i$ (without relabelling) and a locally absolutely continuous limit curve $\bgamma \colon [0, \infty) \to \LtwoR$ such that for all $t \geq 0$, we have that $\bar{U}_{\tau_i}(t) \to \bgamma_t$ in $\law^{-1}(\sigma_{\text{narrow}})$, where $\bar{U}_{\tau_i}$ denotes the piecewise constant interpolation. 
  Since the $\bar{U}_{\tau_i}(t)$ have uniformly integrable second moments by \cref{lemma:uniformly_integrable_sec_moments}, we have in addition that $\bar{U}_{\tau_i}(t) \to \bgamma_t$ in $\law^{-1}(\sigma_{W_2})$. 
  Using a standard estimate for the distance of the De Giorgi interpolant $\tilde{U}_{\tau_i}(t)$ from the piecewise constant interpolation (\cite[Lemma 3.2.2]{ambrosioGradientFlows2008}), we see that for all $t \in [0, T]$, $\|\tilde{U}_{\tau_i}(t) - \bar{U}_{\tau_i}(t)\|^2_{\LtwoR} \leq C \tau_i$, where $C$ is independent of $i$ and $t \in [0, T]$. 
  Thus, we even have $\tilde{U}_{\tau_i}(t) \to \bgamma_t$ in $\law^{-1}(\sigma_{W_2})$, and we can pass to the limit (using \cref{rmk:continuity_of_functionals}) in the estimate (see \cite[Lemma 3.2.2, eq. 3.2.3]{ambrosioGradientFlows2008})
  \begin{equation*}
    \frac{1}{2} \int_{0}^{\tau_i \cdot n} |U'_{\tau_i}|^2(t) \, \diff t + \frac{1}{2} \int_{0}^{\tau_i \cdot n} |\partial \tilde{E}|^2(\tilde{U}_{\tau_i}(t)) \, \diff t \leq \tilde{E}(X_0) - \tilde{E}(U_{\tau_i}^n),
  \end{equation*}
  (with $U'_{\tau_i}$ defined as in \cite[eq. 2.3.4]{ambrosioGradientFlows2008}) to obtain
  \begin{equation} \label{eq:limit_inequality}
    \frac{1}{2} \int_{0}^{T} \|\dot{\bgamma}_t\|_{\LtwoR}^2 \, \diff t + \frac{1}{2} \int_{0}^{T} \mathrm{lsc}(|\partial \tilde{E}|)^2(\bgamma_t) \, \diff t \leq \tilde{E}(\bgamma_0) - \tilde{E}(\bgamma_T),
  \end{equation}
  for any $T > 0$, where we also used the continuity of $E$ w.r.t. $W_2$. By \cref{lem:main_lemma}, $\mathrm{lsc}(|\partial \tilde{E}|)$ is a strong upper gradient, and hence we can deduce from \cref{eq:limit_inequality} that $\bgamma$ is a curve of maximal slope with respect to $\mathrm{lsc}(|\partial \tilde{E}|)$. Thus, we can finish the proof by invoking \cref{lem:main_lemma} once more.
\end{proof}

\section{Newtonian Dynamics in the Wasserstein Space}  \label{sec:newtonian_curves}

While the gradient flow theory in $(\Pp_2,W_2)$ provided an explanation of numerous dissipative PDEs (see, e.g., \cite{OttoGeometry2001,carrilloKineticEquilibrationRates2003}), it remained elusive the relation with other important evolutions, such as Vlasov-type equations driven by suitable potential energies, whose solutions are probability distributions over the phase space.
In the present section, we reveal how to connect such equations with proper notions of second-order dynamics in $\ful{\Pp}_2 \coloneqq \Pp_2(\R^{d}\times \R^{d})$, which we refer to as \emph{Newtonian motion}.
Our results stem from an abstract theory of Newtonian motion in metric spaces, which is developed in the companion paper \cite{fornasierNewtonianMotionCompanion} and is aimed to generalize to metric spaces the notion of solution of the Newton equation $\ddot \gamma = -\nabla E(\gamma)$ (classically phrased leveraging a Hilbert/differentiable structure).
For this reason, we 
advise the reader to proceed the reading by familiarizing with the notions and the tools from \cite{fornasierNewtonianMotionCompanion}, which we concisely recall below, in order to better highlight the challenges posed by the extension of that theory to Wasserstein spaces. 

Given a metric space $(\mathcal S,d_{\mathcal S})$, the construction of Newtonian motion provided in \cite{fornasierNewtonianMotionCompanion} rests on three pillars: \\

\textbf{Pillar 1: A geometric-variational characterization of second-order absolutely continuous curves in metric spaces.} We introduce the space of geodesics\footnote{By \emph{geodesics in $X$}, we mean continuous curves $\xi\colon [a, b]\to \mathcal S$ (with $a<b$) such that, for every $s, t \in [a, b]$, it holds that $	d_{\mathcal S}(\xi_s, \xi_t) = |s - t| \frac{1}{b-a} d_{\mathcal S}(\xi_a, \xi_b)$.}
\begin{equation} \label{eq:recall_geod_space}
	\Tang[\varepsilon] \coloneqq \{\xi \in AC([t, t + \varepsilon]; \mathcal S) \mid \xi \text{ is a constant-speed geodesic in } \mathcal S, t \in \R\}, \quad \varepsilon > 0.
\end{equation}
Then, we say that a curve $\gamma\colon (a, b) \to \mathcal S$ is \emph{\ACt} if there exist $\varepsilon_\gamma > 0$, $\mathfrak{E}_\gamma \geq 0$, $m \in L^2((a, b); \R)$, and a map $V^\gamma\colon (a, b)\to \Tang[\varepsilon_\gamma]$ with $V^\gamma_t  \in AC([t, t + \varepsilon_\gamma]; \mathcal S)$ and $V^\gamma_t(t) = \gamma_t$ such that 
\begin{equation} \label{eq:recall_ac2}
			d_{\mathcal S} \big( V_s^\gamma(T), V_t^\gamma(T) \big) \leq \int_s^t (T-r) \, m(r) \big(1 + \mathfrak{E}_\gamma  (T-r)^2\big) \, \diff r
\end{equation}
for every $s, t \in (a, b)$ and $T \in \R$ with $s \leq t \leq T \leq s + \varepsilon_\gamma$. 
For almost every time $t\in (a, b)$, $V^\gamma_t$ should be interpreted as the geodesic with scalar speed $|\gamma'_t|$ and \emph{tangent} to the curve $\gamma$ at $\gamma_t$. Notice that, for $t\in (a, b)$, the geodesic $V^\gamma_t$ departing from $\gamma_t$ and with scalar speed $|\gamma'_t|$ is (required to be) defined on the domain $[t,t+\varepsilon_\gamma]$, whose length $\varepsilon_\gamma$ does not depend on $t$.
Morally, inequality \eqref{eq:recall_ac2} should be interpreted as ensuring that the ``tangent geodesics'' $V_t^\gamma$ change absolutely continuously over $t \in (a,b)$. 
To this end, in order to fairly compare these changing geodesics, we require a uniform time horizon $\varepsilon_\gamma$ for the time-length of their domains. 
In rough spaces, we might fail to find such a $\varepsilon_\gamma$, even for innocent-looking curves\footnote{As an example, consider the curve $(-1, 0) \ni t \mapsto \frac{1}{2} \delta_{-t} + \frac{1}{2} \delta_{t} \in \Pp_2$. Even though this curve is a geodesic in $(\Pp_2, W_2)$, it fails to pass the definition of an \ACt-curve.} Thus, in order to ensure that the definition of \ACt-regularity is not overly restrictive, the ``extendability of geodesics'' is a crucial assumption. This assumption plays also a relevant role in the formulation of a discrete-time predictor-corrector scheme, see Pillar 3 and \eqref{eq:definition_scheme_recall}.
We mention that $\mathfrak{E}_\gamma > 0$ accounts for a possible curvature in the space $\mathcal S$.
Finally, similarly as in the first-order $AC$-case, the quantity
\begin{equation} \label{eq:recall_second_derivative_bound}
		|\gamma''|(s) \coloneqq \limsupTauTS\left(\frac{d_{\mathcal S} \big(V^\gamma_s(T), V^\gamma_t(T)\big)}{(T-s)(t-s)}\right)
\end{equation}
plays the role of a \emph{metric acceleration}, and it turns out that $|\gamma''|(s) \leq m(s)$ for almost every $s \in (a, b)$.\\

\textbf{Pillar 2: A geometric, variational generalization of Newtonian flows.}
We notice that in a Hilbert space $(\mathcal{H}, |\cdot|)$ and for a smooth energy $E \colon \mathcal{H} \to \R$, owing to the Young inequality yielding $-\langle \ddot \gamma_r, \nabla E(\gamma_r) \rangle \leq \frac{1}{2} |\ddot \gamma_r|^2 + \frac{1}{2} |\nabla E(\gamma_r)|^2$, the relation
\begin{equation} \label{eq:recall_inequality}
	\langle \dot \gamma_s, \nabla E(\gamma_s) \rangle - \langle \dot \gamma_t, \nabla E(\gamma_t) \rangle \leq \int_s^t \frac{1}{2} |\ddot \gamma_r|^2  + \frac{1}{2} |\nabla E(\gamma_r)|^2 - \nabla^2 E(\gamma_r)(\dot \gamma_r, \dot \gamma_r) \, \diff r
\end{equation}
is fulfilled automatically for every $\gamma \in W^{2,2}((a,b);\mathcal{H})$ and every $s,t\in(a,b)$ with $s\leq t$. Moreover, equality in \cref{eq:recall_inequality} implies that $\gamma$ satisfies $\ddot \gamma = - \nabla E(\gamma)$. 
The extension to metric spaces requires the following replacements:
\begin{equation*}
    \begin{split}
        \langle \dot\gamma_s, \nabla E(\gamma_s)\rangle &\longleftrightarrow \frac{\diff}{\diff \sigma} E \big(V_s^\gamma(\sigma) \big)\Big|_{\sigma = s}, \\
         |\ddot \gamma_t| &\longleftrightarrow |\gamma''|(t) \qquad \mbox{ (see \cref{eq:recall_second_derivative_bound})},
    \end{split}
\end{equation*}
and finally we replace the pair $\big(|\nabla E(\gamma_t)|, \nabla^2 E(\gamma_t)(\dot \gamma_t, \dot \gamma_t)\big)$ by a proper pair $(g, H)$, where $g\colon \mathcal S \to [0, +\infty)$ and $H\colon 
\bigcup_{\varepsilon > 0}\Tang[\varepsilon] \to \R$ (see \cref{eq:recall_geod_space} for the definition of $\Tang[\varepsilon]$). 
We name such a replacement an \emph{upper gradient-Hessian pair} if the resulting inequality is fulfilled automatically for any $\eta \in \ACt((a,b);\mathcal S)$, and we call an \ACt-curve $\gamma$ a \emph{curve of Newtonian motion} if equality is achieved. 
We thus characterize curves of Newtonian motion as those curves $\gamma \in \ACt((a,b);\mathcal S)$ that satisfy
\begin{equation} \label{eq:recall_curves_of_newtonian_motion}
		\frac{\diff}{\diff \sigma} E\big( V_s^\gamma(\sigma) \big) \Big|_{\sigma = s} - \frac{\diff}{\diff \sigma} E\big( V_t^\gamma(\sigma) \big) \Big|_{\sigma = t} 
        \geq \int_s^t \frac{1}{2}  |\gamma''|^2(r)\, \diff r + \int_s^t \frac{1}{2} g^2(\gamma_r)\, \diff r - \int_s^t H(V_r^\gamma)\, \diff r
\end{equation}
for every $s,t \in(a,b)$ with $s\leq t$.
For regular enough energies\footnote{We will see an example of sufficient regularity conditions in \cref{prop:upper_gradient_hessian}}, a natural candidate for an upper gradient-Hessian pair is the pair $(|\partial E|, H_0^E)$, where $|\partial E|$ denotes the metric slope and $H_0^E$ is defined as
  \begin{equation*}
    H_0^{E}(\xi) \coloneqq \limsup_{h \searrow 0} \frac{\tilde{E}(\xi_t) - 2 \tilde{E}(\xi_{t + h}) + \tilde{E}(\xi_{t + 2h})}{h^2}, \qquad \text{for } \ \xi \in \mathcal T_\varepsilon \mathcal{S}, \quad  \xi\colon [t,t+\varepsilon] \to \mathcal{S}.
  \end{equation*}
In a Hilbert space $\mathcal{H}$ and for a smooth energy $E\colon \mathcal{H} \to \R$, we can provide a more explicit description of the upper gradient-Hessian pair.
Namely, as every geodesic $\xi \in \bigcup_{\varepsilon > 0}\mathcal T_\varepsilon \mathcal{H}$ has the form $[t_\xi, t_\xi+\varepsilon] \ni \sigma \mapsto x_\xi + (\sigma - t_\xi) v_\xi$ for some $x_\xi , v_\xi \in \mathcal{H}$ and some $\varepsilon>0$, the functions
\begin{equation}
	\Big(g\colon x \mapsto |\nabla E(x)|, \quad H\colon \xi \mapsto \nabla^2 E(x_\xi) [v_\xi,v_\xi] \Big)
\end{equation}
are an upper gradient-Hessian pair, and a Sobolev curve $\gamma \in W^{2,2}((a,b); \mathcal{H})$ is a curve of Newtonian motion if and only if it fulfills $\ddot \gamma = - \nabla E(\gamma)$.
This explains the relevance of the synthetic metric variational characterization \eqref{eq:recall_curves_of_newtonian_motion} based on purely metric quantitites, yet able to fully capture vectorial objects.

\textbf{Pillar 3: A discrete, variational scheme to construct and exhibit curves of Newtonian motion}. 
We introduce a two-step predictor-corrector scheme. 
Given a stepsize $\tau > 0$ and the points $U_{\tau, i-1}, U_{\tau, i} \in \mathcal S$, we recursively define
\begin{equation} \label{eq:definition_scheme_recall}
	\begin{split}
	&Z_{\tau, i} \coloneqq \argmin_{z \in \mathcal S} \Big\{2d_{\mathcal S}^2(U_{\tau, i}, z) - d_{\mathcal S}^2(U_{\tau, i-1}, z) \Big\} \qquad \mbox{(predictor)}, \\
	&U_{\tau, i+1} \in \argmin_{z \in \mathcal S} \Big\{E(z) + \frac{1}{2\tau^2} d_{\mathcal S}^2(Z_{\tau, i}, z) \Big\} \qquad\qquad \mbox{(corrector)}.
	\end{split}
\end{equation}
We notice that the temporary point $Z_{\tau, i}$ is the metric analogue of $U_{\tau, i} + \big( U_{\tau, i} - U_{\tau, i-1} \big)$, and the extendability of geodesics is crucial to obtain that, if $U_{\tau, i}, U_{\tau, i-1}$ are sufficiently close, then $Z_{\tau, i}$ is such that
\begin{equation*}
    d_{\mathcal S}(U_{\tau, i-1}, U_{\tau, i}) = d_{\mathcal S}(U_{\tau, i}, Z_{\tau, i}) = \frac12 d_{\mathcal S}(U_{\tau, i-1}, Z_{\tau, i}).
\end{equation*}
We propose piecewise constant interpolations $u_\tau$ between the points $(U_{\tau, i})_{i \in \NN}$ and ``approximate tangents'' $\mathcal{V}^{\tau}$, and we prove that they fulfill approximate versions of the inequalities in \cref{eq:recall_ac2,eq:recall_curves_of_newtonian_motion}. \\

In the Wasserstein space $(\Pp_2, W_2)$ the crucial condition about the extendability of geodesics is not fulfilled\footnote{E.g., for all $\varepsilon > 0$, the geodesic $\xi\colon [-\varepsilon, 0] \to \Pp_2(\R)$ defined by $\xi(t) = \frac{1}{2}\delta_{t} + \frac{1}{2}\delta_{-t}$ cannot be extended beyond $0$.}. 
This obstruction yields major difficulties to the definition of a proper notion of \ACt-curves as well as to the implementation of the predictor-corrector scheme to construct Newtonian curves.
Fortunately, the Lions lift to the space $L^2$ dissolves the issue about both the notion of \ACt-curve---as we can rely on the Hilbert space structure and the usual notion of acceleration---and the implementation of the predictor-corrector scheme.
Moreover, the compactness and the use of Arzelà--Ascoli arguments can be carried out by taking advantage of the results of \cref{sec:wasserstein_as_quotient}, and we are able to describe curves of Newtonian motion $t\mapsto \bgamma_t$ in the lifted space $L^2$. By considering the curve \emph{in the phase space} $\ful{\bgamma}\colon t\mapsto (\bgamma_t, \dot \bgamma_t)$ and by projecting it through the joint law $\ful{\mathcal{L}}\colon L^2\times L^2\to \Pp_2(\R^{d} \times \R^d)$, we are finally able to describe curves of Newtonian motion in $\ful{\Pp}_2 \coloneqq \Pp_2(\R^{d} \times \R^d)$.

Incidentally, \cite{brigatiKineticOptimalTransport2025} introduced an intrinsic notion of ``physical curves'' on $\ful{\Pp}_2$, which corresponds naturally to projections through  $\ful{\mathcal{L}}$ of $W^{2,2}$-curves in $L^2$, as we prove in \cref{cor:lifting_physical_curves} in \cref{sec:characterizing_wasserstein_newtonian_curves}. 
In the same section, we furthermore see that all the quantities
appearing in the variational relation \eqref{eq:recall_curves_of_newtonian_motion} for a $\LtwoR$-valued $W^{2,2}$-regular curve $\bgamma_\bullet$ can be completely described in terms of the projected curve $\ful{\mu}_\bullet = \ful{\law} \circ \ful{\bgamma}$.
We thus obtain a characterization---through a variational inequality in $\ful{\mathcal{P}}_2$---of those curves $\ful{\mu}_\bullet\colon [0, T] \to \ful{\mathcal{P}}_2$ that arise as compositions of curves of Newtonian motion---in the phase space $\ful{L}^2$---with $\ful{\law}$. 

In \cref{sec:constructing_newtonian_curves_in_l2}, we focus on the construction of such curves of Newtonian motion: We first apply the compactness arguments from \cref{sec:a-a-readiness-of-l2,sec:a-a-w-o-lsc} to the predictor-corrector scheme~\eqref{eq:definition_scheme_recall} in $\LtwoR$, obtaining curves of Newtonian motion through a limiting argument. We furthermore see that we can run the predictor-corrector scheme in $\ful{\mathcal{P}}_2$ itself, without ever considering the space $\LtwoR$. 
Thus, while the space $\LtwoR$ exists in the background---providing both the motivation of the definition, as well as the space where our limiting argument takes place---we can carry out all constructions completely in $\ful{\mathcal{P}}_2$. In this sense, $\ful{\mathcal{P}}_2$ captures ``enough information'' to both characterize and construct curves of Newtonian motion for those energies  that factor through $\law\colon \LtwoR\to \Pp_2$.

Next, in \cref{sec:identifying_newtonian}, we show that the image under $\law$ of such Newtonian curves are solutions to certain Vlasov equations.

As an application of this program, in \cref{sec:example_newtonian} we present the construction and the study of Newtonian trajectories driven by interaction energies of the form
\begin{equation*}
  E\colon \Pp_2 \to \R, \qquad \nu \mapsto \iint K(x, x') \, \diff \nu(x) \, \diff \nu(x') + \int V(x) \, \diff \nu(x).
\end{equation*}
In this case, we obtain that the projected Newtonian trajectories fulfill the Vlasov equation
\begin{align*}
    \partial_t \ful{\mu}_t + v \cdot \nabla_x \ful{\mu}_t - \nabla_v\bigg(\Big(2 \int \nabla_{x_1} K(\cdot, x') \, \diff \ful{\mu}_t(x', v')\Big) \ful{\mu}_t\bigg) - \nabla_v\Big((\nabla V) \ful{\mu}_t\Big) = 0.
\end{align*}

\subsection{Projecting Newtonian Dynamics from \texorpdfstring{$L^2$}{L²} to \texorpdfstring{$\ful{\Pp}_2$}{P̲₂}} \label{sec:characterizing_wasserstein_newtonian_curves}

Applying Pillar 2 from above to the choice\footnote{We recall that we assume that $(\Omega, \mathbb{P})$ to be a probability space which admits a random variable $Z$ such that $\law(Z) = \Leb|_{[0, 1]}$, and set $L^2 \coloneqq L^2(\Omega; \R^d)$. The results hold for any such choice of $(\Omega, \mathbb{P})$; one possible choice is $L^2 = \LtwoR$.} $(\mathcal{S}, d_{\mathcal{S}}) = (L^2, \|\cdot\|_{L^2})$, we obtain that for regular enough\footnote{We single out a regularity regime where this holds in \cref{prop:upper_gradient_hessian}; here, we implicitly mean curves of Newtonian motion \emph{with respect to the upper gradient-Hessian pair $(|\partial \tilde{E}|, H_0^{\tilde{E}})$}} energies $\tilde{E}\colon L^2 \to \R \cup \{+\infty\}$, curves of Newtonian motion in with respect to $\tilde{E}$ are those curves $\bgamma\in W^{2,2}((a, b); L^2)$ which fulfill the variational inequality
	\begin{equation} \label{eq:def_variational_Newtonian_in_intro_of_char}
		\frac{\diff}{\diff \sigma} \tilde E(V_s^{\bgamma}(\sigma)) \Big|_{\sigma = s} 
        - \frac{\diff}{\diff \sigma}  \tilde E(V_t^{\bgamma}(\sigma)) \Big|_{\sigma = t} 
        \geq \int_s^t \frac{1}{2}  \|\ddot{\bgamma}_r\|_{L^2}^2\, \diff r + \int_s^t \frac{1}{2} |\partial \tilde{E}|^2(\bgamma_r)\, \diff r - \int_s^t H_0^{\tilde{E}}(V_r^{\bgamma})\, \diff r,
	\end{equation}
  for almost every $s , t \in (a, b)$  with $s<t$.
  This is in analogy to the first-order case of curves of maximal slope in $L^2$, which are, for regular enough energies $\tilde{E}$, defined to those curves $\bgamma \in W^{1,2}((a, b); L^2)$ that fulfill the variational inequality
  \begin{equation}\label{eq:def_variational_coms_in_intro_of_char}
    \tilde{E}(\bgamma_s) - \tilde{E}(\bgamma_t) \geq \int_s^t \frac{1}{2} \|\dot{\bgamma}_r\|^2_{L^2} \, \diff r + \int_s^t \frac{1}{2} |\partial \tilde{E}|^2(\bgamma_r)\, \diff r
  \end{equation}
  for every $s , t \in (a, b)$  with $s<t$. In that first-order case, we saw in \cref{sec:lifting_gradient_problems} that, whenever $\tilde{E}$ is of the form $\tilde{E} = E \circ \law$ for some energy $E\colon \Pp_2 \to \R$ and furthermore $\bgamma$ is the horizontal lift\footnote{I.e., $\|\ddot{\bgamma}_r\|_{L^2} = |\mu'|(r)$ for almost every $r \in (a, b)$, where $|\mu'|(r)$ denotes the metric derivative of $\mu_\bullet$ in $(\Pp_2, W_2)$. For a further discussion, we refer to \cref{sec:lifting_gradient_problems}.} of the $(\Pp_2, W_2)$-valued  $\mathrm{AC}^2$-curve $\mu_\bullet \coloneqq \law \circ \bgamma$, then each of the terms in \cref{eq:def_variational_coms_in_intro_of_char} can be obtained from $E$ and $\mu_\bullet$ directly. 
  More precisely, we have that
  \begin{equation*}
    |\partial \tilde{E}|(\bgamma_r) = |\partial E|(\mu_r), \qquad \|\dot{\bgamma}_r\|_{L^2} = |\mu'|(r), \qquad \tilde{E}(\bgamma_t) = E(\mu_t),
  \end{equation*}
  for almost every $r \in (a, b)$ and for every $t \in (a, b)$. Thus, in this case, $\bgamma$ fulfills \cref{eq:def_variational_coms_in_intro_of_char} whenever $\mu_\bullet$ fulfills
  \begin{equation}\label{eq:characterization_coms_p2_intro_of_char}
    E(\mu_s) - E(\mu_t) \geq \int_s^t \frac{1}{2} |\mu'|^2(r) \, \diff r + \int_s^t \frac{1}{2} |\partial E|^2(\mu_r)\, \diff r
  \end{equation}
  for every $s , t \in (a, b)$  with $s<t$, which is precisely the definition of $\mu_\bullet$ being a curve of maximal slope for $E$\footnote{We avoid in this outline the notion of upper gradients; for details, see \cref{lem:summarize_sec_1}}. 
  This observation can be summarized in the correspondence
  \begin{equation} \label{eq:correspondence_in_intro_char}
  \left\{\begin{array}{c}\text{Horizontal lifts of curves of }\\
  \text{ maximal slope of } E \text{ in } \Pp_2\end{array}\right\}
  \leftrightarrow
  \left\{\begin{array}{c}\text{Curves of maximal slope}\\
  \text{of } E \circ \law \text{ in } L^2\end{array}\right\}.
  \end{equation}
  
  The task of this section is precisely the second-order analogue of equivalence~\eqref{eq:correspondence_in_intro_char}: 
  Working in the phase space $\ful{\Pp}_2 = \Pp_2(\R^{2d})$, we aim at obtaining all the quantitites involved in \cref{eq:def_variational_Newtonian_in_intro_of_char} directly from the phase-space curve $\ful{\mu}_\bullet = \ful{\law} \circ (\bgamma, \dot{\bgamma})$ and the energy $E$, when the latter has the form $\tilde{E} = E \circ \law$ (see \cref{tab:second_order_equivalences}). 
  Thus, the same way we could translate the first-order inequality \eqref{eq:def_variational_coms_in_intro_of_char} on $L^2$ into the inequality \eqref{eq:characterization_coms_p2_intro_of_char} on $\Pp_2$, we can translate the second-order variational inequality~\eqref{eq:def_variational_Newtonian_in_intro_of_char} into the inequality
	\begin{equation} \label{eq:variational_characterization_newtonian_p2_intro_char}
    \begin{split}
    \frac{\diff}{\diff \sigma} E(V_s^{\ful{\mu}}(\sigma))\Big|_{\sigma = s} - \frac{\diff}{\diff \sigma} E(V_t^{\ful{\mu}}(\sigma))\Big|_{\sigma = t} &\geq \int_s^t \frac{1}{2} \big|\ful{\mu}''\big|^2(r) \, \diff r + 
    \int_s^t \frac{1}{2} g^2\big([(x, v) \mapsto x]_\# \ful{\mu}_r \big)\, \diff r 
    \\& \qquad  - \int_s^t H_0^E(\ful{\mu}_r) \, \diff r.
    \end{split}
	\end{equation}
  We call the class of curves fulfilling the variational inequality~\eqref{eq:variational_characterization_newtonian_p2_intro_char} \emph{Wasserstein-Newtonian curves} (see \cref{def:wasserstein_newtonian_curve}), and, with this definition, we obtain the second-order version of correspondence~\eqref{eq:correspondence_in_intro_char}:
    \begin{equation} \label{eq:correspondence_second_order_in_intro_char}
\left\{\begin{array}{c}\text{Horizontal lifts of Wasserstein-}\\
\text{ Newtonian curves in } \ful{\Pp}_2\end{array}\right\} 
  \leftrightarrow
  \left\{\begin{array}{c}\text{Curves of Newtonian motion}\\
  \text{of } E \circ \law \text{ in } L^2\end{array}\right\},
  \end{equation}
  where horizontal lifts of physical curves in $\ful{\Pp}_2$ are defined in \cref{def:horizontal_lift_physical}.
  The first term in \cref{eq:def_variational_Newtonian_in_intro_of_char} that we tackle is $\|\ddot{\bgamma}\|_{L^2}$. To this end, we recall the the notion of \emph{physical curves} from \cite{brigatiKineticOptimalTransport2025}, which turns out to be precisely the Wasserstein avatar of $L^2$-valued $W^{2,2}$-curves.
  \begin{table}[ht]
    \centering
    \small
    \renewcommand{\arraystretch}{1.4}
    \begin{tabular}{>{\raggedright\arraybackslash}p{0.28\textwidth} >{\raggedright\arraybackslash}p{0.30\textwidth} >{\raggedright\arraybackslash}p{0.34\textwidth}}
      \toprule
      Notions in $L^2$& Equivalences in metric spaces $(\mathcal{S}, d_{\mathcal{S}})$ & Avatars in $\ful{\Pp}_2$ \\
      \midrule
      $W^{2,2}((a, b); L^2)$
      & $\ACt((a, b); \mathcal{S})$ \newline (\cite[\cref*{M-def:AC2_curve}]{fornasierNewtonianMotionCompanion})
      & $\mathrm{PC}^2((a, b); \ful{\Pp}_2)$ \newline (\cref{def:physical_curves} cited from \cite[Definition 1.3]{brigatiKineticOptimalTransport2025}, \cref{cor:lifting_physical_curves}) \\
      $V_t^{\bgamma} \coloneqq \big[\sigma \mapsto \bgamma_t + (\sigma - t) \dot{\bgamma}_t\big]$ \newline (\cref{rmk:recall_gradient-Hessian})
      & $V_t^u$ \newline (\cite[\cref*{M-def:AC2_curve}]{fornasierNewtonianMotionCompanion})
      & $V_t^{\ful{\mu}}(\sigma) \coloneqq \big[(x, v) \mapsto x + (\sigma - t) v\big]_\# \ful{\mu}_t$ \newline (\cref{def:tangent_physical}, \cref{rmk:tangent_definition_passes_downstairs}) \\
      $H_0^{\tilde{E}}(V_t^{\bgamma})$ \newline (\cref{eq:hessian_0})
      & $H_0^F(V_t^u)$ \newline (\cite[\cref*{M-def:slopes_and_H}]{fornasierNewtonianMotionCompanion})
      & $H_0^E(\ful{\mu}_t)$ \newline (\cref{def:tangent_physical}, \cref{rmk:tangent_definition_passes_downstairs}) \\
      \bottomrule
    \end{tabular}
    \caption{Second-order notions in $L^2$, their metric counterparts from \cite{fornasierNewtonianMotionCompanion}, and their avatars in the phase space $\ful{\Pp}_2$. We denote curves in $L^2$ by $\protect\bgamma$, in $\mathcal{S}$ by $u$ and in $\ful{\Pp}_2$ by $\ful{\mu}_\bullet$. Energies are denoted by $\tilde{E}$ in $L^2$, by $F$ in $\mathcal{S}$ and by $E$ in $\ful{\Pp}_2$.}
    \label{tab:second_order_equivalences}
  \end{table}

\subsubsection{Physical Curves} \label{sec:physical_connection}
In~\cite{brigatiKineticOptimalTransport2025}, the authors work towards kinetic optimal transport, starting with the key quantity $\ful{d}_T$ on $\ful{\Pp}_2$ (for consistency, we adapt their notation slightly):
\begin{equation*}
  \ful{d}^2_T(\ful{\mu}, \ful{\nu}) \coloneqq \inf_{\pi \in \Pi(\ful{\mu}, \ful{\nu})} \left( \int 12 \left | \frac{y-x}{T} - \frac{v + w}{2} \right |^2 + |w - v|^2 \, \diff \pi\big( (x,v),(y,w) \big)  \right), \qquad \text{for } \ful{\mu}, \ful{\nu} \in \ful{\Pp}_2.
\end{equation*}
We recall that the expression inside the infimum in the formula above comes as the solution of an optimal control problem with fixed end-points in the phase space. 
More precisely, given fixed and terminal points $(x, v),  (y, w) \in \R^d \times \R^d$ and $T > 0$, and setting $H = \{\gamma \in W^{2,2}([0, T]; \R^d) \mid \gamma_0 = x, \dot\gamma_0 = v, \gamma_T = y, \dot{\gamma}_T = v\}$, we have that
\begin{equation*}
  \inf_{\gamma \in H} \left\{\int_0^T |\ddot{\gamma}_r|^2 \, \diff r \right\} = 12 \left | \frac{y-x}{T} - \frac{v + w}{2} \right |^2 + |w - v|^2.
\end{equation*}

We start by collecting some of the definitions and results from \cite{brigatiKineticOptimalTransport2025}.
\begin{definition}[\cite{brigatiKineticOptimalTransport2025}, Definition 1.3] \label{def:physical_curves}
  Let $\ful{\mu}_\bullet\colon (a, b) \to (\ful{\Pp}, W_2)$ be an absolutely continuous curve. We say that $\ful{\mu}_t$ is physical if, in addition, there exists a function $m \in L^2(a, b)$ such that, for all $s, t \in (a, b)$ with $s < t$,  it holds true that
  \begin{equation*}
    \ful{d}_{t-s}(\ful{\mu}_s, \ful{\mu}_t) \leq \int_s^t m(r) \, \diff r.
  \end{equation*}
  We denote the class of curves $\ful{\mu}_\bullet\colon (a, b) \to \ful{\Pp}_2$ which are physical curves by $\mathrm{PC}^2((a, b); \ful{\Pp}_2)$.
\end{definition}

\begin{theorem}[\cite{brigatiKineticOptimalTransport2025}, Theorem 1.4]\label{theorem:cont_equation_physical}
  The following hold true.
  \begin{enumerate}
    \item If $(\ful{\mu}_t, F)$ is a weak solution to the Vlasov equation
    \begin{equation}\label{eq:vlasov}
      \partial_t \ful{\mu}_t + v \cdot \nabla_x \ful{\mu}_t + \nabla_v \cdot (F_t \ful{\mu}_t) = 0
    \end{equation}
    on $(a, b)$ for some force field $F\colon (a, b) \times \R^{2d} \to \R^{d}$ such that
    \begin{equation*}
      \int_a^b \big(\norm{v}_{L^2(\ful{\mu}_t)}^2 + \norm{F_t}^2_{L^2(\ful{\mu}_t)} \big) \, \diff t < \infty
    \end{equation*}
    then the curve $\ful{\mu}_\bullet$ is physical with
    \begin{equation*}
      m(t) = 2 \norm{F_t}_{L^2(\ful{\mu}_t)}.
    \end{equation*}
    \item Assume that $\ful{\mu}_t$ is a physical curve. Then, there exists a vector field $F\colon (a, b) \times \R^{2d} \to \R^{d}$ with $\norm{F_t}_{L^2(\ful{\mu}_t)} \leq m(t)$ for a.e. $t \in (a, b)$, such that $(\ful{\mu}_t, F)$ is a weak solution to Vlasov's equation~\eqref{eq:vlasov}, and we have the limit
    \begin{equation*}
      \lim_{h \downarrow 0} \frac{\ful{d}_h(\ful{\mu}_t, \ful{\mu}_{t + h})}{h} = \norm{F_t}_{L^2(\ful{\mu}_t)}
    \end{equation*}
    for a.e. $t \in (a, b)$.
  \end{enumerate}
\end{theorem}
\begin{proposition}[\cite{brigatiKineticOptimalTransport2025}, Proposition 4.8; cited from~\cite{ambrosioGradientFlows2008}, 8.2.1]\label{prop:superpos_vlasov}
	Let $\ful{\mu}_\bullet\colon (a, b) \to \ful{\Pp}_2$ be a narrowly continuous curve and let $F\colon (a, b) \times \R^{2d} \to \R^d$ be a vector field such that $(\ful{\mu}_t, F)$ is a weak solution to~\eqref{eq:vlasov} on $(a,b)$, and \begin{equation}
		\int_a^b \int_{\R^{2d}} \bigl( \, |v|^2 + |F_t|^2 \bigr) \diff \ful{\mu}_t \, \diff t <\infty.
	\end{equation}
	Then, there exists $\boldsymbol{\eta} \in \mathcal{P}\big(\R^{2d} \times W^{2,2}(a,b; \R^d)\big)$  
	such that 
	\begin{enumerate}
		\item\label{cond1} the measure $\boldsymbol \eta$ is supported on triples $(x,v,\gamma)$ such that $(\gamma,\dot\gamma)$ is an absolutely continuous curve solving
    \begin{equation*}
      (\gamma_a,\dot\gamma_a) = (x, v), \qquad \partial_t (\gamma_t,\dot\gamma_t) = (\dot\gamma_t, F_t(\gamma_t, \dot\gamma_t)).
    \end{equation*}
		\item we have 
		\begin{equation}\label{eq:superpos_vlasov_eq_1}
			\big(\ful{e}_t\big)_\# \boldsymbol \eta = \ful{\mu}_t,  \qquad \text{for all } t \in (a,b),
		\end{equation}
    where $\ful{e}_t\colon W^{2,2} \to \R^{2d}$ is defined by $\gamma \mapsto (\gamma_t, \dot{\gamma}_t)$.
	\end{enumerate}
	Conversely, any $\boldsymbol \eta \in \mathcal{P}\big(\R^{2d} \times W^{2,2}(a,b;\R^d)\big)$ satisfying Condition~\ref{cond1} and
	\begin{equation}
		\int_a^b \int_{ W^{2,2}(a,b;\R^d)} \big(\, |\dot\gamma_t|^2 + |F_t(\gamma_t,\dot\gamma_t)|^2 \big) \diff \boldsymbol \eta \, \diff t <\infty
	\end{equation}
	induces a solution to~\eqref{eq:vlasov} via
	\begin{equation}
		\ful{\mu}_t \coloneqq  \big(\ful{e}_t\big)_\# \boldsymbol \eta, \qquad  t \in (a,b).
	\end{equation}
\end{proposition}
Using the above theorems, it is easy to see that physical curves are exactly those descending from $W^{2,2}$-curves in $L^2$.
\begin{corollary} \label{cor:lifting_physical_curves}
  Let $(\Omega, \mathcal{F}, \mathbb{P})$ be a probability space which admits a random variable $Z$ such that $\law(Z) = \Leb|_{[0, 1]}$, and set $L^2 \coloneqq L^2(\Omega; \R^d)$. An absolutely continuous curve $\ful{\mu}_\bullet \colon (a, b) \to \ful{\Pp}_2$ is physical in the sense of~\cite{brigatiKineticOptimalTransport2025} if and only if it arises as the composition
  \begin{equation}\label{eq:lifting_equation}
  \ful{\mu}_\bullet = \ful{\law} \circ (\bgamma, \dot{\bgamma}),
  \end{equation}
  where $\bgamma$ is a $W^{2,2}$-curve in $L^2$. 
  More precisely, for each physical curve $\ful{\mu}_\bullet$, there exists a $W^{2,2}$-curve $\bgamma$ in $L^2$ fulfilling \eqref{eq:lifting_equation} such that 
  $\|\ddot{\bgamma}_t\|_{L^2} = \lim_{h \downarrow 0} \frac{\ful{d}_h(\ful{\mu}_t, \ful{\mu}_{t + h})}{h}$ for almost every $t \in (a, b)$. \\
  Conversely, for each $W^{2,2}$-curve $\bgamma$ in $L^2$, the composition $\ful{\mu}_t$ defined as in \eqref{eq:lifting_equation} is physical and $\lim_{h \downarrow 0} \frac{\ful{d}_h(\ful{\mu}_t, \ful{\mu}_{t + h})}{h} \leq \|\ddot{\bgamma}_t\|_{L^2}$, for almost every $t$ in $(a, b)$.
\end{corollary}
\begin{proof}
  To see that each physical curve has a lift as a $W^{2,2}$-curve in $L^2$, let us fix some physical curve $\ful{\mu}_\bullet$. 
  Then, by the second point in~\cref{theorem:cont_equation_physical}, there exists a vector field $F$ such that $(\ful{\mu}_t, F)$ is a weak solution to~\eqref{eq:vlasov}, $\norm{F_t}_{L^2(\ful{\mu}_t)} \leq m(t)$ and such that
    \begin{equation*}
      \lim_{h \downarrow 0} \frac{\ful{d}_h(\ful{\mu}_t, \ful{\mu}_{t + h})}{h} = \norm{F_t}_{L^2(\ful{\mu}_t)}
    \end{equation*}
    for a.e.~$t \in (a, b)$. 
    Now, picking $\boldsymbol{\eta}$ as in \cref{prop:superpos_vlasov}, we can set $X$ to be a measurable function from $\Omega$ to $W^{2,2}((a, b); \R^d)$ such that $\law(X) = \boldsymbol{\eta}$---which we can do, as $W^{2,2}((a, b); \R^d)$ is a standard Borel space\footref{footnote:realize_prop_measure}. 
    We define $\bgamma\colon (a, b) \to L^2$ as $\bgamma_t \coloneqq X(\cdot)(t)$. 
    Then, we first calculate, for $s, t \in (a, b)$ such that $s < t$:
    \begin{equation}\label{eq:second_order_first_minkowski}
      \begin{split}
      \bigg(\int |X(\omega)'(t) - X(\omega)'(s)|^2  &\, \diff \mathbb{P}(\omega)\bigg)^{\frac{1}{2}} \leq \bigg(\int \left| \int_s^t |X(\omega)''(r)|\, \diff r\right|^2 \, \diff \mathbb{P}(\omega)\bigg)^{\frac{1}{2}} \\
      &\overset{(*)}{\leq} \int_s^t\left( \int |X(\omega)''(r)|^2 \, \diff \mathbb{P}(\omega)\right)^{\frac{1}{2}} \, \diff r = \int_s^t \|F_r\|_{L^2(\ful{\mu}_r)} \, \diff r,
      \end{split}
    \end{equation}
    where we use Minkowski's integral inequality in $(*)$ and justify measurability as in the construction of horizontal lifts of Wasserstein-$AC$ curves in the beginning of \cref{sec:lifting_gradient_problems}. Then, we use \cref{eq:second_order_first_minkowski} to show that $\dot{\bgamma}_t = X(\cdot)'(t)$ is in $L^2$.
    Namely, using Minkowski again, we get
    \begin{equation*}
    \begin{split}
      \bigg(\int \bigg| & \frac{\bgamma_{t+h}(\omega) -   \bgamma_t(\omega)}{h}  - X(\omega)'(t)\bigg|^2  \, \diff \mathbb{P}(\omega) \bigg)^{\frac{1}{2}} 
      =  \frac{1}{h} \bigg( \int \bigg| \int_t^{t+h} X(\omega)'(r) - X(\omega)'(t) \, \diff r \bigg|^2 \, \diff \mathbb{P}(\omega)  \bigg)^{\frac{1}{2}} \\
      &\leq \frac{1}{h}  \int_t^{t+h} \left(\int |X(\omega)'(r) - X(\omega)'(t)|^2 \, \diff \mathbb{P}(\omega) \right)^{\frac{1}{2}} \, \diff r \leq \int_t^{t+h} \|F_r\|_{L^2(\ful{\mu}_r)} \, \diff r \to 0.
      \end{split}
    \end{equation*}
    Thus, by~\eqref{eq:superpos_vlasov_eq_1}, we have that
    \begin{align*}
    \ful{\law}\big((\bgamma_t, \dot{\bgamma}_t)\big) = \law\big((X(\cdot)(t), X(\cdot)'(t))\big) = \big(\ful{e}_t\big)_\# \boldsymbol \eta = \ful{\mu}_t,
    \end{align*}
    where $\ful{e}_t\colon W^{2,2} \to \R^{2d}$ is defined by $\gamma \mapsto (\gamma_t, \dot{\gamma}_t)$.
    Next, we use \cref{eq:second_order_first_minkowski} again to see that for $s, t \in (a, b)$ such that $s < t$:
    \begin{equation}\label{eq:minkowski_integral_ac2}
      \begin{split}
      \|\dot{\bgamma}_t - \dot{\bgamma}_s\|_{L^2} &= \left(\int |X(\omega)'(t) - X(\omega)'(s)|^2 \, \diff \mathbb{P}(\omega)\right)^{\frac{1}{2}} \leq \int_s^t \|F_r\|_{L^2(\ful{\mu}_r)} \, \diff r.
      \end{split}
    \end{equation}
    We deduce that $\bgamma$ is a $W^{2,2}$-curve with 
    $$\|\ddot \bgamma_t \|_{L^2} \leq \|F_t\|_{L^2(\ful{\mu}_t)} = \lim_{h \downarrow 0} \frac{\ful{d}_h(\ful{\mu}_t, \ful{\mu}_{t + h})}{h}$$
    for almost every $t \in (a, b)$. The converse inequality follows from the arguments detailed below. \\
    For the converse statement, let us fix some $W^{2,2}$-curve $\bgamma$ in $L^2$. 
    Setting $\ful{\mu}_\bullet \coloneqq \ful{\law} \circ (\bgamma, \dot{\bgamma})$, we clearly have
    \begin{equation} \label{eq:estimate_of_second_derivative}
      \ful{d}_h^2(\ful{\mu}_t, \ful{\mu}_{t+h}) \leq 12 \norm{\frac{\bgamma_{t + h} - \bgamma_{t}}{h} - \frac{\dot{\bgamma}_t + \dot{\bgamma}_{t + h}}{2}}^2_{L^2} + \norm{\dot{\bgamma}_{t + h} - \dot{\bgamma}_t}^2_{L^2}.
    \end{equation}
    For the first summand on the right hand side, we calculate
    \begin{equation} \label{eq:physical_curve_estimate}
      \begin{split}
      \bigg{\|}\frac{\bgamma_{t + h} - \bgamma_{t}}{h} - \frac{\dot{\bgamma}_t + \dot{\bgamma}_{t + h}}{2} \bigg{\|}_{L^2} &= 
      \frac{1}{h} \bigg{\|} \int_{t}^{t+\frac{h}{2}} \int_t^s \ddot{\bgamma}_\tau \, \diff \tau \, \diff s - \int_{t + \frac{h}{2}}^{t + h} \int_s^{t + h} \ddot{\bgamma}_\tau \, \diff \tau \, \diff s \, \bigg{\|}_{L^2} \\
      &=  \frac{1}{h} \bigg{\|} \int_t^{t+h} \left(t + \frac{h}{2} - s\right) \ddot{\bgamma}_s \diff s \, \bigg{\|}_{L^2} \\
      &\leq \frac{1}{2} \int_t^{t+h} \|\ddot{\bgamma}_s \|_{L^2} \, \diff s.
      \end{split}
    \end{equation}
    Hence, we have
    \begin{equation*}
      \ful{d}_h^2(\ful{\mu}_t, \ful{\mu}_{t+h}) \leq
      3 \bigg(\int_t^{t+h} \|\ddot{\bgamma}_s\|_{L^2} \, \diff s \bigg)^2 + \bigg(\int_t^{t+h} \|\ddot{\bgamma}_s\|_{L^2} \, \diff s \bigg)^2 \leq 4 \bigg(\int_t^{t+h} \|\ddot{\bgamma}_s\|_{L^2} \, \diff s \bigg)^2.
    \end{equation*}
    If we take the square root in the above equation, we deduce that $\ful{\mu}_\bullet$ is a physical curve with $m(t) = 2 \|\ddot{\bgamma}_t\|_{L^2}$ ($\ful{\law}$ being $1$-Lipschitz means that $\ful{\mu}_\bullet$ is also $W_2$-absolutely continuous). 
    Furthermore for any Lebesgue point $t$ of $\ddot{\bgamma}$, we use the first two equations from \cref{eq:physical_curve_estimate} to derive
    \begin{align*}
      \frac{1}{h} \bigg{\|}\frac{\bgamma_{t + h} - \bgamma_{t}}{h} - \frac{\dot{\bgamma}_t + \dot{\bgamma}_{t + h}}{2} \bigg{\|}_{L^2} 
      &= \frac{1}{h^2} \bigg{\|} \int_t^{t+h} \left(t + \frac{h}{2} - s\right) \left(\ddot{\bgamma}_s - \ddot{\bgamma}_t\right)\diff s  \, \bigg{\|}_{L^2}\\
      &\leq \frac{1}{2h} \int_t^{t+h} \|\ddot{\bgamma}_s - \ddot{\bgamma}_t\|_{L^2} \, \diff s \xrightarrow{h \to 0} 0.
    \end{align*}
    Thus, dividing by $h^2$ in \cref{eq:estimate_of_second_derivative} and taking the limit $h \to 0$, we obtain the claim.
\end{proof}
For easing the notation, we use the following shorthand:
\begin{definition}\label{def:horizontal_lift_physical}
  Let $\ful{\mu}_\bullet\colon (a, b) \to \ful{\Pp}_2$ be a physical curve. We define $\big|\ful{\mu}''\big|\colon (a, b) \to [0, +\infty]$ as
  \begin{equation*}
    \big|\ful{\mu}''\big|(t) \coloneqq \limsup_{h \downarrow 0} \frac{\ful{d}_h(\ful{\mu}_t, \ful{\mu}_{t + h})}{h}.
  \end{equation*}
  Furthermore, we say that a $W^{2,2}$-curve $\bgamma$ is a \emph{horizontal lift} of $\ful{\mu}_t$ if $\ful{\mu}_t = \ful{\law} \circ (\bgamma_t, \dot{\bgamma}_t)$ and $\|\ddot{\bgamma}_t\|_{L^2} = |\ful{\mu}''|(t)$ for almost every $t \in (a, b)$ (consistently with the first part of \cref{cor:lifting_physical_curves}).
\end{definition}

\subsubsection{Detailing the Projection of Newtonian Dynamics to \texorpdfstring{$\ful{\Pp}_2$}{P̲₂}}\label{sec:wasserstein-newtonian-curves}

Having recalled the definition of physical curves of \cite{brigatiKineticOptimalTransport2025} and defined the related notion of horizontal lifts of physical curves in \cref{def:horizontal_lift_physical}, we have available the avatar of the class $W^{2,2}$ in $\ful{P}_2$ as well as a notion of \emph{acceleration} of such curves, intrinsically to $\ful{P}_2$.
In this section, we will find avatars of the other terms involved in \cref{eq:def_variational_Newtonian_in_intro_of_char}, to finally arrive at the correspondence~\eqref{eq:correspondence_second_order_in_intro_char}. We start by recalling two key notions from \cite{fornasierNewtonianMotionCompanion}, within the present Hilbert space setting. 
For this, we first note that in the Hilbert space setting $(\mathcal{S}, d_{\mathcal{S}}) = (L^2, d_{\|\cdot\|})$, the notion of an $\ACt$-curve from Pillar 1 above coincides with the notion of a $W^{2,2}$ curve, and for a $W^{2,2}$-curve $\bgamma$ in $L^2$, the tangents $V_t^{\bgamma}$ are necessarily of the form $\sigma \mapsto \bgamma_t + (\sigma - t) \dot{\bgamma}_t$.

\begin{definition}[{C.f. \cite[\cref*{M-def:upper_gradient}]{fornasierNewtonianMotionCompanion}}]\label{rmk:recall_gradient-Hessian}
	A pair of functions $\tilde g\colon L^2 \to [0, +\infty]$, $\tilde H\colon \mathcal{T} L^2 \to [-\infty,+\infty]$, is called an \emph{upper gradient-Hessian pair} for the functional $\tilde E\colon L^2 \to (-\infty,+\infty]$ if, for every $W^{2,2}$-curve $\bgamma\colon (a, b)\to L^2$ such that
    \begin{itemize}
        \item we set $V_r^{\bgamma}\colon \sigma \mapsto \bgamma_r + (\sigma - r) \dot \bgamma_r$  for every $r\in (a,b)$,
        \item $t \mapsto \frac{\diff}{\diff \sigma} \tilde E(V_t^{\bgamma}(\sigma))\big|_{\sigma=t}$ exists as a finite real number for all $t \in (a, b)$ and has finite total variation on every compact subinterval of $(a, b)$,
        \item $t \mapsto \tilde H(V_t^{\bgamma})$ is in $L^1_{\mathrm{loc}}((a, b), \R)$,
    \end{itemize}
     we have that $\tilde g \circ \bgamma \colon (a,b) \to [0,+\infty]$ is Borel-measurable, and that, for almost every $s , t \in (a, b)$  with $s<t$, the following inequality holds
	\begin{equation} \label{eq:def_upper_gradient}
		\frac{\diff}{\diff \sigma} \tilde E(V_s^{\bgamma}(\sigma)) \Big|_{\sigma = s} 
        - \frac{\diff}{\diff \sigma}  \tilde E(V_t^{\bgamma}(\sigma)) \Big|_{\sigma = t} 
        \leq \int_s^t \frac{1}{2}  \|\ddot{\bgamma}_r\|_{L^2}^2\, \diff r + \int_s^t \frac{1}{2} \tilde g^2(\bgamma_r)\, \diff r - \int_s^t \tilde H(V_r^{\bgamma})\, \diff r,
	\end{equation}
	where $\ddot \bgamma$ is the second derivative (in the Sobolev sense) of $t\mapsto \bgamma_t$.
\end{definition}

\begin{definition}[{C.f. \cite[\cref*{M-def:curve_newtonian_motion}]{fornasierNewtonianMotionCompanion}}] \label{rmk:recall_curve _Newtonian}
	Let $\tilde E\colon L^2 \to (-\infty,+\infty]$ be a functional and let $\tilde g\colon L^2 \to [0, +\infty]$, $\tilde H\colon \mathcal{T} L^2 \to [-\infty,+\infty]$ be an \emph{upper gradient-Hessian pair} for $\tilde E$ (as per \cref{rmk:recall_gradient-Hessian}).
    Then, a $W^{2,2}$-curve $\bgamma\colon (a, b)\to L^2$ is called a \emph{curve of Newtonian motion} if, setting $V_r^{\bgamma}\colon \sigma \mapsto \bgamma_r + (\sigma - r) \dot \bgamma_r$  for every $r\in (a,b)$, we have
    \begin{itemize}
        \item $t \mapsto \frac{\diff}{\diff \sigma} \tilde E(V_t^{\bgamma}(\sigma))\big|_{\sigma=t}$ exists as a finite real number for all $t \in (a, b)$ and has finite total variation on every compact subinterval of $(a, b)$,
        \item $t \mapsto \tilde H(V_t^{\bgamma})$ is in $L^1_{\mathrm{loc}}((a, b), \R)$,
        \item for almost every $s , t \in (a, b)$  with $s<t$, the following inequality holds
	\begin{equation} \label{eq:def_variational_Newtonian}
		\frac{\diff}{\diff \sigma} \tilde E(V_s^{\bgamma}(\sigma)) \Big|_{\sigma = s} 
        - \frac{\diff}{\diff \sigma}  \tilde E(V_t^{\bgamma}(\sigma)) \Big|_{\sigma = t} 
        \geq \int_s^t \frac{1}{2}  \|\ddot{\bgamma}_r\|_{L^2}^2\, \diff r + \int_s^t \frac{1}{2} \tilde g^2(\bgamma_r)\, \diff r - \int_s^t \tilde H(V_r^{\bgamma})\, \diff r,
	\end{equation}
	where $\ddot \bgamma$ is the second derivative (in the Sobolev sense) of $t\mapsto \bgamma_t$.
    \end{itemize}
\end{definition}

Lastly, we recall the definition of the quantity
  \begin{equation} \label{eq:hessian_0}
    H_0^{\tilde{E}}(\xi) \coloneqq \limsup_{h \searrow 0} \frac{\tilde{E}(\xi_t) - 2 \tilde{E}(\xi_{t + h}) + \tilde{E}(\xi_{t + 2h})}{h^2}, \qquad \text{for } \ \xi \in \mathcal T_\varepsilon L^2, \quad  \xi\colon [t,t+\varepsilon] \to L^2.
  \end{equation}

We are now in a position to define the avatars in $\ful{\Pp}_2$ of the other terms involved in \cref{eq:def_variational_Newtonian_in_intro_of_char}.

\begin{definition}\label{def:tangent_physical}
  Given a physical curve $\ful{\mu}_\bullet$, the \emph{tangent to $\ful{\mu}_\bullet$} at time $t \in (a, b)$ is the absolutely continuous curve $V_t^{\ful{\mu}}\colon \R \to \Pp_2$ defined by
  \begin{equation*}
    V_t^{\ful{\mu}}(\sigma) \coloneqq \big[(x, v) \mapsto x + (\sigma - t) v\big]_\# \ful{\mu}_t.
  \end{equation*}
  Furthermore, for a functional $E\colon \Pp_2 \to \R$, and any $\ful{\nu} \in \ful{\Pp}_2$ we define
  \begin{equation}
    H_0^E(\ful{\nu}) = \frac{\diff^2}{\diff \sigma^2} E\Big(\big[(x, v) \mapsto x + \sigma v\big]_\# \ful{\nu}\Big) \Big{|}_{\sigma=0}
  \end{equation}
\end{definition}

\begin{remark} \label{rmk:tangent_definition_passes_downstairs}
  Given any $W^{2,2}$-curve $\bgamma$ in $L^2$ and setting $\ful{\mu}_\bullet \coloneqq \ful{\law} \circ (\bgamma, \dot{\bgamma})$, we have that $\law\big(V_t^{\bgamma}(\sigma)\big) = V_t^{\ful{\mu}}(\sigma)$, where $V_t^{\bgamma}$ is defined as $\sigma \mapsto \bgamma_t + (\sigma - t) \dot{\bgamma}_t$.
Furthermore, for any functional $E\colon \Pp_2 \to \R$, we have that $H_0^{E \circ \law}(V_t^{\bgamma}) = H_0^E(\ful{\mu}_t)$.
\end{remark}

Next, we intrinsically define in $\Pp_2$ what an upper gradient-Hessian pair is.

\begin{definition}\label{def:upg_wasserstein}
  A pair of functions $g\colon \Pp_2 \to [0, +\infty]$, $H\colon \ful{\Pp}_2 \to [-\infty, +\infty]$ is called an \emph{upper gradient-Hessian pair} for a functional $E\colon \Pp_2 \to \R$ if, for any physical curve $\ful{\mu}_\bullet\colon (a, b) \to \ful{\Pp}_2$ such that 
  \begin{itemize}
    \item $t \mapsto \frac{\diff}{\diff \sigma} E(V_t^{\ful{\mu}}(\sigma))\Big|_{\sigma = t}$ exists as a finite real number for all $t \in (a, b)$ and has finite total variation on every compact subinterval of $(a, b)$ (see \cref{rmk:tangent_definition_passes_downstairs} for the definition of $V_t^{\ful{\mu}}$);
    \item $t \mapsto H(\ful{\mu}_t)$ is in $L^1_{\mathrm{loc}}((a, b))$;
  \end{itemize}
   we have that $t \mapsto g(\ful{\mu}^x_t)$ is Borel and that, for almost all $s < t \in (a, b)$,
  \begin{equation*}
    \frac{\diff}{\diff \sigma} E(V_s^{\ful{\mu}}(\sigma))\Big|_{\sigma = s} - \frac{\diff}{\diff \sigma} E(V_t^{\ful{\mu}}(\sigma))\Big|_{\sigma = t} \leq \int_s^t \frac{1}{2} \big|\ful{\mu}''\big|^2(r) \, \diff r + \int_s^t \frac{1}{2} g^2(\ful{\mu}^x_r)\, \diff r - \int_s^t H(\ful{\mu}_r) \, \diff r,
  \end{equation*}
  where we denote by $\ful{\mu}^x_r$ the state marginal of $\ful{\mu}_r$ (i.e., $\ful{\mu}^x_r \coloneqq [(x, v) \mapsto x]_\# \ful{\mu}_r)$.
\end{definition}

\begin{remark}\label{rmk:ugs_passes_upstairs}
  By \cref{cor:lifting_physical_curves,rmk:tangent_definition_passes_downstairs}, identifying $\mathcal{T}(L^2)$ with $L^2 \times L^2$, a pair of functions $(g, H)$ as above is an upper gradient-Hessian pair for a functional $E$ in the sense of \cref{def:upg_wasserstein} if and only if $(g \circ \law, \tilde{H})$ is an upper gradient-Hessian pair for $E \circ \law$ in the sense of~\cite[\cref*{M-def:upper_gradient}]{fornasierNewtonianMotionCompanion},
where $\tilde{H}\colon \xi \mapsto H\big(\ful{\law}(\xi_t, \xi'_t)\big)$, for every geodesic $\xi \colon [t, t + \varepsilon) \to L^2$.
\end{remark}

\begin{definition}\label{def:wasserstein_newtonian_curve}
  Let $E\colon \Pp_2 \to \R$ be a functional and let $(g, H)$ be an upper gradient-Hessian pair in the sense of \cref{def:upg_wasserstein} for $E$. A physical curve $\ful{\mu}_\bullet\colon (a, b) \to \ful{\Pp}_2$ is called a \emph{Wasserstein-Newtonian curve} for the functional $E$ with respect to the upper gradient-Hessian pair $(g, H)$ if we have
  \begin{itemize}
    \item $t \mapsto \frac{\diff}{\diff \sigma} E(V_t^{\ful{\mu}}(\sigma))\Big|_{\sigma = t}$ exists as a finite real number for every $t \in (a, b)$ and has finite total variation on every compact subinterval of $(a, b)$ (see \cref{rmk:tangent_definition_passes_downstairs} for the definition of $V_t^{\ful{\mu}}$);
    \item $t \mapsto H(\ful{\mu}_t)$ is in $L^1_{\mathrm{loc}}((a, b))$;
    \item for almost all $s < t \in (a, b)$, the following inequality holds
  \begin{equation*}
    \frac{\diff}{\diff \sigma} E(V_s^{\ful{\mu}}(\sigma))\Big|_{\sigma = s} - \frac{\diff}{\diff \sigma} E(V_t^{\ful{\mu}}(\sigma))\Big|_{\sigma = t} \geq \int_s^t \frac{1}{2} \big|\ful{\mu}''\big|^2(r) \, \diff r + \int_s^t \frac{1}{2} g^2(\ful{\mu}^x_r)\, \diff r - \int_s^t H(\ful{\mu}_r) \, \diff r.
  \end{equation*}
  \end{itemize}\
\end{definition}

  We now have the tools to make precise the correspondence~\eqref{eq:correspondence_second_order_in_intro_char}.

\begin{proposition} \label{rmk:wasserstein_newtonian_curve_lifts}
  For a physical curve $\ful{\mu}_\bullet$ and a horizontal lift $\bgamma$ as defined in \cref{def:horizontal_lift_physical}, $\ful{\mu}_\bullet$ is a Wasserstein-Newtonian curve for the functional $E$ with respect to the upper gradient-Hessian pair $(g, H)$ as per \cref{def:wasserstein_newtonian_curve} if and only if $\bgamma$ is a curve of Newtonian motion for the functional $E \circ \law$ with respect to the upper gradient-Hessian pair $(g \circ \law, \tilde{H})$ as per \cref{rmk:recall_curve _Newtonian}, where $\tilde{H}$ is defined as in \cref{rmk:ugs_passes_upstairs}. \\
  Moreover, for any curve of Newtonian motion $\bgamma$ for the functional $E \circ \law$ with respect to the upper gradient-Hessian pair $(g \circ \law, \tilde{H})$, the composition $\ful{\mu}_\bullet \coloneqq \ful{\law} \circ (\bgamma, \dot{\bgamma})$ is a Wasserstein-Newtonian curve for the functional $E$.
\end{proposition}
\begin{proof}
  This follows directly from \cref{cor:lifting_physical_curves,rmk:tangent_definition_passes_downstairs,rmk:ugs_passes_upstairs}.
\end{proof}

\subsection{Constructing Curves of Newtonian Motion in \texorpdfstring{$\LtwoR$}{L²\_ℝ}} \label{sec:constructing_newtonian_curves_in_l2}

In this section, we fix a potential energy $\tilde{E} = E \circ \law\colon \LtwoR \to \R$, and---under suitable assumptions---we construct curves of Newtonian motion (or \emph{Newtonian flows}) $u \colon [0, +\infty) \to \LtwoR$.
For ``Newtonian flows'' we mean second-order dynamics, which fulfill a suitable variational characterization (see \cref{eq:recall_curves_of_newtonian_motion} and the related discussion).
While for first-order gradient flows, this characterization has been fully understood \cite{ambrosioGradientFlows2008}, a proper notion of second order dynamics of Newtonian flows has remained so far open. 
In the companion paper \cite{fornasierNewtonianMotionCompanion} we investigated the construction of Newtonian flows in metric spaces. 
Here, we specify its results to the case of functionals over $L^2_{\mathbb R}$, i.e., the space obtained through the Lions lift of $\Pp_2$.

Before describing the compactness procedure in $\LtwoR$, we recall the machinery introduced in~\cite{fornasierNewtonianMotionCompanion}. To this end, for a fixed $U_0 \in \LtwoR$, $U'_0 \in \LtwoR$ and $\tau > 0$, we say that $(U_{\tau, k})_{k=-1}^\infty$ is a \emph{Newtonian sequence} whenever we have, for all $k \in \NN$,
\begin{equation} \label{eq:newtonian_sequence}
  U_{\tau, k+1} \in J_{\tau^2, \LtwoR}[2 U_{\tau, k} - U_{\tau, k-1}], \qquad J_{\tau^2, \LtwoR}[X] \coloneqq \argmin_{Y \in \LtwoR} \Big{\{}\tilde{E}(Y) + \frac{1}{2\tau^2}\|X - Y\|_{\LtwoR}^2\Big{\}}.
\end{equation}
For such a sequence, we define the velocities $(U'_{\tau, k})_{k=0}^{\infty}$ as
\begin{equation}\label{eq:L2_newtonian_seq_velocities}
  U'_{\tau, k} \coloneqq \frac{U_{\tau, k} - U_{\tau, k-1}}{\tau} \in \LtwoR.
\end{equation}
We provide an explicit description of some of the objects defined in~\cite{fornasierNewtonianMotionCompanion} in our present framework. 
We reiterate that some of the objects below---in particular, the point $Z_{\tau, k}$ in \cref{eq:wasserstein_u} and the approximating tangents in \cref{eq:v_wasserstein,eq:tilde_v_wasserstein,eq:bar_v_wasserstein}---are well-defined due to the extendability of geodesics in $\LtwoR$; this would not be the case when working in $\Pp_2$.
\begin{remark} \label{rmk:recall_objects_Newtonian}
  In our linear setting, we have
  \begin{align}
    U_{\tau, k+1} &\in J_{\tau^2, \LtwoR}[Z_{\tau, k}],  \qquad \text{where } Z_{\tau, k} \coloneqq U_{\tau, k} + \tau U'_{\tau, k}, \label{eq:wasserstein_u} \\
    U'_{\tau, k+1} &=  \frac{1}{\tau}\Big(U_{\tau, k+1} - U_{\tau, k}\Big) = U'_{\tau, k} + \frac{1}{\tau} \Big(U_{\tau, k+1} - Z_{\tau, k}\Big).\label{eq:wasserstein_uprime}
  \end{align}
  We recall that curves $u_\tau$ and $\tilde u_\tau$
\begin{align}
		u_\tau & \colon [0, +\infty) \to \LtwoR,\qquad u_\tau(\tau i + \delta) = U_{\tau, i} \qquad & i \in \NN, \delta \in [0, \tau),  \label{eq:piecewise_constant_interpolant}\\
		\tilde{u}_\tau & \colon [0, +\infty) \to \LtwoR,\qquad \tilde{u}_\tau(\tau i + \delta) \in J_{\delta \cdot \tau, \LtwoR}[Z_{\tau,i}], \qquad & i \in \NN, \delta \in [0, \tau). \label{eq:de_giorgi_interpolant}
\end{align}
and also define
\begin{equation}
  \dot{u}_\tau \colon [0, +\infty) \to \LtwoR, \qquad \dot{u}_\tau(\tau i + \delta) = U'_{\tau, i} = \frac{1}{\tau}\big(U_{\tau, i} - U_{\tau, i-1}\big), \qquad i \in \NN, \delta \in [0, \tau). \label{eq:piecewise_interpolants_velocity}
\end{equation}
In the present framework, we observe that the space of geodesics  $\mathcal T L^2_{\mathbb{R}}$ can be described as $\mathcal T L^2_{\mathbb{R}} \coloneqq \bigcup_{\varepsilon>0}\mathcal T_{\varepsilon} L^2_{\mathbb{R}}$, where
\begin{equation*}
    \mathcal T_{\varepsilon} L^2_{\mathbb{R}} \coloneqq
    \big\{
    \xi\colon [t,t+\varepsilon] \to  L^2_{\mathbb{R}} \mid \exists t\in \R, \, u_0,v_0\in L^2_{\mathbb{R}} : \xi_s = u_0 + s\, v_0 \ \forall s\in [t,t+\varepsilon]
    \big\}.
\end{equation*}
Moreover, regarding the family of ``approximate tangents'' $\mathcal{V}^{\tau, \varepsilon}, \tilde{\mathcal{V}}^{\tau, \varepsilon}, \bar{\mathcal{V}}^{\tau, \varepsilon}$ originally defined in \cite{fornasierNewtonianMotionCompanion},
we have
	\begin{alignat}{3}
		&\mathcal{V}^{\tau, \varepsilon} \colon [0, +\infty)  \to \mathcal T_{\varepsilon} L^2_{\mathbb{R}},  \quad \mathcal{V}^{\tau, \varepsilon}_t\colon [t,t+\varepsilon]\to L^2_{\mathbb{R}}, \quad  t \in [0, +\infty), \nonumber \\  
        &\mathcal{V}^{\tau, \varepsilon}_{\tau i + \delta}  \coloneqq \bigg[\tau i + \delta + \sigma \mapsto U_{\tau,i-1} + \big(2\delta + \sigma\big) U'_{\tau,i} \bigg]\quad & i \in \NN, \, \delta \in [0, \tau), \, \sigma \in [0, \varepsilon], \label{eq:v_wasserstein} \\
        &\tilde{\mathcal{V}}^{\tau, \varepsilon} \colon [0, +\infty)  \to \mathcal T_{\varepsilon} L^2_{\mathbb{R}}, \quad \tilde{\mathcal{V}}^{\tau, \varepsilon}_t\colon [t,t+\varepsilon]\to L^2_{\mathbb{R}}, \quad  t \in [0, +\infty),\nonumber \\  
        &\tilde{\mathcal{V}}^{\tau, \varepsilon}_{\tau i + \delta}  \coloneqq \bigg[\tau i + \delta + \sigma \mapsto U_{\tau,i-1} + \big(2\psi(\delta) + \sigma\big) U'_{\tau,i} \bigg]\quad & i \in \NN,\, \delta \in [0, \tau),\, \sigma \in [0, \varepsilon],  \label{eq:tilde_v_wasserstein}\\
        &\bar{\mathcal{V}}^{\tau, \varepsilon} \colon [0, +\infty)  \to \mathcal T_{\varepsilon} L^2_{\mathbb{R}}, \quad \bar{\mathcal{V}}^{\tau, \varepsilon}_t  \colon[t,t+\varepsilon] \to  L^2_{\mathbb{R}}, \quad  t \in [0, +\infty),\nonumber \\  
        &\bar{\mathcal{V}}^{\tau, \varepsilon}_{\tau i + \delta}  \coloneqq \bigg[\tau i + \delta + \sigma \mapsto U_{\tau,i-1} + \sigma U'_{\tau, i} \bigg]\quad & i \in \NN, \, \delta \in [0, \tau), \, \sigma \in [0, \varepsilon], \label{eq:bar_v_wasserstein}
	\end{alignat}
  where we use the following definition of $\psi$:
  	\begin{equation} \label{eq:def_psi}
		\psi(\delta) \coloneqq \begin{cases}
			\phantom{\tau - {}}\sqrt{\frac{\delta \, \tau}{2}}, \qquad & \delta \leq \frac{\tau}{2}, \\
			\tau - \sqrt{\frac{\tau^2 - \delta \, \tau}{2}}, \qquad & \delta > \frac{\tau}{2},
		\end{cases}
		\quad \text{so that } \quad \psi'(\delta) = \begin{cases}
			\frac{\tau}{4} \frac{1}{\psi(\delta)}, \qquad & \delta \leq \frac{\tau}{2}, \\
			\frac{\tau}{4}\frac{1}{\tau - \psi(\delta)}, \qquad & \delta > \frac{\tau}{2}.
		\end{cases}
	\end{equation}
  We recall that in the general construction as detailed in \cite{fornasierNewtonianMotionCompanion}, the extendability of geodesics is necessary to ensure the existence of a uniform $\varepsilon > 0$ such that the interpolating tangents $\mathcal{V}^{\tau, \varepsilon}$, $\tilde{\mathcal{V}}^{\tau, \varepsilon}$ and $\bar{\mathcal{V}}^{\tau, \varepsilon}$ are well-defined. In $\LtwoR$, geodesics can be extended to arbitrary length, unlocking the possibility of setting $\varepsilon = +\infty$.
  Moreover, we recall that
    \begin{equation*}
        N_\tau(x, y, z)  \coloneqq \frac{\|z - 2y +x\|_{\LtwoR}}{\tau^2}, \qquad  x, y, z \in \LtwoR, \  \tau > 0, 
    \end{equation*}
    and that, for every $i \in \NN$ and $ \delta \in [0, \tau)$,
    	\begin{align}
		\big|u_\tau^{(1)}\big| & \colon [0, +\infty) \to \R_+, \qquad \big|u_\tau^{(1)}\big|(\tau i + \delta) = \|U'_{\tau, i+1}\|_{\LtwoR},  \notag \\
		\big|u_\tau^{(2)}\big| & \colon [0, +\infty) \to \R_+, \qquad \big|u_\tau^{(2)}\big|(\tau i + \delta) = N_\tau(U_{\tau, i-1}, U_{\tau, i}, U_{\tau, i+1})= \frac{\|U'_{\tau, i+1} - U'_{\tau, i}\|_{\LtwoR}}{\tau},   \label{eq:second_derivative_metric}\\
		d(\tilde{E}u_\tau) & \colon [0, +\infty) \to \R, \qquad d(\tilde{E}u_\tau)(\tau i + \delta) = \frac{1}{\tau} (\tilde{E}(U_{\tau, i}) - \tilde{E}(U_{\tau, i-1})). \label{eq:dE_def}
  \end{align}
  For every $i \in \NN$ and $\delta \in (0, \tau]$, we have
  \begin{equation} \label{eq:G_tau}
		|G_\tau|  \colon [0, +\infty) \to \R_+, \qquad |G_\tau|(\tau i + \delta) = \frac{\sup_{Y \in J_{\delta \cdot \tau, \LtwoR}[Z_{\tau,i}]} \| Y - Z_{\tau, i} \|_{\LtwoR} }{\delta \tau}. 
	\end{equation}
  Here, as in~\cite{fornasierNewtonianMotionCompanion}, we set $|G_\tau|(0) = 0$, and we notice that only its equivalence class in $L^1$ (and not its pointwise values) matters.
\end{remark}

Once we are given a Newtonian sequence in the sense of \cref{eq:newtonian_sequence}, we would like to take the limit as $\tau \to 0$ to obtain a limit curve. This is the content of the following lemma:

\begin{lemma} \label{lem:compactness_newtonian}
  Let $(\tau_i)_{i=1}^\infty$ be a sequence of step sizes such that $\tau_i \to 0$ as $i \to \infty$, and let us assume that, for every $i$, $(U_{\tau_i, n})_{n=-1}^{\infty}$ is a Newtonian sequence. 
  Let us furthermore assume that for every $T > 0$, there exists some $C_T > 0$ such that for all $i \in \NN$ and all $0 \leq n \leq T/\tau_i$, we have that
  \begin{equation}\label{eq:required_bounds_for_compactness}
    \|U_{\tau_i, n}\|_{\LtwoR} \leq C_T, \quad \|U'_{\tau_i, n}\|_{\LtwoR} \leq C_T, \quad \int_0^T |u_{\tau_i}^{(2)}|^2(t) \, \diff t \leq C_T.
  \end{equation}
  Then, there exists a subsequence (not relabeled) $(\tau_i)_{i=1}^{\infty}$, a $W^{2,2}_{\mathrm{loc}}$-curve $\bgamma\colon [0, +\infty) \to \LtwoR$ and a function $A \in L^2_{loc}([0, +\infty); \R)$ such that
  \begin{gather}
    \Big(u_{\tau_i}(t), \dot{u}_{\tau_i}(t)\Big) \xrightarrow{\ful{\law}^{-1}(\sigma_{\text{narrow}})} \Big(\bgamma_t, \dot{\bgamma}_t\Big), \qquad \text{for all } t \in [0, +\infty) ,\label{eq:pointwise_convergence_wasserstein}\\
    \big|u^{(2)}_{\tau_i}\big| \rightharpoonup A, \quad A(t) \geq \big\|\ddot{\bgamma}_t\big\|_{\LtwoR}, \qquad \text{for a.e. } t \in [0, +\infty), \label{eq:limit_is_ac2_wasserstein}
  \end{gather}
  where the convergence in~\eqref{eq:limit_is_ac2_wasserstein} is understood in the sense of the weak topology of $L^2_{\mathrm{loc}}([0, +\infty))$.
  Furthermore, if there exist $U_0, U'_0 \in \LtwoR$ such that for all $\tau_i$, $U_{\tau_i, 0} = U_0$ and $U_{\tau_i, -1} = U_0 - \tau_i U'_0$, then $\bgamma_0 = U_0$ and $\dot{\bgamma}_0 = U'_0$.
\end{lemma}
\begin{proof}
  By the assumptions, we can find a subsequence (not relabeled) $(\tau_i)_{i=1}^{\infty}$ and functions $A, B \in L^2_{loc}([0, +\infty); \R)$ such that $|u_{\tau_i}^{(2)}| \rightharpoonup A$ and $|u_{\tau_i}^{(1)}| \rightharpoonup B$ weakly in $L^2_{loc}$. 
  Next, we use \cref{prop:l2_is_aa_prep} and a continuous extension argument as in \cref{prop:aa-ready-allows-aa}---note that we have
  \begin{equation*}
      (U_{\tau_i, n}, U'_{\tau_i, n}) \in \overline{B}_{\sqrt{2}C_T}(0), \quad \|(U_{\tau_i, n}, U'_{\tau_i, n}) - (U_{\tau_i, n+1}, U'_{\tau_i, n+1})\|_{\LtwoR} \leq \tau_i (|u_{\tau_i}^{(1)}|(\tau_i n) + |u_{\tau_i}^{(2)}|(\tau_i n)),
  \end{equation*}
  so that we can use the modulus $\delta(s, t) = \big|\int_s^t A(\sigma) + B(\sigma) \, \diff \sigma\big|$---to choose another subsequence and continuous limiting curves $\bgamma$, $\dot{\bgamma}$ such that
  \begin{enumerate}
    \item \cref{eq:pointwise_convergence_wasserstein} is fulfilled.
    \item $(\bgamma, \dot{\bgamma})$ is an absolutely continuous curve in $\LtwoR \times \LtwoR$.
    \item For any $q_1, \dots, q_n \in \Q$ such that $q_i \geq 0$ for all $i \in [n]$, we have
    \begin{equation}\label{eq:joint_law_wasserstein}
      \law\Big(\big(u_{\tau_i}(q_1), \dot{u}_{\tau_i}(q_1), \dots, u_{\tau_i}(q_n), \dot{u}_{\tau_i}(q_n)\big)\Big) \to \law\Big(\big(\bgamma_{q_1}, \dot{\bgamma}_{q_1}, \dots, \bgamma_{q_n}, \dot{\bgamma}_{q_n}\big)\Big)
    \end{equation}
    narrowly in $\Pp_2\big((\R^{2d})^n\big)$,
    i.e., for any bounded continuous function $f \colon (\R^{2d})^n \to \R$, we have that
    \begin{equation*} 
      \EE\bigg[f\Big(\big(u_{\tau_i}(q_1), \dot{u}_{\tau_i}(q_1)\big), \dots, \big(u_{\tau_i}(q_n), \dot{u}_{\tau_i}(q_n)\big)\Big)\bigg] \to \EE\bigg[f\Big(\big(\bgamma_{q_1}, \dot{\bgamma}_{q_1}\big), \dots, \big(\bgamma_{q_n}, \dot{\bgamma}_{q_n}\big)\Big)\bigg].
    \end{equation*}
  \end{enumerate}
  By using the absolute continuity of $(\bgamma, \dot{\bgamma})$ and the modulus of continuity $\delta(s, t)$, \cref{eq:joint_law_wasserstein} extends to arbitrary $t_1, \ldots, t_n \in [0, +\infty)$ (as in the proof of \cref{prop:aa-ready-allows-aa}).
  If we apply it to 
  \begin{equation*}
      f \colon \big( (x_1, v_1), (x_2, v_2) \big) \mapsto |x_1 + (T-s) v_1 - (x_2 + (T-t) v_2)|^2
  \end{equation*}
  and we employ a truncation argument as in the proof of \cref{cor:l2_is_aa}, we see that, for every $s, t, T \in [0, +\infty)$ with $s \leq t \leq T \leq s + \varepsilon$, we have
  \begin{equation} \label{eq:liminf_vel_Newtonian}
    \liminf_{\tau \to 0} \|\bar{\mathcal{V}}^{\tau, \varepsilon}_s(T) - \bar{\mathcal{V}}^{\tau, \varepsilon}_t(T)\|_{\LtwoR} \geq \|\bar{\mathcal{V}}_s(T) - \bar{\mathcal{V}}_t(T)\|_{\LtwoR},
  \end{equation}
  where we recall that $\bar{\mathcal{V}}^{\tau, \varepsilon}$ is defined in \cref{eq:bar_v_wasserstein}, and where we introduce the geodesics in $\LtwoR$
  \begin{equation} \label{eq:wasserstein_mathcal_v_def}
    \bar{\mathcal{V}}_s \coloneqq [\sigma \mapsto \bgamma_s + (\sigma - s) \dot{\bgamma}_s].
  \end{equation}
  By combining \cref{eq:liminf_vel_Newtonian} with the relation
  \begin{equation} \label{eq:newtonian_sequences_are_pseudo_ac}
		\big\| \bar{\mathcal{V}}^{\tau, \varepsilon}_{t^i_\tau}(T) - \bar{\mathcal{V}}^{\tau, \varepsilon}_{t^j_\tau}(T) \big\|_{\LtwoR} \leq \int_{t^i_\tau}^{t^j_\tau} (T - r) |u_\tau^{(2)}|(r) \, \diff r
	\end{equation}
  from~\cite[\cref*{M-prop:newtonian_sequences_are_proto_ac}]{fornasierNewtonianMotionCompanion}\footnote{While that proposition required further bounds, the proof reveals that those are only required if $\epsilon < +\infty$, i.e., if geodesics cannot be extended infinitely.} (here, $t_\tau^i = i \, \tau$), if we take the limit $\tau \to 0$, we obtain the inequality
  \begin{equation} \label{eq:AC2_relation_wasserstein}
      \|\bar{\mathcal{V}}_s(T) - \bar{\mathcal{V}}_t(T)\|_{\LtwoR} \leq \int_s^t (T - r) A(r) \, \diff r,
  \end{equation}
  which implies that $t\mapsto \bgamma_t$ is a $W^{2,2}$-curve (see \cite[\cref*{M-lem:euclidean_ac2_equivalence}]{fornasierNewtonianMotionCompanion})\footnote{More precisely, \cref{eq:AC2_relation_wasserstein} implies that $t\mapsto \bgamma_t$ is an \ACt-curve as per \cite[\cref*{M-def:AC2_curve}]{fornasierNewtonianMotionCompanion}, and \cite[\cref*{M-lem:euclidean_ac2_equivalence}]{fornasierNewtonianMotionCompanion} establishes the equivalence between the notion of $W^{2,2}$- and \ACt-curve.}.
  The last relation also implies the inequality in \cref{eq:limit_is_ac2_wasserstein} by virtue of \cite[\cref*{M-lem:second_derivative_bound}]{fornasierNewtonianMotionCompanion}. 
  Finally, the last part of the statement follows from the fact that the limit sequence of \cref{prop:l2_is_aa_prep} respects strong limits. 
\end{proof}

In regular enough cases---we will see an example in \cref{sec:example_newtonian}---, we can use the above compactness to obtain curves of Newtonian motion.

\begin{lemma} \label{lemma:limit_curve_is_newtonian}
  Let $E\colon \Pp_2\to \Rext$ and $|\partial E|\colon \Pp_2\to [0,+\infty]$ be, respectively, an energy and its metric slope, and let $\tilde{E} = E \circ \law \colon \LtwoR\to \Rext$ and $|\partial \tilde E|\colon \LtwoR \to [0,+\infty]$ be the respective lifts.
  Let us assume
  \begin{enumerate}
      \item $(|\partial \tilde{E}|, H_0^{\tilde{E}})$ is an upper gradient-Hessian pair (see \cref{eq:hessian_0} and \cref{rmk:recall_gradient-Hessian});
      \item There exists $\delta_0>0$ such that, for $\delta \in (0,\delta_0)$, $J_{\delta, \LtwoR}[X]$ is nonempty, for every $X \in \LtwoR$;
  \end{enumerate}
  we further assume the following regularity assumptions:
  \begin{enumerate}[label=R\arabic*.]
      \item For every $X, Y \in \LtwoR$, $s \mapsto \tilde{E}(X + s Y)$ is in $W_{\text{loc}}^{2,2}$;
      \item $\tilde{E}$ is strictly Hadamard differentiable\footnote{By strict Hadamard differentiability, we mean that for each $X_0$, there exists a continuous linear functional $\nabla \tilde{E}(X_0)$ on $\LtwoR$ such that $\lim_{X \to X_0, \tau \to 0} \frac{\tilde{E}(X + \tau Y) - \tilde{E}(X)}{\tau} = \nabla \tilde{E}(X_0)(Y)$, uniformly over $Y$ within compacts in $\LtwoR$.} (with respect to $\LtwoR$);
      \item $|\partial \tilde{E}|$ and $H_0^{\tilde{E}}$ are $\LtwoR$-continuous\footnote{Here, we identify $\mathcal{T} \LtwoR$ with $\LtwoR \times \LtwoR = \LtwoRful$, and interpret $H_0$, defined in \cref{eq:hessian_0}, as a function on $\LtwoRful$.}. \label{reg:R3}
  \end{enumerate}
   Then, in the situation of \cref{lem:compactness_newtonian}, whenever for each $t$, $(u_{\tau_i}(t))_{i \in \NN}$ and $(\dot u_{\tau_i}(t))_{i \in \NN}$ have uniformly integrable second moments, the limit curve $\bgamma$ in \cref{lem:compactness_newtonian} 
   is a curve of Newtonian motion (see \cref{rmk:recall_curve _Newtonian}).
\end{lemma}
\begin{proof}
  By \cite[\cref*{M-prop:discrete_newtonian_motion}]{fornasierNewtonianMotionCompanion}, we have the discrete identity
  \begin{equation} \label{eq:discrete_newtonian_motion}
		d(\tilde E u_{\tau_i})(t_{\tau_i}^{j_1}) - d(\tilde E u_{\tau_i})(t_{\tau_i}^{j_2}) = \frac{1}{2}\int_{t_{\tau_i}^{j_1}}^{t_{\tau_i}^{j_2}}  |u_{\tau_i}^{(2)}|^2(r)\, \diff r + \frac{1}{2} \int_{t_{\tau_i}^{j_1}}^{t_{\tau_i}^{j_2}} |G_{\tau_i}|^2(r)\, \diff r - \int_{t_{\tau_i}^{j_1}}^{t_{\tau_i}^{j_2}} H^{\tilde E}_{0}(\tilde{\mathcal{V}}^{{\tau_i}, \varepsilon}_{r}) \, \diff r,
	\end{equation}
  for all $i, j_1, j_2 \in \NN$, setting $t_{\tau_i}^j = \tau_i \, j$ with $j_1 \leq j_2$ (see \cref{eq:second_derivative_metric,eq:G_tau,eq:hessian_0} for the definitions of the objects involved).
  Our goal is to take the limit in the above equation as $\tau_i\to 0$.
  First, we observe that, if we consider the $\LtwoR$-geodesic $ \xi\colon [t,t+\varepsilon]\ni\sigma \mapsto X + (\sigma - t) Y$, we can write
  \begin{equation*}
      H_0^{\tilde{E}}(\xi) = H_0^{\tilde{E}}(X, Y) = \limsup_{h \searrow 0} \frac{\tilde{E}(X) - 2 \tilde{E}(X + hY) + \tilde{E}(X + 2hY)}{h^2} 
  \end{equation*}
  and noting that $\law(X + h Y) = ((x, y) \mapsto x + h y)_\#(\law((X, Y)))$, we see that $H_0^{\tilde{E}}$ factors through $\law$ via some $G\colon \Pp_2(\R^d \times \R^d) \to \R$; likewise, we saw in \cref{sec:lifting_gradient_problems} that $|\partial \tilde{E}| = |\partial E| \circ \law$. 
  By \cref{rmk:continuity_of_functionals}, the $\LtwoR$-continuity of $H_0^{\tilde{E}}$ and $|\partial \tilde{E}|$ are equivalent to $W_2$-continuity of $G$ and $|\partial E|$, respectively. 
  Lastly, in the same spirit, the function $D^\text{fin}\colon \LtwoR \times \LtwoR \times \R^+ \mapsto \R$ defined by taking finite differences 
  \begin{equation*}
      D^\text{fin} (X, Y, \tau) \coloneqq \frac{\tilde E(X + \tau Y) - \tilde E(X)}{\tau}
  \end{equation*}
  can be decomposed as $D^\text{fin} = D^\text{fin}_{\Pp_2} \circ (\law, \mathrm{Id})$, for some function $D^\text{fin}_{\Pp_2} \colon \Pp_2(\R^d \times \R^d) \times \R^+ \to \Rext$.
  Moreover, we claim that the strict Hadamard differentiability implies that $D_{\Pp_2}^\text{fin}(\ful{\mu}_i, \tau_i) \to \nabla \tilde{E}(X)(Y)$ whenever $\ful{\mu}_i \to \law((X, Y))$ in $W_2$ and $\tau_i \to 0$. 
  To see this, we use \cref{prop:wasserstein_distance_as_inf_of_l2_distance} to lift $\ful{\mu}_i$ to $(X_i, Y_i)$ such that $\law((X_i, Y_i)) = \ful{\mu}_i$ and such that $(X_i, Y_i) \to (X, Y)$ in $\LtwoR$. 
  Then, $K \coloneqq \{Y_i\}_{i \in \NN} \cup \{Y\}$ is a compact subset of $\LtwoR$ and by 
  the strict Hadamard differentiability, applied uniformly to directions in $K$, we get
  \begin{equation*}
      D^\text{fin} (X_i, Y_i, \tau) - \nabla \tilde{E}(X)(Y_i)
      = \frac{\tilde{E}(X_i + \tau_i Y_i) - \tilde{E}(X_i)}{\tau_i} - \nabla \tilde{E}(X)(Y_i) \to 0
  \end{equation*}
  as $i\to \infty$. Since $\nabla \tilde{E}(X)$ is a continuous linear functional and $Y_i \to Y$, the claim follows.
  By the uniformly integrable second moments as in the assumption, the convergence in \cref{eq:pointwise_convergence_wasserstein} gets upgraded to convergence in $W_2$.
  Furthermore, as $\tau_i \|\dot{u}_{\tau_i}(t)\|_{\LtwoR} \leq \tau_i C_T$, we also have the $W_2$-convergence $\ful{\law}((u_{\tau_i}(t) - \tau_i \dot{u}_{\tau_i}(t), \dot{u}_{\tau_i}(t))) \to \ful{\law}((\bgamma_t, \dot{\bgamma}_t))$. 
  Thus we get, for all $t \in [0, +\infty)$,
  \begin{equation} \label{eq:discrete_identity}
  \begin{split}
    &d(\tilde E u_{\tau_i})(t) \to \nabla \tilde{E}(\bgamma_t)(\dot{\bgamma}_t),  \qquad
    |\partial \tilde E|(u_{\tau_i}(t)) \to |\partial \tilde E|(\bgamma_t),\\
    & \qquad \qquad  \qquad  H_0^{\tilde E}(u_{\tau_i}(t), \dot{u}_{\tau_i}(t)) \to H_0^{\tilde E}(\bgamma_t, \dot{\bgamma}_t),
  \end{split}
  \end{equation}
  where we used that $d(\tilde E u_{\tau_i})(t) = D^\text{fin}\big(u_{\tau_i}(t) - \tau_i \dot{u}_{\tau_i}(t), \dot{u}_{\tau_i}(t), \tau_i\big)$.
  Next, we want to show a uniform bound on $H_0^{\tilde E}(\tilde{\mathcal{V}}^{\tau_i, \varepsilon}_{s})$. To this end, first note that by assumption, the map $(X, Y) \mapsto H_0^{\tilde{E}}(X, Y)$ is jointly continuous in $X$ and $Y$, as well as $2$-homogeneous in $Y$. 
  Thus, we have that $G\big(\ful{\law}((X, Y))\big)$ is uniformly bounded for $X$ in a neighborhood around $\bgamma_t$ and $Y \in B_{C_T}(0)$. 
  By compactness of $[0, t]$, we can even find such a bound uniformly over $X$ in some neighborhood around $\mathrm{Im}(\bgamma|_{[0, t]})$ and $Y \in B_{C_T}(0)$. 
  Hence, again considering the map $G\colon \ful{\Pp}_2 \to \R$ such that $H_0^{\tilde{E}} = G \circ \ful{\law}$, we have that $G$ is bounded uniformly on a $W_2$-neighborhood around $\mathrm{Im}(\law \circ \bgamma|_{[0, t]})$ and $\nu \in B_{C_T}(\delta_0)$ (i.e., the Wasserstein ball around a Dirac distribution centered at zero). 
  Finally, the uniformly integrable second moments of $(u_{\tau_i})|_{[0, t]}$ upgrade the convergence in  \cref{eq:pointwise_convergence_wasserstein} to pointwise $W_2$-convergence. Since $W_2 \leq \|\cdot\|_{\LtwoR}$, the bound $\|U'_{\tau_i, n}\|_{\LtwoR} \leq C_T$ in \cref{eq:required_bounds_for_compactness} ensures that this convergence is uniform on $[0, t]$. 
  Thus, for $i$ large, $\mathrm{Im}(\law \circ (u_{\tau_i})|_{[0, t]})$ is $W_2$-close to $\mathrm{Im}(\law \circ \bgamma|_{[0, t]})$; by \cref{eq:required_bounds_for_compactness,eq:tilde_v_wasserstein}, the same holds for the base point of $\tilde{\mathcal{V}}^{\tau_i, \varepsilon}_{s}$, which shows that $H_0^{\tilde E}(\tilde{\mathcal{V}}^{\tau_i, \varepsilon}_{s})$ is uniformly bounded by some $C > 0$. \\
  Next, by assumption, the set $J_{\delta \cdot \tau, \LtwoR}[Z_{\tau, i}]$ involved in the interpolation introduced in \cref{eq:de_giorgi_interpolant} is nonempty for $\tau$ small enough; in particular, $\tilde{u}_{\tau_i}$ is well-defined, and this is the condition required for \cref{eq:discrete_newtonian_motion} in~\cite[\cref*{M-prop:discrete_newtonian_motion}]{fornasierNewtonianMotionCompanion}. 
  We claim that $\law\big(\tilde{u}_{\tau_i}(s)\big) \to \law(\bgamma_s)$ in $W_2$.
  To see this, by a standard a priori estimate for the proximal map (as done in the proof of \cref{thm:cusp}), using the coercivity of $\tilde{E}$ (as per \cite[\cref*{M-ass:coercivity}]{fornasierNewtonianMotionCompanion}, which follows from $J_{\delta, \LtwoR}[X]$ being nonempty) and the boundedness of $\tilde{E}(Z_{\tau_i, n}) = E\big(\law(Z_{\tau_i, n})\big)$ (as $\law(Z_{\tau_i, n})$ converges in $W_2$ and $E$ is $W_2$-continuous by \cref{rmk:continuity_of_functionals}), together with the bound on $\|U'_{\tau_i, n}\|$ in \cref{eq:required_bounds_for_compactness}, we have that $\tilde{u}_{\tau_i}(s) - u_{\tau_i}(s) \to 0$ in $\LtwoR$, as $\tau_i \to 0$ and for all $s \in [0, +\infty)$; together with the $W_2$ convergence noted before, we have $\law\big(\tilde{u}_{\tau_i}(s)\big) \to \law(\bgamma_s)$ in $W_2$. 
  Moreover, $|G_\tau|$ as in \cref{eq:G_tau} is well-defined and by \cite[\cref*{M-prop:discrete_newtonian_motion}]{fornasierNewtonianMotionCompanion} we have $|\partial \tilde{E}|(\tilde{u}_{\tau_i}(t)) \leq |G_{\tau_i}|(t)$ for all $i \in \NN$ and almost all $t \in [0, +\infty)$. Combining this estimate with the aforementioned $W_2$-convergence, and using that $|\partial \tilde{E}| = |\partial E| \circ \law$ and $|\partial^- \tilde{E}| = \mathrm{lsc}(|\partial E|) \circ \law$ by \cref{rmk:continuity_of_functionals}, we get 
  \begin{equation}\label{eq:g_equality_metric_slope}
    \liminf_{i \to \infty}|G_{\tau_i}|(t) \geq |\partial^- \tilde{E}|(\bgamma_t) = |\partial \tilde{E}|(\bgamma_t), \qquad \text{for almost every } t \in [0, +\infty),
  \end{equation}
  where the last equality holds as $|\partial \tilde{E}|$ is continuous.
  Thus, we have, for any $s', s, t', t \in [0, +\infty)$ such that $s' < s < t < t'$:
  \begin{equation}\label{eq:estimate_newtonian_curve}
  \begin{split}
    	 \int_s^t \frac{1}{2} & \|\ddot{\bgamma}_r\|^2_{\LtwoR} \diff r + \int_s^t \frac{1}{2} |\partial \tilde E|^2(\bgamma_r) \diff r - \int_s^t H_0^{\tilde E}(V^{\bgamma}_{r}) \diff r \\
			& \overset{(1)}{\leq} \liminf_{i \to \infty} \bigg(\int_s^t \frac{1}{2} |u_{\tau_i}^{(2)}|^2(r) \diff r\bigg) + \int_s^t \liminf_{i \to \infty} \frac{1}{2} |G_{\tau_i}|^2(r) \diff r - \int_s^t \lim_{i \to \infty} H_0^{\tilde E}(\tilde{\mathcal{V}}^{\tau_i, \varepsilon}_{r}) \diff r \\
			& \overset{(2)}{\leq} \liminf_{i \to \infty} \bigg(\int_s^t \frac{1}{2} |u_{\tau_i}^{(2)}|^2(r) \diff r + \int_s^t \frac{1}{2} |G_{\tau_i}|^2(r) \diff r - \int_s^t H_0^{\tilde E}(\tilde{\mathcal{V}}^{\tau_i, \varepsilon}_{r}) \diff r\bigg) \\
			& \overset{(3)}{\leq} \liminf_{i \to \infty} \big(d(\tilde Eu_{\tau_i})(s') - d(\tilde Eu_{\tau_i})(t')\big) + C (t' - t) + C (s - s') \\
			 & \overset{(4)}{=}  \frac{\diff}{\diff \sigma} \tilde E(V^{\bgamma}_{s'}(\sigma)) \Big|_{\sigma = s'} - \frac{\diff}{\diff \sigma} \tilde E(V^{\bgamma}_{t'}(\sigma)) \Big|_{\sigma = t'} + C (t' - t) + C (s - s').
    \end{split}
  \end{equation}
  Here, we justify each step as follows:
  \begin{enumerate}[label=(\arabic*)]
  \item \cref{eq:limit_is_ac2_wasserstein}, \cref{eq:g_equality_metric_slope}, and the fact that $H_0^{\tilde E}(\tilde{\mathcal{V}}^{\tau_i, \varepsilon}_{r}) \to H_0^{\tilde E}(V^{\bgamma}_{r})$ for every $r \in [s, t]$ by the $W_2$-continuity of $G$;
  \item Fatou's lemma---note that $|G_{\tau_i}|^2(r) \geq 0$ for all $r$ and $i$---to exchange the $\liminf$ with the integral in the middle term, and dominated convergence---using the uniform bound on $H_0^{\tilde E}(\tilde{\mathcal{V}}^{\tau_i, \varepsilon}_{r})$---to exchange the limit with the integral in the last term; we also use superadditivity of the $\liminf$;
  \item First, $d(\tilde Eu_{\tau_i})(s') = d(\tilde Eu_{\tau_i})(\tau_i \lfloor s' / \tau_i \rfloor)$ and likewise for $t'$ by \cref{eq:dE_def}. For $\tau_i$ small enough, $\tau_i \lfloor t' / \tau_i \rfloor > t$. By applying \cref{eq:discrete_newtonian_motion} at the points $\tau_i \lfloor s' / \tau_i \rfloor$ and $\tau_i \lfloor t' / \tau_i \rfloor$, what is left is to estimate from below the integral terms on the intervals $[\tau_i \lfloor s' / \tau_i \rfloor, s]$ and $[t, \tau_i \lfloor t' / \tau_i \rfloor]$. The first two integrands are nonnegative, and the last is bounded by $C$. Thus, the error is (asymptotically) bounded by $C (t' - t) + C (s - s')$;
  \item \cref{eq:discrete_identity}.
\end{enumerate}
  Next, the uniform bound over $H_0^{\tilde{E}}(\tilde{\mathcal{V}}^{\tau_i, \varepsilon}_{r})$ and \cref{eq:discrete_newtonian_motion} implies that the functions $d(\tilde E u_{\tau_i})$ are the sum of a uniformly bounded monotone non-decreasing and a monotone non-increasing function; hence, also the increase of the limit function $s \mapsto \frac{\diff}{\diff \sigma} \tilde E(V^{\bgamma}_s)(\sigma)$ is bounded, and we deduce that it has finite variation. Thus, for almost every $s, t \in [0, +\infty)$, $s \mapsto \frac{\diff}{\diff \sigma} \tilde E(V^{\bgamma}_s)(\sigma)$ is continuous in $s$ and $t$, which allows us to take the limit in \cref{eq:estimate_newtonian_curve}, and, as $(|\partial \tilde{E}|, H_0^{\tilde{E}})$ is an upper gradient-Hessian pair by assumption, we finally deduce that the limit curve $\bgamma$ is a curve of Newtonian motion.
\end{proof}

\begin{remark}\label{remark:regularity_implies_locally_lipschitz}
  The regularity assumptions in \cref{lemma:limit_curve_is_newtonian} are mildly \emph{stronger} than demanding $\nabla \tilde{E}$ to be locally Lipschitz. In particular, in the regularity regime of \cref{lemma:limit_curve_is_newtonian}, standard Cauchy-Lipschitz arguments ensure the local existence and uniqueness of solutions to the ODE $\ddot \bgamma = -\nabla \tilde{E}$.
  Rather, relying on such mildly stronger regularity assumptions, we gain a new perspective on Newtonian dynamics (more concretely on solutions of certain Vlasov kinetic PDEs) through the variational characterization \cref{eq:def_variational_Newtonian}. While \cref{lemma:limit_curve_is_newtonian} aims at encompassing a very broad spectrum of smooth energies, one can leverage more specific structural properties of energies and trade them with less regularity. We refer specifically to \cref{ex:newtonian_mmd} in \cref{sec:generalized_newtonian}, where \ref{reg:R3} in \cref{lemma:limit_curve_is_newtonian} does not hold, yet the variational characterization and the discrete-time construction hold nevertheless.
\end{remark}

Having shown that, for regular enough energies $\tilde{E}$, the limit curve $\bgamma$ in \cref{lem:compactness_newtonian} is a curve of Newtonian motion \emph{provided that $(|\partial \tilde{E}|, H_0^{\tilde{E}})$ is a upper gradient-Hessian pair}, we next want to show that for regular enough energies $\tilde{E}$, this is in fact the case.
Motivated by a concrete example (see \cref{sec:example_newtonian}) whose energy is not such regular, it would be desirable to establish this result without invoking $C^2$ regularity as done in the very general \cite[\cref*{M-ex:hilbert_space_newtonian_motion}]{fornasierNewtonianMotionCompanion}. 
We aim therefore to lesser regularity as follows.

\begin{proposition} \label{prop:upper_gradient_hessian}
  Let $E\colon \Pp_2 \to \R$ be such that $\tilde{E} \coloneqq E \circ \law$ fulfills 
  \begin{enumerate}
    \item $\tilde{E}$ is Fréchet differentiable and $\nabla \tilde{E}$ is globally Lipschitz.
    \item $\tilde{E}$ is twice Gâteaux differentiable, and $(X,W,W')\mapsto\nabla^2\tilde E(X)[W,W']$ is jointly continuous.
  \end{enumerate}
  Then, the pair $(|\partial \tilde{E}|, H_0^{\tilde{E}})$ is an upper gradient-Hessian pair for $\tilde{E}$ (see \cref{rmk:recall_gradient-Hessian}). 
  Furthermore, a $W^{2,2}$-curve $\bgamma\colon (a, b) \to \LtwoR$ is a curve of Newtonian motion for $\tilde{E}$ with respect to $(|\partial \tilde{E}|, H_0^{\tilde{E}})$ (see \cref{rmk:recall_curve _Newtonian})
  if and only if it fulfills, for almost every $r \in (a, b)$, the equality
  \begin{equation}\label{eq:ex_simple_upper_gradient_hessian_0}
    \ddot \bgamma_r = - \nabla \tilde{E}(\bgamma_r).
  \end{equation}
\end{proposition}

  In the following proof, for a $W^{2,2}$-curve $\bgamma$ we use freely the following identifications (see \cite[\cref*{M-lem:euclidean_ac2_equivalence}]{fornasierNewtonianMotionCompanion} for the first two, and~\cite[\cref*{M-ex:hilbert_space_newtonian_motion}]{fornasierNewtonianMotionCompanion} for the others):
  \begin{gather*}
    V^{\bgamma}_r(\sigma) = \bgamma_r + (\sigma - r) \dot\bgamma_r, \qquad |\bgamma''|(r) = \|\ddot\bgamma_r\|_{\LtwoR}, \qquad \frac{\diff}{\diff \sigma} \tilde{E}(V^{\bgamma}_r(\sigma))\Big|_{\sigma = r} = \langle \nabla \tilde{E}(\bgamma_r), \dot\bgamma_r \rangle, \\
    |\partial \tilde{E}|(\bgamma_r) = |\nabla \tilde{E}(\bgamma_r)|, \qquad H_0^{\tilde{E}}(V^{\bgamma}_r) = \nabla^2 \tilde{E}(\bgamma_r)[\dot\bgamma_r, \dot\bgamma_r],
  \end{gather*}
  which hold for almost every $r \in (a, b)$.

\begin{proof}
  We fix a $W^{2,2}$-curve $\bgamma$ and show the following three claims, which together will show the proposition:
  \begin{claimlist}
    \item The map $t \mapsto \langle \nabla \tilde{E}(\bgamma_t), \dot\bgamma_t \rangle$ is locally absolutely continuous and thus has locally finite total variation. \label{claim:ex_simple_upper_gradient_hessian_1}
    \item The map $t \mapsto \nabla^2 \tilde{E}(\bgamma_t)[\dot\bgamma_t, \dot\bgamma_t]$ is in $L^1_{\mathrm{loc}}$. \label{claim:ex_simple_upper_gradient_hessian_2}
    \item For every $s, t \in (a, b)$ such that $s < t$, we have that
    \begin{equation} \label{eq:ex_simple_upper_gradient_hessian_1}
      \langle \nabla \tilde{E}(\bgamma_s), \dot\bgamma_s \rangle - \langle \nabla \tilde{E}(\bgamma_t), \dot\bgamma_t\rangle \leq \int_s^t \frac{1}{2} \|\ddot \bgamma_r\|^2_{\LtwoR} + \frac{1}{2} \|\nabla \tilde{E}(\bgamma_r)\|^2_{\LtwoR} - \nabla^2 \tilde{E}(\bgamma_r)[\dot\bgamma_r, \dot\bgamma_r] \, \diff r,
    \end{equation}
    with equality if and only if \cref{eq:ex_simple_upper_gradient_hessian_0} is fulfilled for almost any $r \in (s, t)$. \label{claim:ex_simple_upper_gradient_hessian_3}
  \end{claimlist}
  Let us first show how \cref{prop:upper_gradient_hessian} follows from the claims above:
  By~\eqref{eq:ex_simple_upper_gradient_hessian_1}, $(|\partial \tilde{E}|, H_0^{\tilde{E}})$ is an upper gradient-Hessian pair for $\tilde{E}$.
  The fact that equality in \eqref{eq:ex_simple_upper_gradient_hessian_1} implies \cref{eq:ex_simple_upper_gradient_hessian_0} for almost every $r \in (s, t)$ shows that any curve of Newtonian motion (in the sense of \cref{rmk:recall_curve _Newtonian}) must fulfill \cref{eq:ex_simple_upper_gradient_hessian_0}. Moreover, \cref{claim:ex_simple_upper_gradient_hessian_1,claim:ex_simple_upper_gradient_hessian_2} together with the fact that \cref{eq:ex_simple_upper_gradient_hessian_0} yields equality in \eqref{eq:ex_simple_upper_gradient_hessian_1} show that any solution of \cref{eq:ex_simple_upper_gradient_hessian_0} is a curve of Newtonian motion. 
  We note that $\bgamma$ being a $W^{2,2}$-curve implies that $\bgamma \in C^1((a, b); \LtwoR)$ and that $\dot\bgamma$ is an $\mathrm{AC}^2$-curve in the sense of \cite[Definition 1.1.1]{ambrosioGradientFlows2008} (i.e., $\dot\bgamma \in \mathrm{W}^{1,2}((a, b); \LtwoR)$).\\
  \textbf{\cref{claim:ex_simple_upper_gradient_hessian_1}:} We see that
  \begin{equation*}
    |\langle \nabla \tilde{E}(\bgamma_s), \dot\bgamma_s \rangle - \langle \nabla \tilde{E}(\bgamma_t), \dot\bgamma_t \rangle| 
    \leq C_{\mathrm{Lip}} \|\dot\bgamma_s\|_{\LtwoR} \int_s^t \|\dot\bgamma_r\|_{\LtwoR} \, \diff r + \|\nabla \tilde{E}(\bgamma_t)\|_{\LtwoR} \int_s^t \|\ddot\bgamma_r\|_{\LtwoR} \, \diff r,
  \end{equation*}
  where $C_{\mathrm{Lip}}$ is a Lipschitz constant of $\nabla \tilde{E}$. 
  Because $\dot\bgamma$ and $\nabla \tilde{E} \circ \bgamma$ are continuous, $\|\dot\bgamma_s\|_{\LtwoR}$ and $\|\nabla \tilde{E}(\bgamma_t)\|_{\LtwoR}$ are bounded uniformly over $s, t$ in any compact subinterval of $(a, b)$.
  Since furthermore $t \mapsto \|\dot\bgamma_t\|_{\LtwoR}$ and $t \mapsto \|\ddot\bgamma_t\|_{\LtwoR}$ are in $L^1_{\mathrm{loc}}((a, b))$, local absolute continuity follows. \\ 
  \textbf{\cref{claim:ex_simple_upper_gradient_hessian_2}:} This claim follows directly from the fact that $t \mapsto \nabla^2 \tilde{E}(\bgamma_t)[\dot\bgamma_t, \dot\bgamma_t]$ is continuous by continuity of the map $X, W, W' \mapsto \nabla^2 \tilde{E}(X)[W, W']$ and the map $t \mapsto (\bgamma_t, \dot\bgamma_t, \dot\bgamma_t)$.\\ 
  \textbf{\cref{claim:ex_simple_upper_gradient_hessian_3}:} 
  To show this claim, we first want to show that for almost all $t \in (a, b)$, the composition $\nabla \tilde{E} \circ \bgamma$ is differentiable, and we have $\frac{\diff}{\diff r} (\nabla \tilde{E} \circ \bgamma)(r) = \nabla^2\tilde{E}(\bgamma_r) \dot\bgamma_r$. 
  To see that it is differentiable almost everywhere, we use that $\nabla \tilde{E} \circ \bgamma$ is (as the composition of a Lipschitz and an absolutely continous map) absolutely continuous, and that $L^2$ is a Hilbert space, thus reflexive, and thus we can apply \cite[Remark 1.1.3]{ambrosioGradientFlows2008}. 
  To identify $\frac{\diff}{\diff r} (\nabla \tilde{E} \circ \bgamma)(r)$, we fix $r \in (a, b)$ and some element $Y \in \LtwoR$. Then, by the mean value theorem, for each $h \in (0, +\infty)$, there exists some $k_h \in [0, 1]$ such that
  \begin{equation*}
      \begin{split}
          \frac{1}{h} \langle \nabla \tilde{E}(\bgamma_{r+h}) - \nabla \tilde{E}(\bgamma_r), Y \rangle &= 
          \frac{1}{h} \frac{\diff}{\diff \sigma} \langle \nabla 
          \tilde{E}\big(\bgamma_{r} + \sigma (\bgamma_{r+h} - \bgamma_{r})  \big), Y \rangle \Big{|}_{\sigma=k_h} \\
            &=    \left\langle \nabla^2 \tilde{E}\big(\bgamma_{r} + k_h (\bgamma_{r+h} - \bgamma_{r})  \big) \, \frac{\bgamma_{r+h} - \bgamma_{r}}{h}, Y \right\rangle
      \end{split}
  \end{equation*}
  where we emphasize that we only used \emph{Gateaux} differentiability of $\nabla \tilde{E}$ in the above. We also observe that $\frac1h (\bgamma_{r+h} - \bgamma_{r}) = \frac{1}{h}\int_r^{r+h} \dot \bgamma_q \, \diff q$, and if $r \in (a, b)$ is a Lebesgue point of $\dot \bgamma$, we have that $\frac{1}{h}\int_r^{r+h} \dot \bgamma_q \, \diff q \to \dot\bgamma_r$ as $h \to 0$. 
  Furthermore we always have $\bgamma_{r} + k_h (\bgamma_{r+h} - \bgamma_{r}) \to \bgamma_r$ as $h \to 0$. 
  Thus, by joint continuity of $\nabla^2 \tilde{E}$ in its three arguments, we get that 
  $$\left\langle \frac{\diff}{\diff r} (\nabla \tilde{E} \circ \bgamma)(r), Y \right\rangle = \nabla^2 \tilde{E}(\bgamma_r)[\dot\bgamma_r, Y].$$ 
  Since $Y$ was arbitrary, this proves that for almost every $r \in (a, b)$ (namely for those points for which $\nabla \tilde{E} \circ \bgamma$ is differentiable and which are Lebesgue points of $\dot\bgamma$), we have
  $$\frac{\diff}{\diff r} (\nabla \tilde{E} \circ \bgamma)(r) = \nabla^2\tilde{E}(\bgamma_r) \dot\bgamma_r.$$ 
  Then, using \cref{claim:ex_simple_upper_gradient_hessian_1}, we have
  \begin{equation*}
      \begin{split}
          \langle \nabla \tilde{E}(\bgamma_s), \dot\bgamma_s \rangle - \langle \nabla \tilde{E}(\bgamma_t), \dot\bgamma_t\rangle &=  
          - \int_s^t\frac{\diff}{\diff r} \langle \nabla \tilde{E}(\bgamma_r), \dot\bgamma_r \rangle \, \diff r\\
      &= - \int_s^t \langle \nabla^2 \tilde{E}(\bgamma_r)\dot\bgamma_r, \dot\bgamma_r \rangle + \langle \nabla \tilde{E}(\bgamma_r), \ddot \bgamma_r \rangle \, \diff r \\
      &\leq \int_s^t \frac{1}{2} \|\ddot \bgamma_r\|^2_{\LtwoR} + \frac{1}{2} \|\nabla \tilde{E}(\bgamma_r)\|^2_{\LtwoR} - \nabla^2 \tilde{E}(\bgamma_r)[\dot\bgamma_r, \dot\bgamma_r] \, \diff r,
      \end{split}
  \end{equation*}
  with equality if and only if $\ddot \bgamma_r = -\nabla \tilde{E}(\bgamma_r)$ for almost every $r \in (s, t)$.
\end{proof}

We summarize by unifying \cref{lem:compactness_newtonian,lemma:limit_curve_is_newtonian,prop:upper_gradient_hessian}.

\begin{theorem}\label{thm:main_newtonian}
  Let $\tilde{E} \colon \LtwoR \to \R$ fulfill the following regularity conditions:
  \begin{enumerate} 
    \item\label{condition:nonempty_proximal_map} For $\delta$ small enough, $J_{\delta, \LtwoR}[X]$ is nonempty for all $X \in \LtwoR$.
    \item\label{condition:frechet_lipschitz_derivative} $\tilde E$ is Fréchet differentiable and $\nabla\tilde E$ is Lipschitz.
    \item\label{condition:hessian_continuos}$\tilde E$ is twice Gâteaux differentiable, and $(X,W,W')\mapsto\nabla^2\tilde E(X)[W,W']$ is jointly continuous.
   \end{enumerate}
   Let furthermore $U_0, U_0' \in \LtwoR$, let $(\tau_i)_{i\in\NN}$ be a vanishing sequence and let, for each $\tau_i$, let $(U_{\tau_i, n})_{n = -1}^\infty$ be a Newtonian sequence such that
   \begin{itemize}
    \item\label{condition:ui} For each $t$, $(u_{\tau_i}(t))_{i \in \NN}$ and $(\dot{u}_{\tau_i}(t))_{i \in \NN}$ have uniformly integrable second moments.
    \item $U_{\tau_i, 0} = U_0$ and $U_{\tau_i, -1} = U_0 - \tau U'_0$, for every $i \in \NN$.
   \end{itemize}
   Then, there exists a subsequence (not relabeled) $(\tau_i)_{i=1}^{\infty}$, a $W^{2,2}_{\mathrm{loc}}$-curve $\bgamma\colon [0, +\infty) \to \LtwoR$
  \begin{gather}\label{eq:limit_curve_newtonian}
    \Big(u_{\tau_i}(t), \dot{u}_{\tau_i}(t)\Big) \xrightarrow{\ful{\law}^{-1}(\sigma_{\text{narrow}})} \Big(\bgamma_t, \dot{\bgamma}_t\Big), \qquad \text{for all } t \in [0, +\infty),
  \end{gather}
  where $u_{\tau_i}$ and $\dot{u}_{\tau_i}$ are the piecewise constant interpolant for the sequence $(U_{\tau_i, j})_{j=-1}^\infty$, as in \cref{eq:piecewise_constant_interpolant,eq:piecewise_interpolants_velocity}. Furthermore, the limit curve $\bgamma$ is a curve of Newtonian motion with $\bgamma_0 = U_0$ and $\dot{\bgamma}_0 = U'_0$, and it fulfills, for almost every $r \in [0, +\infty)$, the equality
  \begin{equation}\label{eq:ex_simple_upper_gradient_hessian_0_in_main_thm}
    \ddot \bgamma_r = - \nabla \tilde{E}(\bgamma_r).
  \end{equation}
\end{theorem}
\begin{proof}
    First, we apply \cref{prop:upper_gradient_hessian} to deduce that $(|\partial \tilde{E}|, H_0^{\tilde{E}})$ is an upper gradient-Hessian pair. Next, we use \cref{lem:compactness_newtonian} to construct the limit curve $\bgamma$ such that~\eqref{eq:limit_curve_newtonian} holds.
    To this end, we note that by the Lipschitz continuity of $\nabla \tilde{E}$, the energy $\tilde{E}$ fulfills \cite[\cref*{M-ass:coercivity,M-ass:lipschitz}]{fornasierNewtonianMotionCompanion}. Thus, we may invoke the first part of~\cite[\cref*{M-lem:apriori_bounds}]{fornasierNewtonianMotionCompanion} to see that the bounds demanded by \cref{lem:compactness_newtonian} hold.
    Finally, owing to \cref{lemma:limit_curve_is_newtonian}, the limit curve is a curve of Newtonian motion.
    Indeed, let us check the regularity conditions demanded by \cref{lemma:limit_curve_is_newtonian}:
    \begin{enumerate}[label=R\arabic*., ref=R\arabic*]
      \item For every $X, Y \in \LtwoR$, $s \mapsto \tilde{E}(X + s Y)$ is in $W_{\text{loc}}^{2,2}$;\label{condition:W_22}
      \item $\tilde{E}$ is strictly Hadamard differentiable (with respect to $\LtwoR$) \label{condition:hadamard};
      \item $|\partial \tilde{E}|$ and $H_0^{\tilde{E}}$ are $\LtwoR$-continuous. \label{condition:continuity}
  \end{enumerate}
  \Cref{condition:W_22} follows directly from the Lipschitz continuity of $\nabla\tilde{E}$, which gives Lipschitz continuity of the time derivative of the function $s \mapsto \tilde{E}(X + s Y)$ for every $X, Y \in \LtwoR$. 
  \Cref{condition:hadamard} is fulfilled since $\tilde{E}$ is continuously Fréchet differentiable by \cref{lem:energy_lipschitz_cont}, and hence it is in particular strictly Hadamard differentiable. 
  Likewise, \Cref{condition:continuity} is also fulfilled by \cref{lem:energy_lipschitz_cont}, using that $H_0^{\tilde{E}}(X, Y) = \nabla^2 \tilde{E}(X)[Y, Y]$ by \cref{eq:H_0_in_wasserstein_example}. Finally, invoking again \cref{prop:upper_gradient_hessian}, we obtain that $\bgamma$ fulfills \cref{eq:ex_simple_upper_gradient_hessian_0_in_main_thm}.
\end{proof}

\subsubsection*{Constructing Newtonian Sequences in \texorpdfstring{$\ful{\Pp}_2$}{P̲₂}}\label{sec:constructing_newtonian_seq_in_p2}

So far, we showed that one can construct Newtonian sequences in $\LtwoR$, and that those converge---up to subsequence---to curves of Newtonian motion in $\LtwoR$, which correspond to Wasserstein-Newtonian curves in $\ful{\Pp}_2$. 
While this procedure is analogue to the case of curves of maximal slope in $\Pp_2$---where such curves can be constructed by taking the limit of the MMS scheme in $\LtwoR$, see \cref{sec:wasserstein_as_quotient}---, there remains a key difference: In the case of curves of maximal slope, we can apply the MMS scheme in \emph{either $\Pp_2$ or $\LtwoR$} (see \cref{rmk:mms_downstairs_upstairs}), this leading to the same outcome. 
In the case of Wasserstein-Newtonian curves, the construction detailed so far forces us to construct solutions to the predictor-corrector scheme \eqref{eq:newtonian_sequence} in $\LtwoR$, and to take the limit there. 
In this section, we address this unsettling lack of symmetry. 
More precisely, we show that \cref{alg:wasserstein-newtonian-seq} allows us to construct a sequence $(\ful{M}_{\tau, i})_{i \in \NN}$ in $\ful{\Pp}_2$ such that there exists a Newtonian sequence $(U_{\tau, i})_{i \in \NN}$ for which $\ful{\law}\big((U_{\tau, i}, U'_{\tau, i})\big) = \ful{M}_{\tau, i}$ for all $i \in \NN$, where we define $U'_{\tau, i}$ as in \cref{eq:wasserstein_uprime}. 
Defining the piecewise constant interpolation $\ful{\mu}_\tau$ such that $\ful{\mu}_\tau(t) = \ful{M}_{\tau, i}$ for $t \in [\tau i, \tau (i+1))$ and under suitable a priori bounds, we may apply \cref{lem:compactness_newtonian} to see that $\ful{\mu}_\tau(t) \to \ful{\bar{\mu}}_t$ narrowly, for some $\ful{\bar{\mu}}_\bullet\colon [0, +\infty) \to \ful{\Pp}_2$ (up to subsequence) as $\tau \to 0$, and then proceed as in \cref{prop:limiting_newtonian_wasserstein} to see that $\ful{\bar{\mu}}_\bullet$ is a Wasserstein-Newtonian curve.

\begin{theorem}\label{thm:algorithm_in_p_2}
  Let $E\colon \Pp_2 \to \R$ be an energy\footnote{We emphasise that \cref{alg:wasserstein-newtonian-seq} might fail for energies where the proximal point minimizing set is empty.}, and let $\tilde{E} \coloneqq E\circ \law\colon \LtwoR \to \R$ be the correspondent lifted functional. 
  Let us fix $U_0, U'_0 \in \LtwoR$ and set $\ful{M}_0 \coloneqq \ful{\law} \big((U_0, U'_0)\big)$. Then, for each sequence $(\ful{M}_{\tau, i})_{i \in \NN}$ constructed according to \cref{alg:wasserstein-newtonian-seq}, there exists a Newtonian sequence $(U_{\tau, i})_{i = -1}^\infty$ in $\LtwoR$ for $\tilde{E}$ such that, setting $U'_{\tau, i} \coloneqq \frac{U_{\tau, i} - U_{\tau, i-1}}{\tau}$ as in \cref{eq:wasserstein_uprime}, we have that $\ful{\law}\big((U_{\tau, i}, U'_{\tau, i})\big) = \ful{M}_{\tau, i}$ for all $i \in \NN$.
\end{theorem}
\begin{proof}
In this proof, we assume that $\tau$ is fixed and drop all indices $\tau$.
  We show the claim recursively in $i$; for $i=0$, we get equality by setting $U_{\tau, 0} = U_0$ and $U_{\tau, -1} = U_0 -\tau U'_0$. We have to show that for $U_i$, $U'_i$ and $\ful{M}_i$ such that $\ful{\law}\big((U_i, U'_i)\big) = \ful{M}_i$ and for $\ful{M}_{i+1}$ constructed according to a single step in \cref{alg:wasserstein-newtonian-seq}, there exists $U_{i+1}$ such that, setting $U'_{i+1} \coloneqq \frac{U_{i+1} - U_i}{\tau}$, we have that 
  \begin{equation}\label{eq:wasserstein_alg_u_next_gives_right_distribution}
    \ful{\law}\big((U_{i+1}, U'_{i+1})\big) = \ful{M}_{i+1}
  \end{equation} and such that $U_{i+1}$ solves \cref{eq:newtonian_sequence}, which, in our specific example, is reduced to
  \begin{equation} \label{eq:wasserstein_alg_u_next_is_minimizer}
    U_{i+1} \in \argmin_{x \in \LtwoR} \bigg\{\frac{\|x - Z_i\|^2_{\LtwoR}}{2 \tau^2} + E\big(\law(x)\big) \bigg\},
  \end{equation}
  where we set
  \begin{equation*}
    Z_i \coloneqq U_i + \tau U'_i.
  \end{equation*}
  Reusing the notation of \cref{alg:wasserstein-newtonian-seq}, we first note that $\law(Z_i) = \nu_i$. We apply \cref{lem:lift_coupling} with $X = Z_i$, the collection $\{U_i\}$ (i.e., $I$ a singleton) and $\pi$ as in \cref{alg:wasserstein-newtonian-seq}, to obtain $Y = U_{i+1}$ such that
  \begin{enumerate}
    \item $\law\big((Z_i, U_{i+1})\big) = \pi$.\label{enum:law}
    \item $\mathbb{P}_\R(U_{i+1} \in \cdot \mid U_i=x, Z_i=z) = \mathbb{P}_\R(U_{i+1} \in \cdot \mid Z_i=z) = \pi_z$ for $\law\big((U_i, Z_i)\big)$-almost every $(x, z) \in \R^d \times \R^d$. \label{enum:disintegrations_match}
  \end{enumerate}
  From \cref{enum:law}, we immediately get that $\|Z_i - U_{i+1}\|_{\LtwoR} = W_2(M_{i+1}, \nu_i)$ and $\law(U_{i+1}) = M_{i+1}$ and thus that $U_{i+1}$ fulfills \cref{eq:wasserstein_alg_u_next_is_minimizer}. We are left to show \cref{eq:wasserstein_alg_u_next_gives_right_distribution}. To this end, we calculate
\begin{align*}
    \law\big((U_{i+1}, U'_{i+1})\big) &= \law\bigg(\Big(U_{i+1}, \frac{1}{\tau}(U_{i+1} - U_i)\Big)\bigg) \\
    &= \bigg[(x, x', z) \mapsto \Big(x', \frac{1}{\tau} (x' - x)\Big)\bigg]_\# \law\big((U_i, U_{i+1}, Z_i)\big) \\
    &= \bigg[(x, x', z) \mapsto \Big(x', \frac{1}{\tau} (x' - x)\Big)\bigg]_\# \int_{\R^{2d}} \delta_{x} \times \pi_{(x, z)} \times \delta_{z} \, \diff \law\big(U_i, Z_i\big)(x, z) \nonumber\\
    &\overset{(*)}{=} \bigg[(x, x', z) \mapsto \Big(x', \frac{1}{\tau} (x' - x)\Big)\bigg]_\# \int_{\R^{2d}} \delta_{x} \times \pi_{z} \times \delta_{z} \, \diff \law\big(U_i, Z_i\big)(x, z) \\
    &= \bigg[(x, x') \mapsto \Big(x', \frac{1}{\tau} (x' - x)\Big)\bigg]_\# \int_{\R^{2d}} \delta_{x} \times \pi_{x + \tau v} \, \diff \law\big(U_i, U'_i\big)(x, v) \\
    &= \int_{\R^{2d}} \int_{\R^d} \delta_{\left(x', \frac{x' - x}{\tau}\right)} \, \diff \pi_{x + \tau v}(x') \, \diff \ful{M}_i(x, v),
  \end{align*}
  where in $(*)$ we used that \cref{enum:disintegrations_match} is saying precisely that $\pi_{z} = \pi_{(x, z)}$ holds $\law(U_i, Z_i)$-almost everywhere.
\end{proof}

\begin{algorithm}
\caption{Constructing Wasserstein-Newtonian Sequences}\label{alg:wasserstein-newtonian-seq}
\KwIn{$\ful{M}_0 \in \ful{\Pp}_2$}
\For{$i = 0, 1, 2, \ldots$}{
    Set $\nu_i \coloneqq [(x, v) \mapsto x + \tau v]_\# \ful{M}_i$\;
    Choose $M_{i+1} \in \argmin_{\mu \in \Pp_2} \left\{ \dfrac{1}{2 \tau^2} W_2^2(\mu, \nu_i) + E(\mu) \right\}$\;
    Choose an optimal plan $\pi \in \Pi(\nu_i, M_{i+1})$ with disintegration $\{\pi_x\}_{x \in \R^d}$\;
    Set\addtocounter{AlgoLine}{-1}\refstepcounter{AlgoLine}\label{alg:line:define_next_m}
    \[
        \ful{M}_{i+1} \coloneqq \int_{\R^{2d}} \int_{\R^d} \delta_{\left(x',\, \frac{x' - x}{\tau}\right)} \, \diff \pi_{x + \tau v}(x') \, \diff \ful{M}_i(x, v)
    \]
}
\kern2pt\hrule\kern2pt\vspace{0.5em}
\nonl{\footnotesize The measure $\ful{M}_{i+1}$ defined in \cref{alg:line:define_next_m} in density form can also be written in evaluation form as $A \subseteq \R^{2d} \mapsto \int_{\R^{2d}} \int_{\R^d} \chi_A\!\left(x', \frac{x' - x}{\tau}\right) \diff \pi_{x + \tau v}(x') \, \diff \ful{M}_i(x, v)$.}
\end{algorithm}

\subsection{Identifying PDEs of Newtonian Dynamics in \texorpdfstring{$L^2$}{L²}}\label{sec:identifying_newtonian}

  In \cref{sec:physical_connection}, we cite from \cite{brigatiKineticOptimalTransport2025} (see \cref{theorem:cont_equation_physical}) that every physical curve $\ful{\mu}_\bullet\colon (a, b) \to \ful{\Pp}_2$ fulfills (weakly) the Vlasov equation
  \begin{equation}\label{eq:vlasov_in_intro_ident_pdf}
    \partial_t \ful{\mu}_t + v \cdot \nabla_x \ful{\mu}_t + \nabla_v \cdot (F_t \ful{\mu}_t) = 0, \qquad \text{in } \R^{2d} \times (a, b),
  \end{equation}
  for some force field $F\colon (a, b) \times \R^{2d} \to \R^d$ with $\int_a^b \norm{F_t}^2_{L^2(\ful{\mu}_t)} \, \diff t < +\infty$.
  We interpret \cref{eq:vlasov_in_intro_ident_pdf} as the second-order analogue of the continuity equation \eqref{eq:CE}. 
  By \cref{theorem:cont_equation_physical,cor:lifting_physical_curves}, taking any $W^{2,2}$-curve $\bgamma\colon(a, b) \to \LtwoR$, the curve $\ful{\mu}_\bullet \coloneqq \ful{\law} \circ (\bgamma, \dot{\bgamma})$ fulfills \cref{eq:vlasov_in_intro_ident_pdf} for some $F\colon (a, b) \times \R^{2d} \to \R^d$ with $\norm{F_t}_{L^2(\ful{\mu}_t)} = |\ful{\mu}''|(t) \leq \norm{\ddot{\bgamma}_t}_{\LtwoR}$ for almost every $t \in (a, b)$. 
  At this stage, we would like to relate the vector field $F$ to the curve $\bgamma$. In particular, one may ask what form $F$ has in the case of $\bgamma$ being a curve of Newtonian motion. 
  In \cref{prop:abstract_pde_newtonian}, we answer this question, and through correspondence~\eqref{eq:correspondence_second_order_in_intro_char}, we conclude the form of $F$ in \cref{eq:vlasov_in_intro_ident_pdf} for $\ful{\mu}_{\bullet}$ a Wasserstein-Newtonian curve as in \cref{def:wasserstein_newtonian_curve} (see \cref{cor:vlasov_wasserstein_newtonian}).
 
To identify the PDEs related to Wasserstein-Newtonian curves, we will work in a setup similar to \cref{fig:all_spaces}, but using the phase space instead of the state space (see \cref{fig:all_spaces_phase_space}).
\begin{figure}[H]
\begin{center}
\begin{tikzcd}
                                     &                                                                      & \mathbb{R}                                        &                                                                     &                                           &                                     \\
                                     &                                                                      & \Pp_2 \arrow[u, "E"]                           &                                                                     &                                           &                                     \\
{[0, \infty)} \arrow[r, "{(\bgamma, \dot{\bgamma})}"] & \ful{L}^2 \arrow[r, "\ful{\law}", two heads] \arrow[d, no head, dashed] \arrow[rrr, "\ful{G}", bend left] & \ful{\Pp}_2 \arrow[r, "\ful{b}"] \arrow[u, "(p_x)_\#"{fill=white, inner sep=2pt}, crossing over] & \ful{\mathcal{M}} \arrow[r, "\ful{r}"] \arrow[d, no head, dashed] & (\ful{C}_0^2)' \arrow[d, no head, dashed] &                                     \\
                                     & \ful{L}^2                                                            &                                                   & \ful{C}_0                                                         & \ful{C}_0^2 \arrow[l, hook']              & {[0, \infty)} \arrow[l, "\varphi"']
\end{tikzcd}
\end{center}
\caption{The relations of the relevant spaces, similar to \cref{fig:all_spaces}. Here, the vertical dashed lines denote duality pairings. Moreover, $\ful{C}_0$ denotes the space of continuous functions on the phase space $\R^{2d}$ vanishing at infinity, and $\ful{C}_0^2 \coloneqq \{f \in C^2(\R^{2d}) \mid f, \nabla f, \nabla^2 f \in \ful{C}_0\}$; $\ful{\mathcal{M}}$ denotes the space of finite signed Radon measures on $\R^{2d}$, $\ful{b}$ is the canonical inclusion and $\ful{r}$ the restriction of test functions. While we will not make use of the intermediate spaces $\ful{\mathcal{M}}$ and $\ful{C}_0$, we still show them for the complete conceptual picture. The map $p_x$ in the pushforward $(p_x)_\#$ is defined by the projection $(x, v) \mapsto x$.}\label{fig:all_spaces_phase_space}
\end{figure}

Using the above setup and \cref{prop:frechet-derivative-L2-W2}, we derive the analogue result of \cref{prop:weak-formulation-abstract}, where we identify PDEs that correspond to Wasserstein-Newtonian curves.

We recall that under the assumption that $\tilde{E}$ is Fréchet differentiable at some $X \in L^2(\Omega; \R^d)$, $\nabla \tilde{E}(X)$ can be written as $\nabla \tilde{E}(X) = \eta \circ X$ for some $\eta \in L^2_{\law(X)}(\R^d; \R^d)$. Moreover, for $X' \in L^2(\Omega; \R^d)$ such that $\law(X) = \law(X') \eqqcolon \nu$, also $\nabla \tilde{E}(X') = \eta \circ X'$. This $\eta$ is called the L-derivative of $E$ at $\nu$, and its existence is given by \cite[Proposition 5.24]{carmonaProbabilisticTheoryMean2018} for an atomless standard Borel probability space $\Omega$; yet the result extends to $\Omega_{\R} = [0, 1]^\R$ by factoring through countably many coordinates as in the proof of \cref{prop:l2_is_aa_prep}).

\begin{proposition} \label{prop:abstract_pde_newtonian}
  Let $E\colon \Pp_2 \to \R$ be an energy such that the regularity conditions of \cref{thm:main_newtonian} are fulfilled. Furthermore, let $\bgamma\colon [0, +\infty) \to \LtwoR$ be a $W^{2,2}$-curve of Newtonian motion with respect to the upper gradient-Hessian pair $(|\partial \tilde{E}|, H_0^{\tilde{E}})$, and let $\varphi \in C_c^\infty(\R^{2d + 1})$. For any $t \in [0, \infty)$, we define the measure $\ful{\mu}_t \coloneqq (\bgamma_t, \dot{\bgamma}_t)_\#\mathbb{P} = \ful{\law}\big((\bgamma_t, \dot{\bgamma}_t)\big)$. 
  Then, for a.e.~$t \in [0, \infty)$, we have
  \begin{equation}\label{eq:moment-vlasov-momentum-form}
    \begin{split}
    \frac{\diff}{\diff t} \int \varphi_t(x, v) \, \diff \ful{\mu}_t(x, v) &= \int \partial_t \varphi_t(x, v) \, \diff \ful{\mu}_t(x, v) \\
    &+ \int v \cdot \nabla_x \varphi_t(x, v) \, \diff \ful{\mu}_t(x, v) -\int \eta_{(p_x)_\#\ful{\mu}_t}(x) \cdot \nabla_v \varphi_t(x, v) \, \diff \ful{\mu}_t(x, v),
    \end{split}
  \end{equation}
  where $p_x\colon \R^{2d} \to \R^d$ denotes the projection onto the state space and $\eta_\nu$ denotes the L-derivative of $E$ at a probability measure $\nu \in \Pp_2(\R^d)$. 
\end{proposition}
\begin{proof}
  As in \cref{prop:weak-formulation-abstract}, we use the chain rule; in our case, we apply the chain rule inside $(\ful{C}_0^2)'$. To this end, we define the map $\ful{G}\colon \ful{L}^2 \to (\ful{C}_0^2)'$ by $\ful{G}((X, Y))(\varphi) = \int \varphi(X(\omega), Y(\omega)) \, \diff \mathbb{P}(\omega)$. Analogously to the map $G$ treated in \cref{prop:frechet-derivative-L2-W2}, $\ful{G}$ is Fréchet differentiable, and its derivative has the preadjoint
  \begin{equation*}
    D\ful{G}((X, Y))_*\colon \ful{C}_0^2 \to \ful{L}^2, \qquad D\ful{G}((X, Y))_*[\varphi] = \Big[\omega \mapsto (\nabla_x \varphi(X(\omega), Y(\omega)), \nabla_v \varphi(X(\omega), Y(\omega)))\Big].
  \end{equation*}
  We write (only in this proof) $\ful{\bgamma}_t \coloneqq (\bgamma_t, \dot{\bgamma}_t)$, so that $\dot{\ful{\bgamma}}_t = (\dot{\bgamma}_t, \ddot{\bgamma}_t)$ and $\ful{\mu}_t = \ful{\law}(\ful{\bgamma}_t)$, and calculate, for a.e.~$t$:
  \begin{align*}
    \frac{\diff}{\diff t} &\int \varphi_t(x, v) \, \diff \ful{\mu}_t(x, v) = \frac{\diff}{\diff t} \big\langle \ful{G}(\ful{\bgamma}_t), \varphi_t \big\rangle_{(\ful{C}_0^2)' \times \ful{C}_0^2} \\
    &= \big\langle D\ful{G}(\ful{\bgamma}_t)[\dot{\ful{\bgamma}}_t], \varphi_t \big\rangle_{(\ful{C}_0^2)' \times \ful{C}_0^2} + \big\langle \ful{G}(\ful{\bgamma}_t), \partial_t \varphi_t \big\rangle_{(\ful{C}_0^2)' \times \ful{C}_0^2} \\
    &= \big\langle \dot{\ful{\bgamma}}_t, D\ful{G}(\ful{\bgamma}_t)_*[\varphi_t] \big\rangle_{\ful{L}^2} + \int \partial_t \varphi_t(x, v) \, \diff \ful{\mu}_t(x, v) \\
    &= \big\langle \dot{\bgamma}_t, \nabla_x \varphi_t(\ful{\bgamma}_t) \big\rangle_{L^2} + \big\langle \ddot{\bgamma}_t, \nabla_v \varphi_t(\ful{\bgamma}_t) \big\rangle_{L^2} + \int \partial_t \varphi_t(x, v) \, \diff \ful{\mu}_t(x, v) \\
    &\overset{(*)}{=} \int v \cdot \nabla_x \varphi_t(x, v) \, \diff \ful{\mu}_t(x, v) - \big\langle \nabla \tilde{E}(\bgamma_t), \nabla_v \varphi_t(\ful{\bgamma}_t) \big\rangle_{L^2} + \int \partial_t \varphi_t(x, v) \, \diff \ful{\mu}_t(x, v),
  \end{align*}
  where we used in $(*)$ that $\bgamma$ fulfills \cref{eq:ex_simple_upper_gradient_hessian_0}.
  Leveraging that $\eta$ is the L-derivative of $\tilde{E}$, i.e., $\nabla \tilde{E}(\bgamma_t)(\omega) = \eta_{\law(\bgamma_t)}(\bgamma_t(\omega))$, we rewrite the center term as
  \begin{equation*}
    \begin{split}
       \langle \nabla \tilde E(\bgamma_t), \nabla_v \varphi_t(\ful{\bgamma}_t) \rangle_{L^2} &=  \int_{\Omega} \nabla \tilde{E}(\bgamma_t)(\omega) \cdot \nabla_v \varphi_t(\ful{\bgamma}_t(\omega)) \, \diff \mathbb{P}(\omega) \\
      &= \int_{\Omega} \eta_{\law(\bgamma_t)}(\bgamma_t(\omega)) \cdot \nabla_v \varphi_t(\ful{\bgamma}_t(\omega)) \, \diff \mathbb{P}(\omega) \\
      &= \int \eta_{(p_x)_\#\ful{\mu}_t}(x) \cdot \nabla_v \varphi_t(x, v) \, \diff \ful{\mu}_t(x, v),
    \end{split}
  \end{equation*}
  and this concludes the proof.
\end{proof}

\begin{corollary}\label{cor:vlasov_wasserstein_newtonian}
  Let $E\colon \Pp_2 \to \R$ be a functional such that
  \begin{enumerate}
    \item $\tilde E \coloneqq E \circ \law$ is Fréchet differentiable and $\nabla \tilde{E}$ is globally Lipschitz.
    \item $\tilde{E}$ is twice Gâteaux differentiable, and $(X,W,W')\mapsto\nabla^2\tilde E(X)[W,W']$ is jointly continuous.
  \end{enumerate}
  Furthermore, let $\ful{\mu}_\bullet\colon (a, b) \to \ful{\Pp}_2$ be a Wasserstein-Newtonian curve for $E$ with respect to the upper gradient-Hessian pair $(|\partial E|, H_0^E)$ (which is an upper gradient-Hessian pair by \cref{prop:upper_gradient_hessian,rmk:ugs_passes_upstairs,lem:summarize_sec_1,rmk:tangent_definition_passes_downstairs}). Then, for any $\varphi \in C_c^\infty(\R^{2d + 1})$, and for a.e. $t \in (a, b)$, we have
  \begin{equation}\label{eq:momentum-vlasov-equation}
    \begin{split}
    \frac{\diff}{\diff t} \int \varphi_t(x, v) \, \diff \ful{\mu}_t(x, v) &= \int \partial_t \varphi_t(x, v) \, \diff \ful{\mu}_t(x, v) \\
    &+ \int v \cdot \nabla_x \varphi_t(x, v) \, \diff \ful{\mu}_t(x, v) - \int \eta_{(p_x)_\#\ful{\mu}_t}(x) \cdot \nabla_v \varphi_t(x, v) \, \diff \ful{\mu}_t(x, v),
    \end{split}
  \end{equation}
  where $p_x\colon \R^{2d} \to \R^d$ denotes the projection onto the state space and $\eta_\nu$ denotes the L-derivative of $E$ at a distribution $\nu \in \Pp_2(\R^d)$. 
\end{corollary}
\begin{proof}
  We pick a horizontal lift $\bgamma\colon (a, b) \to L^2$ of $\ful{\mu}_\bullet$ as in \cref{rmk:wasserstein_newtonian_curve_lifts}.
  Then, by \cref{lem:summarize_sec_1,rmk:tangent_definition_passes_downstairs,rmk:wasserstein_newtonian_curve_lifts,rmk:ugs_passes_upstairs}, we can apply \cref{prop:abstract_pde_newtonian} to $\bgamma$ to see that $\ful{\mu}_\bullet = \ful{\law} \circ (\bgamma, \dot\bgamma)$ fulfills the claim.
\end{proof}

\subsection{Wasserstein-Newtonian Motion: A Regular Interaction Kernel}\label{sec:example_newtonian}

In this section, we will apply our theory about Wasserstein-Newtonian curves to a simple example driving energy. In this setting, we allow the driving energy to be a combination of an interaction energy and a potential energy
\begin{equation} \label{eq:newtonian_potential_definition}
  E\colon \Pp_2 \to \R, \qquad \nu \mapsto \iint K(x, x') \, \diff \nu(x) \, \diff \nu(x') + \int V(x) \, \diff \nu(x),
\end{equation}
where we put the following assumptions on $K$ and $V$.

\begin{assumption}\label{ass:kernel_and_potential}
  $K \in C^2(\R^d \times \R^d, \R)$ is a symmetric interaction kernel and $V \in C^2(\R^d)$ a global potential such that $|\nabla^2_{x_1} K|$, $|\nabla_{x_2} \nabla_{x_1} K|$ and $|\nabla^2 V|$ are globally bounded.
\end{assumption}

The goal in this section is to show that for each $\ful{\mu}_0 \in \ful{\Pp}_2$, there exists a Wasserstein-Newtonian curve $\ful{\mu}_\bullet$ with respect to $E$, and that $\ful{\mu}_\bullet$ weakly fulfills the PDE
  \begin{equation} \label{eq:example_newtonian_pde}
    \partial_t \ful{\mu}_t + v \cdot \nabla_x \ful{\mu}_t - \nabla_v\bigg(\Big(2 \int \nabla_{x_1} K(\cdot, x') \, \diff \ful{\mu}_t(x', v')\Big) \ful{\mu}_t\bigg) - \nabla_v\Big((\nabla V) \ful{\mu}_t\Big) = 0.
  \end{equation}
    More precisely, we will show the following theorem at the end of this section: 
    
\begin{proposition} \label{prop:limiting_newtonian_wasserstein}
  Let $E\colon \Pp\to \Rext$ be given as in \cref{eq:newtonian_potential_definition} and fulfill \cref{ass:kernel_and_potential}, let $\tilde{E} \coloneqq E\circ \law\colon \LtwoR \to \R$ be the lift of $E$, and let us fix $U_0, U'_0 \in \LtwoR$. 
  Furthermore, for $\tau > 0$, let $(U_{\tau, j})_{j=-1}^\infty$ be a Newtonian sequence constructed according to \cref{eq:newtonian_sequence} for the energy $\tilde{E}$ with $U_{\tau, 0} = U_0$ and $U_{\tau, -1} = U_0 - \tau U'_0$. 
  Then, for every sequence of step sizes $(\tau_i)_{i=1}^\infty$ such that $\tau_i \to 0$ as $i \to \infty$, there exists a subsequence (not relabeled) $(\tau_i)_{i=1}^{\infty}$, a $W^{2,2}_{\mathrm{loc}}$-curve $\bgamma\colon [0, +\infty) \to \LtwoR$
  \begin{gather} \label{eq:limit_newtonian_example}
    \Big(u_{\tau_i}(t), \dot{u}_{\tau_i}(t)\Big) \xrightarrow{\ful{\law}^{-1}(\sigma_{\text{narrow}})} \Big(\bgamma_t, \dot{\bgamma}_t\Big), \qquad \text{for all } t \in [0, +\infty),
  \end{gather}
  where $u_{\tau_i}$ and $\dot{u}_{\tau_i}$ are the piecewise constant interpolant for the sequence $(U_{\tau_i, j})_{j=-1}^\infty$, as in \cref{eq:piecewise_constant_interpolant,eq:piecewise_interpolants_velocity}.
  Furthermore, the limit curve $\bgamma$ is a curve of Newtonian motion with $\bgamma_0 = U_0$ and $\dot{\bgamma}_0 = U'_0$, and the phase space law $\ful{\mu}_t = \ful{\law}(\bgamma_t, \dot{\bgamma}_t)$ fulfills the PDE~\eqref{eq:example_newtonian_pde}.
\end{proposition}

  \begin{remark}
  In \cite{ambrosio2008hamiltonian}, the authors use an explicit discrete-time scheme along Hamiltonian systems to construct solutions to Vlasov equations; in the Hamiltonian system corresponding to Newtonian dynamics along interaction energies, they thus also construct solutions to the Vlasov PDEs of the form~\eqref{eq:example_newtonian_pde}. Comparing with \cite{ambrosio2008hamiltonian}, admittedly, our theory applies to very smooth potentials and interaction kernels; in this sense, we trade more regularity with the variational characterization. As always done, it is possible through ad-hoc arguments---for instance, smoothing, approximations, and compactness---to extend the application of the theory to less regular models, examples are shown in \cref{sec:generalized_newtonian}.
  \end{remark}

We start by showing that $E$ is sufficiently regular; the proof of the following is more tedious than insightful, and we thus postpone it to \cref{sec:proof-of-lipschitz-cont}.
\begin{lemma} \label{lem:energy_lipschitz_cont}
  Let $E$ be given as in \cref{eq:newtonian_potential_definition} and fulfill \cref{ass:kernel_and_potential}. Then, $\tilde{E} \coloneqq E\circ \law\colon \LtwoR \to \R$ is Fréchet differentiable, and $\nabla\tilde{E}$ is given by
  \begin{equation} \label{eq:frechet_derivative_energy_wasserstein}
    \nabla\tilde{E}(X) = \bigg[\omega \mapsto 2 \int \nabla_{x_1} K\Big(X(\omega), X(\omega')\Big)\,\diff\mathbb{P}(\omega') + \nabla V\Big(X(\omega)\Big)\bigg] \in \LtwoR,
  \end{equation}
  and is globally Lipschitz continuous. Furthermore, $\tilde{E}$ is twice Gateaux differentiable, the second derivative is given by
  \begin{equation} \label{eq:second_derivative}
    \begin{split}
    \nabla^2\tilde E(X)W = \Big[\omega\mapsto &2\int \nabla^2_{x_1}K\big(X(\omega),X(\omega')\big)W(\omega) + \nabla_{x_2}\nabla_{x_1}K\big(X(\omega),X(\omega')\big)W(\omega')\,\diff\mathbb P(\omega') \\
    & \qquad + \nabla^2V\big(X(\omega)\big)W(\omega)\Big],
    \end{split}
  \end{equation}
  and finally, the map $(X, W, W') \mapsto \nabla^2 \tilde{E}(X)[W, W']$ is continuous on $(\LtwoR)^3$.
\end{lemma}
\begin{corollary}
  Let $E$ be given as in \cref{eq:newtonian_potential_definition} and fulfill \cref{ass:kernel_and_potential}, and let $\tilde{E} \coloneqq E\circ \law\colon \LtwoR \to \R$. Furthermore, let $X, W \in \LtwoR$ and let the geodesic $\xi\colon \R \to \LtwoR$ be defined by $t \mapsto X + t W$. Then we have
  \begin{align}
    \frac{\diff}{\diff t} \tilde{E}(\xi_t)_{t=0} &= \int \bigg{[}2 \int \nabla_{x_1} K\Big(X(\omega), X(\omega')\Big) \,\diff \mathbb{P}(\omega') \cdot W(\omega) + \nabla V\Big(X(\omega)\Big) \cdot W(\omega) \bigg{]}\, \diff \mathbb{P}(\omega) \label{eq:dE_in_wasserstein_example}\\
    H_0^{\tilde{E}}(\xi) &= \int \bigg{[}2 W(\omega) \cdot \int \nabla^2_{x_1} K\Big(X(\omega), X(\omega')\Big) W(\omega) + \nabla_{x_2} \nabla_{x_1} K\Big(X(\omega), X(\omega')\Big) W(\omega') \, \diff \mathbb{P}(\omega') \nonumber\\
    &\qquad\qquad\qquad\qquad+ W(\omega) \cdot \nabla^2 V\Big(X(\omega)\Big) \cdot W(\omega)\bigg{]} \, \diff \mathbb{P}(\omega). \label{eq:H_0_in_wasserstein_example}
  \end{align}
\end{corollary}
\begin{proof}
  \Cref{eq:dE_in_wasserstein_example} descends directly from \cref{eq:frechet_derivative_energy_wasserstein}, while \cref{eq:H_0_in_wasserstein_example} follows directly from \cref{eq:second_derivative}.
\end{proof}

In the next lemma, we establish the uniformly integrability of second moments of $(U_{\tau, i})_{i=-1}^{\lfloor\frac{T}{\tau} \rfloor}$ and $ (U_{\tau, i})_{i=-1}^{\lfloor\frac{T}{\tau} \rfloor}$, when  $(U_{\tau, i})_{i=-1}^\infty$ is a Newtonian sequence and $\tau>0$ is sufficiently small.

\begin{lemma}\label{lem:uniformly_integrable_wasserstein}
  Let $E\colon \Pp\to \Rext$ be given as in \cref{eq:newtonian_potential_definition} and fulfill \cref{ass:kernel_and_potential}, let $\tilde{E} \coloneqq E\circ \law\colon \LtwoR \to \R$ be the lift of $E$, and let us fix $U_0, U'_0 \in \LtwoR$. 
  Furthermore, for $\tau > 0$, let $(U_{\tau, i})_{i=-1}^\infty$ be a Newtonian sequence constructed according to \cref{eq:newtonian_sequence} for the energy $\tilde{E}$ with $U_{\tau, 0} = U_0$ and $U_{\tau, -1} = U_0 - \tau U'_0$. 
  Then, for $\tau > 0$ small enough, and $T > 0$ fixed, $U_{\tau, i}$ and $U'_{\tau, i}$ are uniformly bounded in $\LtwoR$ and have uniformly integrable second moments over $\tau$ and $i =-1,0,\ldots ,\lfloor\frac{T}{\tau} \rfloor$.
\end{lemma}
\begin{proof}
  We start by fixing $\tau > 0$ and $\omega \in \Omega$ and estimating $x_i \coloneqq U_{\tau, i}(\omega)$ and $v_i \coloneqq U'_{\tau, i}(\omega)$ by $U_0(\omega)$ and $U'_0(\omega)$. 
  To this end, by \cref{eq:wasserstein_uprime,eq:wasserstein_u}---noting that $J_{\tau^2, \LtwoR}[Z_{\tau, i}] + \tau^2 \nabla\tilde{E}(J_{\tau^2, \LtwoR}[Z_{\tau, i}]) = Z_{\tau, i}$---, we have that
  \begin{align*}
    U'_{\tau, i+1} &= U'_{\tau, i} - \tau \nabla\tilde{E}(U_{\tau, i+1}), \\
    U_{\tau, i+1} &= U_{\tau, i} + \tau U'_{\tau, i+1}.
  \end{align*}
  Evaluating these equations pointwise and using the explicit formula for $\nabla\tilde{E}$ in \cref{eq:frechet_derivative_energy_wasserstein}, we get
  \begin{align}
    v_{i+1} &= v_i - \tau \Big(2 \int \nabla_{x_1} K(x_{i+1}, U_{\tau, i+1}(\omega')) \, \diff \mathbb{P}(\omega') + \nabla V(x_{i+1})\Big), \label{eq:equiintegrable_wasserstein_1}\\
    x_{i+1} &= x_i + \tau v_{i+1}; \label{eq:equiintegrable_wasserstein_2}
  \end{align}
  furthermore, using the fact that according to \cref{ass:kernel_and_potential}, $\nabla_x K$ and $\nabla V$ are $C_{\mathrm{Lip}}$-Lipschitz for some $C_{\mathrm{Lip}} > 0$, and setting $C \coloneqq \max\{|\nabla_x K(0,0)|,|\nabla V(0)|\}$, we can estimate from \cref{eq:equiintegrable_wasserstein_1}
  \begin{equation}\label{eq:equiintegrable_wasserstein_3}
    |v_{i+1}| \leq |v_i| + \tau \Big( 2 \big(C + C_{\mathrm{Lip}} |x_{i+1}| + C_{\mathrm{Lip}} \|U_{\tau, i+1}\|_{L^2}\big) + C + C_{\mathrm{Lip}} |x_{i+1}| \Big). 
  \end{equation}
  Considering the sup-norm on the vector 
  $$ \vec{s}_i \coloneqq \begin{pmatrix}
    |x_i| \\
    |v_i|
  \end{pmatrix},
  $$ 
  we can infer from \cref{eq:equiintegrable_wasserstein_2,eq:equiintegrable_wasserstein_3} that $|\vec{s}_{i}|_{l^\infty} \leq |\vec{s}_{i-1}|_{l^\infty} + \tau (1 + 2C_{\mathrm{Lip}}) |\vec{s}_{i}|_{l^\infty} + \tau(3C + 2C_{\mathrm{lip}} \|U_{\tau, i+1}\|_{L^2})$
  \begin{equation}\label{eq:somewhere_in_the_uniform_integrability_proof}
  |\vec{s}_{i}|_{l^\infty} - |\vec{s}_{i-1}|_{l^\infty} \leq \tau (1 + 2C_{\mathrm{Lip}}) |\vec{s}_{i}|_{l^\infty} + \tau(3C + 2C_{\mathrm{lip}} \|U_{\tau, i+1}\|_{L^2}).
  \end{equation}
  Next, note that $\tilde{E}$ fulfills \cite[\cref*{M-ass:coercivity}]{fornasierNewtonianMotionCompanion} because 
  $$\tilde{E}(X) \geq \tilde{E}(0) - |\nabla \tilde{E}(0)| \|X\|_{L^1} - \frac{5 C_{\mathrm{Lip}}}{2} \|X\|^2_{L^2}$$
  (as $\nabla \tilde{E}$ is $5 C_{\mathrm{Lip}}$-Lipschitz by the proof of \cref{lem:energy_lipschitz_cont}) and $\|X\|_{L^1} \leq \|X\|_{L^2}$.
  Moreover, $\tilde{E}$ fulfills as well \cite[\cref*{M-ass:lipschitz}]{fornasierNewtonianMotionCompanion} by virtue of \cref{lem:energy_lipschitz_cont}. 
  Thus, we may invoke the first part of~\cite[\cref*{M-lem:apriori_bounds}]{fornasierNewtonianMotionCompanion} to see that there exists some $C'$ independent of $\tau$ such that $\sup_{i \leq \frac{T}{\tau}}\{\|U_{\tau, i+1}\|_{L^2}\} \leq C'$ for $\tau$ small enough. Setting $C'' \coloneqq (3C + 2C_{\mathrm{lip}} C')$, we thus get from \cref{eq:somewhere_in_the_uniform_integrability_proof} the inequality
  \begin{equation}
  |\vec{s}_{i}|_{l^\infty} - |\vec{s}_{i-1}|_{l^\infty} \leq \tau (1 + 2C_{\mathrm{Lip}}) |\vec{s}_{i}|_{l^\infty} + \tau C''.
  \end{equation}
  Applying a discrete Gronwall lemma (see \cite[Lemma 3.2.4]{ambrosioGradientFlows2008}), we get, for $\tau i \leq T$,
  \begin{equation}\label{eq:wasserstein_exponential}
    |\vec{s}_{i}|_{l^\infty} \leq \underbrace{\frac{1}{1 - \tau (1 + 2 C_{\mathrm{Lip}})}  e^{\frac{1}{1 - \tau (1 + 2 C_{\mathrm{Lip}})} (1 + 2 C_{\mathrm{Lip}}) T}}_{\substack{\leq C''' \text{ for some } C'''\\\text{and } \tau \text{ small enough}}} |\vec{s}_0|_{l^\infty} C'' T.
  \end{equation}
  We now deduce from~\cref{eq:wasserstein_exponential} that both $U_{\tau, i}$ and $U'_{\tau, i}$ are dominated by $C''' C'' T \max\{|U_0|, |U'_0|\}$, which has integrable second moments, and thus are uniformly integrable.
\end{proof}
We are now ready to prove \cref{prop:limiting_newtonian_wasserstein}.
\begin{proof}[Proof of \cref{prop:limiting_newtonian_wasserstein}]
  We preliminary observe that for $\delta$ small enough, $J_{\delta, \LtwoR}[X]$ is nonempty for all $X \in \LtwoR$. Indeed, this follows directly from the fact that $\nabla \tilde{E}$ is Lipschitz continuous, and so 
  $$Y \mapsto \tilde{E}(Y) + \frac{\|X - Y\|^2_{\LtwoR}}{2\delta}$$
  is continuous and strongly convex for $\delta$ small enough, for every $X \in \LtwoR$. 
  Hence, we can choose a vanishing sequence $(\tau_i)_{i \in \NN}$ and, for each $i \in \NN$, we construct a Newtonian sequence $(U_{\tau_i, n})_{n=-1}^\infty$ such that $U_{\tau_i, 0} = U_0$ and such that $U_{-1} = U_0 - \tau U'_0$.
  \newcommand{\crossout}[1]{%
  \par\noindent
  \begin{tikzpicture}
    \node[inner sep=0pt] (box) {\parbox{\linewidth}{#1}};
    \draw[red] (box.south west) -- (box.north east);
    \draw[red] (box.north west) -- (box.south east);
  \end{tikzpicture}\par}
  As $(u_{\tau_i}(t))_{i \in \NN}$ and $(\dot{u}_{\tau_i}(t))_{i \in \NN}$ (see, respectively, \cref{eq:piecewise_constant_interpolant,eq:piecewise_interpolants_velocity} for their definitions) have uniformly integrable second moments by \cref{lem:uniformly_integrable_wasserstein}, we are in a position to invoke \cref{thm:main_newtonian}, and we obtain the existence of a subsequence (without relabelling) $(\tau_i)_{i\in \NN}$ and a limit curve of Newtonian motion $\bgamma$ such that the convergence~\eqref{eq:limit_newtonian_example} holds. Finally, we can invoke \cref{prop:abstract_pde_newtonian} to deduce that $\ful{\mu}_\bullet \coloneqq \ful{\law} \circ (\bgamma, \dot{\bgamma})$ fulfills, for any $\varphi \in C^\infty_c(\R^{2d+1})$ and for a.e.~$t \in [0, +\infty)$, the equation
  \begin{equation*}
  \begin{split}
      \frac{\diff}{\diff t} \int \varphi_t(x, v) \, \diff \ful{\mu}_t(x, v) &= \int \partial_t \varphi_t(x, v) \, \diff \ful{\mu}_t(x, v) \\
    &+ \int v \cdot \nabla_x \varphi_t(x, v) \, \diff \ful{\mu}_t(x, v) -\langle \nabla\tilde E(\bgamma_t), \nabla_v \varphi_t(\bgamma_t, \dot{\bgamma}_t) \rangle_{\LtwoR}.
  \end{split}
  \end{equation*}
  Substituting $\nabla \tilde{E}$ using \cref{eq:frechet_derivative_energy_wasserstein} in the last term, we are left precisely with the weak formulation of \cref{eq:example_newtonian_pde}.
\end{proof}

\newpage
\printbibliography
\bigskip
\begin{center}
  \FundingLogos
  
  \vspace{0.5em}
  \begin{tcolorbox}\centering\small
   
    Funded by the European Union. Views and opinions expressed are however those of the author(s) only and do not necessarily reflect those of the European Union or the European Research Council Executive Agency. Neither the European Union nor the granting authority can be held responsible for them. This project has received funding from the European Research Council (ERC) under the European Union’s Horizon Europe research and innovation programme (grant agreement No. 101198055, project acronym NEITALG).
    
  \end{tcolorbox}
\end{center}
\begin{table}[ht]
\centering
\caption{Disclosure of AI use in the preparation of this work.}
\label{tab:ai-disclosure}

\begingroup
\renewcommand{\arraystretch}{1.2}
\begin{tabular}{@{}lc@{}}
\toprule
\textbf{Activity} & \textbf{AI used} \\
\midrule
Asking the research question & No \\
Coming up with the approach  & No \\
Proof development            & No$^*$ \\
Writing and exposition       & No \\
Checking for bugs            & Yes \\
Other supporting tasks       & Yes \\
\bottomrule
\end{tabular}
\endgroup

\captionsetup{width=0.8\textwidth}
\caption*{\textbf{Details.}
We employed Claude Opus~5.5 to find bugs in the arguments and notation inconsistencies, and for the tikz code of the figures.\\
\noindent
$^*$: \Cref{prop:chain_rule_generalized} was discovered with the assistance of an AI system (Claude~5.5 Opus). Let us stress that this is the only contribution in this paper obtained through interaction with AI. We deliberately wish to keep it isolated and we place it in the Appendix as a side contribution.}
\end{table}
\appendix
\pdfbookmark[section]{Appendix}{apx:top}
\bookmarksetupnext{level=subsection}
\section{Wasserstein-Newtonian Curves on Nonsmooth Kernels}\label{sec:generalized_newtonian}
The goal of this section is to extend the results from \cref{sec:newtonian_curves} beyond potential- and interaction energies for which the gradient of the Lions lift is globally Lipschitz continuous. 
To this end, we will strongly exploit the specific form of the energy $E\colon \Pp_2\to \Rext$, which is given as  
\begin{equation} \label{eq:potential_interaction_energy}
  E\colon \mu \mapsto \int V(x) \, \diff \mu(x) + \iint \mathcal{K}(x, y) \, \diff \mu(x)\, \diff \mu(y),
\end{equation}
where  $V\colon \R^d \to \R$ is a potential and $\mathcal{K}\colon (\R^d)^2 \to \R$ is a kernel.
The key argument for such an extension to less regular energies is \cref{prop:chain_rule_generalized}, which was discovered with the assistance of an AI system (Claude~5.5 Opus). Let us stress that this is the only contribution in this paper obtained through interaction with AI. We deliberately wish to keep it isolated and we place it in the Appendix as a side contribution.

We recall that for constructing curves of Newtonian motion in the framework of \cite{fornasierNewtonianMotionCompanion}, there are three main challenges, and to tackle them we relied on the strong assumptions formulated in \cref{sec:newtonian_curves}:
\begin{enumerate}
  \item \textbf{Upper gradient-Hessian pair.} In \cref{prop:upper_gradient_hessian}, we relied on global Lipschitz continuity of the gradient of the Lions lift $\tilde{E}$, together with continuity of its (Gâteaux) Hessian, to show that the pair $(|\partial \tilde{E}|, H_0^{\tilde{E}})$ is an upper gradient-Hessian pair for $\tilde{E}$. 
  We proved this for every energy with this regularity, and we did not restrict ourselves to the potential/interaction form \eqref{eq:potential_interaction_energy}.
  \item \textbf{A priori bounds.} In \cref{lem:uniformly_integrable_wasserstein}, we considered energies of the form \eqref{eq:potential_interaction_energy}, and we relied on both $\mathcal{K}$ and $V$ being globally $C^2$, with globally bounded Hessians.
  \item \textbf{Compactness.} 
  When taking the limit of discrete Newtonian sequences, the variational inequality~\eqref{eq:def_variational_Newtonian} that characterizes the curves of Newtonian motion is established by studying the limit of the discrete identity in \cref{eq:discrete_newtonian_motion}.
  In this procedure, among all the objects involved in those variational relations, the biggest difficulty is posed by the Hessian term $H^E_0(\cdot)$. 
  For this reason, in \cref{lemma:limit_curve_is_newtonian} we relied on strong regularity assumptions on the lifted energy $\tilde{E}$ (including continuity of the Hessian term).
\end{enumerate}
Eventually we proved in \cref{prop:limiting_newtonian_wasserstein} that, for energies of the form \eqref{eq:potential_interaction_energy} with suitably regular $V, \mathcal{K}$, we can use both the compactness and the upper gradient-Hessian pair results to construct curves of Newtonian motion. 
However, as mentioned above, we restricted the discussion to a regularity regime on $V$ and $\mathcal{K}$ that implies that the gradient $\nabla \tilde{E}$ is locally Lipschitz (see \cref{remark:regularity_implies_locally_lipschitz}). 
In this section, we go beyond this regime, and we replace the proofs for upper gradient-Hessian pair (i.e., \cref{prop:upper_gradient_hessian}), a-priori bounds (i.e., \cref{lem:uniformly_integrable_wasserstein}), and compactness (i.e., \cref{lemma:limit_curve_is_newtonian}), by proofs which work under weaker assumptions on $V$ and $\mathcal{K}$. 
In what follows, we make strong use of the specific form \eqref{eq:potential_interaction_energy} of the energy, and the proofs will often rely on per-agent computations, i.e., fixing $\omega \in \Omega_\R$ and considering the curve $t \mapsto X_t(\omega)$ for some process $X\colon [a, b] \to L^2(\Omega_\R; \R^d)$.
\begin{remark}[Notation]
  For a random variable $X\colon \Omega_\R \to \R^d$, we denote by $X \otimes X\colon \Omega^2_\R \to (\R^d)^2$ the random variable that maps $(\omega, \omega')$ to $(X(\omega), X(\omega'))$.
\end{remark}
In the rest of this section, we will use the following assumption:
\begin{assumption}\label{ass:growth_condition}
  We have $V \in C^1(\R^d)$ and $\mathcal{K} \in C^1((\R^d)^2)$, and there exists a constant $C$ such that the following two inequalities hold for all $x, y \in \R^d$:
  \begin{align*}
    |\nabla V(x)| &\leq C (1 + |x|), \\
    |\nabla \mathcal{K}(x, y)| &\leq C(1 + |x| + |y|).
  \end{align*}
\end{assumption}
\bookmarksetupnext{level=subsubsection}

\subsection{Upper Gradient-Hessian Pair and Chain Rule}
The main insight of this section is that, leveraging the identification $$W^{2, 2}([a, b]; L^2\big(\Omega_\R; \R^d)\big) \simeq L^2\big(\Omega_\R; W^{2, 2}([a, b]; \R^d)\big),$$ 
one can consider, for $F\colon \R^d \to \R$ regular enough and for $\gamma \in W^{2,2}([a,b]; \R^d)$, the (formal) chain rule
\begin{equation*}
  \langle \nabla F(\gamma_s), \dot\gamma_s \rangle_{\R^d} - \langle \nabla F(\gamma_t), \dot\gamma_t \rangle_{\R^d} = \int_s^t -\langle \nabla F(\gamma_r), \ddot\gamma_r \rangle_{\R^d} - \nabla^2 F(\gamma_r)[\dot\gamma_r, \dot\gamma_r] \, \diff r, \quad \forall s \leq t \in [a, b],
\end{equation*}
which holds for a.e.~particle (or a.e.~particle pair, in the case of interaction energies).
As made precise in the forthcoming result, under the correct integrability assumptions, we can employ Fubini Theorem to integrate the previous identity with respect to $\mathbb{P}_{\R}$ (or $\mathbb{P}_{\R} \otimes \mathbb{P}_{\R}$ in case of interaction energies).

\begin{proposition} \label{prop:chain_rule_generalized}
 Let $E\colon \Pp\to \Rext$ be given as in \cref{eq:potential_interaction_energy} and fulfill \cref{ass:growth_condition}, and let $\tilde{E} \coloneqq E\circ \law\colon \LtwoR \to \R$. 
 Let $H_V\colon \R^d \times \R^d \to \R$ and $H_{\mathcal{K}}\colon (\R^d)^2 \times (\R^d)^2 \to \R$ be such that $V$ and $\mathcal{K}$ fulfill a chain rule along $H_V$ and $H_{\mathcal{K}}$, i.e., such that, for all $W^{2,2}$-curves $\gamma^V\colon [a, b] \to \R^d$, $\gamma^{\mathcal{K}}\colon [a, b] \to (\R^d)^2$, we have, for all $s, t \in [a, b]$ such that $s \leq t$:
  \begin{align}
  \langle \nabla V(\gamma^V_s), \dot\gamma^V_s \rangle_{\R^d} - \langle \nabla V(\gamma^V_t), \dot\gamma^V_t \rangle_{\R^d} &= 
  \int_s^t -\langle \nabla V(\gamma^V_r), \ddot\gamma^V_r \rangle_{\R^d} - H_V(\gamma^V_r, \dot\gamma^V_r) \, \diff r, \label{eq:pointwise_chain_rule_v} \\
  \langle \nabla \mathcal{K}(\gamma^{\mathcal{K}}_s), \dot\gamma^{\mathcal{K}}_s \rangle_{\R^d} - \langle \nabla \mathcal{K}(\gamma^{\mathcal{K}}_t), \dot\gamma^{\mathcal{K}}_t \rangle_{\R^d} &=
  \int_s^t -\langle \nabla \mathcal{K}(\gamma^{\mathcal{K}}_r), \ddot\gamma^{\mathcal{K}}_r \rangle_{\R^d} - H_{\mathcal{K}}(\gamma^{\mathcal{K}}_r, \dot\gamma^{\mathcal{K}}_r) \, \diff r. \label{eq:chain_rule_kernel}
  \end{align}
  Furthermore, let us assume that there exist $C_1, C_2 \in \R$ and $\mathfrak{s}_1, \mathfrak{s}_2 \in \{-1, 1\}$ such that, for all $x_1, y_1 \in \R^d$ and $x_2, y_2 \in \R^{2d}$, we have that
  \begin{enumerate}
    \item $\mathfrak{s}_1 H_V(x_1, y_1) \leq C_1 (1 + |y_1|)^2$;
    \item $\mathfrak{s}_2 H_{\mathcal{K}}(x_2, y_2) \leq C_2 (1 + |y_2|)^2$.
  \end{enumerate}
  Then, letting $\bgamma\colon [a, b] \to \LtwoR$ be a curve $W^{2, 2}$-regular, for almost all $s, t \in [a, b]$ such that $s \leq t$ we have the equality
  \begin{equation} \label{eq:chain_rule_integrated}
    \langle \nabla \tilde{E}(\bgamma_s), \dot\bgamma_s \rangle_{\LtwoR} - \langle \nabla \tilde{E}(\bgamma_t), \dot\bgamma_t \rangle_{\LtwoR} = \int_s^t -\langle \nabla \tilde{E}(\bgamma_r), \ddot\bgamma_r \rangle_{\LtwoR} - \mathbf{H}(\bgamma_r, \dot\bgamma_r) \, \diff r,
  \end{equation}
  where $\mathbf{H}\colon \LtwoR \times \LtwoR \to \R \cup \{\pm \infty\}$ is defined as
  \begin{equation*}
      (X, W) \mapsto \int_{\Omega} H_V\big(X,  W\big) \, \diff \mathbb{P}_{\mathbb{R}} + \iint_{\Omega^2} H_{\mathcal{K}}\big(X \otimes X, W \otimes W\big) \, \diff \mathbb{P}_{\mathbb{R}}\otimes\mathbb{P}_{\mathbb{R}}.
  \end{equation*}
\end{proposition}
\begin{proof}
  First, we note that under \cref{ass:growth_condition}, it is well-known that $\tilde{E} \in C^1$, with the Fréchet differential given as
  \begin{equation}\label{eq:frechet_differential}
    \nabla \tilde{E}(X) = \bigg[\omega \mapsto \nabla V(X(\omega)) + \int_{\Omega} \nabla_{x_1} \mathcal{K}(X(\omega), X(\omega')) + \nabla_{x_2} \mathcal{K}(X(\omega'), X(\omega)) \, \diff \mathbb{P}(\omega')\bigg].
  \end{equation}
  Next, we see that by additivity of all the terms involved in \cref{eq:chain_rule_integrated}, it suffices to treat the two cases where either $V=0$ or $\mathcal{K}=0$. However, the kernel case (i.e., $V=0$) is just the potential case (i.e., $\mathcal{K}=0$) in disguise: If the proposition held in the potential case, we could apply it to the kernel case by replacing $d$ with $2d$, setting $V'\colon \R^{2d} \to \R$ with $V'(x, y) \coloneqq \mathcal{K}(x, y)$, $\Omega'_{\R} \coloneqq \Omega_{\R} \times \Omega_{\R}$ and $\bgamma' = \bgamma \otimes \bgamma$. 
  Thus, it suffices to treat the case where $\mathcal{K}=0$.

  For the case where $\mathcal{K}=0$, we identify $W^{2,2}\big([a, b]; L^2(\Omega_{\R}; \R^d)\big)$ with $L^2\big( \Omega_{\R}; W^{2,2}([a, b]; \R^d)\big)$. 
  We denote the image of $\bgamma\colon [a,b]\to \LtwoR$ under this identification by $X\colon \Omega_\R\to W^{2,2}([a, b]; \R^d)$. 
  Then, the terms involved in \cref{eq:chain_rule_integrated} behave as follows under this identification:
  \begin{align*}
    &\langle \nabla \tilde{E}(\bgamma_s), \dot{\bgamma}_s \rangle_{\LtwoR} = 
    \int_{\Omega} \langle \nabla V(X(\omega)_s), (\dot X(\omega))_s \rangle_{\R^d} \, \diff \mathbb{P}_{\R}(\omega), \\
    &\langle \nabla \tilde{E}(\bgamma_t), \dot{\bgamma}_t \rangle_{\LtwoR}
    = \int_{\Omega} \langle \nabla V(X(\omega)_t), (\dot X(\omega))_t \rangle \, \diff \mathbb{P}_{\R}(\omega), \\
    &\langle \nabla \tilde{E}(\bgamma_r), \ddot{\bgamma}_r \rangle_{\LtwoR} 
    = \int_{\Omega} \langle \nabla V(X(\omega)_r), (\ddot X(\omega))_r \rangle_{\R^d} \, \diff \mathbb{P}_{\R}(\omega), \\
    &\mathbf{H}(\bgamma_r, \dot{\bgamma}_r) = \int_{\Omega} H_V(X(\omega)_r, (\dot X(\omega))_r) \, \diff \mathbb{P}_{\R}(\omega).
  \end{align*}
  Thus, by integrating the pointwise chain rule \cref{eq:pointwise_chain_rule_v}, we only need to make sure that we can apply Fubini-Tonelli over the integration domain $[a, b] \times \Omega_{\R}$ and to the two maps 
  $$
  (r, \omega) \mapsto \langle \nabla V(\bgamma_r(\omega)), \ddot{\bgamma}_r(\omega) \rangle_{\R^d}, \qquad (r, \omega) \mapsto H_V (\bgamma_r(\omega), \dot{\bgamma}_r(\omega)), $$
  in order to to arrive at \cref{eq:chain_rule_integrated} (using again the identification provided above). 
  For the former map, Fubini-Tonelli is justified by the fact that, owing to the (at most) linear growth of $\nabla V$ guaranteed by \cref{ass:growth_condition}, the map $(x, y) \mapsto \langle \nabla V(x), y \rangle$ has (at most) quadratic growth, and---since $\ddot{\bgamma}$ and $\bgamma$ are both in $L^2([a, b]; L^2(\Omega_{\R}; \R^d))$---the integrability follows. 
  Likewise, for the latter map, the quadratic bound on $H_V$ together with $\dot{\bgamma} \in L^2([a, b]; L^2(\Omega_{\R}; \R^d))$ justifies the usage of Tonelli; a posteriori, the integral is finite.
\end{proof}

\begin{remark}
  In the proof above, the trick of reducing to potential energies works analogously for kernels related to three-body interaction energies (or more generally, $n$-body interactions). 
  For brevity, we focus here on binary interaction energies.
\end{remark}

\begin{corollary} \label{cor:upper_gradient_hessian}
  Let $E\colon \Pp\to \Rext$ be given as in \cref{eq:potential_interaction_energy} and fulfill \cref{ass:growth_condition}, and let $\tilde{E} \coloneqq E\circ \law\colon \LtwoR \to \R$.
  Let $H_V$ and $H_{\mathcal{K}}$ be as in \cref{prop:chain_rule_generalized}, so that, in particular, they fulfill the quadratic bounds in \cref{prop:chain_rule_generalized}. Then, the pair $(|\nabla \tilde{E}|, \mathbf{H})$ is an upper gradient-Hessian pair, and any curve of Newtonian motion $\bgamma\colon [a, b] \to \LtwoR$ satisfies, for almost every $r \in [a, b]$, the ODE
  \begin{equation}
    \ddot{\bgamma}_r = -\nabla \tilde{E}(\bgamma_r).
  \end{equation}
\end{corollary}

\begin{proof}
  From \cref{eq:chain_rule_integrated}, we get, for any $W^{2,2}$-curve $\bgamma\colon [a, b] \to \LtwoR$, the inequality
  \begin{align*}
    \langle \nabla \tilde{E}(\bgamma_s), \dot\bgamma_s \rangle - \langle \nabla \tilde{E}(\bgamma_t), \dot\bgamma_t \rangle &= \int_s^t -\langle \nabla \tilde{E}(\bgamma_r), \ddot\bgamma_r \rangle - \mathbf{H}(\bgamma_r, \dot\bgamma_r) \, \diff r, \\
    &\overset{(*)}{\leq} \int_s^t \frac{1}{2} |\ddot \bgamma_r|^2  + \frac{1}{2} |\nabla \tilde{E}(\bgamma_r)|^2 - \mathbf{H}(\bgamma_r, \dot\bgamma_r) \, \diff r,
  \end{align*}
  which shows that $(|\nabla \tilde{E}|, \mathbf{H})$ is an upper gradient-Hessian pair. Furthermore, if $\bgamma$ is a curve of Newtonian motion, it attains equality in $(*)$ above, which implies that $\ddot{\bgamma}_r = -\nabla \tilde{E}(\bgamma_r)$ for almost every $r \in (a, b)$.
\end{proof}

\begin{corollary}
  Let $E\colon \Pp\to \Rext$ be given as in \cref{eq:potential_interaction_energy} and fulfill \cref{ass:growth_condition}, and let $\tilde{E} \coloneqq E\circ \law\colon \LtwoR \to \R$. 
  Let us furthermore assume that $V \in C^2(\R^d; \R)$ and $\mathcal{K} \in C^2((\R^d)^2; \R)$, and that both $V$ and $\mathcal{K}$ are semiconvex or semiconcave.
  Then, setting
  \begin{equation*}
    \mathbf{H}\colon (X, W) \mapsto \int_{\Omega} \nabla^2 V(X)[W, W] \, \diff \mathbb{P}_\R + \iint_{\Omega^2} \nabla^2 \mathcal{K}(X \otimes X)[W \otimes W, W \otimes W] \, \diff \mathbb{P}_\R \otimes \mathbb{P}_\R,
  \end{equation*}
  we have that $(|\nabla \tilde{E}|, \mathbf{H})$ is an upper gradient-Hessian pair for $\tilde{E}$, and any curve of Newtonian motion $\bgamma\colon [a, b] \to \LtwoR$ satisfies, for almost every $r \in [a, b]$, the ODE
  \begin{equation*}
    \ddot{\bgamma}_r = -\nabla \tilde{E}(\bgamma_r).
  \end{equation*}
\end{corollary}

\begin{proof}
  Since both $V$ and $\mathcal{K}$ are continuously twice differentiable, the choices $H_V\colon (x, v) \mapsto \nabla^2 V(x)[v, v]$ and $H_{\mathcal{K}}\colon (x, v) \mapsto \nabla^2 \mathcal{K}(x)[v, v]$ fulfill the chain rules \cref{eq:pointwise_chain_rule_v,eq:chain_rule_kernel}.
  Furthermore, semiconvexity or semiconcavity imply the quadratic bounds in \cref{prop:chain_rule_generalized}. 
  Thus, the assumptions put us in a position where we can apply \cref{cor:upper_gradient_hessian} with the above choices of $H_V$ and $H_{\mathcal{K}}$.
\end{proof}

\bookmarksetupnext{level=subsubsection}
\subsection{A Priori Bounds}
To show the necessary a priori bounds, we need only \cref{ass:growth_condition}.
\begin{proposition} \label{prop:a_priori}
  Let $E\colon \Pp\to \Rext$ be given as in \cref{eq:potential_interaction_energy} and fulfill \cref{ass:growth_condition}, and let $\tilde{E} \coloneqq E\circ \law\colon \LtwoR \to \R$.
  Furthermore, let us consider $U_0, U'_0 \in \LtwoR$ and $T > 0$. 
  Then, there exists $\Theta \in \LtwoR$ such that the following holds: For $\tau>0$ small enough, for every $i \in \NN$ such that $i \leq T/\tau$, and for every Newtonian sequence $(U_{\tau, i})_{i \geq -1}$ for $\tilde{E}$ with $U_{\tau, 0} = U_0$ and $U_{\tau, -1} = U_0 - \tau U'_0$, we have that
  \begin{equation*}
  |U_{\tau, i}(\omega)| \leq \Theta(\omega), \qquad |U'_{\tau, i}(\omega)| \leq \Theta(\omega), \qquad |\nabla \tilde{E}(U_{\tau, i})(\omega)| \leq \Theta(\omega), \qquad \text{for a.e.~} \omega \in \Omega_\R,
  \end{equation*}
  where we set $U'_{\tau, i} \coloneqq \frac{1}{\tau} (U_{\tau, i} - U_{\tau, i-1})$.
\end{proposition}
\begin{proof}
  By \cref{eq:wasserstein_uprime,eq:wasserstein_u} and using the Fr\'echet differentiability of $\tilde{E}$, we have that, for all $i \in \NN$, $U_{\tau, i + 1}' = U_{\tau, i}' - \tau \nabla \tilde{E}(U_{\tau, i+1})$. 
  Thus, using \cref{eq:frechet_differential}, we have that, for almost all $\omega \in \Omega$,
  \begin{equation} \label{eq:domination_proof_1}
      \begin{split}
    &\Big|\nabla \tilde{E}(U_{\tau, i+1}) (\omega)\Big| \leq \,\Big|\nabla V(U_{\tau, i+1})(\omega)\Big| \\
    &  \quad \qquad\qquad \qquad + \Big|\int_\Omega \nabla_{x_1} \mathcal{K}(U_{\tau, i+1}(\omega), U_{\tau, i+1}(\omega')) + \nabla_{x_2} \mathcal{K}(U_{\tau, i+1}(\omega'), U_{\tau, i+1}(\omega)) \, \diff \mathbb{P}_\R (\omega') \Big| \\
    & \qquad\qquad \qquad \leq  3C + 3C \big|U_{\tau, i+1}(\omega)\big| + 2C\|U_{\tau, i+1}\|_{L^1}, 
    \end{split}
  \end{equation}
  and hence, recalling that $U_{\tau, i + 1}' = U_{\tau, i}' - \tau \nabla \tilde{E}(U_{\tau, i+1})$, by the triangle inequality we further have
  \begin{equation*}
    \Big|U_{\tau, i + 1}'(\omega)\Big| 
    \leq \big|U_{\tau, i}'(\omega)\big| + \tau \Big[3C + 3C \big|U_{\tau, i+1}(\omega)\big| + 2C\|U_{\tau, i+1}\|_{L^1}\Big]. 
 \end{equation*}
  Integrating over $\Omega$, we get $\|U_{\tau, i+1}'\|_{L^1} \leq \|U_{\tau, i}'\|_{L^1} + \tau \big[3C + 5C \|U_{\tau, i+1}\|_{L^1}\big]$; furthermore, we have the two estimates $\|U_{\tau, i+1}\|_{L^1} \leq \|U_{\tau, i}\|_{L^1} + \tau \|U_{\tau, i+1}'\|_{L^1}$ and $|U_{\tau, i+1}(\omega)| \leq |U_{\tau, i}(\omega)| + \tau |U'_{\tau, i+1}(\omega)|$. 
  By combining the last four inequalities and considering the vector
  \begin{equation*}
    P_i \coloneqq \big(|U_{\tau, i}'(\omega)|, |U_{\tau, i}(\omega)|, \|U'_{\tau, i}\|_{L^1}, \|U_{\tau, i}\|_{L^1} \big)^T,
  \end{equation*}
  we deduce that $|P_{i+1}|_{\infty} \leq |P_i|_{\infty} + \tau [3C + C' |P_{i+1}|_{\infty}]$, where $C' \coloneqq \max\{5C, 1\}$. 
  Using a discrete Gronwall lemma (\cite[Lemma 3.2.4]{ambrosioGradientFlows2008}), we thus see that there exists $D>0$ only depending on $C$ and $T$, such that for all $\tau$ small enough and all $i \in \NN$ such that $\tau i \leq T$, we have that $|P_i|_{\infty} \leq D (1 + |P_0|_{\infty})$. Since this held for any fixed $\omega \in \Omega_\R$, we have that, defining
  \begin{equation*}
    \omega \mapsto \bar{\Theta}(\omega) \coloneqq D \Big(1 + \max \big\{|U_0'(\omega)|, |U_0(\omega)|, \|U'_0\|_{L^1}, \|U_0\|_{L^1} \big\}\Big) \in \LtwoR,
  \end{equation*}
  $\bar{\Theta}$ dominates $U_{\tau, i}$ and $U'_{\tau, i}$ as required.
  Finally, by setting $\Theta\coloneqq 3C + 3C\bar{\Theta} + 2C\|\bar{\Theta}\|_{L^1} + \bar{\Theta}$, we have by \cref{eq:domination_proof_1} that $\Theta$ also dominates $\nabla \tilde{E}(U_{\tau, i})$ as required.
\end{proof}

\begin{proposition} \label{prop:discrete_newtonain_curve_equation}
  Let $E\colon \Pp\to \Rext$ be given as in \cref{eq:potential_interaction_energy} and fulfill \cref{ass:growth_condition}, and let $\tilde{E} \coloneqq E\circ \law\colon \LtwoR \to \R$.
  Let us in addition assume that there exists $\tau^*$ such that $J_{\tilde{\tau}, \LtwoR}[x]$ is nonempty for $\tilde{\tau} \leq \tau^*$. 
  Furthermore, let $\tau \in (0, \sqrt{\tau^*}]$ and $(U_{\tau, j})_{j \geq -1}$ be a Newtonian sequence for $\tilde{E}$, and let $\mathbf{H} \colon \LtwoR \times \LtwoR \to \R \cup \{\pm \infty\}$ be such that the chain rule \eqref{eq:chain_rule_integrated} is fulfilled along all geodesics.
  Then, setting $t_\tau^j \coloneqq \tau j$, we have, for all $i, j \in \NN$ such that $i < j$,
  \begin{equation}\label{eq:discrete_newtonian_pointwise}
    d(\tilde E u_{\tau})(t_{\tau}^{i}) - d(\tilde E u_{\tau})(t_{\tau}^{j}) = \frac{1}{2}\int_{t_{\tau}^{i}}^{t_{\tau}^{j}} |u_{\tau}^{(2)}|^2(r)\, \diff r + \frac{1}{2} \int_{t_{\tau}^{i}}^{t_{\tau}^{j}} |G_{\tau}|^2(r)\, \diff r - \int_{t_{\tau}^{i}}^{t_{\tau}^{j}} \mathbf{H}(\tilde{\mathcal{V}}^{\tau, \varepsilon}_{r}) \, \diff r,
  \end{equation}
  where $d(\tilde E u_\tau)$ is defined in \cref{eq:dE_def}, $|u_\tau^{(2)}|$ in \cref{eq:second_derivative_metric}, $|G_\tau|$ in \cref{eq:G_tau}, and $\tilde{\mathcal{V}}^{\tau, \varepsilon}$ in \cref{eq:tilde_v_wasserstein}, and where $\mathbf{H}(\tilde{\mathcal{V}}^{\tau, \varepsilon}_r)$ denotes $\mathbf{H}$ evaluated at the base point and velocity of the geodesic $\tilde{\mathcal{V}}^{\tau, \varepsilon}_r$.
\end{proposition}
\begin{proof}
  Following \cite[\cref*{M-prop:discrete_newtonian_motion}]{fornasierNewtonianMotionCompanion} with $\mathbf{H}$ in place of $H_0^{\tilde{E}}$, the single equation that we have to prove is
  \begin{equation}\label{eq:discrete_formula_intermediate}
    \int_0^\tau \frac{1}{\psi'(\psi^{-1}(\delta))} H_0^{\tilde{E}}(\mathcal{V}^{\tau, \varepsilon}_{t_{\tau}^i + \delta}) \, \diff \delta = \int_0^\tau \frac{1}{\psi'(\psi^{-1}(\delta))} \mathbf{H}(\mathcal{V}^{\tau, \varepsilon}_{t_{\tau}^i + \delta}) \, \diff \delta,
  \end{equation}
  where $\psi$ is defined in \cref{eq:def_psi}
  and $\mathcal{V}^{\tau, \varepsilon}_{t_{\tau}^i + \delta} = \Big[t_{\tau}^i + \delta + \sigma \mapsto \xi_{2\delta + \sigma}\Big]$, where $\xi$ is the geodesic with $\xi_0 = U_{\tau, i-1}$ and $\xi_\tau = U_{\tau, i}$, as defined in \cref{eq:v_wasserstein}.
  However, we see from the definitions of $H_0^{\tilde{E}}$ in \cref{eq:hessian_0}
  and from $\mathbf{H}$ fulfilling the chain rule, that $H_0^{\tilde{E}}(\mathcal{V}^{\tau, \varepsilon}_{t_{\tau}^i + \delta}) = \frac{\diff^2}{\diff s^2} \tilde{E}(\xi_s)|_{s=2\delta} = \mathbf{H}(\mathcal{V}^{\tau, \varepsilon}_{t_{\tau}^i + \delta})$ for almost every $\delta \in [0, \tau]$, from which \cref{eq:discrete_formula_intermediate} follows.
\end{proof}

\bookmarksetupnext{level=subsubsection}
\subsection{Compactness}
For the compactness step, we need another piece of information, i.e., the fact that the Hessian $\mathbf{H}$ converges.

\begin{proposition} \label{prop:compactness}
  Let $E\colon \Pp\to \Rext$ be given as in \cref{eq:potential_interaction_energy} and fulfill \cref{ass:growth_condition}, and let $\tilde{E} \coloneqq E\circ \law\colon \LtwoR \to \R$.
  Let us in addition assume that there exists $\tau^*$ such that $J_{\tilde{\tau}, \LtwoR}[x]$ is nonempty for $\tilde{\tau} \leq \tau^*$, and that $\mathbf{H}\colon \LtwoR\times \LtwoR \to \R \cup \{\pm \infty\}$ is such that the chain rule \cref{eq:chain_rule_integrated} is fulfilled along every $W^{2,2}$-curve. 
  Let us further consider $U_0, U'_0 \in \LtwoR$. 
  Then, there exist a vanishing sequence $(\tau_i)_{i\in \NN}$, and, for every $i \in \NN$, a Newtonian sequence $(U_{\tau_i, j})_{j=-1}^{\infty}$ with $U_{\tau_i, 0} = U_0$ and $U_{\tau_i, -1} = U_0 - \tau_i U'_0$, such that, for all $t \in [0, +\infty)$,
  \begin{equation} \label{eq:compactness_convergence}
    \big(u_{\tau_i}(t), \dot{u}_{\tau_i}(t)\big) \to \big(\bar{u}_t, \dot{\bar{u}}_t\big), \qquad \text{as } i \to \infty, \text{ in } \ful{\law}^{-1}(W_2),
  \end{equation}
  where $\bar{u}\colon [0, +\infty) \to \LtwoR$ is a $W^{2,2}$-regular limit curve with $\bar{u}_0 = U_0$ and $\dot{\bar{u}}_0 = U'_0$, and 
  where $u_\tau$ and $\dot{u}_\tau$ are the piecewise constant interpolations introduced, respectively, in \cref{eq:piecewise_constant_interpolant,eq:piecewise_interpolants_velocity}.
  Moreover, whenever $\limsup_{i \to \infty} \int_s^t\mathbf{H}(\tilde{\mathcal{V}}^{\tau_i, \varepsilon}_{r})\,\diff r \leq \int_s^t\mathbf{H}(\bar{u}_r, \dot{\bar{u}}_r)\, \diff r$ for every $s, t \in [0, +\infty)$ such that $s < t$, with $\tilde{\mathcal{V}}^{\tau, \varepsilon}$ and the evaluation of $\mathbf{H}$ as in \cref{eq:discrete_newtonian_pointwise}, then the limit curve $\bar{u}$ is a curve of Newtonian motion with respect to $(|\nabla \tilde{E}|, \mathbf{H})$.
\end{proposition}
\begin{remark}
  The proof below follows the general strategy of the proofs of \cref{lem:compactness_newtonian,lemma:limit_curve_is_newtonian}. 
  One important distinction is that, in \cref{eq:estimate_newtonian_curve}, we were able to estimate the errors that came from shifting from the times $s'$ to $s$ and from $t'$ to $t$ by $C(s-s')$ and, respectively, uniformly over $\tau_i$. 
  This was needed because the discrete jump points $\tau_i j$ (for $j \in \NN$) did not agree between different values for $\tau_i$, so we chose a small intervals of positive length $[s', s]$, and $[t', t]$, which are guaranteed to contain one discrete jump point for all $\tau_i$ small enough. 
  In the setting of \cref{prop:compactness}, we are not able to control the error uniformly over $\tau_i$, so we have to resort to another strategy: We consciously pick $\tau_i = 2^{-M-i}$ for some $M \in \NN$. 
  That way, we can first show that the limit curve fulfills the inequality~\eqref{eq:recall_curves_of_newtonian_motion} (which characterizes curves of Newtonian motion) for all dyadic times $s, t \in \{k 2^{-j} \mid k, j \in \NN\}$, and then generalize by continuity of both sides.
\end{remark}

\begin{proof}
  First, we choose the vanishing sequence $(\tau_i)_{i\in \NN} = (2^{-M-i})_{i \in \NN}$ (for some $M \in \NN$ such that $(2^{-M})^2 \leq \tau^*$), and, for each $i \in \NN$, we construct a Newtonian sequence $(U_{\tau_i, j})_{j=-1}^{\infty}$ with $U_{\tau_i, 0} = U_0$ and $U_{\tau_i, -1} = U_0 - \tau_i U'_0$. 
  For every $T \in \NN$, we can choose $\Theta_T$ as in \cref{prop:a_priori}.
  Then, we see from \cref{eq:second_derivative_metric} that 
  $$|u^{(2)}_{\tau_i}|(t) = \| \nabla \tilde{E}(U_{\tau_i, \lfloor t / \tau_i \rfloor + 1}) \|_{\LtwoR} \leq \|\Theta_T\|_{\LtwoR}$$
  for all $t \in [0, T-\tau_i]$ and for all $i$ large enough. 
  Thus, after passing to a subsequence (without relabeling), we can choose $A$ such that $|u_{\tau_i}^{(2)}| \rightharpoonup A$ weakly in $L^2_{loc}([0, +\infty))$. 
  Likewise, we have that
  $$|u^{(1)}_{\tau_i}|(t) = \|U'_{\tau_i, \lfloor t / \tau_i \rfloor + 1}\|_{L^2} \leq \|\Theta_T\|_{L^2},$$ 
  and so, again, we can choose $B$ such that $|u_{\tau_i}^{(1)}| \rightharpoonup B$ weakly in $L^2_{loc}([0, +\infty))$. These estimates guarantee that $(u_{\tau_i}, \dot{u}_{\tau_i})$ fulfill the conditions of the Arzelà--Ascoli argument in \cref{prop:l2_is_aa_prep,prop:aa-ready-allows-aa} (as in the proof of \cref{lem:compactness_newtonian}) with $A + B$ as the modulus of continuity. Thus, after passing to a further subsequence, we can find continuous curves $\bar{u}, \dot{\bar{u}}\colon [0, +\infty) \to \LtwoR$ such that, for all $t \in [0, +\infty)$, $(u_{\tau_i}(t), \dot{u}_{\tau_i}(t)) \to (\bar{u}_t, \dot{\bar{u}}_t)$ in $\ful{\law}^{-1}(\sigma_{\text{narrow}})$, and such that also the finite-dimensional joint laws converge as in \cref{eq:joint_law_wasserstein}.
  It remains to show that $\dot{\bar{u}}$ is indeed the derivative of $\bar{u}$. Exactly as in the proof of \cref{lem:compactness_newtonian}, applying the convergence of the joint laws to 
  $$f \colon \big( (x_1, v_1), (x_2, v_2) \big) \mapsto |x_1 + (T-s) v_1 - (x_2 + (T-t) v_2)|^2$$ and using a truncation argument as in the proof of \cref{cor:l2_is_aa}, we deduce that, for any $s, t, T \in [0, +\infty)$ with $s \leq t \leq T$, we have
  \begin{equation*}
    \liminf_{i \to \infty} \|\bar{\mathcal{V}}^{\tau_i, \varepsilon}_s(T) - \bar{\mathcal{V}}^{\tau_i, \varepsilon}_t(T)\|_{L^2} \geq \|\bar{\mathcal{V}}_s(T) - \bar{\mathcal{V}}_t(T)\|_{L^2}, \qquad \bar{\mathcal{V}}_s \coloneqq [\sigma \mapsto \bar{u}_s + (\sigma - s) \dot{\bar{u}}_s],
  \end{equation*}
  where $\bar{\mathcal{V}}^{\tau, \varepsilon}$ is defined in \cref{eq:bar_v_wasserstein}. Combining this inequality with \cref{eq:newtonian_sequences_are_pseudo_ac} and taking the limit $i \to \infty$, we obtain $\|\bar{\mathcal{V}}_s(T) - \bar{\mathcal{V}}_t(T)\|_{L^2} \leq \int_s^t (T - r) A(r) \, \diff r$. 
  This implies that $t\mapsto \bar u(t)$ is a $W^{2,2}$-curve (see \cite[\cref*{M-lem:euclidean_ac2_equivalence}]{fornasierNewtonianMotionCompanion}) with $|\ddot{\bar{u}}(t)| \leq A$ for almost all $t \in [0, +\infty)$\footnote{More precisely, the above inequality implies that $t\mapsto \bar u(t)$ is an \ACt-curve as per \cite[\cref*{M-def:AC2_curve}]{fornasierNewtonianMotionCompanion}, and \cite[\cref*{M-lem:euclidean_ac2_equivalence}]{fornasierNewtonianMotionCompanion} establishes the equivalence between the notion of $W^{2,2}$- and \ACt-curve.}, so, $\bar{u} \in C^1$. 
  Again owing to \cite[\cref*{M-lem:euclidean_ac2_equivalence}]{fornasierNewtonianMotionCompanion} (more precisely, its proof), we have that 
  $$\frac{\diff}{\diff s} \bar{u}_s = \frac{\diff}{\diff \sigma} \bar{\mathcal{V}}_s(\sigma) = \dot{\bar{u}}_s$$
  for almost every $s$, and hence, as both sides are continuous in $s$, for every $s \in [0, +\infty)$. In particular, this gives that $\dot{\bar{u}}$ is in fact the derivative of $\bar{u}$. 
  Finally, $\bar{u}_0 = U_0$ and $\dot{\bar{u}}_0 = U'_0$ follow as in the last line of the proof of \cref{lem:compactness_newtonian}.\\
  Next, we assume that $\limsup_{i \to \infty} \int_s^t\mathbf{H}(\tilde{\mathcal{V}}^{\tau_i, \varepsilon}_{r})\,\diff r \leq \int_s^t\mathbf{H}(\bar{u}_r, \dot{\bar{u}}_r)\, \diff r$ for every $s, t \in [0, +\infty)$ such that $s < t$ and we show that the limit curve $\bar{u}$ is a curve of Newtonian motion. To this end, by~\cref{eq:discrete_newtonian_pointwise}, we have the discrete identity
  \begin{equation}\label{eq:discrete_newtonian_motion_generalized}
		d(\tilde E u_{\tau_i})(t_{\tau_i}^{j_1}) - d(\tilde E u_{\tau_i})(t_{\tau_i}^{j_2}) = \frac{1}{2}\int_{t_{\tau_i}^{j_1}}^{t_{\tau_i}^{j_2}}  |u_{\tau_i}^{(2)}|^2(r)\, \diff r + \frac{1}{2} \int_{t_{\tau_i}^{j_1}}^{t_{\tau_i}^{j_2}} |G_{\tau_i}|^2(r)\, \diff r - \int_{t_{\tau_i}^{j_1}}^{t_{\tau_i}^{j_2}} \mathbf{H}(\tilde{\mathcal{V}}^{{\tau_i}, \varepsilon}_{r}) \, \diff r,
	\end{equation}
  for all $i, j_1, j_2 \in \NN$, setting $t_{\tau_i}^j = \tau_i \, j$ with $j_1 \leq j_2$. We recall that $d(\tilde E u_\tau)(t_\tau^i)$ is defined as in \cref{eq:dE_def}. Here, our goal is to take the limit in the above equation.
  First, note that, as we saw in \cref{sec:lifting_gradient_problems}, $|\partial \tilde{E}| = |\partial E| \circ \law$.
  By \cref{rmk:continuity_of_functionals}, the $\LtwoR$-continuity of $|\partial \tilde{E}|$ is equivalent to $W_2$-continuity of $|\partial E|$; here, $|\partial \tilde{E}| = \|\nabla \tilde{E}\|_{L^2}$ is $\LtwoR$-continuous since $\tilde{E} \in C^1$.
  By the uniformly integrable second moments (which follow from the domination by $\Theta_T$ in \cref{prop:a_priori}), the convergence in $\ful{\law}^{-1}(\sigma_{\text{narrow}})$ shown above gets upgraded to convergence in $W_2$.
  Furthermore, as $\tau_i \|\dot{u}_{\tau_i}(t)\|_{\LtwoR} \leq \tau_i \|\Theta_T\|_{L^2}$, we also have the $W_2$-convergence $\ful{\law}((u_{\tau_i}(t) - \tau_i \dot{u}_{\tau_i}(t), \dot{u}_{\tau_i}(t))) \to \ful{\law}((\bar{u}(t), \dot{\bar{u}}(t)))$. Thus we get, for all $t \in [0, +\infty)$,
  \begin{align}\label{eq:discrete_identity_generalized}
    d(\tilde E u_{\tau_i})(t) \to \nabla \tilde{E}(\bar{u}(t))(\dot{\bar{u}}(t)), \quad
    |\partial \tilde E|(u_{\tau_i}(t)) \to |\partial \tilde E|(\bar{u}(t)),
  \end{align}
  For the first convergence, we set $(X_i, Y_i) \coloneqq (u_{\tau_i}(t) - \tau_i \dot{u}_{\tau_i}(t), \dot{u}_{\tau_i}(t))$ and note that, as $\tilde{E} \in C^1$,
  \begin{equation*}
    d(\tilde E u_{\tau_i})(t) = \frac{1}{\tau_i}\big(\tilde{E}(X_i + \tau_i Y_i) - \tilde{E}(X_i)\big) = \int_0^1 \big\langle \nabla \tilde{E}(X_i + \sigma \tau_i Y_i), Y_i \big\rangle \, \diff \sigma.
  \end{equation*}
  The integrand equals $\frac{\diff}{\diff h} \tilde{E}(X_i + h Y_i)|_{h = \sigma \tau_i}$ and thus only depends on $\ful{\law}((X_i, Y_i))$; hence, by \cref{eq:wasserstein_distance_as_inf_of_l2_distance}, we may replace $(X_i, Y_i)$ by a pair with the same joint law such that $(X_i, Y_i) \to (\bar{u}(t), \dot{\bar{u}}(t))$ in $\LtwoR$. Then $X_i + \sigma \tau_i Y_i \to \bar{u}(t)$ in $\LtwoR$, uniformly in $\sigma \in [0, 1]$, and the continuity of $\nabla \tilde{E}$ together with the boundedness of $(Y_i)_{i \in \NN}$ in $\LtwoR$ yields 
  $$d(\tilde E u_{\tau_i})(t) \to \langle \nabla \tilde{E}(\bar{u}(t)), \dot{\bar{u}}(t) \rangle = \nabla \tilde{E}(\bar{u}(t))(\dot{\bar{u}}(t)).$$
  Next, by assumption, $J_{\delta \cdot \tau}[Z_{\tau, i}]$ as in \cref{eq:de_giorgi_interpolant} is nonempty for $\tau$ small enough; in particular, $\tilde{u}_{\tau_i}$ is well-defined, and this provides the nonemptiness required for \cref{eq:discrete_newtonian_motion_generalized} in~\cite[\cref*{M-prop:discrete_newtonian_motion}]{fornasierNewtonianMotionCompanion}.
  We first claim that $\law\big(\tilde{u}_{\tau_i}(s)\big) \to \law(\bar{u}(s))$ in $W_2$. 
  To see this, we combine the following:
  \begin{itemize}
      \item a standard a priori estimate for the proximal map (as in the proof of \cref{thm:cusp});
      \item the coercivity of $\tilde{E}$ (as per \cite[\cref*{M-ass:coercivity}]{fornasierNewtonianMotionCompanion}, which follows from the nonemptiness of $J_\delta[X]$);
      \item the boundedness of $\tilde{E}(Z_{\tau_i, n}) = E\big(\law(Z_{\tau_i, n})\big)$ (as $\law(Z_{\tau_i, n})$ converges in $W_2$ and $E$ is $W_2$-continuous by \cref{rmk:continuity_of_functionals});
      \item the bound on $\|U'_{\tau_i, n}\|$ in \cref{prop:a_priori};
  \end{itemize}
  and we deduce that $\tilde{u}_{\tau_i}(s) - u_{\tau_i}(s) \to 0$ in $\LtwoR$ as $\tau_i \to 0$ and for all $s \in [0, +\infty)$. 
  Together with the $W_2$ convergence noted before, we have $\law\big(\tilde{u}_{\tau_i}(s)\big) \to \law(\bar{u}(s))$ in $W_2$.\\
  Moreover, $|G_\tau|$ as in \cref{eq:G_tau} is well-defined and by \cite[\cref*{M-prop:discrete_newtonian_motion}]{fornasierNewtonianMotionCompanion} we have $|\partial \tilde{E}|(\tilde{u}_{\tau_i}(t)) \leq |G_{\tau_i}|(t)$ for all $i \in \NN$ and almost all $t \in [0, +\infty)$. 
  Combining this estimate with the aforementioned $W_2$-convergence, and using that $|\partial \tilde{E}| = |\partial E| \circ \law$ and $|\partial^- \tilde{E}| = \mathrm{lsc}(|\partial E|) \circ \law$ by \cref{rmk:continuity_of_functionals}, we get 
  $$\liminf_{i \to \infty}|G_{\tau_i}|(t) \geq |\partial^- \tilde{E}|(\bar{u}(t)) = |\partial \tilde{E}|(\bar{u}(t))$$
  for almost all $t$, where the last equality holds as $|\partial \tilde{E}|$ is continuous.
  Thus, we have, for any $s, t \in [0, +\infty)$ with $s < t$ of the form $s=k_s 2^{-m}$, $t = k_t 2^{-m}$, that
  \begin{equation} \label{eq:estimate_newtonian_curve_2}
  \begin{split}
    	 \int_s^t \frac{1}{2} |\bar{u}''&|^2(r)  \, \diff r + \int_s^t \frac{1}{2} |\partial \tilde E|^2(\bar{u}(r)) \diff r - \int_s^t \mathbf{H}(V^{\bar{u}}_{r}) \, \diff r \\
			& \leq \liminf_{i \to \infty} \bigg(\int_s^t \frac{1}{2} |u_{\tau_i}^{(2)}|^2(r) \, \diff r\bigg) + \int_s^t \liminf_{i \to \infty} \frac{1}{2} |G_{\tau_i}|^2(r) \diff r - \limsup_{i \to \infty} \int_s^t \mathbf{H}(\tilde{\mathcal{V}}^{\tau_i, \varepsilon}_{r}) \, \diff r \\
			& \leq \liminf_{i \to \infty} \bigg(\int_s^t \frac{1}{2} |u_{\tau_i}^{(2)}|^2(r) \,  \diff r + \int_s^t \frac{1}{2} |G_{\tau_i}|^2(r) \, \diff r - \int_s^t \mathbf{H}(\tilde{\mathcal{V}}^{\tau_i, \varepsilon}_{r}) \, \diff r\bigg) \\
			& \leq \liminf_{i \to \infty} \big(d(\tilde Eu_{\tau_i})(s) - d(\tilde Eu_{\tau_i})(t)\big)=  \frac{\diff}{\diff \sigma} \tilde E(V^{\bar{u}}_{s}(\sigma)) \Big|_{\sigma = s} - \frac{\diff}{\diff \sigma} \tilde E(V^{\bar{u}}_{t}(\sigma)) \Big|_{\sigma = t},
  \end{split}
  \end{equation}
  	where, in the first step, we used that $\limsup_{i \to \infty} \int_s^t \mathbf{H}(\tilde{\mathcal{V}}^{\tau_i, \varepsilon}_{r}) \, \diff r \leq \int_s^t \mathbf{H}(V^{\bar{u}}_{r}) \, \diff r$ by assumption, we used Fatou's lemma in the second step---note that $|G_{\tau_i}|^2(r) \geq 0$ for all $r$ and $i$, and that $\limsup_{i \to \infty} \int_s^t \mathbf{H}(\tilde{\mathcal{V}}^{\tau_i, \varepsilon}_{r}) \, \diff r < +\infty$ by assumption---, and used \cref{prop:discrete_newtonain_curve_equation} in the third, and \cref{eq:discrete_identity_generalized} in the last step. 
    By continuity of both sides---the left-hand side is continuous in $(s, t)$ since $|\bar{u}''| \leq A \in L^2_{loc}$, since $|\partial \tilde{E}|(\bar{u}_\bullet)$ is continuous, and since $\mathbf{H}(V^{\bar{u}}_\bullet) \in L^1_{loc}$ by the chain rule along the $W^{2,2}$-curve $\bar{u}$; the right-hand side equals $\langle \nabla \tilde{E}(\bar{u}_s), \dot{\bar{u}}_s \rangle - \langle \nabla \tilde{E}(\bar{u}_t), \dot{\bar{u}}_t \rangle$, which is continuous as $\bar{u} \in C^1$ and $\nabla \tilde{E}$ is continuous---, we can expand \cref{eq:estimate_newtonian_curve_2} to all $s, t \in [0, +\infty)$ such that $s < t$.


\end{proof}

\bookmarksetupnext{level=subsubsection}
\subsection{Examples}

The first example is the case of concave $C^2$-regular potentials and kernels.
\begin{example} \label{ex:C2}
  Let $V \in C^2(\R^d)$ and $\mathcal{K} \in C^2(\R^d \times \R^d)$, fulfilling \cref{ass:growth_condition}, with subquadratic growth and with Hessians uniformly bounded from above (i.e., there exists $C \in \R$ such that for all $x_1, y_1 \in \R^d$, $x_2, y_2 \in \R^{2d}$, we have that $\nabla^2 V(x_1)[y_1, y_1] \leq C|y_1|^2$ and $\nabla^2 \mathcal{K}(x_2)[y_2, y_2] \leq C|y_2|^2$). 
  Then, we define $\mathbf{H}\colon \LtwoR \times  \LtwoR \to \R \cup \{\pm \infty\}$ as
  \begin{equation} \label{eq:integral_h_def_in_ex}
      \mathbf{H}(X, W) \coloneqq \int_{\Omega_\R} \nabla^2 V(X)[W, W] \, \diff \mathbb{P}_\R + \iint_{\Omega_\R^2} \nabla^2 \mathcal{K}\big(X \otimes X\big)\big[W \otimes W, W \otimes W] \, \diff \mathbb{P}_\R^2,
  \end{equation}
  we have that the pair $(|\nabla \tilde{E}|, \mathbf{H})$ is an upper gradient-Hessian pair by \cref{cor:upper_gradient_hessian}. 
  Furthermore, $\mathbf{H}$ is strongly upper semicontinuous: Strong $L^2$ convergence implies pointwise convergence of the arguments; $V$ and $\mathcal{K}$ being in $C^2$ then implies pointwise convergence of the integrands in \cref{eq:integral_h_def_in_ex}; the reverse Fatou then gives upper semicontinuity of $\mathbf{H}$ with regards to the strong topology in $L^2$; by \cref{rmk:continuity_of_functionals} (applied in $L^2(\Omega_\R; \R^{2d}) \cong L^2(\Omega_\R; \R^d)^2$) and since $\mathbf{H}$ factors through $\law$, this is equivalent to upper semicontinuity with regards to $\law^{-1}(W_2)$. 
  Moreover, because $V$ and $\mathcal{K}$ have subquadratic growth, we can argue as in the proof of \cref{thm:cusp} to see that for $\tau$ small enough, the proximal point operator $J_{\tau, \LtwoR}[X]$ is nonempty for all $X \in \LtwoR$.
  Thus, for any $X, Y \in L^2(\Omega_\R; \R^d)$, there exists a curve of Newtonian motion $\bar{u}_{\bullet}\colon [0, +\infty) \to \LtwoR$ such that $\bar{u}_0 = X$ and $\dot{\bar{u}}_0 = Y$ by \cref{prop:compactness}. 
  Again, by \cref{cor:upper_gradient_hessian}, this curve of Newtonian motion fulfills the ODE $\ddot{\bar{u}}_t = -\nabla \tilde{E}(\bar{u}_t)$ for almost every $t \in [0, +\infty)$. 
  If we consider the function $G\colon \LtwoR \to (C_0^2(\R^d))'$ defined as
  $$\langle G(Z),  \psi \rangle \coloneqq  \int_{\R^d} \psi \, \diff \law(Z)$$
  for every $Z\in \LtwoR$ and for every $\psi \in C_0^2(\R^d)$,
  we have seen in \cref{prop:frechet-derivative-L2-W2} that $D^*G[Z] = [\psi \mapsto \nabla \psi \circ Z]$. In particular, by applying the analogous statement for the phase-space map $\ful{G}$ as in the proof of \cref{prop:abstract_pde_newtonian}, the phase-space law $\ful{\mu}_\bullet \coloneqq \ful{\law} \circ (\bar{u}_\bullet, \dot{\bar{u}}_\bullet)$ weakly fulfills the Vlasov equation
  \begin{equation*}
    \partial_t \ful{\mu}_t + v \cdot \nabla_x \ful{\mu}_t - \nabla_v \cdot \bigg(\Big(\nabla V + \int \big(\nabla_{x_1} \mathcal{K}(\cdot, x') + \nabla_{x_2} \mathcal{K}(x', \cdot)\big) \, \diff \ful{\mu}_t(x', v')\Big) \ful{\mu}_t\bigg) = 0,
  \end{equation*}
  i.e., 
  \begin{equation*}
    \begin{split}
    \frac{\diff}{\diff t} \int \varphi_t \, \diff \ful{\mu}_t = \int \partial_t \varphi_t \, \diff \ful{\mu}_t + \int v \cdot \nabla_x \varphi_t(x, v) \, \diff \ful{\mu}_t(x, v) - \langle \nabla \tilde{E}(\bar{u}_t), \nabla_v \varphi_t(\bar{u}_t, \dot{\bar{u}}_t) \rangle_{L^2}
    \end{split}
  \end{equation*}
  for every $\varphi \in C_c^\infty(\R^{2d+1})$ and almost every $t \in [0, +\infty)$
\end{example}

The kernels we considered in \cref{sec:ex_nonlocal}
are not $C^2$ and hence are not covered by the previous example. 
This is investigated in the next example.

\begin{example}[{Nonsingular attraction-repulsion, cf.~\cref{sec:ex_nonlocal}}]\label{ex:newtonian_mmd}
  Let $K_s(z) = |z|^{s+1}$ for $z \in \R^d$ and $s \in (0, 1)$\footnote{We emphasise that we do \emph{not} consider the singular case $s=0$ here.}, let $\mathcal{K}_s(x, y) = -\frac{1}{2} K_s(x - y)$, and let $V_{s, \nu} = x \mapsto \int_{\R^d} K_s(x - y) \, \diff \nu(y)$, for a reference signed measure $\nu = f \, \Leb^d$ with $f \in L^\infty(\R^d) \cap L^1(\R^d)$ and such that $\int_{\R^d} |y|^{s+1} |f(y)|\, \diff y < +\infty$. Then, up to an additive constant, the energy \eqref{eq:potential_interaction_energy} with $V = V_{s, \nu}$ and $\mathcal{K} = \mathcal{K}_s$ is the energy in \cref{eq:definition_of_potential}. Furthermore, $V_{s, \nu} \in W^{2, \infty}_{loc}$ with $\nabla^2 V_{s, \nu} \in L^\infty$ by \cref{lem:linear_potential_globally_lipschitz}.
  In particular, \cref{ass:growth_condition} is met. The main obstruction is that $\mathcal{K}_s$ is not $C^2$, since $\nabla^2 \mathcal{K}_s$ explodes around the diagonal $D = \{(x, x) \mid x \in \R^d\} \subseteq (\R^d)^2$. However, extending $\nabla^2 \mathcal{K}_s$ onto $D$ by setting
  \begin{equation*}
    \nabla^2 \mathcal{K}_s(x, x)[y, y] = \begin{cases}
      0 \qquad &\text{if } y \in D, \\
      -\infty \qquad &\text{otherwise,}
    \end{cases}
  \end{equation*}
  we have that $\nabla^2 \mathcal{K}_s$ still fulfills a chain rule as in \cref{eq:chain_rule_kernel}: We pick any $W^{2,2}$-curve $\gamma^{\mathcal{K}}\colon [a, b] \to \R^{2d}$ and note that by a mollification argument and because $\mathcal{K}_s$ is concave, we have that, for all $t_0, t_1 \in [a, b]$ such that $t_0 < t_1$,
  \begin{equation}
    \langle \nabla \mathcal{K}_s(\gamma^{\mathcal{K}}_{t_0}), \dot\gamma^{\mathcal{K}}_{t_0} \rangle - \langle \nabla \mathcal{K}_s(\gamma^{\mathcal{K}}_{t_1}), \dot\gamma^{\mathcal{K}}_{t_1} \rangle \geq \int_{t_0}^{t_1} -\langle \nabla \mathcal{K}_s(\gamma^{\mathcal{K}}_r), \ddot\gamma^{\mathcal{K}}_r \rangle \, \diff r.
  \end{equation}
  Since the integrand on the right hand side is in $L^1$, we obtain that the function $f\colon r \mapsto \langle \nabla \mathcal{K}_s(\gamma^{\mathcal{K}}_{r}), \dot\gamma^{\mathcal{K}}_{r} \rangle$ has finite variation on $[a, b]$. Note that $f$ also maps any nullset $N \subseteq [a, b]$ to a nullset: The image of $[a, b] \cap (\gamma^{\mathcal{K}})^{-1}(D)$ under $f$ is $\{0\}$, since $\nabla \mathcal{K}_s = 0$ on $D$; $[a, b] \setminus (\gamma^{\mathcal{K}})^{-1}(D)$ consists of a countable union of relatively open intervals, and on each of those intervals, $f$ is locally absolutely continuous, and so maps nullsets to nullsets. 
  Thus, $f$ also being continuous, it is absolutely continuous by the Banach-Zarecki theorem (see, e.g., \cite[Theorem~6.3.1]{heil2019introduction}). 
  Hence, we have to show that 
  \begin{equation}\label{eq:differential_equality}
    \frac{\diff}{\diff \sigma}\langle \nabla \mathcal{K}_s(\gamma^{\mathcal{K}}_\sigma), \dot\gamma^{\mathcal{K}}_{\sigma} \rangle|_{\sigma=r} = \langle \nabla \mathcal{K}_s(\gamma^{\mathcal{K}}_r), \ddot\gamma^{\mathcal{K}}_r \rangle + \nabla^2 \mathcal{K}_s(\gamma^{\mathcal{K}}_r)[\dot{\gamma}^{\mathcal{K}}_r, \dot{\gamma}^{\mathcal{K}}_r],
  \end{equation}
  for almost every $r \in [a, b]$; then it follows that the chain rule holds.
  Note that \eqref{eq:differential_equality} holds outside of $(\gamma^{\mathcal{K}})^{-1}(D)$. For $r \in (\gamma^{\mathcal{K}})^{-1}(D)$, it holds whenever $\dot\gamma^{\mathcal{K}}_r \in D$ \emph{and} the left hand side in \eqref{eq:differential_equality} is zero. Thus, the only points where it might not hold are those where either (1) $\langle \nabla \mathcal{K}_s(\gamma^{\mathcal{K}}_\sigma), \dot\gamma^{\mathcal{K}}_{\sigma} \rangle|_{\sigma=r} = 0$ (since $\gamma^{\mathcal{K}}_r \in D$), but its derivative is not, or (2) $\gamma^{\mathcal{K}}_r \in D$, but $\dot\gamma^{\mathcal{K}}_r \not\in D$. Both of those conditions are only fulfilled for countably many $r$ (since for any absolutely continuous curve, points with non-vanishing derivative are isolated in their level sets). Thus, the chain rule \cref{eq:chain_rule_kernel} holds. Next, we note that $V_{s, \nu} \in C^2(\R^d)$: To see this, we split $\nabla^2 V_{s, \nu}$ as in the proof of \cref{lem:linear_potential_globally_lipschitz}, and with,
  for $x, z \in \R^d$ close:
  \begin{equation*}
    \nabla^2 V_{s, \nu}(x) = \int_{B_1(z)} \underbrace{(\nabla^2 K_s)(x - y)}_{\in L^1(B_1(z))} f(y)\, \diff y + \int_{B_1^c(z)} (\nabla^2 K_s)(x-y) f(y)\, \diff y.
  \end{equation*}
  This expression is continuous in $x$, since translation of a function in $L^1$ is continuous (for the first term) and by dominated convergence (for the second term, where $|(\nabla^2 K_s)(x-y)| \leq (s+1)(2-s)\left|\frac{1}{2}\right|^{s-1}$ if $|z-x| \leq \frac{1}{2}$). Thus, also the pointwise chain rule \cref{eq:pointwise_chain_rule_v} holds with $H_V = \nabla^2 V_{s, \nu}$.
  Since $\nabla^2 V_{s, \nu} \in L^\infty$, and since $\mathcal{K}_s$ is concave, both pointwise Hessians are quadratically bounded from above as required by \cref{prop:chain_rule_generalized}. Thus, by \cref{cor:upper_gradient_hessian}, with the pointwise Hessian $\mathbf{H}$ defined as in \cref{eq:integral_h_def_in_ex}, the pair $(|\nabla \tilde{E}|, \mathbf{H})$ is an upper gradient-Hessian pair.\\
    As in \cref{ex:C2}, because $V$ and $\mathcal{K}_s$ have subquadratic growth, we can argue as in the proof of \cref{thm:cusp} to see that for $\tau$ small enough, the proximal point operator $J_{\tau, \LtwoR}[X]$ is nonempty for all $X \in \LtwoR$. Furthermore, the upper semicontinuity on $\mathbf{H}$ demanded by \cref{prop:compactness} is fulfilled: Indeed, the extended $\nabla^2 \mathcal{K}_s$ is pointwise upper semicontinuous and the map $(x, y) \mapsto \nabla^2 \mathcal{K}_s(x)[y, y]$ has quadratic upper growth in $x, y$ (and likewise for $\nabla^2 V_{s, \nu}$), and hence we can invoke \cite[Lemma 5.1.7, Lemma 7.1.5]{ambrosioGradientFlows2008} to see that $\mathbf{H}$ is $W_2$-upper semicontinuous. 
    Combining this upper semicontinuity with the convergence~\eqref{eq:compactness_convergence}, we see that for each fixed $t \in [0, +\infty)$, it holds that $\limsup_{i \to \infty} \mathbf{H}(\tilde{\mathcal{V}}^{\tau_i, \varepsilon}_{t}) \leq \mathbf{H}(\bar{u}_t, \dot{\bar{u}}_t)$. 
    By reverse Fatou and the uniform (in time) upper bound on $\mathbf{H}(\tilde{\mathcal{V}}^{\tau_i, \varepsilon}_{t})$ following from quadratic upper growth of the pointwise Hessians and \cref{prop:a_priori}, we even have $\limsup_{i \to \infty} \int_s^t \mathbf{H}(\tilde{\mathcal{V}}^{\tau_i, \varepsilon}_{r})\, \diff r \leq \int_s^t \mathbf{H}(\bar{u}_{r}, \dot{\bar{u}}_r)\, \diff r$, for $s, t \in [0, +\infty)$ such that $s < t$.  
    Thus, we may invoke \cref{prop:compactness} to see that for any $X, Y \in \LtwoR$ there exists a curve of Newtonian motion $\bar{u}_\bullet\colon [0, +\infty) \to \LtwoR$ such that $\bar{u}_0 = X$ and $\dot{\bar{u}}_0 = Y$, and by \cref{cor:upper_gradient_hessian}, the curve $\bar{u}$ fulfills the ODE $\ddot{\bar{u}}_t = - \nabla \tilde{E}(\bar{u}_t)$ for almost every $t \in [0, +\infty)$. 
    Lastly, as in the previous example, we deduce that the phase-space law $\ful{\mu}_\bullet \coloneqq \ful{\law} \circ (\bar{u}_\bullet, \dot{\bar{u}}_\bullet)$ fulfills weakly in the sense of distributions the Vlasov equation
    \begin{equation*}
      \partial_t \ful{\mu}_t + v \cdot \nabla_x \ful{\mu}_t - \nabla_v \cdot \bigg(\Big(\nabla V_{s, \nu} - \int \nabla K_s(\cdot - x') \, \diff \ful{\mu}_t(x', v')\Big) \ful{\mu}_t\bigg) = 0.
  \end{equation*}
\end{example}

\bookmarksetupnext{level=subsection}
\section{Weak Solutions to PDEs without Using Curves of Maximal Slope}\label{sec:weak_solutions_without_coms}

In \cref{sec:wasserstein_as_quotient,sec:PDE_gradient_flows}, we provided two different tools:
\begin{enumerate}
  \item A relaxed form of Arzelà--Ascoli, which allows us to construct curves of maximal slope directly in $L^2$, without passing to $(\Pp_2, W_2)$ (see \cref{rmk:generalized_aa_yields_coms}).
  \item A weak description of the differentials of curves of maximal slope, which allows us to identify which PDEs are weakly fulfilled by curves of maximal slope in $L^2$.
\end{enumerate}
Taking those arguments together, this allows one to construct weak solutions to PDEs in the following way\footnote{This is a simplified picture of the situation; in practice and depending on the functional, one might use a slight variation of the below picture. Furthermore, here we assume that $|\partial^-E| = |\partial E|$, i.e., that $|\partial E|$ is weakly lsc.}: \\
\begin{tikzcd}[arrows=Rightarrow, column sep=small,font=\small]
  & \fbox{\Centerstack[l]{A priori bounds\\\cite[L. 3.2.2]{ambrosioGradientFlows2008}\\\cite[R. 3.2.3]{ambrosioGradientFlows2008}}} \arrow{d}  & \fbox{\Centerstack[l]{$|\partial^- E|$ is a\\strong upper\\gradient}} \arrow{dr}\\
  \fbox{\Centerstack[l]{MMS scheme}} \arrow{ur} \arrow{dr} & \fbox{\Centerstack[l]{Arzelà--Ascoli\\\cref{prop:aa-ready-allows-aa}\\\cref{cor:l2_is_aa}}} \arrow{r} & \fbox{\Centerstack[l]{Limit curve}} \arrow{d} & \fbox{\Centerstack[l]{Curve of \\maximal\\slope}} \arrow{r} & \fbox{\Centerstack[l]{Integral\\
  equality on \\
  test functions\\\cref{prop:weak-formulation-abstract}}}\\
  & \fbox{\Centerstack[l]{Integral inequality I\\\cite[Eq.3.2.4]{ambrosioGradientFlows2008}}} \arrow{r} & \fbox{\Centerstack[l]{Integral inequality II\\\cite[Eq. 3.4.4]{ambrosioGradientFlows2008}}} \arrow{ur}
\end{tikzcd}
In the above picture, \emph{Integral inequality I} is an integral inequality that holds on certain interpolations of the discrete MMS solutions; in passing through the limit, it yields \emph{Integral inequality II}.

Using the machinery we developed so far, we are able to identify PDEs which are fulfilled by curves of maximal slope---the last step in the above picture---by only testing the derivative of $E$ into certain directions---those corresponding to test functions---without requiring a general understanding of $E$. However, \emph{constructing} curves of maximal slope---all the other steps in the above picture---still requires a strong understanding of $E$, e.g., one needs to prove that $|\partial^-E|$ is a strong upper gradient.

On the other hand, the above picture immediately suggest a shortcut often taken in practice\footnote{In fact, such a ``shortcut'' was used in~\cite{jordan1998variational}, where the scheme~\eqref{eq:mms} was used in Wasserstein space to construct solutions to the Fokker-Planck equation. This was done before \cite{ambrosioGradientFlows2008} made available the theory of curves of maximal slope in Wasserstein space to analyze PDEs, and led to the name ``JKO scheme'' for the scheme~\eqref{eq:mms} in the case $(\mathcal{S}, d_{\mathcal{S}}) = (\Pp_2, W_2)$.}: In our current approach, we have an integral inequality on the MMS interpolations---\emph{integral inequality I}---which passes to an integral inequality of the limit curve---\emph{integral inequality II}. In the last step, we prove integral equalities of the limit curve. Instead, we might prove integral equalities on the MMS interpolations, and pass them to the limit, obtaining integral equalities on the limit curve directly, bypassing completely the necessity of the limit curve to be a curve of maximal slope: \\
\begin{center}
\begin{tikzcd}[arrows=Rightarrow, column sep=small,font=\small]
  & \fbox{\Centerstack[l]{A priori bounds\\\cite[L. 3.2.2]{ambrosioGradientFlows2008}\\\cite[R. 3.2.3]{ambrosioGradientFlows2008}}} \arrow{d}\\
  \fbox{\Centerstack[l]{MMS scheme}} \arrow{ur} \arrow{dr} & \fbox{\Centerstack[l]{Arzelà--Ascoli\\\cref{prop:aa-ready-allows-aa}\\\cref{cor:l2_is_aa}}} \arrow{r} & \fbox{\Centerstack[l]{Limit curve}} \arrow{d}\\
  & \fbox{\Centerstack[l]{Integral\\ equalities on\\test functions I}} \arrow{r} & \fbox{\Centerstack[l]{Integral\\
  equalities on \\
  test functions}} 
\end{tikzcd}
\end{center}
Let us informally show what we hope to get for the \emph{integral equalities on test functions I} in the special case of $E$ being the Boltzmann entropy as in \cref{ex:bolzmann-coms}. We choose a test function $\varphi$ and $\tau$ small. We denote by $(U_\tau^n)_{n=0}^{+\infty}$ a solution to the MMS scheme on the Boltzmann entropy, by $\bar{U}_\tau$ the corresponding piecewise constant interpolation satisfying $\bar{U}_\tau(t) = U_\tau^{\lceil t / \tau \rceil}$, and we set $t_i \coloneqq \tau i$. We then hope to formalize the following informal computation:
\begin{align}
  \langle G(\bar{U}_\tau(T)), \varphi_T \rangle_{C_0^2} &- \langle G(\bar{U}_\tau(0)), \varphi_0 \rangle_{C_0^2} = \sum_{i=1}^{\lceil T / \tau \rceil} \langle G(U_\tau^i), \varphi_{t_i} \rangle - \langle G(U_\tau^{i-1}), \varphi_{t_{i-1}}  \rangle \\
  &= \sum_{i=1}^{\lceil T / \tau \rceil} \langle G(U_\tau^i) - G(U_\tau^{i-1}), \varphi_{t_i} \rangle + \langle G(U_\tau^{i-1}), \varphi_{t_i} - \varphi_{t_{i-1}}  \rangle \\
  &\approx \sum_{i=1}^{\lceil T / \tau \rceil} \langle U_\tau^i - U_\tau^{i-1}, \nabla \varphi_{t_i}(U_\tau^i) \rangle + \tau \langle G(U_\tau^{i-1}), \partial_t \varphi_{t_{i-1}} \rangle \label{eq:approximate-g-differential-discrete}\\
  &= \sum_{i=1}^{\lceil T / \tau \rceil} \tau \bigg(-\frac{\diff}{\diff s} E\big(\law(U_\tau^i + s \nabla \varphi_{t_i}(U_\tau^i))\big)\Big|_{s=0}  + \langle G(U_\tau^{i-1}), \partial_t \varphi_{t_{i-1}} \rangle \bigg) \label{eq:approximate-coms-technical-discrete} \\
  &= \sum_{i=1}^{\lceil T / \tau \rceil} \tau \bigg( \langle G(U_\tau^i), \Delta \varphi_{t_i}\rangle  + \langle G(U_\tau^{i-1}), \partial_t \varphi_{t_{i-1}} \rangle \bigg) \\
  &\approx \int_0^T \int \Delta \varphi_t + \partial_t \varphi_t \diff (\bar{U}_\tau(t)_\# \mathbb{P}) \diff t
\end{align}
If we are able to control all errors by controlling $\tau$ we are able to show that any limit curve---in the appropriate sense---of piecewise interpolations $\bar{U}_\tau$ as $\tau \to 0$ is a weak solution to the heat equation. For the approximation in \cref{eq:approximate-g-differential-discrete}, we need an approximate, discrete version of \cref{prop:frechet-derivative-L2-W2}; for the equality in \cref{eq:approximate-coms-technical-discrete}, we need a discrete version of \cref{lem:technical-coms}. We start with the latter, where the proof is much easier than that of \cref{lem:technical-coms}.
\begin{lemma} \label{lem:untechnical-coms}
  Let $H$ be a Hilbert space, let $F\colon H \to \R$ be a functional, and let $x \in H$. Furthermore, let some
  \begin{equation*}
    y \in \argmin_{z \in H} F(z) + \frac{|x - z|^2}{2\tau}
  \end{equation*}
  be given. Then, for any $v \in H$, we have that
  \begin{equation*}
    \langle y - x, v \rangle = - \frac{\diff}{\diff s} F(y + s v) \bigg{|}_{s=0}\cdot \tau,
  \end{equation*}
  whenever the differential on the right-hand side exists.
\end{lemma}
\begin{proof}
  Let us assume that the statement does not hold; then by replacing $v$ with $-v$ if necessary, we may assume that $\langle y - x, v \rangle < - \frac{\diff}{\diff s} F(y + s v)\big{|}_{s=0} \cdot \tau$. Setting $w^s \coloneqq y  + s  v - x$, we have that
  \begin{equation*}
    \frac{\diff}{\diff s} \frac{|w^s|^2}{2} \bigg{|}_{s=0} = \langle y - x, v \rangle < - \frac{\diff}{\diff s} F(y + s v) \bigg{|}_{s=0}\cdot \tau.
  \end{equation*}
  Thus, for $s$ small enough, we have that
  \begin{equation*}
    \frac{|w^s|^2}{2\tau} - \frac{|w^0|^2}{2\tau} < F(x + w^0) - F(x + w^s),
  \end{equation*}
  which, after rewriting and substituting the definition of $w^s$, yields
  \begin{equation*}
    \frac{|y  + s  v - x|^2}{2\tau} + F(y + sv)  < \frac{|y - x|^2}{2\tau} + F(y),
  \end{equation*}
  which contradicts the optimality of $y$.
\end{proof}
Next, let us continue with an approximate, discrete variant of \cref{prop:frechet-derivative-L2-W2}.
\begin{lemma}\label{lemma:approximate_discrete_frechet_derivative_l2_w2}
  Let us assume the situation of \cref{fig:all_spaces}. Let $X, Y \in L^2$, set $G \coloneqq r \circ b \circ \law$ and assume $\psi \in C_0^2$. Then, we have that
  \begin{equation*}
    \big| \langle G(Y) - G(X), \psi \rangle_{C_0^2} - \langle Y - X, \nabla \psi(Y) \rangle_{L^2}\big| \leq \frac{1}{2} \|Y - X\|_{L^2}^2 \cdot \|\nabla^2 \psi\|_{L^\infty}
  \end{equation*}
\end{lemma}
\begin{proof}
  We calculate
  \begin{align*}
    &\big| \langle G(Y) - G(X), \psi \rangle_{C_0^2} - \langle Y - X, \nabla \psi(Y) \rangle_{L^2}\big| \\
    =\, &\big| \int \psi(Y) - \psi(X) - \langle Y - X, \nabla \psi(Y) \rangle \, \diff \mathbb{P} \big| \\
    =\, &\big| \int - \int_0^1 (1-s)\langle Y - X, \nabla^2 \psi \big(Y + s (X - Y)\big) (Y - X) \rangle \, \diff s\, \diff \mathbb{P}| \\
    \leq\, &\frac{1}{2} \int |Y - X|^2 \cdot \|\nabla^2 \psi\|_{L^\infty} \, \diff \mathbb{P} = \frac{1}{2}\|Y - X\|_{L^2}^2 \cdot \|\nabla^2 \psi\|_{L^\infty}
  \end{align*}
\end{proof}
Before we write down the main result of this section, let us formally gather the assumptions we put on $E$.
\begin{assumption}\label{ass:psi_good_in_good_directions}
  The functional $E$ is narrowly lower semicontinuous. For any function $\psi \in C_c^\infty$ and measure $\mu \in D(E)$, the limit $\frac{\diff}{\diff s} E\big((T_s^\psi)_\#\mu\big)\big|_{s=0}$ exists. Furthermore, for fixed $\psi \in C_c^\infty$, the quantity $\frac{\diff}{\diff s} E\big((T_s^\psi)_\#\mu\big)\big|_{s=0}$ is uniformly continuous in $\mu$ on sublevels of $E$ with respect to the bounded-Lipschitz distance on $\Pp_2$.
\end{assumption}
\begin{theorem}\label{thm:weak-pde-wo-coms}
  Let $E\colon \Pp_2 \to \Rext$ be a functional bounded from below fulfilling \cref{ass:psi_good_in_good_directions}. Furthermore, let $(\tau_i)_{i \in \NN}$ be a nullsequence and for each $\tau_i$, assume that we are given a MMS sequence $(U_{\tau_i}^n)_{n \in \NN}$ with respect to $E \circ \law$ and a constant $S > 0$ such that $E\big(\law(U_{\tau_i}^0)\big) \leq S$ and $\|U_{\tau_i}^0\|_{L^2} \leq S$ for all $i \in \NN$. Defining, for each $i \in \NN$, $\bar{U}_{\tau_i}(t) \coloneqq U_{\tau_i}^{\lceil t / \tau_i \rceil}$, we furthermore assume that there exists a limit curve $u$ such that $\law\big(\bar{U}_{\tau_i}(t)\big) \rightharpoonup \law\big(u(t)\big)$ narrowly in $\Pp_2$, uniformly on $[0, T]$ with respect to the bounded-Lipschitz distance for each $T > 0$. Then, for each $C_c^\infty(\R^{d + 1}) \ni \varphi\colon (t, x) \mapsto \varphi_t(x)$ and each $T > 0$, the law of the limit curve $\mu_t \coloneqq \law\big(u(t)\big)$ fulfills the integral equality 
  \begin{equation}\label{eq:abstract-PDE-strong}
    \int \varphi_T \, \diff \mu_T - \int \varphi_0 \, \diff \mu_0 = \int_0^T \bigg( - \frac{\diff}{\diff s} E\big((T_s^{\varphi_t})_\#\mu_t\big)\Big|_{s=0} + \int \partial_t \varphi_t \, \diff \mu_t \bigg) \, \diff t
  \end{equation}
\end{theorem}
\begin{proof}
  We fix $\varphi$ and $T$ and show~\cref{eq:abstract-PDE-strong}. By~\cite[Lemma 3.2.2]{ambrosioGradientFlows2008}, we have that there exists $C > 0$ such that we have, up to discarding finitely many $i$,
  \begin{equation}\label{eq:integral_speed_estimate}
    \sum_{j=0}^{\lceil T / \tau_i \rceil - 1} \frac{|U_{\tau_i}^j- U_{\tau_i}^{j+1}|^2}{2\tau_i} \leq C, \qquad \forall i \in \NN,
  \end{equation}
  and
  \begin{equation}\label{eq:pointwise_speed_estimate}
    |U_{\tau_i}^j - U_{\tau_i}^{j+1}|^2 \leq C \tau_i \qquad \forall i \in \NN, \, j \leq \lceil T / \tau_i \rceil,
  \end{equation}
  and also $|U_{\tau_i}^j| \leq C$. We set, for each $i \in \NN$, $V_{\tau_i}(t) \coloneqq \frac{U_{\tau_i}^{\lceil t / \tau_i \rceil} - U_{\tau_i}^{\lfloor t / \tau_i \rfloor}}{\tau_i}$, $\tilde{V}_{\tau_i}(t) \coloneqq \frac{G(U_{\tau_i}^{\lceil t / \tau_i \rceil}) - G(U_{\tau_i}^{\lfloor t / \tau_i \rfloor})}{\tau_i}$ and calculate
  \begin{align}
    \int \varphi_T \, \diff \mu_T &- \int \varphi_0 \, \diff \mu_0 = \lim_{i \to \infty} \int \varphi_T \, \diff (\bar{U}_{\tau_i}(T)_\# \mathbb{P}) - \int \varphi_0 \, \diff (\bar{U}_{\tau_i}(0)_\# \mathbb{P}) \\
    &= \lim_{i \to \infty} \sum_{j=0}^{\lceil T / \tau_i \rceil - 1} \langle G(U_{\tau_i}^{j+1}), \varphi_{(j+1) \cdot \tau_i} \rangle - \langle G(U_{\tau_i}^{j}), \varphi_{(j) \cdot \tau_i} \rangle \label{eq:limit_proof_1}\\
    &= \lim_{i \to \infty} \sum_{j=0}^{\lceil T / \tau_i \rceil - 1} \langle G(U_{\tau_i}^{j+1}) - G(U_{\tau_i}^{j}), \varphi_{(j+1) \cdot \tau_i} \rangle + \langle G(U_{\tau_i}^{j}), \varphi_{(j+1) \cdot \tau_i} - \varphi_{(j) \cdot \tau_i} \rangle  \label{eq:limit_proof_2} \\
    &= \lim_{i \to \infty} \int_{0}^T \langle \tilde{V}_{\tau_i}(t), \varphi_t \rangle \diff t + \sum_{j=0}^{\lceil T / \tau_i \rceil - 1}  \tau_i\langle G(U_{\tau_i}^{j}), \partial_t \varphi_{j \cdot \tau_i} \rangle \label{eq:limit_proof_3} \\
    &= \lim_{i \to \infty} \int_{0}^T \Big( \langle \tilde{V}_{\tau_i}(t), \varphi_t \rangle + \langle G(\bar{U}_{\tau_i}(t)), \partial_t \varphi_t \rangle \Big) \, \diff t \label{eq:limit_proof_4} \\
    &= \lim_{i \to \infty} \int_{0}^T \Big( \langle V_{\tau_i}(t), \nabla \varphi_t(\bar{U}_{\tau_i}(t)) \rangle + \langle G(\bar{U}_{\tau_i}(t)), \partial_t \varphi_t \rangle \Big) \, \diff t \label{eq:limit_proof_5} \\
    &= \lim_{i \to \infty} \int_{0}^T \bigg( - \frac{\diff}{\diff s} E\big(\law(\bar{U}_{\tau_i}(t) + s \nabla \varphi_t(\bar{U}_{\tau_i}(t)))\big)\Big|_{s=0} + \langle G(\bar{U}_{\tau_i}(t)), \partial_t \varphi_t \rangle \bigg) \, \diff t \label{eq:limit_proof_6} \\
    &= \lim_{i \to \infty} \int_{0}^T \bigg( - \frac{\diff}{\diff s} E\Big((T_s^{\varphi_t})_\#(\bar{U}_{\tau_i}(t)_\#\mathbb{P})\Big)\Big|_{s=0} + \langle G(\bar{U}_{\tau_i}(t)), \partial_t \varphi_t \rangle \bigg) \, \diff t \label{eq:limit_proof_7} \\
    &= \int_{0}^T \lim_{i \to \infty} \bigg( - \frac{\diff}{\diff s} E\Big((T_s^{\varphi_t})_\#(\bar{U}_{\tau_i}(t)_\#\mathbb{P})\Big)\Big|_{s=0} + \langle G(\bar{U}_{\tau_i}(t)), \partial_t \varphi_t \rangle \bigg) \, \diff t \label{eq:limit_proof_8} \\
    &=\int_{0}^T \bigg( - \frac{\diff}{\diff s} E\big((T_s^{\varphi_t})_\#\mu_t\big)\Big|_{s=0} + \int \partial_t \varphi_t \,\diff \mu_t \bigg) \, \diff t. \label{eq:limit_proof_9}
  \end{align}
  Let us justify each step of the calculation above: In~\eqref{eq:limit_proof_1}, the terms differ by at most $\tau_i \|\partial_t \varphi\|_\infty$, as $0 \leq \lceil T/\tau_i \rceil \tau_i - T < \tau_i$; in~\eqref{eq:limit_proof_2}, we have an exact equality in each term of the sequence. In~\eqref{eq:limit_proof_3}, we used (1.) that $|\tilde{V}_{\tau_i}(t)| \cdot \tau_i$ goes to zero uniformly as $i \to \infty$---by~\cref{eq:pointwise_speed_estimate} and using that $G$ is 1-Lipschitz by~\cref{prop:frechet-derivative-L2-W2}---to transform the first part of the sum to an integral and (2.) that $\frac{1}{\tau_i} \langle \mu, \varphi_{(j+1) \cdot \tau_i} - \varphi_{(j) \cdot \tau_i} \rangle$ converges to $\langle \mu, \partial_t \varphi_{j \cdot \tau_i} \rangle$, uniformly over $\mu \in \Pp_2$ to transform the second part of the sum. In~\eqref{eq:limit_proof_4}, we used that $\langle G(\mu), \partial_t \varphi_t \rangle$ is bounded uniformly over $\mu \in \Pp_2$ and $t \in \R$, and again~\cref{eq:pointwise_speed_estimate} together with the 1-Lipschitz continuity of $G$. Further, in~\eqref{eq:limit_proof_5}, we first apply~\cref{lemma:approximate_discrete_frechet_derivative_l2_w2} to obtain the error bound $|\langle \tilde{V}_{\tau_i}(t), \varphi_t \rangle - \langle V_{\tau_i}(t), \nabla \varphi_t(\bar{U}_{\tau_i}(t)) \rangle| \leq \frac{\tau_i}{2} \|V_{\tau_i}(t)\|_{L^2}^2 \cdot \|\nabla^2 \varphi\|_{L^\infty}$ and then~\cref{eq:integral_speed_estimate} to bound the global error by $\tau_i C \|\nabla^2 \varphi\|_{L^\infty}$. In~\eqref{eq:limit_proof_6}, we again have exact equality under the integral using~\cref{lem:untechnical-coms}; likewise, in~\eqref{eq:limit_proof_7}, we simply rewrite the limit. Finally, in~\eqref{eq:limit_proof_8} and~\eqref{eq:limit_proof_9}, we use the uniform convergence of $\law(\bar{U}_{\tau_i}(t)) \rightharpoonup \law(u(t))$ in $[0, T]$ in the assumptions of the theorem and the uniform continuity in $\mu$ of $\frac{\diff}{\diff s} E\big((T_s^\varphi)_\#\mu\big)\big|_{s=0}$ on sublevels of $E$ (note that $E(\law(\bar{U}_{\tau_i}(t))) \leq S$ along the MMS, and hence $E(\mu_t) \leq S$ by lower semicontinuity), which together yield uniform convergence of the integrand and thus allow for dominated convergence.
\end{proof}
\begin{remark}
  The required uniform continuity in $\mu$ of $\frac{\diff}{\diff s} E\big((T_s^\varphi)_\#\mu\big)\big|_{s=0}$ on sublevels of $E$ is required for the limiting passages in \cref{eq:limit_proof_8,eq:limit_proof_9}, and may look too strong for certain applications of interest. However, in practice, these passages are instead often proven ad-hoc, by means of problem-specific analysis tools.
\end{remark}
\bookmarksetupnext{level=subsection}
\section{\texorpdfstring{Proof of \cref{prop:frechet-derivative-L2-W2}}{Proof of Lemma~\getrefnumber{prop:frechet-derivative-L2-W2}}}\label{sec:proof-of-frechet-derivative-l2-w2}
\begin{proof}[Proof of \cref{prop:frechet-derivative-L2-W2}]
  To warm up, we start by showing that $G$ is $1$-Lipschitz. To this end, we have to show that for every $X, Y \in L^2$ and every $\psi \in C_0^2$ such that $\|\psi\|_{C^2} \leq 1$, we have
  \begin{equation}
    \Big| \int \psi(x) \, \diff (\law(X))(x) - \int \psi(x) \, \diff (\law(Y))(x) \Big| \leq \|X - Y\|_{L^2}.
  \end{equation}
  We calculate
  \begin{align}
    &\Big| \int \psi(x) \, \diff (\law(X))(x) - \int \psi(x) \, \diff (\law(Y))(x) \Big| = \Big| \int \psi(X(\omega)) - \psi(Y(\omega)) \, \diff \mathbb{P}(\omega) \Big| \\
    &\leq \int |\psi(X(\omega)) - \psi(Y(\omega))| \, \diff \mathbb{P}(\omega) \leq \int \|\nabla \psi\|_\infty |X(\omega) - Y(\omega)|_{\R^d} \, \diff \mathbb{P}(\omega) \\
    &\leq \int |X(\omega) - Y(\omega)|_{\R^d} \,\diff \mathbb{P}(\omega) \leq  \|X - Y\|_{L^2}.
  \end{align}
  Next, we show the Fréchet differentiability. First, we show the equality in \cref{eq:frechet-derivative-L2-W2}. Starting from the right-hand side, we have
  \begin{align}
    \int \big\langle \nabla \psi\big(X(\omega)\big), H(\omega) \big\rangle \, \diff \mathbb{P}(\omega) &= \int \big\langle \nabla \psi\big(X(\omega)\big), \mathbb{E}[H | X](\omega) \big\rangle \, \diff \mathbb{P}(\omega) \\
    &= \int \big\langle \nabla \psi\big(X(\omega)\big), Z\big(X(\omega)\big) \big\rangle \, \diff \mathbb{P}(\omega) \\
    &= \int \big\langle \nabla \psi(x), Z(x) \big\rangle \diff (\law(X))(x).
  \end{align}
Next, note that the map $(H, \psi) \mapsto DG(X)[H](\psi)$ as defined in \cref{eq:frechet-derivative-L2-W2} is bilinear and continuous in both components:
\begin{equation*}
  \left| \int \big\langle \nabla \psi\big(X(\omega)\big), H(\omega) \big\rangle \, \diff \mathbb{P}(\omega) \right| \leq \|\nabla \psi\|_\infty \|H\|_{L^1} \leq \|\psi\|_{C^2} \|H\|_{L^2}.
\end{equation*}
What is left to show is that $DG(X)$ is indeed the Fréchet derivative of $G$. To this end, we calculate for $X, H \in L^2$ and $\psi \in C_0^2$:
\begin{align}
  &\Big| G(X + H)(\psi) - G(X)(\psi) - DG(X)[H](\psi) \Big| \\
  &= \Big| \int \psi\big(X(\omega) + H(\omega)\big) - \psi\big(X(\omega)\big) - \langle \nabla \psi\big(X(\omega)\big), H(\omega) \rangle \, \diff \mathbb{P}(\omega) \Big| \\
  &\leq \frac12 \int \|\nabla^2 \psi\|_\infty |H(\omega)|_{\R^d}^2 \diff \mathbb{P}(\omega) \leq  \frac12 \|\psi\|_{C^2} \|H\|_{L^2}^2,
\end{align}
  and thus
  \begin{equation}
    \lim_{H \to 0} \frac{\|G(X + H) - G(X) - DG(X)[H]\|_{(C_0^2)'}}{\|H\|_{L^2}} = 0,
  \end{equation}
  which proves the claim.
\end{proof}

\bookmarksetupnext{level=subsection}

\section{\texorpdfstring{Proof of \cref{lem:energy_lipschitz_cont}}{Proof of Lemma~\getrefnumber{lem:energy_lipschitz_cont}}}\label{sec:proof-of-lipschitz-cont}
\begin{proof}[Proof of \cref{lem:energy_lipschitz_cont}]
  We fix $X, Y \in L^2$, and calculate
  \begin{equation} \label{eq:wasserstein_diff_int_exchange}
    \begin{split}
    &\frac{\diff}{\diff \tau} E\big(\law(X + \tau Y)\big) \\
    =\,& \frac{\diff}{\diff \tau}\bigg[\iint K\big(X(\omega) + \tau Y(\omega), X(\omega') + \tau Y(\omega')\big) \, \diff \mathbb{P}(\omega) \, \diff \mathbb{P}(\omega') + \int V\big(X(\omega) + \tau Y(\omega) \big) \, \diff \mathbb{P}(\omega)\bigg] \\
    \overset{(*)}{=}\,& \iint \frac{\diff}{\diff \tau} K\Big(X(\omega) + \tau Y(\omega), X(\omega') + \tau Y(\omega')\Big)\, \diff \mathbb{P}(\omega) \, \diff \mathbb{P}(\omega') + \int \frac{\diff}{\diff \tau} V\Big(X(\omega) + \tau Y(\omega)\Big)\, \diff \mathbb{P}(\omega)\\
    =\, &2 \iint \nabla_{x_1} K\Big(X(\omega), X(\omega')\Big) \cdot Y(\omega) \, \diff \mathbb{P}(\omega) \, \diff \mathbb{P}(\omega') + \int \nabla V\Big(X(\omega)\Big) \cdot Y(\omega) \, \diff \mathbb{P}(\omega);
  \end{split}
  \end{equation}
  where we justify the exchange of differentiation and integration in $(*)$ as follows: To exchange with the outer integral, we use the global Lipschitz continuity of the derivative (by \cref{ass:kernel_and_potential}) to see that for fixed $\omega'$ and $\tau \in (0, 1)$, and setting $x = X(\omega')$, $y=Y(\omega')$
  \begin{align*}
    &\bigg|\frac{\int K(X(\omega) + \tau Y(\omega), x + \tau y)\,\diff \mathbb{P}(\omega) - \int K(X(\omega), x) \, \diff \mathbb{P}(\omega)}{\tau}\bigg| \\
    &\leq \int \sup_{\sigma \in [0, \tau]} \Big|\nabla_{x_1,x_2} K\Big(X(\omega) + \sigma Y(\omega), x + \sigma y\Big) \cdot (Y(\omega), y)\Big| \, \diff \mathbb{P}(\omega) \\
    &\leq \int \Big(|\nabla_{x_1,x_2} K(0, 0)| + C_{\mathrm{Lip}} \big(|X(\omega)| + |Y(\omega)| + |x| + |y|\big)\Big) \Big(|Y(\omega)| + |y|\Big) \, \diff \mathbb{P}(\omega) \\
    &\leq C (1 + |y|^2 + |x|^2 + \|X\|_{L^2}^2 + \|Y\|_{L^2}^2)
  \end{align*}
  for some $C > 0$ large enough; a similar calculation for the global potential allows us to pick $C' > 0$ large enough such that the finite differences for the whole (outer) integrand are dominated by some function
  \begin{equation}
    \Big[\omega \mapsto C' \big(1 + |X(\omega)|^2 + |Y(\omega)|^2 + \|X\|_{L^2}^2 + \|Y\|_{L^2}^2\big)\Big] \in L^1
  \end{equation}
  for small enough step size, and thus we may exchange the differentiation and the outer integral; by a similar estimate we can move the differentiation into the inner integral of kernel.
  From the above calculation, we see that $\tilde E$ is Gateaux differentiable everywhere and that the Gateaux differential $\nabla \tilde E$ is linear and given by
  \begin{equation*}
    \nabla \tilde E(X) = \bigg[\omega \mapsto 2 \int \nabla_{x_1} K\Big(X(\omega), X(\omega')\Big)\,\diff\mathbb{P}(\omega') + \nabla V\Big(X(\omega)\Big)\bigg] \in L^2.
  \end{equation*}
  Showing that the map $\nabla \tilde E\colon L^2 \to L^2$ is globally Lipschitz will automatically upgrade the Gateaux differential to a Fréchet differential. To this end, we fix $X, X' \in L^2$ and calculate
  \begin{align*}
    & \| \nabla\tilde{E}(X) - \nabla\tilde{E}(X') \|^2_{L^2} \\
    = \,&\bigg\|\bigg[\omega \mapsto 2 \int \nabla_{x_1} K\big(X(\omega), X(\omega')\big) - \nabla_{x_1} K\big(X'(\omega), X'(\omega') \big) \,\diff\mathbb{P}(\omega') + \nabla V\big(X(\omega)\big) - \nabla V\big(X'(\omega)\big)\bigg]\bigg\|^2_{L^2} \\
    \leq \,& \bigg\|\bigg[ \omega \mapsto 2 \int C_{\mathrm{Lip}} \Big(\big|X(\omega') - X'(\omega')\big| + \big|X(\omega) - X'(\omega)\big|\Big) \, \diff \mathbb{P}(\omega') + C_{\mathrm{Lip}} \big|X(\omega) - X'(\omega)\big|\bigg]\bigg\|^2_{L^2} \\
    = \, & \bigg\|\bigg[ \omega \mapsto 2 \int C_{\mathrm{Lip}} \big|X(\omega') - X'(\omega')\big| \, \diff \mathbb{P}(\omega') + 3 C_{\mathrm{Lip}} \big|X(\omega) - X'(\omega)\big|\bigg]\bigg\|^2_{L^2} \\
    = \, & \mathlarger\int \bigg( 2 \int C_{\mathrm{Lip}} \big|X(\omega') - X'(\omega')\big| \, \diff \mathbb{P}(\omega') + 3 C_{\mathrm{Lip}} \big|X(\omega) - X'(\omega)\big| \bigg)^2 \, \diff \mathbb{P} (\omega) \\
    = \, & C_{\mathrm{Lip}}^2 \mathlarger{\int} \bigg(\quad 4 \Big(\int \big|X(\omega') - X'(\omega')\big| \, \diff \mathbb{P}(\omega')\Big)^2 \\
    &\qquad\quad\quad + 12  \int \big|X(\omega') - X'(\omega')\big| \, \diff \mathbb{P}(\omega') \big|X(\omega) - X'(\omega)\big|  + 9 \big|X(\omega) - X'(\omega)\big|^2 \bigg)\, \diff \mathbb{P} (\omega) \\
    \leq \, &C_{\mathrm{Lip}}^2 \mathlarger{\int} \bigg(\quad 4 \int \big|X(\omega') - X'(\omega')\big|^2 \, \diff \mathbb{P}(\omega') \\
    &\qquad\quad\quad + 6 \Big(\int \big|X(\omega') - X'(\omega')\big|^2 \, \diff \mathbb{P}(\omega') + \big|X(\omega) - X'(\omega)\big|^2\Big)  + 9 \big|X(\omega) - X'(\omega)\big|^2 \bigg)  \, \diff \mathbb{P} (\omega)  \\
    = \, &25 C_{\mathrm{Lip}}^2 \|X - X'\|^2_{L^2}.
  \end{align*}
  Next, we calculate the second derivative:
  \begin{align*}
    &\nabla^2 \tilde{E}(X)[W, \cdot] \\
    &= \frac{\diff}{\diff \tau} \bigg[\omega \mapsto 2 \int \nabla_{x_1} K(X(\omega) + \tau W(\omega), X(\omega') + \tau W(\omega')) \, \diff \mathbb{P}(\omega') + \nabla V(X(\omega) + \tau W(\omega)) \bigg] \\
    &= \bigg[ \omega \mapsto 2 \frac{\diff}{\diff \tau} \int \nabla_{x_1} K(X(\omega) + \tau W(\omega), X(\omega') + \tau W(\omega')) \, \diff \mathbb{P}(\omega') + \nabla^2 V(X(\omega))[W(\omega), \cdot] \bigg] \\
    &= \bigg[ \omega \mapsto 2 \int \nabla^2_{x_1} K(X(\omega), X(\omega'))[W(\omega), \cdot] + 2 \nabla_{x_2} \nabla_{x_1} K(X(\omega), X(\omega'))[W(\omega'), \cdot] \, \diff \mathbb{P}(\omega') \\
    &\qquad\qquad\qquad\qquad\qquad+ \nabla^2 V(X(\omega))[W(\omega), \cdot] \bigg].
  \end{align*}
  Here, we justify switching differentiation and integral in the second equality by the fact that the corresponding differential quotient in the inner integral is dominated, since $|\nabla^2_{x_1} K|$ and $|\nabla_{x_2} \nabla_{x_1} K|$ are globally bounded by \cref{ass:kernel_and_potential} and $W \in L^1$; similarly, for the first equality the integrand is dominated by the same argument and because in addition also $|\nabla^2 V|$ is globally bounded (and also $W \in L^2$). Lastly, we have to show that the map $X, W, W' \mapsto \nabla^2 \tilde{E}(X)[W, W']$ is jointly continuous in $X, W, W'$; to this end, we pick a sequence $(X_i, W_i, W'_i) \to (X, W, W')$ and calculate
  \begin{align*}
    \bigg|\int \nabla^2 &V(X_i(\omega))[W_i(\omega), W'_i(\omega)] - \int \nabla^2 V(X(\omega))[W(\omega), W'(\omega)] \bigg| \\
    &\leq M (\|W_i - W\|_{L^2} \, \|W'_i\|_{L^2} + \|W\|_{L^2} \, \|W'_i - W'\|_{L^2}) + \bigg|\int W (\nabla^2 V(X_i) - \nabla^2 V(X)) W' \bigg| \\
    &\to 0
  \end{align*}
  where $M \coloneqq \sup_{x \in \R^d} |\nabla^2 V(x)| < +\infty$ by \cref{ass:kernel_and_potential} and where the last term goes to zero by dominated convergence (with dominating function $2M|W||W'| \in L^1$) along any subsequence along which $X_i \to X$ almost everywhere, and hence, by the subsequence principle, along the whole sequence. Likewise, the difference of the terms involving the interaction kernel converges to zero, which yields the required joint convergence.
\end{proof}

\end{document}